\pdfoutput=1
\documentclass{amsart}
\usepackage{todonotes}
\usepackage{comment}
\usepackage{graphicx, amsmath, amssymb, amsthm,booktabs}
\usepackage{tikz-cd}
\usepackage{spectralsequences,bbm}
\usepackage{hyperref}
\usepackage[T1]{fontenc}
\hypersetup{
  colorlinks   = true,
  urlcolor     = blue,
  linkcolor    = blue,
  citecolor   = red
}
\usepackage{adjustbox}
\usepackage[normalem]{ulem}
\usepackage{contour}
\usepackage{mathtools} 

\contourlength{0.8pt}

\newcommand{\myuline}[1]{%
  \ifmmode
    \uline{\hphantom{#1}}
    \mathllap{\contour{white}{\ensuremath{#1}}}
  \else
    \uline{\phantom{#1}}%
    \llap{\contour{white}{#1}}%
  \fi
}
\usepackage[backend=biber, style=numeric, maxbibnames=99, doi=true,url=false,
    giveninits=true, hyperref]{biblatex}
\renewbibmacro{in:}{}

\newcommand{\ip}[1]{\langle#1\rangle}

\newcommand{\gr}{\mathrm{gr}}
\newcommand{\slice}[1]{\mathrm{Slice}(#1)}

\newcommand{\bprn}[1]{%
  BP_\R\ip{#1}%
}

\newcommand{\bpcm}[1]{BP^{(\!(C_4)\!)}\langle#1\rangle}

\newcommand{\bpgm}[1]{BP^{(\!(G)\!)}\langle#1\rangle}

\newcommand{\tmf}{\mathrm{tmf}}
\newcommand{\TMF}{\mathrm{TMF}}

\newcommand{\Z}{{\mathbb  Z}}

\newcommand{\R}{{\mathbb R}}
\newcommand{\F}{{\mathbb F}}

\DeclareMathOperator{\Hom}{Hom}
\DeclareMathOperator{\Ext}{Ext}
\DeclareMathOperator{\Tor}{Tor}

\DeclareMathOperator{\Res}{Res}

\newcommand{\Sp}{\mathcal Sp}

\mathchardef\mhyphen=45

\numberwithin{equation}{section}

\newtheorem{theorem}{Theorem}[section]
\newtheorem{lemma}[theorem]{Lemma}
\newtheorem{corollary}[theorem]{Corollary}

\newtheorem{proposition}[theorem]{Proposition}

\newtheorem*{theorem*}{Theorem}
\newtheorem*{proposition*}{Proposition}

\newtheorem{MainTheorems}{Theorem}

\theoremstyle{remark}
\newtheorem{remark}[theorem]{Remark}

\theoremstyle{definition}

\newtheorem{definition}[theorem]{Definition}

\newcommand{\Sq}{{\rm Sq}}

\SseqNewClassPattern{hsss}{
  (0.000,0.000);
  (-0.120,0.000)(0.120,0.000);
  (-0.200,0.000)(0.000,0.000)(0.200,0.000);
  (-0.120,0.100)(0.120,0.100)(-0.120,-0.100)(0.120,-0.100);
  (-0.200,0.100)(0.000,0.100)(0.200,0.100)(-0.100,-0.100)(0.100,-0.100);
  (-0.200,0.100)(0.000,0.100)(0.200,0.100)(-0.200,-0.100)(0.000,-0.100)(0.200,-0.100);
  (-0.200,0.200)(0.000,0.200)(0.200,0.200)(-0.200,0.000)(0.000,0.000)(0.200,0.000)(0.000,-0.200);
  (-0.200,0.200)(0.000,0.200)(0.200,0.200)(-0.200,0.000)(0.000,0.000)(0.200,0.000)(-0.100,-0.200)(0.100,-0.200);
  (-0.200,0.200)(0.000,0.200)(0.200,0.200)(-0.200,0.000)(0.000,0.000)(0.200,0.000)(-0.200,-0.200)(0.000,-0.200)(0.200,-0.200);
  (-0.240,0.200)(-0.080,0.200)(0.080,0.200)(0.240,0.200)(-0.240,0.000)(-0.080,0.000)(0.080,0.000)(0.240,0.000)(-0.080,-0.200)(0.080,-0.200);
  (-0.240,0.200)(-0.080,0.200)(0.080,0.200)(0.240,0.200)(-0.240,0.000)(-0.080,0.000)(0.080,0.000)(0.240,0.000)(-0.160,-0.200)(0.000,-0.200)(0.160,-0.200);
  (-0.240,0.200)(-0.080,0.200)(0.080,0.200)(0.240,0.200)(-0.240,0.000)(-0.080,0.000)(0.080,0.000)(0.240,0.000)(-0.240,-0.200)(-0.080,-0.200)(0.080,-0.200)(0.240,-0.200);
  (-0.240,0.240)(-0.080,0.240)(0.080,0.240)(0.240,0.240)(-0.240,0.080)(-0.080,0.080)(0.080,0.080)(0.240,0.080)(-0.240,-0.080)(-0.080,-0.080)(0.080,-0.080)(0.240,-0.080)(0.000,-0.240);
  (-0.240,0.240)(-0.080,0.240)(0.080,0.240)(0.240,0.240)(-0.240,0.080)(-0.080,0.080)(0.080,0.080)(0.240,0.080)(-0.240,-0.080)(-0.080,-0.080)(0.080,-0.080)(0.240,-0.080)(-0.080,-0.240)(0.080,-0.240);
  (-0.240,0.240)(-0.080,0.240)(0.080,0.240)(0.240,0.240)(-0.240,0.080)(-0.080,0.080)(0.080,0.080)(0.240,0.080)(-0.240,-0.080)(-0.080,-0.080)(0.080,-0.080)(0.240,-0.080)(-0.160,-0.240)(0.000,-0.240)(0.160,-0.240);
  (-0.240,0.240)(-0.080,0.240)(0.080,0.240)(0.240,0.240)(-0.240,0.080)(-0.080,0.080)(0.080,0.080)(0.240,0.080)(-0.240,-0.080)(-0.080,-0.080)(0.080,-0.080)(0.240,-0.080)(-0.240,-0.240)(-0.080,-0.240)(0.080,-0.240)(0.240,-0.240);
  (-0.320,0.240)(-0.160,0.240)(0.000,0.240)(0.160,0.240)(0.320,0.240)(-0.320,0.080)(-0.160,0.080)(0.000,0.080)(0.160,0.080)(0.320,0.080)(-0.320,-0.080)(-0.160,-0.080)(0.000,-0.080)(0.160,-0.080)(0.320,-0.080)(-0.080,-0.240)(0.080,-0.240);
  (-0.320,0.240)(-0.160,0.240)(0.000,0.240)(0.160,0.240)(0.320,0.240)(-0.320,0.080)(-0.160,0.080)(0.000,0.080)(0.160,0.080)(0.320,0.080)(-0.320,-0.080)(-0.160,-0.080)(0.000,-0.080)(0.160,-0.080)(0.320,-0.080)(-0.160,-0.240)(0.000,-0.240)(0.160,-0.240);
  (-0.320,0.240)(-0.160,0.240)(0.000,0.240)(0.160,0.240)(0.320,0.240)(-0.320,0.080)(-0.160,0.080)(0.000,0.080)(0.160,0.080)(0.320,0.080)(-0.320,-0.080)(-0.160,-0.080)(0.000,-0.080)(0.160,-0.080)(0.320,-0.080)(-0.240,-0.240)(-0.080,-0.240)(0.080,-0.240)(0.240,-0.240);
  (-0.320,0.240)(-0.160,0.240)(0.000,0.240)(0.160,0.240)(0.320,0.240)(-0.320,0.080)(-0.160,0.080)(0.000,0.080)(0.160,0.080)(0.320,0.080)(-0.320,-0.080)(-0.160,-0.080)(0.000,-0.080)(0.160,-0.080)(0.320,-0.080)(-0.320,-0.240)(-0.160,-0.240)(0.000,-0.240)(0.160,-0.240)(0.320,-0.240);
  (-0.320,0.320)(-0.160,0.320)(0.000,0.320)(0.160,0.320)(0.320,0.320)(-0.320,0.160)(-0.160,0.160)(0.000,0.160)(0.160,0.160)(0.320,0.160)(-0.320,0.000)(-0.160,0.000)(0.000,0.000)(0.160,0.000)(0.320,0.000)(-0.320,-0.160)(-0.160,-0.160)(0.000,-0.160)(0.160,-0.160)(0.320,-0.160)(0.000,-0.320);
  (-0.320,0.320)(-0.160,0.320)(0.000,0.320)(0.160,0.320)(0.320,0.320)(-0.320,0.160)(-0.160,0.160)(0.000,0.160)(0.160,0.160)(0.320,0.160)(-0.320,0.000)(-0.160,0.000)(0.000,0.000)(0.160,0.000)(0.320,0.000)(-0.320,-0.160)(-0.160,-0.160)(0.000,-0.160)(0.160,-0.160)(0.320,-0.160)(-0.080,-0.320)(0.080,-0.320);
  (-0.320,0.320)(-0.160,0.320)(0.000,0.320)(0.160,0.320)(0.320,0.320)(-0.320,0.160)(-0.160,0.160)(0.000,0.160)(0.160,0.160)(0.320,0.160)(-0.320,0.000)(-0.160,0.000)(0.000,0.000)(0.160,0.000)(0.320,0.000)(-0.320,-0.160)(-0.160,-0.160)(0.000,-0.160)(0.160,-0.160)(0.320,-0.160)(-0.160,-0.320)(0.000,-0.320)(0.160,-0.320);
  (-0.320,0.320)(-0.160,0.320)(0.000,0.320)(0.160,0.320)(0.320,0.320)(-0.320,0.160)(-0.160,0.160)(0.000,0.160)(0.160,0.160)(0.320,0.160)(-0.320,0.000)(-0.160,0.000)(0.000,0.000)(0.160,0.000)(0.320,0.000)(-0.320,-0.160)(-0.160,-0.160)(0.000,-0.160)(0.160,-0.160)(0.320,-0.160)(-0.240,-0.320)(-0.080,-0.320)(0.080,-0.320)(0.240,-0.320);
  (-0.320,0.320)(-0.160,0.320)(0.000,0.320)(0.160,0.320)(0.320,0.320)(-0.320,0.160)(-0.160,0.160)(0.000,0.160)(0.160,0.160)(0.320,0.160)(-0.320,0.000)(-0.160,0.000)(0.000,0.000)(0.160,0.000)(0.320,0.000)(-0.320,-0.160)(-0.160,-0.160)(0.000,-0.160)(0.160,-0.160)(0.320,-0.160)(-0.320,-0.320)(-0.160,-0.320)(0.000,-0.320)(0.160,-0.320)(0.320,-0.320);
}

\colorlet{dthree}{blue!75!black}
\colorlet{dfive}{red!75!black}
\colorlet{dseven}{green!50!black}
\colorlet{dnine}{orange!90!black}
\sseqset{
  hsss/.global sseq style = {
    Adams grading, class pattern = hsss,
    classes = {fill, draw = none, inner sep = 0pt, minimum size = 1.6pt},
    differentials = {line width = 0.25pt, shorten >= 0.4pt, shorten <= 0.4pt},
    struct lines = {black!55, line width = 0.3pt},
    xscale = 0.38, yscale = 0.38,
    grid = go, grid color = black!12,
    x tick step = 2, y tick step = 2,
    tick style = {font = \tiny},
    axes type = border,
  }
}

\colorlet{dtwo}{orange!85!black}
\colorlet{dthr}{red!75!black}
\colorlet{dfour}{green!50!black}
\sseqset{
  adams/.global sseq style = {
    Adams grading, class pattern = hsss,
    classes = {fill, draw = none, inner sep = 0pt, minimum size = 1.6pt},
    differentials = {line width = 0.3pt, shorten >= 0.4pt, shorten <= 0.4pt},
    struct lines = {black!55, line width = 0.25pt},
    xscale = 0.5, yscale = 0.45,
    grid = go, grid color = black!12,
    x tick step = 2, y tick step = 2,
    tick style = {font = \tiny},
    axes type = border,
    class labels = {fill = white, fill opacity = 0.85, text opacity = 1, inner sep = 0.3pt},
  },
  adamswide/.global sseq style = {
    Adams grading, class pattern = hsss,
    classes = {fill, draw = none, inner sep = 0pt, minimum size = 1.6pt},
    differentials = {line width = 0.3pt, shorten >= 0.4pt, shorten <= 0.4pt},
    struct lines = {black!55, line width = 0.25pt},
    xscale = 0.95, yscale = 0.55,
    grid = go, grid color = black!12,
    x tick step = 2, y tick step = 2,
    tick style = {font = \tiny},
    axes type = border,
    class labels = {fill = white, fill opacity = 0.85, text opacity = 1, inner sep = 0.3pt},
  },
  tmfclass/.sseq style = {draw = none, fill = blue!75!black, minimum size = 2.1pt},
}
\sseqset{d2/.sseq style = {dtwo}, d3/.sseq style = {dthr}, d4/.sseq style = {dfour}}
\colorlet{hiddencol}{violet!85!black}
\usetikzlibrary{arrows.meta}
\sseqset{
  hiddeneta/.sseq style = {hiddencol, dashed, line width = 0.45pt},
  hiddennu/.sseq style = {hiddencol, dashed, line width = 0.45pt},
}
\tikzset{hiddentwo/.style = {draw = hiddencol, dashed, line width = 0.45pt, -{Latex[length=2.2pt]}, shorten >= 1.2pt, shorten <= 1.2pt}}

\begin{document}

\title{On the Adams spectral sequence of $eo_2(C_4)$}
\author{Christian Carrick}
\address{Mathematisches Institut, Universität Bonn, Bonn, Germany}
\email{carrick@math.uni-bonn.de}
\thanks{The author was supported by the National Science Foundation under Grant No. DMS-2401918.}

\begin{abstract}
   We calculate the mod $2$-homology of the fixed points of the $C_4$-equivariant Hill--Hopkins--Ravenel theory $\bpcm{1}$, as an $\mathcal{A}_*$-comodule algebra. $2$-locally these fixed points are $\tmf_0(5)$ by a theorem of Quinn. We deduce that $\tmf_0(5)$ is an indecomposable $\tmf$-module that is an extension of two $\tmf$-modules of the form $\tmf\otimes K$ for a finite spectrum $K$. The $K(2)$-localization of $\tmf_0(5)$ is $32$-periodic, and we compute the Adams spectral sequence of $\tmf_0(5)$ through the $32$-stem.
\end{abstract}

\maketitle
\tableofcontents
\section{Introduction}\label{sec:intro}
Chromatic homotopy theory organizes the stable homotopy groups of spheres by height. The fundamental object at height $h$ is the $K(h)$-local sphere, and a theorem of Devinatz--Hopkins \cite{devhop} identifies it as a homotopy fixed point spectrum
\[L_{K(h)}\mathbb{S}\simeq E_h^{h\mathbb{G}_h},\]
where $E_h$ is a Lubin--Tate theory of height $h$ and $\mathbb{G}_h$ is the extended Morava stabilizer group. Hopkins and Miller observed that the fixed points $EO_h(G)=E_h^{hG}$ with respect to \emph{finite} subgroups $G\subseteq\mathbb{G}_h$ are tractable approximations to $L_{K(h)}\mathbb{S}$, and that they are higher height analogues of real $K$-theory \cite{hopmiller}: at heights $1$ and $2$ the maximal finite $2$-subgroups of $\mathbb{G}_h$ are $C_2$ and $G_{48}$, the fixed points $E_1^{hC_2}$ are the $2$-completion of $\mathrm{KO}$, and $E_2^{hG_{48}}$ is closely related to $\TMF$ (see \cite{ghmr}). These higher real $K$-theories are a principal source of what is known about $v_h$-periodic phenomena in $\pi_*\mathbb{S}$.

Working with the theories $EO_h(G)$ in practice is limited by the fact that they are not connective, their homotopy groups are not degreewise finitely generated, and their mod $p$ homology vanishes. Much of what is known about $\pi_*\mathbb{S}$ at the prime $2$ has been obtained with the classical Adams spectral sequence. This presents an obstacle to using the $EO_h(G)$'s to understand $\pi_*\mathbb S$, as the Adams spectral sequence of $EO_h(G)$ therefore vanishes.

The classical remedy at heights $1$ and $2$ is to replace $\mathrm{KO}$ and $\TMF$ by the connective theories $\mathrm{ko}$ and $\tmf$. These are connective spectra with strong finiteness properties, from which the periodic theories are recovered by $K(h)$-localization. They are the prototypical examples of fp spectra in the sense of Mahowald--Rezk \cite{MR}: bounded below, $2$-complete spectra whose mod $2$ cohomology is finitely presented over the Steenrod algebra. Their mod $2$ homology is induced up from the finite Hopf algebras $\mathcal{A}(1)$ and $\mathcal{A}(2)$ respectively, which makes their Adams spectral sequences computable, and these Adams spectral sequences are among the most useful tools in computational homotopy theory.

This is especially true of $\tmf$. Its Adams spectral sequence has been studied in great detail by Bruner--Rognes \cite{BR}, and the comparison with it organizes a significant amount of our knowledge about the Adams spectral sequence of the sphere, e.g. as in \cite{iwx} or \cite{bmq}. The unit $\mathbb{S}\to\tmf$ induces a map from the Adams spectral sequence of the sphere to that of $\tmf$, which on $E_2$-pages is the restriction map $\Ext_{\mathcal{A}}(\F_2,\F_2)\to\Ext_{\mathcal{A}(2)}(\F_2,\F_2)$, and since the target spectral sequence is completely understood, naturality carries information back to the source. For example, a class of $\Ext_{\mathcal{A}}(\F_2,\F_2)$ whose image supports a differential in $\tmf$ cannot be a permanent cycle, such as the differential $d_2(e_0)=h_1^2d_0$ in the $17$-stem, which is forced in the sphere by the same differential in $\tmf$. This kind of information requires knowing not only $\pi_*\tmf$, but also the structure of the Adams spectral sequence converging to it. In this way, studying the Adams spectral sequences of \emph{connective models} of $EO_n(G)$'s allows us to pass information from $EO_n(G)$ back to the Adams spectral sequence of the sphere.

This paper explores a new case of this phenomenon, namely by computing the mod $2$ homology and Adams spectral sequence of a connective model of $EO_2(C_4)$, known as $\tmf_0(5)$. We compute homology relative to $\tmf$. For a $\tmf$-module $X$, we write
\[H_*^{\tmf}(X)=\pi_*(\F_2\otimes_{\tmf}X)\]
for the homology of $X$ relative to $\tmf$. This is a comodule over $\mathcal{A}(2)_*$, it determines the mod $2$ homology of $X$ by $H_*(X;\F_2)\cong\mathcal{A}_*\Box_{\mathcal{A}(2)_*}H_*^{\tmf}(X)$, and so the Adams spectral sequence of $X$ has $E_2$-page $\Ext_{\mathcal{A}(2)_*}(\F_2,H_*^{\tmf}(X))$. Our main computational result is the following.

\begin{MainTheorems}\label{thm:introA}
The $\mathcal{A}(2)_*$-comodule algebra $H_*^{\tmf}(\tmf_0(5))$ is $24$-dimensional, concentrated in degrees $0$ through $18$, with Poincar\'e series
\[1+t+t^2+2t^4+t^5+t^6+2t^7+3t^8+t^{10}+2t^{11}+2t^{12}+t^{13}+2t^{14}+t^{15}+t^{16}+t^{17}+t^{18}.\]
Its coaction and multiplication are given explicitly in Theorem \ref{thm:Mprime}.
\end{MainTheorems}

The comodule of Theorem \ref{thm:introA} is of a different nature from the homology of $\tmf_1(3)$, $\tmf_0(3)$, and $\tmf_1(5)$, which are the homologies of finite complexes: as a $\tmf$-module, $\tmf_0(5)$ is built from two such modules, but is not one itself.

\begin{MainTheorems}\label{thm:introB}
Two-locally, there is a cofiber sequence of $\tmf$-modules
\[\tmf\otimes K_1\to\tmf_0(5)\to\tmf\otimes K_2,\]
where $K_1$ and $K_2$ are finite complexes with $12$ cells each, in dimensions $0$ through $8$ and $10$ through $18$, respectively. On the other hand, $\tmf_0(5)$ is not equivalent to $\tmf\otimes K$ for any spectrum $K$, and it is indecomposable as a $\tmf$-module.
\end{MainTheorems}

\tikzset{
  cell/.style={circle, fill=black, inner sep=0pt, minimum size=3pt},
  sqone/.style={line width=0.5pt},
  sqtwo/.style={line width=0.5pt},
  sqfour/.style={line width=0.5pt},
  clab/.style={font=\tiny, inner sep=1.5pt}
}
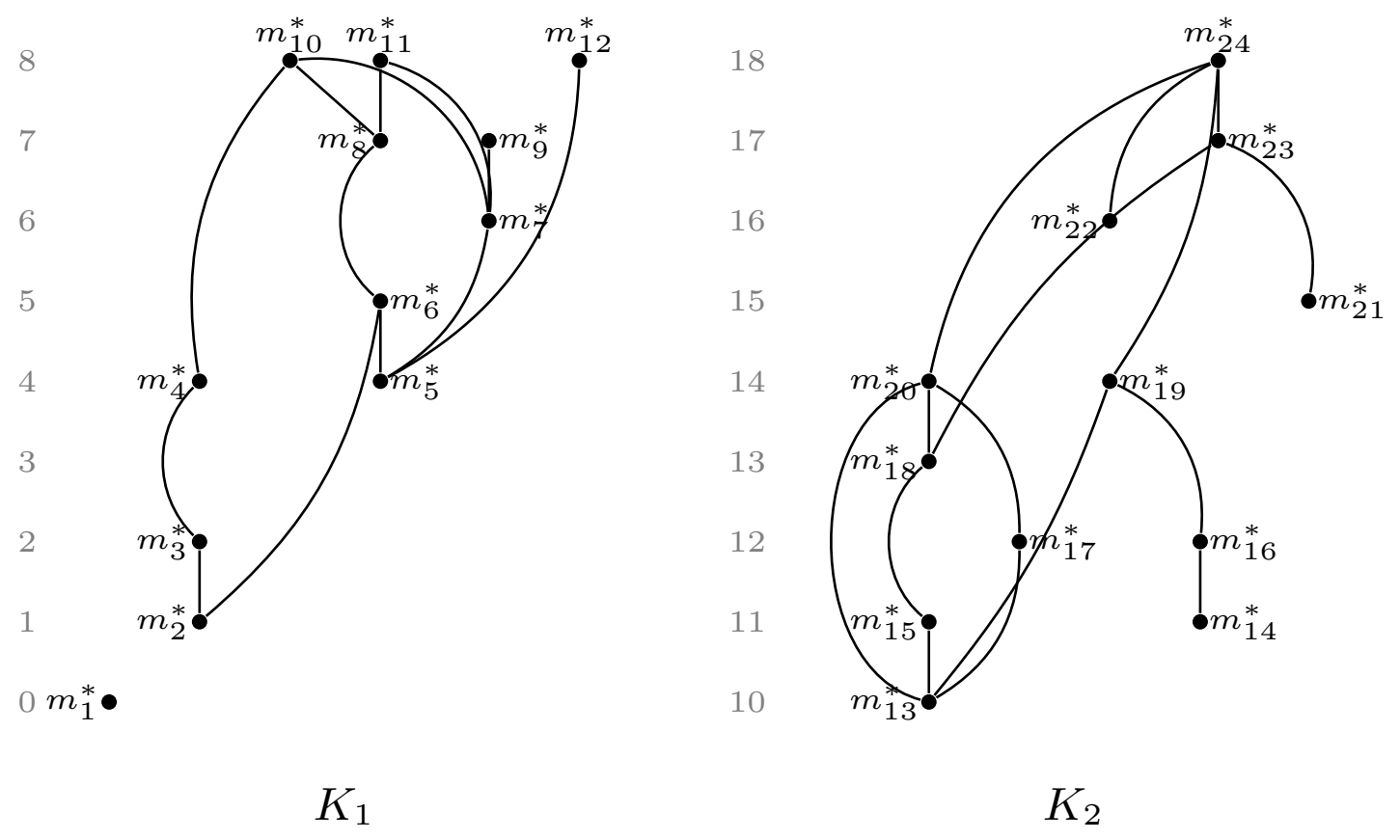
\begin{figure}[htbp]\centering
\begin{tikzpicture}[x=0.62cm,y=0.55cm]
\begin{scope}[xshift=0cm]
\foreach \j in {0,...,8} \node[font=\tiny, anchor=east, gray] at (-0.6,\j-0) {$\j$};
\node[cell] (m1) at (0,0) {}; \node[clab, left] at (m1) {$m_{1}^*$};
\node[cell] (m2) at (1,1) {}; \node[clab, left] at (m2) {$m_{2}^*$};
\node[cell] (m3) at (1,2) {}; \node[clab, left] at (m3) {$m_{3}^*$};
\node[cell] (m4) at (1,4) {}; \node[clab, left] at (m4) {$m_{4}^*$};
\node[cell] (m5) at (3,4) {}; \node[clab, right] at (m5) {$m_{5}^*$};
\node[cell] (m6) at (3,5) {}; \node[clab, right] at (m6) {$m_{6}^*$};
\node[cell] (m7) at (4.2,6) {}; \node[clab, right] at (m7) {$m_{7}^*$};
\node[cell] (m8) at (3,7) {}; \node[clab, left] at (m8) {$m_{8}^*$};
\node[cell] (m9) at (4.2,7) {}; \node[clab, right] at (m9) {$m_{9}^*$};
\node[cell] (m10) at (2,8) {}; \node[clab, above] at (m10) {$m_{10}^*$};
\node[cell] (m11) at (3,8) {}; \node[clab, above] at (m11) {$m_{11}^*$};
\node[cell] (m12) at (5.2,8) {}; \node[clab, above] at (m12) {$m_{12}^*$};
\draw[sqone] (m2)--(m3);
\draw[sqone] (m5)--(m6);
\draw[sqone] (m7)--(m9);
\draw[sqone] (m8)--(m10);
\draw[sqone] (m8)--(m11);
\draw[sqtwo] (m3) to[bend left=45] (m4);
\draw[sqtwo] (m5) to[bend right=25] (m7);
\draw[sqtwo] (m6) to[bend left=50] (m8);
\draw[sqtwo] (m7) to[bend right=45] (m10);
\draw[sqtwo] (m7) to[bend right=40] (m11);
\draw[sqfour] (m2) to[bend right=20] (m6);
\draw[sqfour] (m4) to[bend left=25] (m10);
\draw[sqfour] (m5) to[bend right=30] (m12);
\node[font=\small] at (2.6,-1.3) {$K_1$};
\end{scope}
\begin{scope}[xshift=5.0cm]
\foreach \j in {10,...,18} \node[font=\tiny, anchor=east, gray] at (-0.6,\j-10) {$\j$};
\node[cell] (m13) at (1,0) {}; \node[clab, left] at (m13) {$m_{13}^*$};
\node[cell] (m15) at (1,1) {}; \node[clab, left] at (m15) {$m_{15}^*$};
\node[cell] (m18) at (1,3) {}; \node[clab, left] at (m18) {$m_{18}^*$};
\node[cell] (m20) at (1,4) {}; \node[clab, left] at (m20) {$m_{20}^*$};
\node[cell] (m17) at (2,2) {}; \node[clab, right] at (m17) {$m_{17}^*$};
\node[cell] (m14) at (4,1) {}; \node[clab, right] at (m14) {$m_{14}^*$};
\node[cell] (m16) at (4,2) {}; \node[clab, right] at (m16) {$m_{16}^*$};
\node[cell] (m19) at (3,4) {}; \node[clab, right] at (m19) {$m_{19}^*$};
\node[cell] (m21) at (5.2,5) {}; \node[clab, right] at (m21) {$m_{21}^*$};
\node[cell] (m22) at (3,6) {}; \node[clab, left] at (m22) {$m_{22}^*$};
\node[cell] (m23) at (4.2,7) {}; \node[clab, right] at (m23) {$m_{23}^*$};
\node[cell] (m24) at (4.2,8) {}; \node[clab, above] at (m24) {$m_{24}^*$};
\draw[sqone] (m13)--(m15);
\draw[sqone] (m14)--(m16);
\draw[sqone] (m18)--(m20);
\draw[sqone] (m23)--(m24);
\draw[sqtwo] (m13) to[bend right=30] (m17);
\draw[sqtwo] (m15) to[bend left=50] (m18);
\draw[sqtwo] (m16) to[bend right=35] (m19);
\draw[sqtwo] (m17) to[bend right=30] (m20);
\draw[sqtwo] (m21) to[bend right=40] (m23);
\draw[sqtwo] (m22) to[bend left=30] (m24);
\draw[sqfour] (m13) to[bend right=10] (m19);
\draw[sqfour] (m13) to[bend left=75] (m20);
\draw[sqfour] (m18) to[bend left=15] (m23);
\draw[sqfour] (m19) to[bend right=15] (m24);
\draw[sqfour] (m20) to[bend left=30] (m24);
\node[font=\small] at (2.6,-1.3) {$K_2$};
\end{scope}
\end{tikzpicture}
\caption{Cell diagrams of the $\mathcal{A}(2)$-modules $H^*(K_1)$ and $H^*(K_2)$, read off from Table \ref{tab:Mstar}; the two complexes are the layers of the presentation of $\tmf_0(5)$ in Theorem \ref{thm:introB}.}
\label{fig:cells}
\end{figure}

The existence of these finite complexes $K_i$ is not obvious. There are multiple $\mathcal{A}$-module structures extending the given $\mathcal{A}(2)$-module structures, and none of them can be shown to underlie a finite spectrum directly from Toda obstruction theory. In addition, a cell-by-cell construction would involve delicate computations with Toda brackets. We show they exist by the existence of the finite $\tmf$-module $\tmf_0(5)$ and the fact that the unit map $\mathbb S\to\tmf$ is $7$-connected.

Finally, we run the Adams spectral sequence of $\tmf_0(5)$. Its $E_2$-page is the direct sum of $\Ext_{\mathcal{A}(2)_*}(\F_2,\F_2)$, the $E_2$-page for $\tmf$, and a second summand which is a module over it.

\begin{MainTheorems}\label{thm:introC}
The Adams spectral sequence of $\tmf_0(5)$ is determined through the $32$-stem: its $E_2$-page is the direct sum of $\Ext_{\mathcal{A}(2)_*}(\F_2,\F_2)$ and a module over it with $25$ generators in stems at most $34$, its differentials are listed in Theorem \ref{thm:adamsdiff}, and it collapses at $E_5$. The $E_\infty$-page is shown in Figures \ref{fig:adamsEinfa} and \ref{fig:adamsEinfb}, together with hidden $2$, $\eta$, and $\nu$ extensions.
\end{MainTheorems}

The homotopy groups of $\mathrm{Tmf}_0(5)$ and $\tmf_0(5)=\tau_{\ge0}\mathrm{Tmf}_0(5)$ are well known, from the slice spectral sequence of $\bpcm{1}$ \cite{hhrc4}\cite{hswx}. The homotopy groups of $\TMF_0(5)$ were studied by Behrens--Ormsby \cite{BO} via the descent spectral sequence, and the descent spectral sequence for $\mathrm{Tmf}_0(5)$ is also approachable via the techniques of the author, Davies, and van Nigtevecht in \cite{smf}. 

The goal of Theorem \ref{thm:introC} is not to recompute them. It is the Adams spectral sequence itself, which the slice spectral sequence does not see: how the elements of $\pi_*\tmf_0(5)$ are detected by $\Ext$ over the Steenrod algebra, how the image of $\pi_*\tmf$, and through it the image of $\pi_*\mathbb{S}$, is distributed over the two summands of the $E_2$-page. Since $EO_2(C_4)$ is $32$-periodic, the range of Theorem \ref{thm:introC} displays the basic pattern.

\subsection{Hill--Hopkins--Ravenel connective models}
The relation of a connective model to its periodic theory is like that of an integral model of a variety over a $p$-adic field to its generic fiber: the connective model is an integral lift of the periodic theory, with a nontrivial special fiber, its mod $p$ homology, which is where the arithmetic information lives, and the periodic theory is recovered from it by localization, as the generic fiber is recovered from the integral model by inverting $p$.

There are two known successful approaches to constructing connective models of the theories $EO_h(G)$. The first is through spectral algebraic geometry: $\tmf$ is the connective cover of the global sections of a sheaf of $\mathbb{E}_\infty$-rings on the compactified moduli stack of elliptic curves, and the same construction with level structures produces the connective models $\tmf_1(n)$ and $\tmf_0(n)$ of $\TMF_1(n)$ and $\TMF_0(n)$ \cite{HL}\cite{Mei23}. These theories are integral away from $n$ and come with a great deal of structure and with the control of algebraic geometry, but the construction has only seen success at heights at most $2$.

The second approach is that of Hill--Hopkins--Ravenel through Real bordism. Their theories $\bpgm{m}$ are genuine $G$-spectra for $G=C_{2^k}$, built from the norm $N_{C_2}^{G}BP_\R$ of the Real Brown--Peterson spectrum by killing $G$-orbits of generators \cite{HHR}, and by work of Beaudry--Hill--Shi--Zeng the fixed points $\bpgm{m}^{H}$ are connective models of $E_h^{hH}$ at height $h=2^{k-1}m$ \cite{BHSZ}. This construction works at every height, and it makes the full apparatus of genuine equivariant homotopy theory available: slice filtrations, norms, transfers and $RO(G)$-graded homotopy. In previous work with Hill, the author showed that these theories share the finiteness properties of $\mathrm{ko}$ and $\tmf$.

\begin{theorem*}[\cite{carrickhill}]\label{thm:carrickhill}
Let $G=C_{2^k}$. For every subgroup $H\subseteq G$, the fixed points $\bpgm{m}^H$ are fp spectra of type $2^{k-1}m$.
\end{theorem*}

Following \cite{carrickhill}, we write $eo_h(H)=\bpgm{m}^H$ for these connective higher real $K$-theories, where $h=2^{k-1}m$. The above theorem says that the Adams spectral sequences of the $eo_h(H)$ are tractable: their $E_2$-pages are $\Ext$ groups of finitely presented modules over the Steenrod algebra. The present paper follows up on \cite{carrickhill} by computing such an Adams spectral sequence with $h=2$ and $H=G=C_4$. This Adams spectral sequence lies under that of $\tmf$, so it cannot provide more information about the Adams spectral sequence of the sphere than that of $\tmf$, but this begins the study of the Adams spectral sequences of the above $eo_h$ theories in earnest by developing strategies that work for groups larger than $C_2$.

The first step in such a computation is the $E_2$-page, and hence the mod $2$ homology of $eo_h(H)$ as a comodule over the dual Steenrod algebra. For $G=C_2$ this was studied by the author with Hill and Ravenel by means of the slice filtration \cite{CHR}. At $C_2$, however, one also has the $C_2$-equivariant Adams spectral sequence \cite{HK}\cite{GHIR}, which in many cases is a better behaved replacement for the Adams spectral sequence of the fixed points. For $G=C_4$ no such replacement is known to be available, and the homology of the fixed points of $\bpcm{m}$ has not been computed for any $m\ge1$. This paper treats the first case, $m=1$.

This case is also the one where the two approaches to connective models meet. By a theorem of Quinn in forthcoming work, there is a $C_4$-$\mathbb{E}_\infty$-form of $\bpcm{1}$ whose underlying spectrum is $\tmf_1(5)$, with $C_4\cong(\Z/5)^\times$ acting through the level structure \cite[Theorem B, Corollary 5.3.2]{quinn}, just as $\tmf_1(3)$ is a $C_2$-$\mathbb{E}_\infty$-form of $\bprn{2}$ by Hill--Meier \cite{HM17}. We work with Quinn's form throughout, so that
\[eo_2(C_4)=\tmf_0(5)=\tmf_1(5)^{C_4}=\bpcm{1}^{C_4}\]
is a connective $\mathbb{E}_\infty$-model for $\mathrm{TMF}_0(5)$ and  $EO_2(C_4)$.

We have told the story relative to $\tmf$ for convenience. However: the same spectral sequence, with $\mathcal{A}_*$ in place of $\mathcal{A}(2)_*$, computes the mod $2$ homology of the $C_4$-fixed points of any form of $\bpcm{1}$ as an $\mathcal{A}_*$-comodule algebra, and our arguments apply to it with mostly notational changes.

\subsection{Methods} Apart from the results themselves, the main point of this paper is a set of techniques for computing with the slice filtration in homology beyond the group $C_2$. We describe them here, and refer to the body of the paper for the precise statements.

\subsubsection*{The $E_2$-page} The homological slice spectral sequence of \cite{CHR} is the spectral sequence obtained by tensoring the slice filtration of a $G$-spectrum $X$ with the inflation of $\F_2$. It converges to the $RO(G)$-graded homotopy of $\inf(\F_2)\otimes X$, whose integer-graded part is the mod $2$ homology of the fixed points $X^G$, and it is a spectral sequence of comodule algebras over the dual Steenrod algebra and a spectral sequence of Mackey functors. For $G=C_2$, as in \cite{CHR}, the $E_2$-page is small, the differentials are determined by the images of the classes $\bar v_1$ and $\bar v_2$, and the pages can be computed by algebraic Bockstein spectral sequences.

For $G=C_4$, the slice associated graded of $\tmf_1(5)$ splits into norm summands, generated by the norm $N=N(\bar t_1)$ in degree $\rho$, and transfer summands, which are induced from $C_2$ \cite[Corollary C]{quinn}. The $E_2$-page splits accordingly as $\mathcal{N}\oplus\mathcal{T}$. The transfer part $\mathcal{T}$ is determined by the $C_2$-level. The norm part requires the $RO(C_4)$-graded homotopy of $\inf(\F_2)\otimes_{\inf(\tmf)}\underline{\Z}$, which we compute in the positive cone in terms of the classes $a_\sigma$, $a_\lambda$, $u_\lambda$ and two further classes $x_\sigma$ and $x_\lambda$ in degrees $1-\sigma$ and $1-\lambda$ (Proposition \ref{prop:positiveconehomologyC4}). In integer degrees this gives
\[\mathcal{N}=\mathcal{A}\langle0\rangle[\eta',\nu,e,d]/(\nu^2,\eta'\nu,\eta'd+\nu e),\]
where
\[\eta'=Na_\sigma a_\lambda,\qquad\nu=Na_\sigma u_\lambda,\qquad e=Nx_\sigma a_\lambda,\qquad d=Nx_\sigma u_\lambda\]
lie in bidegrees $(1,3)$, $(3,1)$, $(2,2)$ and $(4,0)$ (Proposition \ref{prop:normpartE2}). The class $\nu$ detects the spherical element of the same name, and $e$ is a fourth root of the class of $\epsilon\in\pi_8\mathbb{S}$. Since $\nu$ and $\epsilon$ map to zero in homology, both $\nu$ and $e^4$ must be killed by differentials, and much of the computation is organized around how this happens.

\subsubsection*{Transfer strata} The spectral sequence has three levels, connected by restriction and transfer, and we never compute a page at the $C_4$-level directly. The image of the transfer on $E_2$ is $\mathcal{T}\oplus\nu\mathcal{A}\langle0\rangle[e,w]$, where $w=d+\zeta_1^2e+\zeta_2\eta'$ (Proposition \ref{prop:imtr}), and we define the \emph{transfer stratum} $\mathcal{T}_r\subseteq E_r$ to be the set of classes admitting a representative in the image of the transfer on $E_2$. By Frobenius reciprocity this is an ideal, and each page sits in a short exact sequence
\[0\to\mathcal{T}_r\to E_r\to E_r/\mathcal{T}_r\to0.\]
We compute the two ends separately. The stratum $\mathcal{T}_r$ is controlled by the $C_2$-level, where everything is known, and the quotient $E_r/\mathcal{T}_r$ is small. For example, on the $E_4$-page the quotient is
\[E_4/\mathcal{T}_4=\mathcal{A}\langle1\rangle[e,w]\oplus\eta'\mathcal{A}\langle0\rangle[\eta',e]\]
(Theorem \ref{thm:E4C4}), on the $E_6$-page it is $\F_2[e,w]\{1,g\}$ together with the eight classes $\eta'c$ for $c\in\mathcal{A}\langle1\rangle$, where $g=w\zeta_1^4+e\zeta_2^2$ (Theorem \ref{thm:E6C4}), and it does not change from $E_6$ to $E_8$ (Theorem \ref{thm:E8C4}).

\subsubsection*{The differentials} There are nonzero differentials $d_3$, $d_5$, $d_7$ and $d_9$, and each is of a different nature.

The $d_3$-differentials are determined by those at the $C_2$-level. Indeed, the restriction is injective on $\mathcal{T}$, which forces $d_3(\zeta_1^2)=\mathrm{tr}(x)$ and $d_3(\zeta_2)=\mathrm{tr}(y)$ (Proposition \ref{prop:d3atC4}). At both levels $d_3$ is then in a certain sense a Koszul differential, and the restriction and transfer are maps of Koszul complexes (Proposition \ref{prop:E2koszul}). The $E_4$-page is then assembled from Koszul homology by the long exact sequences associated to restriction and transfer.

The differential $d_5$ vanishes at the $C_2$-level, and hence on the image of the transfer. On the quotient $E_4/\mathcal{T}_4$ it is given as follows.

\begin{theorem}[Theorem \ref{thm:d5}]
On the summand $\mathcal{A}\langle1\rangle[e,w]$ of the $E_4$-page of the homological slice spectral sequence of $\tmf_1(5)$ at the $C_4$-level, the differential $d_5$ is $\F_2[e,w]$-linear and is given by
\[d_5(c)=\nu\,\partial c,\qquad c\in\mathcal{A}\langle1\rangle=E(\zeta_1^4,\zeta_2^2,\zeta_3),\]
where $\nu$ is the Hurewicz image of the Hopf map and $\partial$ is the derivation given by the cap product with $\Sq^4$.
\end{theorem}

\begin{remark}
This family follows from the fact that $\nu$ must die, which forces $d_5(\zeta_1^4)=\nu$, and from the comodule structure. The $C_4$-level also has differentials which cannot be seen on the integer-graded page at all, and these are new compared with the $C_2$-level. To determine them we use classes in the $RO(C_4)$-graded degrees $n-\sigma$, together with the norm maps from the $C_2$-level. For example, we have the differentials
\begin{align*}
d_5(a_\sigma\zeta_1^2)&=a_\sigma\eta',\\
d_5(a_\sigma\zeta_2)&=a_\sigma e,\\
d_5(a_\sigma\zeta_1^2\zeta_2)&=a_\sigma(\eta'\zeta_2+e\zeta_1^2)+x_\sigma\nu.
\end{align*}
The first two follow from the differential $d_5(u_{2\sigma})=Na_\sigma^3a_\lambda$ of Hill--Hopkins--Ravenel \cite[Theorem 11.13]{hhrc4} (Proposition \ref{prop:ROinput}).
\end{remark}

After multiplying by norm classes, these differentials determine $d_5$ on the summand $\eta'\mathcal{A}\langle0\rangle[\eta',e]$: modulo $\mathcal{T}_4$, the classes $\zeta_1^2$ and $\zeta_2$ act there as $\eta'$ and $e$, so that this summand is a Koszul complex for the regular sequence $\eta',e$, and only the eight classes $\eta'c$ survive. On the summand $\mathcal{A}\langle1\rangle[e,w]$, the family of the theorem leaves the classes of $\F_2[e,w]\{1,g\}$.

The image of the differential $d_7$ lies in the transfer stratum, and we show that it vanishes on the quotient $E_6/\mathcal{T}_6$ (Lemma \ref{lem:imd7}, Proposition \ref{prop:d7bar}). It is therefore the transfer of the $d_7$ at the $C_2$-level, and $\mathcal{T}_8$ is the homology of $\mathrm{tr}(E_6^{C_2})$ with respect to it. 

Finally, $e$ and $w$ are permanent cycles, so the pages from $E_6$ on are modules over $\F_2[e,w]$ and the differentials from $d_7$ on are linear over it (Lemma \ref{lem:ew}). The class $e^4$ must die, and the only possibility is the differential
\[d_9(\nu\zeta_1^4)=e^3,\]
which kills the tower $e^k$ for $k\ge3$. All other $d_9$-differentials are determined by this one together with $\F_2[e,w]$-linearity, naturality and the comodule structure (Theorem \ref{thm:bit}).

\subsubsection*{Collapse from a finiteness argument} The pages of our spectral sequence are infinite, while the target is finite-dimensional, so almost everything has to die, and one has to prove this for infinitely many stems at once. After the $d_9$-differential above we have $e^3=0$, and we conclude with a Hilbert series argument.

\begin{theorem}[Theorem \ref{thm:vanish}, Corollary \ref{cor:collapse}]
The homological slice spectral sequence of $\tmf_1(5)$ at the $C_4$-level collapses at $E_{10}$, and $E_{10}=E_\infty$ is $24$-dimensional, concentrated in stems at most $18$.
\end{theorem}

\begin{remark}
The argument uses only that the pages are finitely generated modules over a polynomial ring and that the differentials are linear. The Hilbert series of such a module gives a linear recursion for its dimensions, and combined with a periodicity in the quotient by $e$ it shows that $e$ acts injectively in high stems. Since $e^3=0$ after the $d_9$, the vanishing of $E_{10}$ in high stems is then a consequence of the acyclicity of $E_8/e$, which is a computation in finitely many stems.
\end{remark}

\subsubsection*{Extension problems} The $E_\infty$-page is only the associated graded of $H^{\tmf}_*(\tmf_0(5))$, and the comodule structure and the products of the latter have additional terms of higher filtration, most of which are not directly visible in integer degrees. Solving these extension problems takes up the longest section of the paper, Section \ref{sec:tmf05}, and three ideas go into it.

The first is to use the Mackey functor structure of integral homology. The $C_4$-spectrum $\inf(H\Z)\otimes_{\inf(\tmf)}\tmf_1(5)$ is a module over the constant Mackey functor $\underline{\Z}$, and this alone bounds the torsion in the integral homology of the fixed points, and hence the number of $\Sq^1$'s in the mod $2$ homology.

The second is to use classes in $RO(C_4)$-graded degrees $n+\sigma$ and $n+\sigma-\lambda$. Their coactions can be computed, and they map to the integer-graded classes under multiplication by the Euler classes $a_\sigma$ and $a_\lambda$, which decides the extensions by $\Sq^2$. The third is an integral computation in the ring of $\Gamma_1(5)$-modular forms, using Quinn's explicit description of $\tmf_1(5)$ with its $C_4$-action, which decides the multiplicative extensions.

\begin{remark}
In each case an extension which is invisible in integer degrees becomes visible after moving to a neighbouring $RO(C_4)$-graded degree, to the integral homology, or to the underlying level, where the answer is known. We expect the methods of Section \ref{sec:tmf05} to apply to $\bpcm{m}$ for $m>1$ and to larger cyclic groups.
\end{remark}

\subsubsection*{From homology to $\tmf$-modules} Theorem \ref{thm:introB} is deduced from Theorem \ref{thm:introA} by two general principles. The first uses that the unit $\mathbb{S}\to\tmf$ is $7$-connected: a connective $\tmf$-module of finite type whose $\tmf$-homology is concentrated in nine consecutive degrees is of the form $\tmf\otimes K$ for a finite spectrum $K$ (Lemma \ref{lem:narrow}). The homology of $\tmf_0(5)$ is the union of two such windows, separated by the gap in degree $9$, which gives the two-layer presentation. For the second, we observe that if $\tmf_0(5)$ were of the form $\tmf\otimes K$, then the $\mathcal{A}(2)$-module dual to $H^{\tmf}_*(\tmf_0(5))$ would extend to a module over the full Steenrod algebra, and the Adem relation $\Sq^4\Sq^6=\Sq^{10}+\Sq^8\Sq^2$ shows that it does not (Proposition \ref{prop:noAlift}). The indecomposability uses the relation $4\nu=0$ in $\pi_3\tmf_0(5)$.

\subsubsection*{The Adams spectral sequence} In Section \ref{sec:adamsss} the $E_2$-page is computed by machine from Theorem \ref{thm:introA}, and we name its generators after the cells of the complexes $K_1$ and $K_2$. The differentials are forced by the orders of the groups $\pi_n\tmf_0(5)$ and by the behaviour of a few elements of $\pi_*\tmf$ under the map to $\pi_*\TMF_0(5)$, which is known from the work of Behrens--Ormsby \cite{BO}. The hidden extensions come from comparing the $E_\infty$-page with the slice spectral sequence.

\subsubsection*{Organization} Sections \ref{sec:hsss}--\ref{sec:C2stuff} recall the spectral sequence and compute its underlying and $C_2$-levels, Sections \ref{sec:C4E2}--\ref{sec:collapse} compute the $C_4$-level, Section \ref{sec:tmf05} solves the extension problems and proves Theorem \ref{thm:introA}, Section \ref{sec:tmfmodule} proves Theorem \ref{thm:introB}, and Section \ref{sec:adamsss} proves Theorem \ref{thm:introC}.

\subsection{Acknowledgments} The author would like to thank Mike Hill and Doug Ravenel; this work builds on ideas developed in our joint work \cite{CHR}.
We thank Ryan Quinn for sharing a draft of \cite{quinn}. We thank Jack Davies for helpful comments on a draft. We used LLM's to produce Sage scripts to carry out straightforward linear algebra computations over $\F_2$. These scripts were verified by the author, and he is happy to share them.

\subsection{Notation and conventions}
\begin{enumerate}
\item We work $2$-locally throughout.
\item $\mathcal{A}_*$ is the dual Steenrod algebra, with Milnor generators $\xi_i$ and their conjugates $\zeta_i=\chi(\xi_i)$, and $\mathcal{A}(2)_*=\mathcal{A}_*/(\zeta_1^8,\zeta_2^4,\zeta_3^2,\zeta_4,\dots)$. We write $\mathcal{A}\langle n\rangle=\mathcal{A}(2)_*\Box_{\mathcal{A}(n)_*}\F_2$, so that
\[\mathcal{A}\langle 0\rangle=\F_2[\zeta_1^2,\zeta_2,\zeta_3]/(\zeta_1^8,\zeta_2^4,\zeta_3^2),\qquad \mathcal{A}\langle 1\rangle=E(\zeta_1^4,\zeta_2^2,\zeta_3).\]
The coaction on $\mathcal{A}\langle n\rangle$ is the coproduct $\Delta(\zeta_k)=\sum_i\zeta_i\otimes\zeta_{k-i}^{2^i}$.
\item For an $\mathcal{A}(2)_*$-comodule $N$ we write the coaction as $\psi(n)=\sum_I\xi^I\otimes\psi_I(n)$, where $I=(i_1,i_2,\dots)$ runs over exponent sequences and $\xi^I=\xi_1^{i_1}\xi_2^{i_2}\cdots$. We also write $\psi_{\xi^I}$ for $\psi_I$, and $\Sq^r_*=\psi_{\xi_1^r}$. The dual $N^*=\Hom(N,\F_2)$ is a left $\mathcal{A}(2)$-module, with $(\Sq^rf)(n)=f(\Sq^r_*n)$.
\item We write $H^e_*$, $H^{C_2}_*$ and $H^{C_4}_*$ for $H^{\tmf}_*$ of $\tmf_1(5)$, $\tmf_1(5)^{C_2}$ and $\tmf_0(5)$.
\item $\sigma$ and $\lambda$ are the sign and rotation representations of $C_4$, and $\rho$ denotes a regular representation. In Section \ref{sec:C2stuff} we write $\sigma$ for the sign representation of $C_2$; when both groups are present we write $\sigma_2$ for it, and $a_{\sigma_2}$ for its Euler class. We write $\Res$ and $\mathrm{tr}$ for restriction and transfer, and $\gamma$ for a generator of $C_4$.
\item Spectral sequences are drawn in Adams grading: a class in bidegree $(t-s,s)$ has stem $t-s$ and filtration $s$, and in the slice spectral sequence it comes from $\pi_{t-s}P^t_t$. We call $t$ its \emph{slice degree}; we refer frequently to this quantity because it is preserved by the coaction on the pages of our spectral sequence.
\end{enumerate}

\section{The slice filtration and homology}\label{sec:hsss} To calculate the homology of the fixed points of a $G$-spectrum $X$, one may use the equivalence
\[\F_2\otimes X^G\simeq(\inf(\F_2)\otimes X)^G,\]
where $\inf$ is the inflation functor $\Sp\to \Sp^G$, the unique symmetric monoidal left adjoint; see \cite[Proposition 2.12]{CHR}. The upshot is that the mod $2$ homology of $X^G$ is the integer-graded part of the $RO(G)$-graded homotopy of the $G$-spectrum $\inf(\F_2)\otimes X$, which we may filter equivariantly. The same holds relative to $\tmf$: if $X$ is a module over $\inf(\tmf)$, then $\F_2\otimes_{\tmf}X^G\simeq(\inf(\F_2)\otimes_{\inf(\tmf)}X)^G$.

Since $(-)^G$ commutes with limits and colimits, and $\inf(\F_2)$ is bounded below, tensoring the slice filtration $\mathrm{Slice}(X)$ of $X$ gives a complete filtration
\[\inf(\F_2)\otimes_{\inf(\tmf)}\mathrm{Slice}(X)\]
of $\inf(\F_2)\otimes_{\inf(\tmf)} X$, and hence a complete filtration of $\F_2\otimes_{\tmf}X^G$. When $X$ admits an $\mathbb E_\infty$-structure, so do these filtrations, and this gives the following (cf. \cite[Section 2.1]{CHR}).

\begin{proposition}\label{prop:hsss}
    For a $G$-ring spectrum $R$ under $\inf(\tmf)$, there is a $RO(G)$-graded spectral sequence
    \[E_2^{s,\bigstar}=\pi_\bigstar(\inf(\F_2)\otimes_{\inf(\tmf)} \gr\slice{R})\implies\pi_\bigstar(\inf(\F_2)\otimes_{\inf(\tmf)}R)\]
    of $\mathcal{A}(2)_*$-comodule algebras. It is compatible with restriction and transfer, and with the action of the Weyl group.
\end{proposition}

Note that by adjunction a $G$-ring spectrum under $\inf(\tmf)$ is the same thing as a $G$-ring spectrum $R$ such that $R^G$ is a ring spectrum under $\tmf$. We apply this to $R=\tmf_1(5)$, which is a $C_4$-$\mathbb{E}_\infty$-ring \cite[Corollary 5.3.2]{quinn} with fixed points the $\mathbb{E}_\infty$-$\tmf$-algebra $\tmf_0(5)$. By the theorem of Quinn recalled in the introduction, $\tmf_1(5)$ is a form of $\bpcm{1}$.

\section{The underlying computation and $\tmf_1(5)$}\label{sec:underlying}

Restricting the filtration $\inf(\F_2)\otimes_{\inf(\tmf)}\mathrm{Slice}(\tmf_1(5))$ to the trivial
subgroup gives the Atiyah--Hirzebruch spectral sequence of the underlying
spectrum, $\tmf_1(5)$.

\begin{proposition}\label{prop:underlyingE2}
The underlying spectral sequence has
\[
E_2^e=\mathcal{A}\langle 0\rangle[t_1,\gamma t_1],
\]
where $|b|=(k,-k)$ for $b\in\mathcal{A}\langle 0\rangle_k$, and $|t_1|=|\gamma t_1|=(2,0)$. The classes $t_1$ and $\gamma t_1$ are comodule primitives.
Its differentials are generated under the Leibniz rule by
\[
d_3(\zeta_2)=t_1+\gamma t_1,\qquad d_7(\zeta_3)=t_1^3 .
\]
\end{proposition}
\begin{proof}
This is the Atiyah--Hirzebruch spectral sequence
\[H_*^{\tmf}(H\Z)\otimes\pi_*\tmf_1(5)\implies H_*^\tmf(\tmf_1(5)),\]
with $H^{\tmf}_*(H\Z)=\mathcal{A}\langle 0\rangle$ and $\pi_*\tmf_1(5)=\Z_{(2)}[t_1,\gamma t_1]$ \cite[Lemma 4.3.3]{quinn}. The formal group law classes $v_1$ and $v_2$ must die, as they are zero in $\F_2$-homology, hence the differentials are forced by \cite[Theorem 1.1]{BHSZ}, which implies here that
\[v_1\equiv t_1+\gamma t_1\bmod 2,\qquad v_2\equiv t_1^3\bmod(2,v_1).\qedhere\]
\end{proof}

There is no room for comodule algebra extension problems on the $E_\infty$-page, and one has the following.

\begin{theorem}\label{thm:underlying}
The underlying spectral sequence collapses at $E_8^e$, and
\[
H^e_*=\F_2[\zeta_1^2,\zeta_2^2]/(\zeta_1^8,\zeta_2^4)\otimes\F_2[t_1]/(t_1^3)
\]
as an $\mathcal{A}(2)_*$-comodule algebra, where $t_1$ is primitive.
\end{theorem}

\begin{sseqdata}[name = Einfe, hsss, x range = {0}{16}, y range = {-12}{1}, y tick step = 2]
\class(0,0)
\class(2,-2)
\class(2,0)
\class(4,-4)
\class(4,-2)
\class(4,0)
\class(6,-6)
\class(6,-6)
\class(6,-4)
\class(6,-2)
\class(8,-8)
\class(8,-6)
\class(8,-6)
\class(8,-4)
\class(10,-10)
\class(10,-8)
\class(10,-6)
\class(10,-6)
\class(12,-12)
\class(12,-10)
\class(12,-8)
\class(14,-12)
\class(14,-10)
\class(16,-12)
\structline(0,0,1)(2,0,1)
\structline(0,0,1)(2,-2,1)
\structline(2,-2,1)(4,-2,1)
\structline(2,-2,1)(4,-4,1)
\structline(2,0,1)(4,0,1)
\structline(2,0,1)(4,-2,1)
\structline(4,-4,1)(6,-4,1)
\structline(4,-4,1)(6,-6,2)
\structline(4,-2,1)(6,-2,1)
\structline(4,-2,1)(6,-4,1)
\structline(4,0,1)(6,-2,1)
\structline(6,-6,1)(8,-6,1)
\structline(6,-6,1)(8,-8,1)
\structline(6,-6,2)(8,-6,2)
\structline(6,-4,1)(8,-4,1)
\structline(6,-4,1)(8,-6,2)
\structline(6,-2,1)(8,-4,1)
\structline(8,-8,1)(10,-8,1)
\structline(8,-8,1)(10,-10,1)
\structline(8,-6,1)(10,-6,1)
\structline(8,-6,1)(10,-8,1)
\structline(8,-6,2)(10,-6,2)
\structline(8,-4,1)(10,-6,2)
\structline(10,-10,1)(12,-10,1)
\structline(10,-10,1)(12,-12,1)
\structline(10,-8,1)(12,-8,1)
\structline(10,-8,1)(12,-10,1)
\structline(10,-6,1)(12,-8,1)
\structline(12,-12,1)(14,-12,1)
\structline(12,-10,1)(14,-10,1)
\structline(12,-10,1)(14,-12,1)
\structline(12,-8,1)(14,-10,1)
\structline(14,-12,1)(16,-12,1)
\structline(14,-10,1)(16,-12,1)
\classoptions["\scriptscriptstyle 1" above, label distance = -1pt](0,0,1)
\classoptions["\scriptscriptstyle \zeta_1^{2}" below left, label distance = -1pt](2,-2,1)
\classoptions["\scriptscriptstyle t_1" above, label distance = -1pt](2,0,1)
\classoptions["\scriptscriptstyle \zeta_1^{4}" below left, label distance = -1pt](4,-4,1)
\classoptions["\scriptscriptstyle t_1^{2}" above, label distance = -1pt](4,0,1)
\classoptions["\scriptscriptstyle \zeta_2^{2}" above left, label distance = -1pt](6,-6,1)
\classoptions["\scriptscriptstyle \zeta_1^{6}" below left, label distance = -1pt](6,-6,2)

\end{sseqdata}

\begin{figure}[htbp]
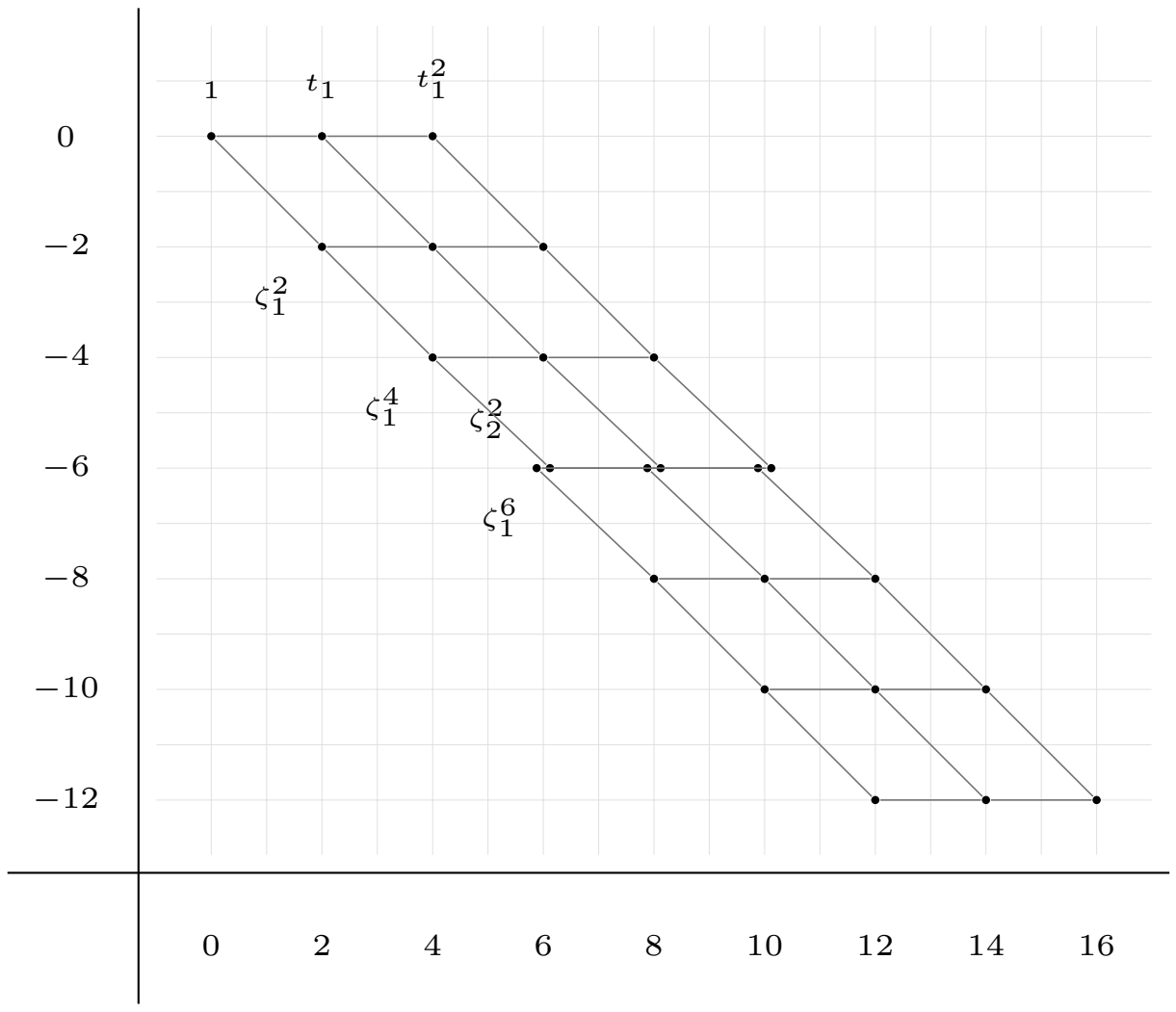

\centering
\printpage[name = Einfe, page = 2]
\caption{The $E_\infty$-page $E_8^e$ of the underlying
spectral sequence. Horizontal lines are multiplication by $t_1$ and diagonal lines
multiplication by $\zeta_1^2$. In filtration $-6$ the classes $\zeta_2^2t_1^k$ are drawn to the left of $\zeta_1^6t_1^k$.}
\label{fig:einfty-e}
\end{figure}

The underlying homology $H_*^e$ carries a $C_4$-action, which we now make explicit. Since $t_1=\gamma t_1$ on $E_\infty^e$, the generator $\gamma$ acts trivially on $E_\infty^e$, so that $\gamma(u)-u$ has strictly larger filtration than $u$ for every $u\in H^e_*$.

\begin{proposition}\label{prop:C4actionunderlying}
    The $C_4$-action on $H^e_*$ satisfies
\[\gamma(t_1)=t_1,\qquad \gamma(\zeta_1^2)=\zeta_1^2+t_1.\]
\end{proposition}
\begin{proof}
For the first claim, $t_1=\gamma t_1$ on $E_\infty$, and both classes are in the highest nonzero filtration in the $2$-stem. For the second claim, note that, by \cite[5.4.1]{HHR}, the $C_4$-module $H_2(\Res^{C_4}_e BP^{((C_4))};\Z)$ is isomorphic to $\mathrm{Ind}_{C_2}^{C_4}(\sigma)$, an induced sign representation. The map $H_2(\Res^{C_4}_e BP^{((C_4))};\Z)\to H_2(\Res^{C_4}_e \bpcm{1};\Z)$ is an isomorphism. It follows that $H_2^e\cong \mathrm{Ind}_{C_2}^{C_4}(1)$, and since $t_1$ is fixed, the class $\zeta_1^2$ cannot also be fixed. The claimed formula for $\gamma(\zeta_1^2)$ is the only possibility.
\end{proof}

\section{The $C_2$-equivariant calculation and $\tmf_1(5)^{C_2}$}\label{sec:C2stuff}
Setting $R=\tmf_1(5)$, the $C_2$-level of the spectral sequence of Proposition \ref{prop:hsss} was essentially computed in \cite[Section 5.2]{CHR}. More precisely, the spectral sequence computed there was that of $R=\bprn{2}$, which is closely related to the $C_2$-equivariant spectrum $\tmf_1(5)$. We briefly review the setup and then run the necessary algebraic Bockstein spectral sequences to calculate the homology.

\subsection{The $E_2$ and $E_4$ pages at $C_2$}

\begin{proposition}\label{prop:C2ROSS}
    The $E_2$-page of the positive cone part of the homological slice spectral sequence for $\Res^{C_4}_{C_2}\bpcm{1}$ is given by
    \[E_2=\mathcal{A}\langle 0\rangle[a_\sigma,x_1,\bar{t}_1,\gamma\bar{t}_1],\]
    where $|b|=(k,-k)$ for $b\in \mathcal{A}\langle 0\rangle_k$, $|a_\sigma|=(-\sigma,1)$, $|x_1|=(1-\sigma,0)$, $|\bar{t}_1|=|\gamma\bar{t}_1|=(\rho,0)$. The classes $a_\sigma$, $\bar{t}_1$, $\gamma\bar{t}_1$ are comodule primitives, and
    \[\psi(x_1)=1\otimes x_1+\xi_1\otimes a_\sigma.\]
    The subalgebra generated by $a_\sigma$, $x_1$, $\bar{t}_1$, $\gamma\bar{t}_1$, $x_2:=a_\sigma\zeta_2+x_1\zeta_1^2$, and $x_3:=a_\sigma\zeta_3+x_1\zeta_2^2+x_2\zeta_1^4$ consists of permanent cycles. The spectral sequence is determined up to use of the Leibniz rule by the differentials
\begin{align*}
d_3(\zeta_1^2)&=a_\sigma(\bar{t}_1+\gamma\bar{t}_1), & d_7(\zeta_1^4)&=a_\sigma^3\bar{t}_1^3,\\
d_3(\zeta_2)&=x_1(\bar{t}_1+\gamma\bar{t}_1), & d_7(\zeta_2^2)&=a_\sigma x_1^2\bar{t}_1^3,\\
&& d_7(\zeta_3)&=(x_1^3+a_\sigma^2x_2)\bar{t}_1^3.
\end{align*}
\end{proposition}
\begin{proof}
    Combine \cite[Proposition 1.4, Proposition 1.6, Theorem 1.7]{CHR}, and \cite[Theorem 1.1]{BHSZ}. Indeed, the latter states that the map $BP_\R\to \Res^{C_4}_{C_2}\bpcm{1}$ sends $\Bar{v}_1\mapsto \Bar{t}_1+\gamma\bar{t}_1\mod 2$ and $\bar{v}_2\mapsto \bar{t}_1^3\mod(2,\bar{v}_1)$.
\end{proof}

In \cite{CHR}, we had to take care to carry out the above computation in the positive cone, and we lifted the calculation to the $\R$-motivic category to do this. Here the calculation is easier, since we can present the $E_2$-page directly in integer degrees.

\begin{proposition}\label{prop:C2ZSS}
    The $E_2$-page of the integer-graded part of the homological slice spectral sequence for $\Res^{C_4}_{C_2}\bpcm{1}$ is given by
    \[E_2=\mathcal{A}\langle 0\rangle[x,\gamma x,y,\gamma y]/(x\gamma y+y\gamma x),\]
    where $x=a_\sigma\bar{t}_1$ and $y=x_1\bar{t}_1$, so that $|x|=|\gamma x|=(1,1)$ and $|y|=|\gamma y|=(2,0)$. The classes $x$ and $\gamma x$ are comodule primitives, while
    \[
        \psi(y)=1\otimes y+\xi_1\otimes x,\qquad
        \psi(\gamma y)=1\otimes \gamma y+\xi_1\otimes \gamma x.
    \]
    The subalgebra generated by $x$, $\gamma x$, $y$, $\gamma y$, $z:=x_2\bar{t}_1=y\zeta_1^2+x\zeta_2$, and $q:=x_3\bar{t}_1=x\zeta_3+y\zeta_2^2+z\zeta_1^4$  consists of permanent cycles. The spectral sequence is determined up to use of the Leibniz rule by the differentials
\begin{align*}
d_3(\zeta_1^2)&=x+\gamma x, & d_7(\zeta_1^4)&=x^3,\\
d_3(\zeta_2)&=y+\gamma y, & d_7(\zeta_2^2)&=xy^2,\\
&& d_7(\zeta_3)&=y^3+x^2z.
\end{align*}
\end{proposition}
\begin{proof}
    This follows from the previous proposition.
\end{proof}

We work in integer degrees from now on. The $E_4$-page follows directly from the $d_3$-differentials.

\begin{corollary}\label{cor:E4pageC2}
    The $E_4$-page of the integer-graded homological slice spectral sequence of $\Res^{C_4}_{C_2}\bpcm{1}$ is given by
    \[E_4=\mathcal{A}\langle 1\rangle[x,y,z]/(z^2+\zeta_1^4y^2+\zeta_2^2x^2).\]
\end{corollary}

\subsection{The $E_\infty$-page at $C_2$} The page $E_8=E_\infty$ is determined algebraically by Proposition \ref{prop:C2ZSS}, by running algebraic Bockstein spectral sequences as in \cite[Section 5]{CHR}.

\begin{theorem}\label{thm:C2Einfty}
    The $E_\infty$-page of the homological slice spectral sequence of the $C_2$-spectrum $\Res^{C_4}_{C_2}\bpcm{1}$ is the following $\mathcal{A}(2)_*$-comodule algebra of dimension $28$:
    \[
E^{C_2}_\infty=\frac{\F_2[x,y,z,q]}
{\left(\begin{gathered}
x^3,\quad xy^2,\quad y^3+x^2z,\quad z^2x,\quad z^2y+x^2q,\quad z^3+y^2q,\\
 z^4,\quad z^2q,\quad q^2
\end{gathered}\right)}
\]
where $|x|=(1,1)$, $|y|=(2,0)$, $|z|=(4,-2)$, $|q|=(8,-6)$, with coactions
\begin{align*}
    \psi(x)&=1\otimes x,\\
    \psi(y)&=1\otimes y+\xi_1\otimes x,\\
    \psi(z)&=1\otimes z+\xi_1^2\otimes y+\xi_2\otimes x,\\
    \psi(q)&=1\otimes q+\xi_1^4\otimes z+\xi_2^2\otimes y+\xi_3\otimes x.
\end{align*}
\end{theorem}
\begin{proof}
The coactions follow from those of $a_\sigma$, $x_1$ and the elements of $\mathcal{A}\langle 0\rangle$. On $E_7=E_4$, the Leibniz rule gives
\begin{align*}
d_7(\zeta_1^4\zeta_2^2)&=z^2x, & d_7(\zeta_1^4\zeta_3)&=z^2y+x^2q, & d_7(\zeta_2^2\zeta_3)&=z^3+y^2q
\end{align*}
together with $d_7(\zeta_1^4\zeta_2^2\zeta_3)=z^2q$, so the claimed relations hold in $E_\infty$; note that $z^4$ and $q^2$ vanish already on $E_4$, as $\zeta_1^8=\zeta_2^4=\zeta_3^2=0$. That $x,y,z,q$ generate and that the relations are complete may be checked with the $x$-adic filtration of the differential graded algebra $(E_7,d_7)$. Its associated graded is
\[E(\zeta_1^4,\zeta_2^2,\zeta_3)[x,y,z]/(z^2+\zeta_1^4y^2)\]
with $d_0(\zeta_3)=y^3$, and since $y^3$ is a nonzerodivisor one gets
\[E_1=E(\zeta_1^4,\zeta_2^2)[x,y,z]/(y^3,z^2+\zeta_1^4y^2).\]
Next $d_1(\zeta_1^4)=0$ and $d_1(\zeta_2^2)=xy^2$, and the annihilator of $xy^2$ in $E_1$ is $(y)$, so by \cite[Lemma 9.5]{carrickguilloupetersen}
\[
E_2
=E(\zeta_1^4)[x,y,z,q]/
\bigl(y^3,\;z^2+\zeta_1^4y^2,\;xy^2,\;
q^2,\;y^2q+\zeta_1^4y^2z\bigr).
\]
Finally $d_3(\zeta_1^4)=x^3$, the same analysis computes $E_4$, and the $x$-adic spectral sequence collapses there.
\end{proof}

\begin{sseqdata}[name = EinfCtwo, hsss, x range = {0}{16}, y range = {-8}{3}, y tick step = 2]
\class(0,0)
\class(1,1)
\class(2,2)
\class(2,0)
\class(3,1)
\class(4,2)
\class(4,0)
\class(4,-2)
\class(5,-1)
\class(6,0)
\class(6,-2)
\class(7,-1)
\class(8,0)
\class(8,-2)
\class(8,-4)
\class(8,-6)
\class(9,-5)
\class(10,-4)
\class(10,-6)
\class(11,-5)
\class(12,-4)
\class(12,-6)
\class(12,-8)
\class(13,-7)
\class(14,-6)
\class(14,-8)
\class(15,-7)
\class(16,-6)
\structline(0,0)(1,1)
\structline(0,0)(2,0)
\structline(1,1)(2,2)
\structline(1,1)(3,1)
\structline(2,2)(4,2)
\structline(2,0)(3,1)
\structline(2,0)(4,0)
\structline(3,1)(4,2)
\structline(4,0)(6,0)
\structline(4,-2)(5,-1)
\structline(4,-2)(6,-2)
\structline(5,-1)(7,-1)
\structline(6,0)(8,0)
\structline(6,-2)(7,-1)
\structline(6,-2)(8,-2)
\structline(8,-4)(10,-4)
\structline(8,-6)(9,-5)
\structline(8,-6)(10,-6)
\structline(9,-5)(11,-5)
\structline(10,-4)(12,-4)
\structline(10,-6)(11,-5)
\structline(12,-6)(14,-6)
\structline(12,-8)(13,-7)
\structline(12,-8)(14,-8)
\structline(13,-7)(15,-7)
\structline(14,-6)(16,-6)
\structline(14,-8)(15,-7)
\classoptions["\scriptscriptstyle 1" above, label distance = -1pt](0,0)
\classoptions["\scriptscriptstyle x" above, label distance = -1pt](1,1)
\classoptions["\scriptscriptstyle y" below, label distance = -1pt](2,0)
\classoptions["\scriptscriptstyle z" below, label distance = -1pt](4,-2)
\classoptions["\scriptscriptstyle q" below, label distance = -1pt](8,-6)
\classoptions["\scriptscriptstyle y^{2}z^{3}" above, label distance = -1pt](16,-6)

\end{sseqdata}

\begin{figure}[htbp]
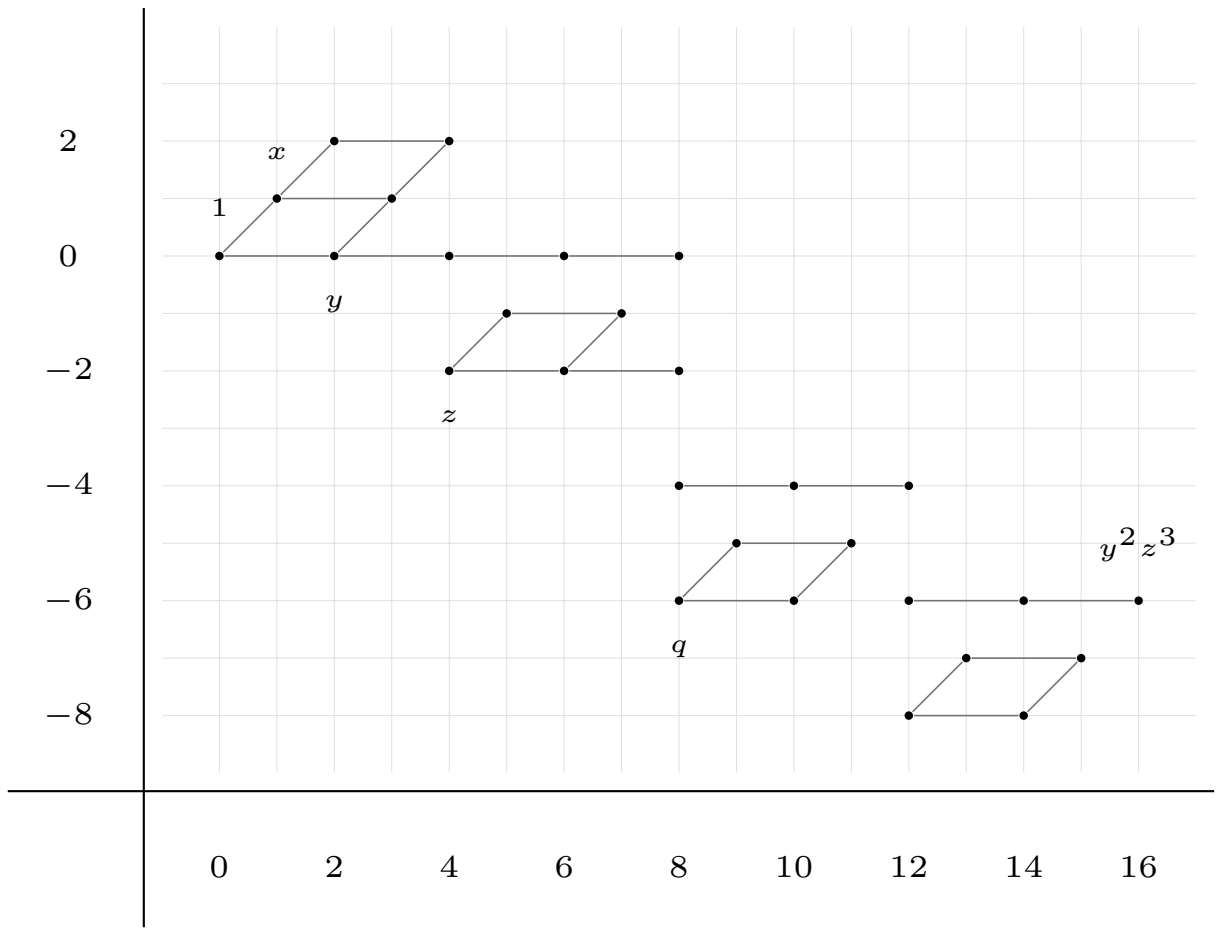

\centering
\printpage[name = EinfCtwo, page = 2]
\caption{The $E_\infty$-page of the $C_2$-equivariant
spectral sequence. Lines of slope $1$ are multiplication by $x$
and horizontal lines are multiplication by $y$.}
\label{fig:einfty-C2}
\end{figure}

\subsection{The comodule algebra $H^{C_2}_*$}\label{sec:C2ext}
We now determine $H^{C_2}_*$ itself.  The following is \cite{CHR}. 

\begin{proposition}\label{prop:C2coactions}
The classes $x_1$ and $x_2$ admit unique lifts to
$\pi_\bigstar(\inf(\F_2)\otimes BP_\R)$, and for any lift of $x_3$,
\[
\psi(x_n)=\sum_{i=0}^{n}\xi_i^{2^{n-i}}\otimes x_{n-i},\qquad x_0=a_\sigma,\qquad n\le 3.
\]
\end{proposition}

Since $a_\sigma$ and $\bar t_1$ are comodule primitives, and the map $BP_\R\to \Res^{C_4}_{C_2}\tmf_1(5)$ is one of ring spectra, the same formulas hold in $\pi_\bigstar(\inf(\F_2)\otimes_{\inf(\tmf)}\Res^{C_4}_{C_2}\tmf_1(5))$ with $\mathcal{A}(2)_*$ in place of $\mathcal{A}_*$, and the classes $x=a_\sigma\bar t_1$, $y=x_1\bar t_1$, $z=x_2\bar t_1$ and $q=x_3\bar t_1$ of $H^{C_2}_*$ have the coactions of Theorem \ref{thm:C2Einfty}. We use these lifts from now on.

\begin{proposition}\label{prop:C2mult}
In $H^{C_2}_*$ the relations $x^3$, $xy^2$, $y^3+x^2z$, $z^2x$, $z^2y+x^2q$ and
$z^3+y^2q$ of Theorem \ref{thm:C2Einfty} hold exactly, and
\[
q^2=c_1\,y^2z^3,\qquad z^2q=c_2\,y^2z^3,\qquad z^4=c_3\,y^2z^3
\]
for some $c_1,c_2,c_3\in\F_2$.
\end{proposition}
\begin{proof}
The relations $x^3$, $xy^2$,
$y^3+x^2z$, $z^2x$ and $z^2y+x^2q$ lie in bidegrees $(3,3)$, $(5,1)$, $(6,0)$,
$(9,-3)$ and $(10,-4)$, and by Figure \ref{fig:einfty-C2} there are no such classes above these bidegrees. The relation $z^3+y^2q$ lies in bidegree $(12,-6)$, so the only possibility
is $z^3=y^2q+\epsilon\,y^2z^2$ with $\epsilon\in\F_2$; and $q^2$, $z^2q$, $z^4$ lie in
bidegrees $(16,-12)$, $(16,-10)$, $(16,-8)$, so each is a multiple of the class $y^2z^3$ in $(16,-6)$.

It remains to prove $\epsilon=0$. We apply $\psi$ and compare coefficients of $\xi_1^4$. The coefficient of $\xi_1^4$ in $\psi(z)^3+\psi(y)^2\psi(q)$ is
$y^2z+y^2z=0$, whereas that in $\psi(y)^2\psi(z)^2$ is $y^4\neq0$. So
$\epsilon=0$.
\end{proof}

\begin{theorem}\label{thm:C2comodule}
The coaction of $H^{C_2}_*$ on the basis of monomials in $x,y,z,q$ is that of $E^{C_2}_\infty$.
\end{theorem}
\begin{proof}
By Proposition \ref{prop:C2mult} the products of
$H^{C_2}_*$ agree with those of $E^{C_2}_\infty$ except possibly in degree $16$. For any monomial $m$, the non primitive terms in $\psi(m)$ are of the form $\xi^I\otimes c$ where $|c|<16$. The result then follows from the fact that the spectral sequence is one of comodule algebras.
\end{proof}

\begin{proposition}\label{prop:C2res}
The restriction $\Res\colon H^{C_2}_*\to H^e_*$ satisfies
\[
\Res(x)=0,\qquad \Res(y)=t_1,\qquad \Res(z)\equiv\zeta_1^2t_1\mod t_1^2,
\]
and the transfer $H^e_{16}\to H^{C_2}_{16}$ is an isomorphism, $\mathrm{tr}(\zeta_1^6\zeta_2^2t_1^2)=y^2z^3$.
\end{proposition}
\begin{proof}
On the $E_\infty$-page the restriction is determined by $\Res(a_\sigma)=0$ and $\Res(x_1)=1$. In degree $2$ there are no classes of larger filtration than $t_1$, and in degree $4$ the only class of larger filtration than $\zeta_1^2t_1$ is $t_1^2$. For the transfer, the cofiber sequence $C_{2+}\to S^0\to S^\sigma$ shows that the image of the transfer is the kernel of $a_\sigma$, which contains $y^2z^3$ for degree reasons. The groups $H^e_{16}$ and $H^{C_2}_{16}$ are spanned by $\zeta_1^6\zeta_2^2t_1^2$ and $y^2z^3$.
\end{proof}

On $E^{C_2}_4$ we have $x+\gamma x=d_3(\zeta_1^2)$ and $y+\gamma y=d_3(\zeta_2)$ by
Proposition \ref{prop:C2ZSS}, so the residual action of $C_4/C_2$ is the identity on $E^{C_2}_\infty$. It is
not the identity on $H^{C_2}_*$, as the underlying action of Proposition \ref{prop:C4actionunderlying} shows.

\begin{proposition}\label{prop:C2gamma}
There is an $\varepsilon\in\F_2$ with
\[
\gamma(x)=x,\qquad
\gamma(y)=y+x^2,\qquad
\gamma(z)=z+y^2+x^2y,\qquad
\gamma(q)=q+z^2+y^2z+\varepsilon\,y^4.
\]
Moreover $z^4=0$ in $H^{C_2}_*$, that is, $c_3=0$ in Proposition \ref{prop:C2mult}.
\end{proposition}
\begin{proof}
The map $\gamma+1$ strictly raises filtration, so by Figure \ref{fig:einfty-C2} we may write $\gamma(x)=x$, $\gamma(y)=y+\epsilon_1\,x^2$,
$\gamma(z)=z+\epsilon_2\,y^2+\epsilon_3\,x^2y$ and
$\gamma(q)=q+\epsilon_4\,z^2+\epsilon_5\,y^2z+\varepsilon\,y^4$.
Since $\gamma$ is a comodule map, applying $1\otimes\gamma$ to
$\psi(z)=1\otimes z+\xi_1^2\otimes y+\xi_2\otimes x$ and comparing with $\psi(\gamma z)$ gives $\epsilon_2=\epsilon_1$, using
$\psi(y^2)=1\otimes y^2+\xi_1^2\otimes x^2$ and $\psi(x^2y)=1\otimes x^2y$. In the same way, the $\xi_1^4$-coefficient of $\psi(\gamma q)$ gives $\epsilon_4=\epsilon_2$ and $\epsilon_5=\epsilon_3$, using $\psi(z^2)=1\otimes z^2+\xi_1^4\otimes y^2+\xi_2^2\otimes x^2$, $\psi(y^2z)=1\otimes y^2z+\xi_1^4\otimes x^2y$ and $\psi(y^4)=1\otimes y^4$.

Restriction is equivariant, and $\gamma$ fixes $t_1^2$, so by Propositions \ref{prop:C4actionunderlying} and \ref{prop:C2res}
\[
\epsilon_2\,t_1^2=\Res\bigl(\gamma(z)+z\bigr)
=\gamma(\zeta_1^2t_1)+\zeta_1^2t_1=t_1^2.
\]
So $\epsilon_1=\epsilon_2=\epsilon_4=1$.

The class $y^2z^3$ spans $H^{C_2}_{16}$, so $\gamma$ fixes it, and in particular $\gamma$ fixes $q^2$ and $z^2q$. On the other hand,
\[
\gamma(q)^2=(q+z^2+\epsilon_5\,y^2z+\varepsilon y^4)^2=q^2+z^4+\epsilon_5\,y^4z^2+\varepsilon\,y^8=q^2+z^4,
\]
since $y^4z^2$ and $y^8$ lie in degree $16$ and filtrations $-4$ and $0$, where $H^{C_2}_*$ vanishes. As $\gamma$ is a ring map, $z^4=0$. Using this, $y^4z^2=y^6z=y^8=0$, and $y^4q=y^2z^3$, one finds
\[
\gamma(z^2q)=(z^2+y^4)(q+z^2+\epsilon_5\,y^2z+\varepsilon y^4)=z^2q+(1+\epsilon_5)\,y^2z^3,
\]
so $\epsilon_5=1$, and therefore $\epsilon_3=1$.
\end{proof}

\begin{theorem}\label{thm:C2abutment}
For some $c_1,c_2\in\F_2$ there is an isomorphism of algebras
\[
H^{C_2}_*=\frac{\F_2[x,y,z,q]}
{\left(\begin{gathered}
x^3,\quad xy^2,\quad y^3+x^2z,\quad z^2x,\quad z^2y+x^2q,\quad z^3+y^2q,\quad z^4,\\
q^2+c_1\,y^2z^3,\quad z^2q+c_2\,y^2z^3
\end{gathered}\right)},
\]
with the coaction of Theorem \ref{thm:C2Einfty} and the
action of $C_4/C_2$ of Proposition \ref{prop:C2gamma}. In particular $H^{C_2}_*\cong E^{C_2}_\infty$ as $\mathcal{A}(2)_*$-comodules.
\end{theorem}

\begin{remark}\label{rmk:c1c2}
The coefficient $c_1$ does not depend on the choice of the lift $q$, since $z^4=0$ and $y^4z^2=0$. The coefficients $c_2$ and $\varepsilon$ do depend on it: replacing $q$ by $q+y^2z$ changes both $c_2$ and $\varepsilon$ by $1$, as one checks using the relations of Theorem \ref{thm:C2abutment}. So we may choose $q$ with $c_2=0$, and then $\varepsilon=1$ by Proposition \ref{prop:abeps} below. Only the sum $c_2+\varepsilon$ enters the calculation of $H^{C_4}_*$, and we do not determine $c_1$ here.
\end{remark}

\section{The $C_4$-equivariant $E_2$- and $E_4$-pages}\label{sec:C4E2}

\subsection{The $E_2$-page at $C_4$} Computing the $E_2$-page of the homological slice spectral sequence of Proposition \ref{prop:hsss} for $\tmf_1(5)$ splits into two pieces, as the slice associated graded of $\tmf_1(5)$ splits into two types of summands, norm and transfer summands \cite[Corollary C]{quinn}. The contribution of transfer summands to our $E_2$-page is determined by the $C_2$-level, i.e. by Proposition \ref{prop:C2ROSS}. The norm summand portion of $\gr\slice{\tmf_1(5)}$ is the polynomial algebra $\underline{\Z}[N(\bar{t}_1)]$, hence its contribution comes from the $RO(C_4)$-graded homotopy groups of $\inf(\F_2)\otimes_{\inf(\tmf)}\underline{\Z}$.

As our applications are to understanding the integer degree portion of homology, we only determine the ``positive cone'' portion of the latter, i.e. those $RO(C_4)$-degrees $a+b\sigma+c\lambda$ where $b,c\le 0$. In the following, recall that at the $C_2$-level, we have
\[\pi_\bigstar^{C_2}(\inf(\F_2)\otimes_{\inf(\tmf)}\underline{\Z})=\mathcal{A}\langle 0\rangle[a_{\sigma_2},x_1]\]
in the positive cone, as in Proposition \ref{prop:C2ROSS}.

\begin{proposition}\label{prop:positiveconehomologyC4}
    The positive cone portion of
    \[\pi_\bigstar^{C_4}(\inf(\F_2)\otimes_{\inf(\tmf)} \underline{\Z})\]
    is isomorphic to the ring
    \[\frac{\mathcal{A}\langle 0\rangle[a_\sigma,a_\lambda,u_\lambda,x_\sigma,x_\lambda]}{(a_\sigma^2u_\lambda,a_\sigma x_\lambda,x_\lambda^2,x_\sigma x_\lambda+a_\sigma u_\lambda)}\]
    with degrees
    \[
        |a_\sigma|=-\sigma,\quad
        |a_\lambda|=-\lambda,\quad
        |u_\lambda|=2-\lambda,\quad
        |x_\sigma|=1-\sigma,\quad
        |x_\lambda|=1-\lambda,
    \]
    where $a_\sigma$, $a_\lambda$, $u_\lambda$, and $x_\lambda$ are comodule primitives, and $\psi(x_\sigma)=1\otimes x_\sigma+\xi_1\otimes a_\sigma$. The restrictions of $a_\sigma$ and $x_\lambda$ are zero, and
    \[
        \Res(a_\lambda)=a_{\sigma_2}^2,\quad
        \Res(u_\lambda)=u_{2\sigma_2}=x_1^2+\zeta_1^2a_{\sigma_2}^2,\quad
        \Res(x_\sigma)=1.
    \]
    The class $x_\lambda$ is the transfer of $x_1a_{\sigma_2}$, and the residual action of $C_4/C_2$ on the positive cone of $\pi_\bigstar^{C_2}(\inf(\F_2)\otimes_{\inf(\tmf)}\underline{\Z})$ is trivial.
\end{proposition}
\begin{proof}
    The positive cone portion of $\pi_\bigstar^{C_4}\underline{\Z}$ is given by the ring
    \[\frac{\Z[a_\sigma,a_\lambda,u_{2\sigma},u_\lambda]}{(2a_\sigma,4a_\lambda,2u_{2\sigma }a_\lambda-u_\lambda a_\sigma^2)},\]
    which gives a splitting
    \begin{align*}
        \bigg(\bigoplus\limits_{b,c\ge0}\Sigma^{b\sigma+c\lambda}\underline{\Z}\bigg)^{C_4}
&\simeq
\bigoplus\limits_{\gamma,\delta\ge 0}\Sigma^{2(\gamma+\delta)}\Z\{u_{2\sigma}^\gamma u_\lambda^\delta\}
\oplus
\bigoplus\limits_{\substack{\beta>0\\ \gamma,\delta\ge0}}\Sigma^{2(\gamma+\delta)}\Z/4\{a_\lambda^\beta u_{2\sigma}^\gamma u_\lambda^\delta\}
\\
&\oplus
\bigoplus\limits_{\substack{\alpha>0,\ \beta,\gamma,\delta\ge0\\ \delta>0\,\Rightarrow\,\alpha=1}}
\Sigma^{2(\gamma+\delta)}\F_2\{a_\sigma^\alpha a_\lambda^\beta u_{2\sigma}^\gamma u_\lambda^\delta\}.
    \end{align*}

    We use this decomposition to define $x_\sigma$ as the class $[\zeta_1a_\sigma]$, i.e. the element $\zeta_1\in \mathcal{A}_1$ corresponding to the summand indexed by $a_\sigma$, and we define $x_\lambda$ as the class $[x\cdot a_\lambda]$, i.e. the nonzero element $x\in H_1(\Z/4;\F_2)$ corresponding to the summand indexed by $a_\lambda$. The claimed values of the restriction map follow from that of $\pi^{C_4}_\bigstar\underline{\Z}$, except for that of $x_\sigma$ and $x_\lambda$. The restriction of $x_\sigma$ must be nonzero, as it is not divisible by $a_\sigma$, as follows for degree reasons. The only possibility is that its restriction is $1$. The restriction of $x_\lambda$ follows from its description as a transfer, which we verify below. The residual action on $\pi_\bigstar^{C_2}(\inf(\F_2)\otimes_{\inf(\tmf)}\underline{\Z})$ must be trivial in the positive cone as $\mathcal{A}\langle 0\rangle$ is in the image of the restriction map and therefore fixed, $a_{\sigma_2}$ is fixed as it comes from the sphere, and $x_1$ is fixed as it is the unique nonzero class in degree $1-\sigma_2$.

    The claimed description now follows, additively, from the claim that
    \[u_{2\sigma}=x_\sigma^2+\zeta_1^2a_\sigma^2.\]
    To see this, note that for degree reasons
    \[x_\sigma^2=\epsilon_1u_{2\sigma}+\epsilon_2\zeta_1^2a_\sigma^2.\]
    The coaction on $x_\sigma$ forces $\epsilon_2=1$, and the restriction on $x_\sigma$ forces $\epsilon_1=1$.

    The relation $a_\sigma^2u_\lambda$ is determined by the gold relation $a_\sigma^2u_\lambda=2u_{2\sigma}a_\lambda$ in $\pi_\bigstar^{C_4}\underline{\Z}$. This implies that $a_\sigma u_\lambda$ is a transfer, and for degree reasons one sees that $a_\sigma u_\lambda$ is the transfer of $x_1a_{\sigma_2}$ in degree $2-\sigma-\lambda$. The relation $a_\sigma x_\lambda=0$ holds for degree reasons and implies that $x_\lambda$ is a transfer; for degree reasons the only possibility is that $x_\lambda$ is the transfer of $x_1a_{\sigma_2}$ in degree $1-\lambda$. The Frobenius relation now implies $x_\sigma x_\lambda=a_\sigma u_\lambda$.
\end{proof}

As at the $C_2$-level, the previous proposition does not account for all of our $RO(C_4)$-graded $E_2$-page, even within the positive cone. As in Proposition \ref{prop:C2ZSS}, we will therefore proceed to give a direct description of the $E_2$-page within integer degrees.

The slice associated graded of $\tmf_1(5)$ has the free $\mathbb E_1$-$\underline{\Z}$-algebra on a class $\bar{t}_1$ in degree $\rho$, denoted $\underline{\Z}[\overline{t}_1]$, as a summand. We refer to the homotopy groups
\[\mathcal{N}:=\pi_*^{C_4}(\inf(\F_2)\otimes_{\inf(\tmf)} \underline{\Z}[\overline{t}_1])\]
as the \emph{norm part} of our $E_2$-page. We write $N=N(\bar t_1)$ for the norm of $\bar t_1$, in degree $\rho=1+\sigma+\lambda$.

\begin{proposition}\label{prop:normpartE2}
 There is an isomorphism of $\mathcal{A}(2)_*$-comodule algebras
 \[
\mathcal{N}\cong\mathcal{A}\langle 0\rangle[\eta',\nu,e,d]/(\nu^2,\eta'\nu,\eta'd+\nu e),
 \]
 where
 \[
     \eta'=Na_\sigma a_\lambda,\qquad
     \nu=Na_\sigma u_\lambda,\qquad
     e=Nx_\sigma a_\lambda,\qquad
     d=Nx_\sigma u_\lambda,
\]
in bidegrees $(1,3)$, $(3,1)$, $(2,2)$ and $(4,0)$. The classes $\eta'$ and $\nu$ are comodule primitives, and
 \[
     \psi(e)=1\otimes e+\xi_1\otimes \eta',\qquad
     \psi(d)=1\otimes d+\xi_1\otimes \nu.
 \]
\end{proposition}
\begin{proof}
    The relations, degrees, and coactions all follow by definition from Proposition \ref{prop:positiveconehomologyC4}. In fact, the presentation of Proposition \ref{prop:positiveconehomologyC4} implies that an $\F_2$-basis is given by the list of monomials of the form
    \[bN^ma_\sigma^\alpha a_\lambda^\beta x_\sigma^\gamma u_\lambda^\delta,\]
    with $b\in\mathcal{A}\langle 0\rangle$, $\alpha+\gamma=\beta+\delta=m$ and $\delta\ge 1\implies \alpha\le 1$, noting that no monomial in $x_\lambda$ appears in integer degrees by the relations $x_\lambda a_\sigma=0$ and $x_\lambda x_\sigma=a_\sigma u_\lambda$. To see that the four generators generate, note that each above monomial with $m> 1$ is divisible by a monomial with $m=1$, and our four generators are precisely those monomials with $m=1$. The three relations enforce precisely the restriction $\delta\ge 1\implies \alpha\le 1$.
\end{proof}

The element $\nu$ detects the spherical element of the same name, and $e$ is a fourth root of the spherical element $\epsilon\in\pi_8\mathbb{S}$ (see \cite[pg. 31]{hswx}). As they are zero in homology, both $\nu$ and $e^4$ must be killed by differentials.

It remains to account for the transfer summands in the slice associated graded of $\tmf_1(5)$. This portion consists of the sum
\[\bigoplus\limits_{i>j\ge0}\mathrm{Ind}_{C_2}^{C_4}\underline{\Z}\{\bar{t}_1^i\gamma\bar{t}_1^j\},\]
and we refer to
\[\mathcal{T}:=\pi_*^{C_4}\bigg(\bigoplus\limits_{i>j\ge0}\inf(\F_2)\otimes_{\inf(\tmf)}\mathrm{Ind}_{C_2}^{C_4}\underline{\Z}\{\bar{t}_1^i\gamma\bar{t}_1^j\}\bigg)\]
as the \emph{transfer part} of $E_2$. In the following, we use the notation of Proposition \ref{prop:C2ZSS}.

\begin{proposition}\label{prop:transferpartE2C4}
    The transfer part $\mathcal{T}$ of $E_2$ is the free $\mathcal{A}\langle 0\rangle$-module with basis the monomials
    \[\mathrm{tr}(x^p\gamma x^q y^r\gamma y^s),\qquad ps=0,\qquad p+r>q+s.\]
   The $\mathcal{N}$-module structure along with the multiplication map $\mathcal{T}\otimes \mathcal{T}\to \mathcal{T}$ is determined by Frobenius reciprocity along with the restrictions
    \[
        \Res(\eta')=0,\qquad
        \Res(\nu)=0,\qquad
        \Res(e)=x\gamma x,\qquad
        \Res(d)=y\gamma y+\zeta_1^2\,x\gamma x.
    \]
\end{proposition}

\begin{proposition}\label{prop:imtr}
The image of the transfer on $E_2$ is
\[\mathrm{im}(\mathrm{tr})=\mathcal{T}\oplus\nu\mathcal{A}\langle 0\rangle[e,w],\qquad w:=d+\zeta_1^2e+\zeta_2\eta'.\]
\end{proposition}
\begin{proof}
The transfer of a monomial whose $\bar{t}_1$- and $\gamma\bar{t}_1$-degrees differ is a basis element of $\mathcal{T}$: after applying $\mathrm{tr}\gamma=\mathrm{tr}$ if need be and rewriting by $x\gamma y=y\gamma x$ to the normal form with $ps=0$, this is Proposition \ref{prop:transferpartE2C4}. A monomial whose two degrees agree is a product of the classes $x\gamma x$, $y\gamma y$ and $x\gamma y$, and $(x\gamma y)^2=(x\gamma x)(y\gamma y)$, so such monomials are precisely
\[(x\gamma x)^k(y\gamma y)^n=\Res(e^kw^n)\qquad\text{and}\qquad (x\gamma x)^k(y\gamma y)^n\,x\gamma y=\Res(e^kw^n)\cdot x\gamma y,\]
since $\Res(w)=y\gamma y$. By Frobenius reciprocity their transfers are
\[
\mathrm{tr}(\Res(e^kw^n))=e^kw^n\cdot\mathrm{tr}(1)=0,\qquad
\mathrm{tr}(\Res(e^kw^n)\,x\gamma y)=e^kw^n\cdot\mathrm{tr}(x\gamma y)=\nu\, e^kw^n,
\]
since $\mathrm{tr}(1)=0$ and $\mathrm{tr}(x\gamma y)=\nu$ by Proposition \ref{prop:positiveconehomologyC4}.
\end{proof}

\begin{lemma}\label{lemma:permanentcyclesE2}
The classes $\eta'$, $\nu$, and $e$ are permanent cycles.
\end{lemma}
\begin{proof}
    The classes $\eta'$ and $\nu$ are the image of permanent classes of the same name in the slice spectral sequence of $\tmf_1(5)$ \cite{hswx}. The classes $N$ and $a_\lambda$ are permanent cycles for the same reason, and the class $x_\sigma$ is a permanent cycle for degree reasons: the $-\sigma$ stem is $0$ on $E_2$ in filtrations $>1$. Together these imply $e$ is permanent.
\end{proof}

\subsection{$d_3$ differentials at $C_4$} The $d_3$-differential in the homological slice spectral sequence for $\tmf_1(5)$ is determined by the $d_3$-differentials on the $C_2$ level of Proposition \ref{prop:C2ZSS} using the restriction and transfer maps. We need the following observation.

\begin{lemma}\label{lem:resinj}
    The restriction map is injective on $\mathcal{T}$.
\end{lemma}
\begin{proof}
    By Proposition \ref{prop:transferpartE2C4} and the $\mathcal{A}\langle 0\rangle$-linearity of $\Res$, it suffices to show that the elements $\Res(\mathrm{tr}(m))=m+\gamma m$, for $m$ running over the monomials $x^p\gamma x^q y^r\gamma y^s$ with $ps=0$ and $p+r>q+s$, are linearly independent over $\mathcal{A}\langle 0\rangle$. Grade $E_2^{C_2}$ by the degrees in $\bar{t}_1$ and $\gamma\bar{t}_1$. Then $m$ has bidegree $(p+r,q+s)$ and $\gamma m$ has bidegree $(q+s,p+r)$, so the component of $m+\gamma m$ in bidegree $(p+r,q+s)$ is $m$, and a linear dependence among the $m+\gamma m$ would give one among the monomials $m$.
\end{proof}

\begin{proposition}\label{prop:d3atC4}
    On the norm part $\mathcal{N}$, the $d_3$-differential is determined by
    \[
        d_3(\zeta_1^2)=\mathrm{tr}(x),\qquad
        d_3(\zeta_2)=\mathrm{tr}(y),\qquad
        d_3(d)=\mathrm{tr}(x^2\gamma x),
    \]
    together with $d_3(\zeta_3)=0$, the fact that $\eta'$, $\nu$ and $e$ are $d_3$-cycles (Lemma \ref{lemma:permanentcyclesE2}), and the Leibniz rule. On the transfer part $\mathcal{T}$ it is given by naturality of the transfer.
    In particular, the class $w=d+\zeta_1^2e+\zeta_2\eta'$ is a $d_3$-cycle.
\end{proposition}
\begin{proof}
    The differentials on $\zeta_1^2$ and $\zeta_2$ follow from the corresponding differentials at the $C_2$-level: since $\zeta_1^2$ and $\zeta_2$ restrict to the classes of the same name that support $d_3$ differentials, they must support differentials whose targets restrict to the targets of the $C_2$-differentials. The claimed differentials are the only possibility, as $\Res$ is injective on $\mathcal{T}$ by Lemma \ref{lem:resinj}.
    The same analysis holds for $d$, as $\Res(d)=y\gamma y+\zeta_1^2x\gamma x$, so that
    \[\Res(d_3(d))=d_3\big(y\gamma y+\zeta_1^2x\gamma x\big)=x^2\gamma x+\gamma\big(x^2\gamma x\big)=\Res(\mathrm{tr}(x^2\gamma x)).\qedhere\]
\end{proof}

Using that $w$ is a cycle, we now give a different presentation of $E_3$ that makes computing its homology much easier.

\begin{proposition}\label{prop:E3asacomplex}
  As an $\mathcal{A}\langle 0\rangle$-module, we have
  \begin{align*}
      \mathcal{N}&=\mathcal{A}\langle 0\rangle[\eta',\nu,e,w]\big/
  \bigl(\nu^2,\ \eta'\nu,\ \eta'w+\nu e+\zeta_1^2\,\eta'e+\zeta_2\,\eta'^2\bigr)\\
      &=\mathcal{A}\langle 0\rangle[e,w]\ \oplus\ \eta'\mathcal{A}\langle 0\rangle[\eta',e]\ \oplus\ \nu\mathcal{A}\langle 0\rangle[e,w],\\
      E_2&=\mathcal{T}\oplus\mathcal{N},
  \end{align*}
  and $d_3$ is $\mathcal{A}\langle 1\rangle$-linear, vanishes on the monomials in $\eta'$, $\nu$, $e$, $w$ and on the basis of $\mathcal{T}$, and is given on the $\mathcal{A}\langle 1\rangle$-basis $\{1,\zeta_1^2,\zeta_2,\zeta_1^2\zeta_2\}$ of $\mathcal{A}\langle 0\rangle$ by $\zeta_1^2\mapsto\mathrm{tr}(x)$, $\zeta_2\mapsto\mathrm{tr}(y)$ and the Leibniz rule.
\end{proposition}
\begin{proof}
    The basis of $\mathcal{N}$ is a translation of the monomial basis of Proposition \ref{prop:normpartE2} under the invertible change of variables $w=d+\zeta_1^2e+\zeta_2\eta'$; the rest is Proposition \ref{prop:d3atC4}, $w$ being a $d_3$-cycle.
\end{proof}

\begin{sseqdata}[name = Etwo, hsss, x range = {0}{8}, y range = {-8}{24}, y tick step = 4, differentials = {dthree}]
\class(0,0)
\class(1,1)
\class(1,3)
\class(2,-2)
\class(2,0)
\class(2,2)
\class(2,2)
\class(2,6)
\class(3,-3)
\class(3,-1)
\class(3,1)
\class(3,1)
\class(3,1)
\class(3,3)
\class(3,3)
\class(3,5)
\class(3,9)
\class(4,-4)
\class(4,-2)
\class(4,-2)
\class(4,0)
\class(4,0)
\class(4,0)
\class(4,0)
\class(4,0)
\class(4,2)
\class(4,2)
\class(4,4)
\class(4,4)
\class(4,4)
\class(4,4)
\class(4,8)
\class(4,12)
\class(5,-5)
\class(5,-3)
\class(5,-3)
\class(5,-1)
\class(5,-1)
\class(5,-1)
\class(5,-1)
\class(5,-1)
\class(5,1)
\class(5,1)
\class(5,1)
\class(5,1)
\class(5,3)
\class(5,3)
\class(5,3)
\class(5,3)
\class(5,3)
\class(5,5)
\class(5,5)
\class(5,5)
\class(5,7)
\class(5,7)
\class(5,11)
\class(5,15)
\class(6,-6)
\class(6,-6)
\class(6,-4)
\class(6,-4)
\class(6,-2)
\class(6,-2)
\class(6,-2)
\class(6,-2)
\class(6,-2)
\class(6,-2)
\class(6,-2)
\class(6,0)
\class(6,0)
\class(6,0)
\class(6,0)
\class(6,0)
\class(6,0)
\class(6,2)
\class(6,2)
\class(6,2)
\class(6,2)
\class(6,2)
\class(6,2)
\class(6,2)
\class(6,2)
\class(6,4)
\class(6,4)
\class(6,4)
\class(6,6)
\class(6,6)
\class(6,6)
\class(6,6)
\class(6,6)
\class(6,6)
\class(6,10)
\class(6,10)
\class(6,14)
\class(6,18)
\class(7,-7)
\class(7,-7)
\class(7,-5)
\class(7,-5)
\class(7,-5)
\class(7,-3)
\class(7,-3)
\class(7,-3)
\class(7,-3)
\class(7,-3)
\class(7,-3)
\class(7,-3)
\class(7,-3)
\class(7,-1)
\class(7,-1)
\class(7,-1)
\class(7,-1)
\class(7,-1)
\class(7,-1)
\class(7,1)
\class(7,1)
\class(7,1)
\class(7,1)
\class(7,1)
\class(7,1)
\class(7,1)
\class(7,1)
\class(7,1)
\class(7,1)
\class(7,1)
\class(7,3)
\class(7,3)
\class(7,3)
\class(7,3)
\class(7,3)
\class(7,3)
\class(7,5)
\class(7,5)
\class(7,5)
\class(7,5)
\class(7,5)
\class(7,5)
\class(7,5)
\class(7,7)
\class(7,7)
\class(7,7)
\class(7,7)
\class(7,9)
\class(7,9)
\class(7,9)
\class(7,13)
\class(7,13)
\class(7,17)
\class(7,21)
\class(8,-8)
\class(8,-6)
\class(8,-6)
\class(8,-6)
\class(8,-6)
\class(8,-4)
\class(8,-4)
\class(8,-4)
\class(8,-4)
\class(8,-4)
\class(8,-4)
\class(8,-4)
\class(8,-4)
\class(8,-4)
\class(8,-4)
\class(8,-2)
\class(8,-2)
\class(8,-2)
\class(8,-2)
\class(8,-2)
\class(8,-2)
\class(8,-2)
\class(8,-2)
\class(8,0)
\class(8,0)
\class(8,0)
\class(8,0)
\class(8,0)
\class(8,0)
\class(8,0)
\class(8,0)
\class(8,0)
\class(8,0)
\class(8,0)
\class(8,0)
\class(8,0)
\class(8,0)
\class(8,0)
\class(8,2)
\class(8,2)
\class(8,2)
\class(8,2)
\class(8,2)
\class(8,2)
\class(8,2)
\class(8,2)
\class(8,2)
\class(8,4)
\class(8,4)
\class(8,4)
\class(8,4)
\class(8,4)
\class(8,4)
\class(8,4)
\class(8,4)
\class(8,4)
\class(8,4)
\class(8,4)
\class(8,6)
\class(8,6)
\class(8,6)
\class(8,6)
\class(8,8)
\class(8,8)
\class(8,8)
\class(8,8)
\class(8,8)
\class(8,8)
\class(8,8)
\class(8,8)
\class(8,12)
\class(8,12)
\class(8,12)
\class(8,16)
\class(8,16)
\class(8,20)
\class(8,24)
\d3(2,-2,1,1)
\d3(3,-3,1,1)
\d3(3,-1,1,2)
\d3(4,-2,1,3)
\d3(4,0,1,1)
\d3(4,0,5,2)
\d3(5,-5,1,2)
\d3(5,-3,1,4)
\d3(5,-1,1,1)
\d3(5,-1,4,2)
\d3(5,1,2,3)
\d3(5,1,4,4)
\d3(6,-6,2,2)
\d3(6,-4,2,5)
\d3(6,-2,1,1)
\d3(6,-2,5,3)
\d3(6,0,2,4)
\d3(6,0,5,5)
\d3(6,2,1,1)
\d3(6,2,6,2)
\d3(6,2,8,3)
\d3(7,-7,2,1)
\d3(7,-5,1,6)
\d3(7,-5,2,7)
\d3(7,-3,1,3)
\d3(7,-3,2,1)
\d3(7,-3,6,4)
\d3(7,-3,8,6)
\d3(7,-1,1,5)
\d3(7,-1,4,7)
\d3(7,1,1,1)
\d3(7,1,7,2)
\d3(7,1,10,3)
\d3(7,3,2,4)
\d3(7,3,4,5)
\d3(7,3,6,6)
\d3(8,-8,1,3)
\d3(8,-6,1,7)
\d3(8,-4,1,3)
\d3(8,-4,8,5)
\d3(8,-4,9,6)
\d3(8,-2,1,6)
\d3(8,-2,4,8)
\d3(8,-2,5,9)
\d3(8,-2,8,11)
\d3(8,0,2,1)
\d3(8,0,9,3)
\d3(8,0,13,5)
\d3(8,2,2,5)
\d3(8,2,5,6)
\d3(8,2,8,7)
\d3(8,4,1,1)
\d3(8,4,7,2)
\d3(8,4,9,3)
\d3(8,4,11,4)
\classoptions["\scriptscriptstyle 1" left, label distance = -1pt](0,0,1)
\classoptions["\scriptscriptstyle \eta'" left, label distance = -1pt](1,3,1)
\classoptions["\scriptscriptstyle e" above, label distance = -1pt](2,2,1)
\classoptions["\scriptscriptstyle \nu" above, label distance = -1pt](3,1,1)
\classoptions["\scriptscriptstyle d" below, label distance = -1pt](4,0,1)
\classoptions["\scriptscriptstyle \zeta_1^2" below left, label distance = -1pt](2,-2,1)
\classoptions["\scriptscriptstyle \zeta_2" below left, label distance = -1pt](3,-3,1)
\classoptions["\scriptscriptstyle \zeta_1^4" below left, label distance = -1pt](4,-4,1)
\classoptions["\scriptscriptstyle \zeta_2^2" below left, label distance = -1pt](6,-6,1)
\classoptions["\scriptscriptstyle \zeta_3" below left, label distance = -1pt](7,-7,1)

\end{sseqdata}

\begin{figure}[htbp]
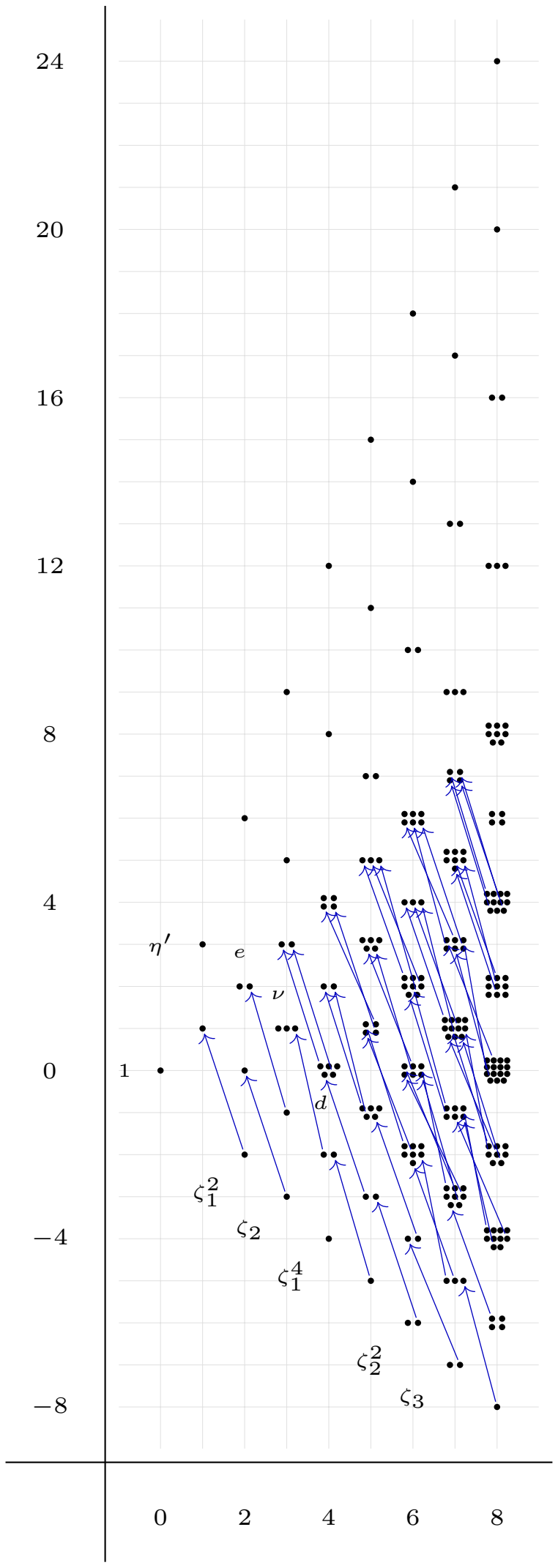
\centering
\printpage[name = Etwo, page = 3]
\caption{The $E_2$-page at $C_4$ in stems $0$ to $8$, with the $d_3$-differentials. A bidegree of dimension $k$ carries $k$ dots, and each arrow is one rank unit of $d_3$.}
\label{fig:E2d3}
\end{figure}

Each of $\eta'$, $\nu$, $e$, $w$ has slice degree $t=4$, each of $x,\gamma x,y,\gamma y$ has slice degree $t=2$, and $\mathcal{A}\langle 0\rangle$ has slice degree $t=0$. So the slice degree is divisible by $4$ on $\mathcal{N}$, while it is $2(p+q+r+s)$ on the transfer basis. Since $d_r$ raises the slice degree by $r-1$, all differentials of even length vanish, and a bidegree of slice degree $t\equiv2\bmod4$ consists of transfers alone.

\begin{definition}\label{def:strata}
For every page, the \emph{transfer stratum} $\mathcal{T}_r\subseteq E_r$ is the set of classes admitting a representative on $E_2$ in $\mathrm{im}(\mathrm{tr})=\mathcal{T}\oplus\nu\mathcal{A}\langle 0\rangle[e,w]$.
\end{definition}

By Frobenius reciprocity, $\mathcal{T}_r$ is an ideal.

\subsection{Koszul complexes}\label{sec:koszul} Every differential computed below is a Koszul differential, and every page is assembled from Koszul homology by long exact sequences. Let $M$ be an abelian group with two commuting endomorphisms $s_1,s_2$. The \emph{Koszul complex} $K(s_1,s_2;M)$ is the chain complex
\[M\{\kappa_1\kappa_2\}\xrightarrow{\ \partial\ }M\{\kappa_1,\kappa_2\}\xrightarrow{\ \partial\ }M\{1\},\qquad \partial(\kappa_i)=s_i.\]
It is graded by the exterior degree $i\in\{0,1,2\}$, which $\partial$ lowers by one, and we write $M^{(i)}$ for the part in exterior degree $i$ and $H_i(M)$ for the homology there. 

We repeatedly use three standard facts about Koszul complexes. First, $H_2(M)=\{m: s_1m=s_2m=0\}$ and $H_0(M)=M/(s_1M+s_2M)$. Second, a short exact sequence of modules compatible with the $s_i$ gives a long exact sequence in Koszul homology, whose connecting map sends a cycle to $\partial$ of any lift. Third, if $s_1,s_2$ is a regular sequence on $M$, then $H_i(M)=0$ for $i>0$.

In the following, let $P\subseteq E_2^{C_2}$ be the $\mathcal{A}\langle 1\rangle$-subalgebra generated by $x,\gamma x,y,\gamma y$, and let $V\subseteq E_2^{C_4}$ be the $\mathcal{A}\langle 1\rangle$-span of the monomials of Proposition \ref{prop:E3asacomplex}.

\begin{proposition}\label{prop:E2koszul}
At both levels, the $d_3$-differential is a Koszul differential in the following sense.
\begin{enumerate}
\item $E_2^{C_2}=P\{1,\zeta_1^2,\zeta_2,\zeta_1^2\zeta_2\}$ and $E_2^{C_4}=V\{1,\zeta_1^2,\zeta_2,\zeta_1^2\zeta_2\}$.
\item With $\zeta_1^2,\zeta_2$ in the role of $\kappa_1,\kappa_2$,
\[(E_2^{C_2},d_3)=K(x+\gamma x,\,y+\gamma y;\,P),\qquad (E_2^{C_4},d_3)=K\bigl(\mathrm{tr}(x),\mathrm{tr}(y);\,V\bigr).\]
\item The maps $\Res$ and $\mathrm{tr}$ carry $V$ into $P$ and $P$ into $V$, and are maps of Koszul complexes.
\end{enumerate}
\end{proposition}
\begin{proof}
We have $\mathcal{A}\langle 0\rangle=\mathcal{A}\langle 1\rangle\{1,\zeta_1^2,\zeta_2,\zeta_1^2\zeta_2\}$, and $d_3$ is $\mathcal{A}\langle 1\rangle$-linear and vanishes on $P$ and $V$, with $d_3(\zeta_1^2)=x+\gamma x$, $d_3(\zeta_2)=y+\gamma y$ at $C_2$ and $d_3(\zeta_1^2)=\mathrm{tr}(x)$, $d_3(\zeta_2)=\mathrm{tr}(y)$ at $C_4$ (Propositions \ref{prop:C2ZSS} and \ref{prop:E3asacomplex}). Since $\Res(e)=x\gamma x$, $\Res(w)=y\gamma y$ and $\Res(\eta')=\Res(\nu)=0$, $\Res$ carries $V$ into $P$, and so by Frobenius reciprocity multiplication by $\mathrm{tr}(x)$ and $\mathrm{tr}(y)$ preserves $V$; the transfer of a monomial of $P$ is a basis element of $V$ or zero by Propositions \ref{prop:transferpartE2C4} and \ref{prop:imtr}. Both maps are $\mathcal{A}\langle 0\rangle$-linear.
\end{proof}

\subsection{The $E_4$-page at $C_4$}\label{sec:E4C4} We compute $E_4^{C_4}$ by comparison with the $C_2$-level, using restriction and transfer.

\begin{lemma}\label{lem:resexact}
The restriction and transfer on $E_2$ have the following properties.
\begin{enumerate}
\item The sequence $E_2^{C_4}\xrightarrow{\Res}E_2^{C_2}\xrightarrow{\mathrm{tr}}E_2^{C_4}$ is exact.
\item $\ker(\Res)=\eta'\mathcal{A}\langle 0\rangle[\eta',e]\oplus\nu\mathcal{A}\langle 0\rangle[e,w]$, and $d_3$ vanishes on $\ker(\Res)$.
\item $\ker(\mathrm{tr})=(1+\gamma)E_2^{C_2}\oplus\mathcal{A}\langle 0\rangle[x\gamma x,y\gamma y]$.
\end{enumerate}
\end{lemma}
\begin{proof}
The monomials in $x$, $\gamma x$, $y$, $\gamma y$ whose degrees in $\bar{t}_1$ and in $\gamma\bar{t}_1$ agree are $(x\gamma x)^k(y\gamma y)^n$ and $(x\gamma x)^k(y\gamma y)^nx\gamma y$, and $\gamma$ fixes these and freely permutes all other monomials. By Propositions \ref{prop:transferpartE2C4} and \ref{prop:imtr}, $\mathrm{tr}$ carries any other monomial $m$ and $\gamma m$ to the same basis element of $\mathcal{T}$, carries $(x\gamma x)^k(y\gamma y)^nx\gamma y$ to $\nu e^kw^n$, and kills $(x\gamma x)^k(y\gamma y)^n$. This gives $\ker(\mathrm{tr})$. On the other side $\Res$ kills $\eta'$ and $\nu$, is injective on $\mathcal{T}$ with image $(1+\gamma)E_2^{C_2}$ by Lemma \ref{lem:resinj}, and carries $\mathcal{A}\langle 0\rangle[e,w]$ isomorphically onto $\mathcal{A}\langle 0\rangle[x\gamma x,y\gamma y]$. Finally $d_3$ vanishes on $\ker(\Res)$ because $\eta',\nu,e$ are $d_3$-cycles and $\mathrm{tr}(x)\eta'=\mathrm{tr}(x)\nu=0$ by Frobenius reciprocity.
\end{proof}

By Proposition \ref{prop:E2koszul} we therefore have a short exact sequence of Koszul complexes
\[
0\longrightarrow\ker(\Res)\longrightarrow E_2^{C_4}\xrightarrow{\ \Res\ }\ker(\mathrm{tr})\longrightarrow0,
\]
the first with zero differential. Its connecting maps vanish, because one checks that no nonzero element of $\ker(\Res)$ is a $d_3$-boundary in $E_2^{C_4}$.
So for each $i$ there is a short exact sequence
\begin{equation}\label{eq:resSESi}
0\longrightarrow\ker(\Res)^{(i)}\longrightarrow E_4^{C_4,(i)}\xrightarrow{\ \Res\ }H_i(\ker(\mathrm{tr}))\longrightarrow0,
\end{equation}
where $(-)^{(i)}$ denotes the part in exterior degree $i$, as in Section \ref{sec:koszul}. Concretely, the exterior degree of an element of $E_2^{C_4}=V\{1,\zeta_1^2,\zeta_2,\zeta_1^2\zeta_2\}$ is the number of the factors $\zeta_1^2,\zeta_2$ it involves, so that $E_2^{C_4,(0)}=V$, $E_2^{C_4,(1)}=V\{\zeta_1^2,\zeta_2\}$ and $E_2^{C_4,(2)}=V\{\zeta_1^2\zeta_2\}$, and $E_4^{C_4,(i)}=H_i(E_2^{C_4})$. Everything therefore comes down to the Koszul homology of $\ker(\mathrm{tr})$, which is a statement at the $C_2$-level.

\begin{proposition}\label{prop:kertr}
We have the following identifications.
\begin{enumerate}
\item $H_2(E_2^{C_2})=H_2(\ker(\mathrm{tr}))=0$, and
\[H_2(\mathrm{im}(\mathrm{tr}))=\nu\zeta_1^2\zeta_2\,\mathcal{A}\langle 1\rangle[e,w].\]
\item The map $H_0(\ker(\mathrm{tr}))\to H_0(E_2^{C_2})=\mathcal{A}\langle 1\rangle[x,y]$ is injective, with image $\mathcal{A}\langle 1\rangle[x^2,y^2]$, so $H_0(\ker(\mathrm{tr}))$ is free over $\mathcal{A}\langle 1\rangle$ on the classes of $(x\gamma x)^k(y\gamma y)^n$.
\item Every $d_3$-cycle of $\ker(\mathrm{tr})^{(1)}$ has the form $\zeta_1^2(1+\gamma)v_1+\zeta_2(1+\gamma)v_2$.
\end{enumerate}
\end{proposition}
\begin{proof}
These are statements about the $C_2$-level page $E_2^{C_2}$ of Proposition \ref{prop:C2ZSS}. The ring $\F_2[x,\gamma x,y,\gamma y]/(x\gamma y+y\gamma x)$ is a domain, $\gamma$ is the identity modulo $(x+\gamma x,y+\gamma y)$, and modulo this ideal the $\gamma$-invariants reduce into $\mathcal{A}\langle 1\rangle[x^2,y^2]\{1,xy\}$; each statement follows by comparing parities of exponents.
\end{proof}

\begin{theorem}\label{thm:E4C4}
There is a short exact sequence of $\mathcal{A}\langle 1\rangle[e,w]$-modules
\[0\to\mathcal{T}_4\to E_4^{C_4}\to \mathcal{A}\langle 1\rangle[e,w]\oplus\eta'\mathcal{A}\langle 0\rangle[\eta',e]\to 0,\]
which is split as a sequence of $\mathcal{A}\langle 1\rangle$-modules. Here $w$ acts on the second summand of the quotient by $w\eta'=\eta'(\zeta_1^2e+\zeta_2\eta')$.
\end{theorem}

The second summand of the quotient has $\F_2$-basis the monomials $b\eta'^ae^k$ with $b\in\mathcal{A}\langle 0\rangle$ and $a\ge1$.

\begin{proof}
We use the short exact sequences \eqref{eq:resSESi} together with Proposition \ref{prop:kertr}. In exterior degree $2$, $H_2(\ker(\mathrm{tr}))=0$, so $E_4^{C_4,(2)}=\ker(\Res)^{(2)}$: the monomials $c\,\eta'^ae^k\zeta_1^2\zeta_2$ ($c\in\mathcal{A}\langle 1\rangle$) and the classes $\nu\zeta_1^2\zeta_2\,c\,e^kw^n$ of $\mathcal{T}_4$. In degree $1$, a class of $H_1(\ker(\mathrm{tr}))$ is $\Res$ of $\zeta_1^2\mathrm{tr}(v_1)+\zeta_2\mathrm{tr}(v_2)$, whose $d_3$ is a boundary lying in $\ker(\Res)$, hence zero; so these classes lie in $\mathcal{T}_4$, and the rest is $\ker(\Res)^{(1)}$, spanned by the monomials $c\,\eta'^ae^k\zeta_1^2$, $c\,\eta'^ae^k\zeta_2$ and classes of $\mathcal{T}_4$. In degree $0$, $H_0(\ker(\mathrm{tr}))$ is the restriction of the cycles $c\,e^kw^n$, and the kernel is spanned by the monomials $c\,\eta'^ae^k$ and classes of $\mathcal{T}_4$.

So $E_4^{C_4}$ is spanned by $\mathcal{T}_4$ and the stated monomials. These monomials are independent modulo $\mathcal{T}_4$. Indeed, $\Res$ kills the $\eta'$-monomials and carries $\mathcal{T}_4$ into $(1+\gamma)E_2^{C_2}$, which is zero in $H_0(\ker(\mathrm{tr}))$ by Proposition \ref{prop:kertr}(2), so a relation has no $e,w$-part. An $\eta'$-combination lying in $\mathcal{T}_4$ lies in $\ker(\Res)\cap\mathrm{im}(\mathrm{tr})=\nu\mathcal{A}\langle 0\rangle[e,w]$, since every $d_3$-boundary is a transfer, and so it is zero. Finally, the module structure follows from Frobenius reciprocity and the identity $w\eta'=\eta'(\zeta_1^2e+\zeta_2\eta')+\nu e$, in which $\nu e\in\mathrm{im}(\mathrm{tr})$.
\end{proof}

\begin{proposition}\label{prop:tronE4}
In $E_4^{C_4}$ one has
\[\mathrm{tr}(1)=\mathrm{tr}(x)=\mathrm{tr}(y)=\mathrm{tr}(z)=0,\qquad\mathrm{tr}(xy)=\nu.\]
\end{proposition}
\begin{proof}
The first vanishes already on $E_2^{C_4}$, and $\mathrm{tr}(x)=d_3(\zeta_1^2)$, $\mathrm{tr}(y)=d_3(\zeta_2)$ by Proposition \ref{prop:d3atC4}; since $z=\zeta_1^2y+\zeta_2x$ this gives $\mathrm{tr}(z)=0$ as well. Finally $y\gamma x=x\gamma y$ and Frobenius reciprocity give
\[\mathrm{tr}(y)\mathrm{tr}(x)=\mathrm{tr}\big(y(x+\gamma x)\big)=\mathrm{tr}(xy)+\nu,\]
and the left-hand side is $d_3\big(\zeta_2\,\mathrm{tr}(x)\big)$.
\end{proof}

The transfer gives a second short exact sequence of Koszul complexes,
\[
0\longrightarrow\ker(\mathrm{tr})\longrightarrow E_2^{C_2}\xrightarrow{\ \mathrm{tr}\ }\mathrm{im}(\mathrm{tr})\longrightarrow0,
\]
and $\mathcal{T}_4$ is by definition the image of $H(\mathrm{im}(\mathrm{tr}))$ in $E_4^{C_4}$. Its long exact sequence measures the difference between $\mathcal{T}_4$ and the image of the transfer on $E_4^{C_4}$.

\begin{proposition}\label{prop:T4}
\[\mathcal{T}_4=\mathrm{im}\bigl(\mathrm{tr}\colon E_4^{C_2}\to E_4^{C_4}\bigr)+\nu\zeta_1^2\zeta_2\,\mathcal{A}\langle 1\rangle[e,w].\]
\end{proposition}
\begin{proof}
In the long exact sequence, the cokernel of $\mathrm{tr}\colon E_4^{C_2}\to H(\mathrm{im}(\mathrm{tr}))$ is all of $H_2(\mathrm{im}(\mathrm{tr}))=\nu\zeta_1^2\zeta_2\,\mathcal{A}\langle 1\rangle[e,w]$ in exterior degree $2$, since $H_2(E_2^{C_2})=0$; it is zero in degree $1$, since $H_0(\ker(\mathrm{tr}))\to H_0(E_2^{C_2})$ is injective; and it is zero in degree $0$. Both facts are Proposition \ref{prop:kertr}.
\end{proof}

\section{The $E_6$- and $E_8$-pages at $C_4$}\label{sec:C4higherpages}
From now on $E_r$ denotes the $C_4$-level page $E_r^{C_4}$, and the $C_2$-level page is always written $E_r^{C_2}$.

\subsection{$d_5$ differentials at $C_4$} The differential $d_5$ vanishes at the $C_2$-level, so $d_5$ vanishes at $C_4$ on the image of the transfer. This does not apply to all of $\mathcal{T}_4$, which by Proposition \ref{prop:T4} also contains the classes $\nu\zeta_1^2\zeta_2\,c\,e^kw^n$ with $c\in\mathcal{A}\langle 1\rangle$.

\begin{proposition}\label{prop:d5zeta}
\[d_5(\zeta_1^4)=\nu,\qquad d_5(\zeta_2^2)=\zeta_1^2\nu,\qquad d_5(\zeta_3)=\zeta_2\nu.\]
\end{proposition}
\begin{proof}
The class $\nu$ must be killed by a differential, and $d_5(\zeta_1^4)$ is the only possibility.
For the other two, we apply $1\otimes d_5$ to
\begin{align*}
\psi(\zeta_2^2)&=1\otimes\zeta_2^2+\zeta_1^2\otimes\zeta_1^4+\zeta_2^2\otimes1,\\
\psi(\zeta_3)&=1\otimes\zeta_3+\zeta_1\otimes\zeta_2^2+\zeta_2\otimes\zeta_1^4+\zeta_3\otimes1.
\end{align*}
The results are the coactions of $\zeta_1^2\nu$ and of $\zeta_2\nu$, so the differentials are as claimed up to a comodule primitive of $E_4$ in bidegree $(5,-1)$, resp. $(6,-2)$. In the first of these, $E_4$ is spanned by $\eta'\zeta_1^4$ and $\zeta_1^2\nu$, whose coactions have nonzero $\xi_1^4\otimes\eta'$ and $\xi_1^2\otimes\nu$ terms, so there is no nonzero primitive; the second is the same.
\end{proof}

Let $\partial\colon\mathcal{A}\langle 1\rangle\to\mathcal{A}\langle 0\rangle$ be the cap product with $\Sq^4$, that is, $\partial c=\Sq^4_*c$ for the coaction of $\mathcal{A}\langle 1\rangle$. Since $\Sq^1$, $\Sq^2$ and $\Sq^3$ act trivially on $\mathcal{A}\langle 1\rangle$, the Cartan formula makes $\partial$ a derivation there, and
\[\partial(\zeta_1^4)=1,\qquad \partial(\zeta_2^2)=\zeta_1^2,\qquad \partial(\zeta_3)=\zeta_2,\]
so that Proposition \ref{prop:d5zeta} reads $d_5(c)=\nu\,\partial c$ for $c\in\mathcal{A}\langle 1\rangle$.

The elements $\zeta_1^2$ and $\zeta_2$ are not classes on $E_4$, but $a_\sigma\zeta_1^2$ and $a_\sigma\zeta_2$ survive and support the following differentials.

\begin{proposition}\label{prop:ROinput}
$d_5(a_\sigma\zeta_1^2)=a_\sigma\eta'$ and $d_5(a_\sigma\zeta_2)=a_\sigma e$.
\end{proposition}
\begin{proof}
In the bidegrees of $d_5(a_\sigma\zeta_1^2)$ and $d_5(a_\sigma\zeta_2)$, the page $E_2$ is spanned by $a_\sigma\eta'$ and $a_\sigma e$ respectively, so it suffices to show that the two values are nonzero. Hill--Hopkins--Ravenel show $d_5(u_{2\sigma})=Na_\sigma^3a_\lambda$ \cite[Theorem 11.13]{hhrc4}, and $u_{2\sigma}=x_\sigma^2+\zeta_1^2a_\sigma^2$ with $x_\sigma$ permanent, so $a_\sigma d_5(a_\sigma\zeta_1^2)\ne0$. Since $\psi(a_\sigma\zeta_2)$ has $\xi_1$-component $a_\sigma\zeta_1^2$, the class $\psi(d_5(a_\sigma\zeta_2))$ has $\xi_1$-component $a_\sigma\eta'\ne0$.
\end{proof}

The behavior of $d_5$ comes down to proving one more differential, namely
\[d_5(a_\sigma\zeta_1^2\zeta_2)=a_\sigma(\eta'\zeta_2+e\zeta_1^2)+x_\sigma\nu.\]
For this we need the following $d_3$-cycle.

\begin{lemma}\label{lem:wlambda}
The class $w_\lambda=u_\lambda+\zeta_1^2a_\lambda$, in bidegree $(2-\lambda,0)$, is a $d_3$-cycle with $\Res(w_\lambda)=x_1^2$, and
\[d_5(w_\lambda)=a_\lambda\eta'.\]
Moreover $w=Nx_\sigma w_\lambda+\zeta_2\eta'$, and $w$ is a $d_5$-cycle.
\end{lemma}
\begin{proof}
In the slice spectral sequence $d_3(u_\lambda)=\mathrm{tr}(a_{\sigma_2}^3\bar t_1)$ \cite[Theorem 11.13]{hhrc4}, \cite[Theorem 3.1]{hswx}, so here $d_3(u_\lambda)=\mathrm{tr}(a_{\sigma_2}^2x)=a_\lambda\mathrm{tr}(x)=d_3(\zeta_1^2a_\lambda)$.

The bidegree of $d_5(w_\lambda)$ is $(1-\lambda,5)$, where $E_2$ is spanned by $Na_\sigma a_\lambda^2=a_\lambda\eta'$. So it suffices to show that $d_5(w_\lambda)\ne0$. Now $a_\sigma w_\lambda=a_\sigma u_\lambda+a_\lambda\cdot a_\sigma\zeta_1^2$, and $a_\sigma u_\lambda$ is a permanent transfer, so $a_\sigma d_5(w_\lambda)=a_\lambda a_\sigma\eta'$ by Proposition \ref{prop:ROinput}. The class $a_\sigma a_\lambda\eta'=Na_\sigma^2a_\lambda^2$ is the only class of $E_2$ in $(1-\sigma-\lambda,6)$ and is hit by nothing.

Finally $Nx_\sigma w_\lambda=d+\zeta_1^2e$, and the differential $d_5(\zeta_2\eta')=d_5(Na_\lambda\cdot a_\sigma\zeta_2)=\eta'e$ cancels $d_5(Nx_\sigma w_\lambda)=e\eta'$.
\end{proof}

\begin{proposition}\label{prop:ROinput2}
$d_5(\eta'\zeta_1^2\zeta_2)=\eta'w$ and $d_5(\nu\zeta_1^2\zeta_2)=\nu d$, where $\nu d=\nu(w+\zeta_1^2e)$.
\end{proposition}
\begin{proof}
Both classes are expressible in terms of the $d_3$-cycle $a_\sigma\zeta_1^2\zeta_2$: on $E_2$,
\begin{align*}
\eta'\zeta_1^2\zeta_2&=Na_\lambda\cdot a_\sigma\zeta_1^2\zeta_2,\\
\nu\zeta_2&=Nw_\lambda\cdot a_\sigma\zeta_2+Na_\lambda\cdot a_\sigma\zeta_1^2\zeta_2,\\
\nu\zeta_1^2\zeta_2&=Nw_\lambda\cdot a_\sigma\zeta_1^2\zeta_2+\eta'\zeta_1^4\zeta_2.
\end{align*}
One checks that $\nu\zeta_2=\mathrm{tr}(z\gamma y)$, so that $d_5(\nu\zeta_2)=0$, and the second identity with Lemma \ref{lem:wlambda}, Proposition \ref{prop:ROinput} and $Nw_\lambda a_\sigma=\nu+\zeta_1^2\eta'$ gives
\[d_5(\eta'\zeta_1^2\zeta_2)=Na_\lambda d_5(a_\sigma\zeta_1^2\zeta_2)=d_5(Nw_\lambda\cdot a_\sigma\zeta_2)=\eta'^2\zeta_2+(\nu+\zeta_1^2\eta')e=\eta'w,\]
the last step by $\eta'd=\nu e$.

For the second value we need $d_5(a_\sigma\zeta_1^2\zeta_2)$ itself. It lies in bidegree $(4-\sigma,1)$, where $E_2$ is spanned by the norm monomials $a_\sigma\eta'\zeta_2$, $a_\sigma e\zeta_1^2$, $x_\sigma\nu$ and the transfer $x_\sigma\mathrm{tr}(xy)$, and $d_3\bigl(x_\sigma\mathrm{tr}(\zeta_2x)\bigr)=x_\sigma\mathrm{tr}(xy)+x_\sigma\nu$ is the only boundary. Proposition \ref{prop:ROinput} gives $a_\sigma d_5(a_\sigma\zeta_1^2\zeta_2)=a_\sigma^2(\eta'\zeta_2+e\zeta_1^2)$, so
\[d_5(a_\sigma\zeta_1^2\zeta_2)=a_\sigma(\eta'\zeta_2+e\zeta_1^2)+\epsilon\,x_\sigma\nu,\qquad \epsilon\in\F_2,\]
and $Na_\lambda d_5(a_\sigma\zeta_1^2\zeta_2)=\eta'^2\zeta_2+\eta'e\zeta_1^2+\epsilon\,\nu e$. Comparing with $d_5(\eta'\zeta_1^2\zeta_2)$ gives $(1+\epsilon)\nu e=0$ on $E_4$, and therefore $\epsilon=1$, as $\nu e\ne0$ in $E_4$. Now the third identity, with $Nw_\lambda x_\sigma=d+\zeta_1^2e$ and $d_5(\eta'\zeta_1^4\zeta_2)=d_5\bigl(Na_\lambda\cdot\zeta_1^4\cdot a_\sigma\zeta_2\bigr)=\eta'e\zeta_1^4$, gives
\[d_5(\nu\zeta_1^2\zeta_2)=(\nu+\zeta_1^2\eta')(\eta'\zeta_2+e\zeta_1^2)+(d+\zeta_1^2e)\nu+\eta'^2\zeta_1^2\zeta_2+\eta'e\zeta_1^4=\nu d.\qedhere\]
\end{proof}

Applying the Leibniz rule now gives the following complete description of $d_5$.

\begin{theorem}\label{thm:d5}
The differential $d_5$ is $\F_2[e,w]$-linear and vanishes on the image of the transfer. On the summand $\mathcal{A}\langle 1\rangle[e,w]$ of Theorem \ref{thm:E4C4} it is
\begin{equation}\label{eq:d5free}
    d_5(e^kw^nc)=e^kw^n\,\nu\,\partial c,\qquad c\in\mathcal{A}\langle 1\rangle;
\end{equation}
on the summand $\eta'\mathcal{A}\langle 0\rangle[\eta',e]$, for $c\in\mathcal{A}\langle 1\rangle$ and $a\ge1$,
\begin{equation}\label{eq:d5eta}
    d_5(\eta'^ae^kb)=\begin{cases}
    0, & b=c,\\
    \eta'^{a+1}e^kc, & b=c\zeta_1^2,\\
    \eta'^ae^{k+1}c, & b=c\zeta_2,\\
    \eta'^ae^k\,w\,c, & b=c\zeta_1^2\zeta_2;
    \end{cases}
\end{equation}
and on the classes $\nu b\,e^kw^n$ of $\mathcal{T}_4$,
\begin{equation}\label{eq:d5nu}
    d_5(e^kw^n\nu b)=\begin{cases}
    0, & b=c,\ c\zeta_1^2\ \text{or}\ c\zeta_2,\\
    e^kw^n\,\nu d\,c, & b=c\zeta_1^2\zeta_2.
    \end{cases}
\end{equation}
Together with Proposition \ref{prop:T4} this determines $d_5$ on all of $E_4$.
\end{theorem}
\begin{proof}
The classes $\eta',\nu,e$ are permanent and $w$ is a $d_5$-cycle by Lemma \ref{lem:wlambda}, which gives linearity, and \eqref{eq:d5free} is Proposition \ref{prop:d5zeta}. For \eqref{eq:d5eta}, write
\[\eta'^ae^kc\,\zeta_1^2=(Na_\lambda\eta'^{a-1}e^kc)\cdot(a_\sigma\zeta_1^2)\]
and apply Proposition \ref{prop:ROinput}, noting that $d_5(c)=\nu\,\partial c$ contributes nothing as $a_\sigma\nu=0$; the cases $b=c$ and $b=c\zeta_2$ are the same, and the case $b=c\zeta_1^2\zeta_2$ follows in the same way for $a\ge2$, and from Proposition \ref{prop:ROinput2} for $a=1$. For \eqref{eq:d5nu}, the classes $\nu c=\mathrm{tr}(c\,x\gamma y)$, $\nu c\zeta_1^2=\mathrm{tr}(c\,x\gamma z)$ and $\nu c\zeta_2=\mathrm{tr}(c\,z\gamma y)$ are transfers of $d_3$-cycles (Proposition \ref{prop:d7nu} below), and the last case is Proposition \ref{prop:ROinput2}.
\end{proof}

\begin{sseqdata}[name = Efour, hsss, x range = {0}{12}, y range = {-12}{24}, y tick step = 4, differentials = {dfive}]
\class(0,0)
\class(1,3)
\class(2,2)
\class(2,6)
\class(3,1)
\class(3,1)
\class(3,5)
\class(3,9)
\class(4,-4)
\class(4,0)
\class(4,0)
\class(4,4)
\class(4,4)
\class(4,8)
\class(4,12)
\class(5,-1)
\class(5,-1)
\class(5,3)
\class(5,3)
\class(5,3)
\class(5,7)
\class(5,7)
\class(5,11)
\class(5,15)
\class(6,-6)
\class(6,-2)
\class(6,-2)
\class(6,-2)
\class(6,2)
\class(6,2)
\class(6,2)
\class(6,6)
\class(6,6)
\class(6,6)
\class(6,10)
\class(6,10)
\class(6,14)
\class(6,18)
\class(7,-7)
\class(7,-3)
\class(7,-3)
\class(7,-3)
\class(7,-1)
\class(7,1)
\class(7,1)
\class(7,1)
\class(7,1)
\class(7,5)
\class(7,5)
\class(7,5)
\class(7,5)
\class(7,9)
\class(7,9)
\class(7,9)
\class(7,13)
\class(7,13)
\class(7,17)
\class(7,21)
\class(8,-4)
\class(8,-4)
\class(8,-4)
\class(8,-4)
\class(8,-4)
\class(8,0)
\class(8,0)
\class(8,0)
\class(8,0)
\class(8,0)
\class(8,0)
\class(8,4)
\class(8,4)
\class(8,4)
\class(8,4)
\class(8,8)
\class(8,8)
\class(8,8)
\class(8,8)
\class(8,12)
\class(8,12)
\class(8,12)
\class(8,16)
\class(8,16)
\class(8,20)
\class(8,24)
\class(9,-5)
\class(9,-5)
\class(9,-5)
\class(9,-5)
\class(9,-1)
\class(9,-1)
\class(9,-1)
\class(9,-1)
\class(9,-1)
\class(9,-1)
\class(9,1)
\class(9,3)
\class(9,3)
\class(9,3)
\class(9,3)
\class(9,3)
\class(9,3)
\class(9,7)
\class(9,7)
\class(9,7)
\class(9,7)
\class(9,7)
\class(9,11)
\class(9,11)
\class(9,11)
\class(9,11)
\class(9,15)
\class(9,15)
\class(9,15)
\class(9,19)
\class(9,19)
\class(9,23)
\class(10,-10)
\class(10,-6)
\class(10,-6)
\class(10,-6)
\class(10,-6)
\class(10,-6)
\class(10,-6)
\class(10,-2)
\class(10,-2)
\class(10,-2)
\class(10,-2)
\class(10,-2)
\class(10,-2)
\class(10,-2)
\class(10,2)
\class(10,2)
\class(10,2)
\class(10,2)
\class(10,2)
\class(10,2)
\class(10,2)
\class(10,2)
\class(10,6)
\class(10,6)
\class(10,6)
\class(10,6)
\class(10,6)
\class(10,6)
\class(10,10)
\class(10,10)
\class(10,10)
\class(10,10)
\class(10,10)
\class(10,14)
\class(10,14)
\class(10,14)
\class(10,14)
\class(10,18)
\class(10,18)
\class(10,18)
\class(10,22)
\class(10,22)
\class(11,-11)
\class(11,-7)
\class(11,-7)
\class(11,-7)
\class(11,-7)
\class(11,-5)
\class(11,-3)
\class(11,-3)
\class(11,-3)
\class(11,-3)
\class(11,-3)
\class(11,-3)
\class(11,-3)
\class(11,-3)
\class(11,-1)
\class(11,1)
\class(11,1)
\class(11,1)
\class(11,1)
\class(11,1)
\class(11,1)
\class(11,1)
\class(11,1)
\class(11,3)
\class(11,5)
\class(11,5)
\class(11,5)
\class(11,5)
\class(11,5)
\class(11,5)
\class(11,5)
\class(11,5)
\class(11,9)
\class(11,9)
\class(11,9)
\class(11,9)
\class(11,9)
\class(11,9)
\class(11,9)
\class(11,13)
\class(11,13)
\class(11,13)
\class(11,13)
\class(11,13)
\class(11,17)
\class(11,17)
\class(11,17)
\class(11,17)
\class(11,21)
\class(11,21)
\class(11,21)
\class(12,-8)
\class(12,-8)
\class(12,-8)
\class(12,-8)
\class(12,-8)
\class(12,-8)
\class(12,-4)
\class(12,-4)
\class(12,-4)
\class(12,-4)
\class(12,-4)
\class(12,-4)
\class(12,-4)
\class(12,-4)
\class(12,-4)
\class(12,-4)
\class(12,0)
\class(12,0)
\class(12,0)
\class(12,0)
\class(12,0)
\class(12,0)
\class(12,0)
\class(12,0)
\class(12,0)
\class(12,0)
\class(12,0)
\class(12,4)
\class(12,4)
\class(12,4)
\class(12,4)
\class(12,4)
\class(12,4)
\class(12,4)
\class(12,4)
\class(12,4)
\class(12,8)
\class(12,8)
\class(12,8)
\class(12,8)
\class(12,8)
\class(12,8)
\class(12,8)
\class(12,8)
\class(12,12)
\class(12,12)
\class(12,12)
\class(12,12)
\class(12,12)
\class(12,12)
\class(12,12)
\class(12,16)
\class(12,16)
\class(12,16)
\class(12,16)
\class(12,16)
\class(12,20)
\class(12,20)
\class(12,20)
\class(12,20)
\class(12,24)
\class(12,24)
\class(12,24)
\d5(3,1,1,1)
\d5(4,-4,1,1)
\d5(4,0,1,1)
\d5(4,4,1,1)
\d5(5,3,1,1)
\d5(5,7,1,1)
\d5(6,-6,1,1)
\d5(6,-2,1,1)
\d5(6,-2,2,2)
\d5(6,2,1,1)
\d5(6,6,1,1)
\d5(6,10,1,1)
\d5(7,-7,1,1)
\d5(7,-3,1,1)
\d5(7,1,1,1)
\d5(7,5,1,1)
\d5(7,9,1,1)
\d5(7,13,1,1)
\d5(8,-4,1,1)
\d5(8,-4,2,2)
\d5(8,-4,3,3)
\d5(8,0,1,1)
\d5(8,0,2,2)
\d5(8,0,3,3)
\d5(8,4,1,1)
\d5(8,4,2,2)
\d5(8,8,1,1)
\d5(8,12,1,1)
\d5(8,16,1,1)
\d5(9,-5,1,1)
\d5(9,-5,2,2)
\d5(9,-1,1,1)
\d5(9,3,1,1)
\d5(9,3,2,2)
\d5(9,7,1,1)
\d5(9,7,2,2)
\d5(9,11,1,1)
\d5(9,15,1,1)
\d5(9,19,1,1)
\d5(10,-10,1,1)
\d5(10,-6,1,1)
\d5(10,-6,2,2)
\d5(10,-6,3,3)
\d5(10,-6,4,4)
\d5(10,-2,1,1)
\d5(10,-2,2,2)
\d5(10,-2,3,3)
\d5(10,-2,4,4)
\d5(10,2,1,1)
\d5(10,2,2,2)
\d5(10,2,3,3)
\d5(10,6,1,1)
\d5(10,6,2,2)
\d5(10,6,3,3)
\d5(10,10,1,1)
\d5(10,10,2,2)
\d5(10,14,1,1)
\d5(10,18,1,1)
\d5(11,-11,1,1)
\d5(11,-7,1,1)
\d5(11,-7,2,2)
\d5(11,-3,1,1)
\d5(11,-3,2,2)
\d5(11,-3,3,3)
\d5(11,-3,4,4)
\d5(11,1,1,1)
\d5(11,1,2,2)
\d5(11,5,1,1)
\d5(11,5,2,2)
\d5(11,9,1,1)
\d5(11,9,2,2)
\d5(11,9,3,3)
\d5(11,13,1,1)
\d5(11,13,2,2)
\d5(11,17,1,1)
\d5(12,-8,1,1)
\d5(12,-8,2,2)
\d5(12,-8,3,3)
\d5(12,-4,1,1)
\d5(12,-4,2,2)
\d5(12,-4,3,3)
\d5(12,-4,4,4)
\d5(12,-4,5,5)
\d5(12,0,1,1)
\d5(12,0,2,2)
\d5(12,0,3,3)
\d5(12,0,4,4)
\d5(12,0,5,5)
\d5(12,0,6,6)
\d5(12,4,1,1)
\d5(12,4,2,2)
\d5(12,4,3,3)
\d5(12,4,4,4)
\d5(12,8,1,1)
\d5(12,8,2,2)
\d5(12,8,3,3)
\d5(12,12,1,1)
\d5(12,12,2,2)
\d5(12,12,3,3)
\d5(12,16,1,1)
\d5(12,16,2,2)

\classoptions["{\scriptscriptstyle 1}" left, label distance = -1pt](0,0,1)
\classoptions["{\scriptscriptstyle \eta^{\prime}}" left, label distance = -1pt](1,3,1)
\classoptions["{\scriptscriptstyle e}" left, label distance = -1pt](2,2,1)
\classoptions["{\scriptscriptstyle \nu,\,\eta^{\prime}\zeta_1^2}" right, label distance = -1pt](3,1,1)
\classoptions["{\scriptscriptstyle w,\,\eta^{\prime}\zeta_2}" right, label distance = -1pt](4,0,1)
\classoptions["{\scriptscriptstyle \zeta_1^4}" left, label distance = -1pt](4,-4,1)
\classoptions["{\scriptscriptstyle \zeta_2^2}" left, label distance = -1pt](6,-6,1)
\classoptions["{\scriptscriptstyle \zeta_3}" left, label distance = -1pt](7,-7,1)
\classoptions["{\scriptscriptstyle \mathrm{tr}(xyz)}" below, label distance = -1pt](7,-1,1)
\end{sseqdata}

\begin{figure}[htbp]
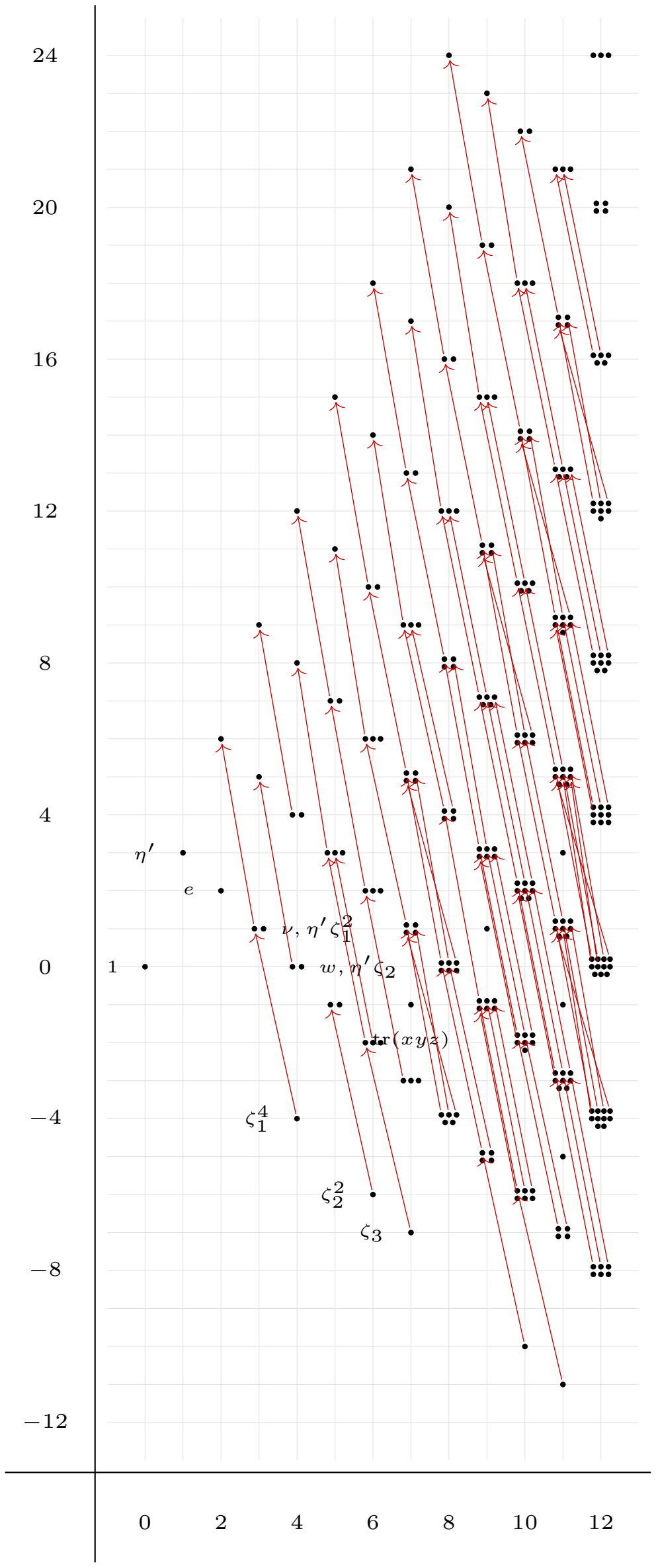
\centering
\printpage[name = Efour, page = 5]
\caption{The page $E_4=E_5$ in stems $0$ to $12$ and filtrations at most $24$, with the $d_5$-differentials. Where several classes share a bidegree, the labels name a basis of it, and an arrow is one unit of rank.}
\label{fig:E4d5}
\end{figure}

\subsection{The $E_6$-page at $C_4$} We compute $E_6$ from the filtration $\mathcal{T}_4\subseteq E_4$. The point is that $\mathcal{T}_4$ is a subcomplex for $d_5$, and that on the quotient $E_4/\mathcal{T}_4$ the two summands of Theorem \ref{thm:E4C4} do not interact: $\mathcal{A}\langle 1\rangle[e,w]$ has zero differential there, and $\eta'\mathcal{A}\langle 0\rangle[\eta',e]$ is a Koszul complex.

\begin{theorem}\label{thm:E6C4}
\begin{enumerate}
\item The quotient $E_6/\mathcal{T}_6$ is
\[\F_2[e,w]\{1,g\}\ \oplus\ \F_2\{\eta'c:c\in\mathcal{A}\langle 1\rangle\},\qquad g:=w\zeta_1^4+e\zeta_2^2,\]
where $g$, in bidegree $(8,-4)$, lifts to the $d_5$-cycle $g+\nu\zeta_1^2\zeta_2$, and the eight classes $\eta'c$, in bidegrees $(1+|c|,3-|c|)$, are annihilated by $\eta'$, $e$ and $w$.
\item $\mathcal{T}_6$ is the image of the transfer $E_6^{C_2}\to E_6$.
\end{enumerate}
\end{theorem}
\begin{proof}
We write $E_4=\mathcal{T}_4\oplus\mathcal{A}\langle 1\rangle[e,w]\oplus\eta'\mathcal{A}\langle 0\rangle[\eta',e]$ by Theorem \ref{thm:E4C4}. By Proposition \ref{prop:T4}, an element of $\mathcal{T}_4$ is $v+\nu\zeta_1^2\zeta_2\beta$ with $v$ in the image of the transfer and $\beta\in\mathcal{A}\langle 1\rangle[e,w]$, so Theorem \ref{thm:d5} gives three facts:
\begin{enumerate}
\item[(i)] $d_5(\mathcal{T}_4)=\nu d\,\mathcal{A}\langle 1\rangle[e,w]\subseteq\mathcal{T}_4$, with $d_5(v+\nu\zeta_1^2\zeta_2\beta)=\nu d\beta$;
\item[(ii)] $d_5$ carries $\mathcal{A}\langle 1\rangle[e,w]$ into $\mathcal{T}_4$;
\item[(iii)] $d_5$ carries $\eta'\mathcal{A}\langle 0\rangle[\eta',e]$ into $\eta'\mathcal{A}\langle 0\rangle[\eta',e]\oplus\mathcal{T}_4$, and the only $\mathcal{T}_4$-component is the term $\nu e^{k+1}c$ of $d_5(\eta'e^kc\zeta_1^2\zeta_2)=\eta'e^kwc=\eta'e^k(\zeta_1^2e+\zeta_2\eta')c+\nu e^{k+1}c$.
\end{enumerate}
So $\mathcal{T}_4$ is a subcomplex, and the quotient complex $E_4/\mathcal{T}_4$ is the direct sum of $\mathcal{A}\langle 1\rangle[e,w]$, with zero differential, and $\eta'\mathcal{A}\langle 0\rangle[\eta',e]$, with the differential $\bar d_5$ induced modulo $\mathcal{T}_4$. The long exact sequence of $0\to\mathcal{T}_4\to E_4\to E_4/\mathcal{T}_4\to0$ gives
\[E_6/\mathcal{T}_6=\ker\Bigl(H(E_4/\mathcal{T}_4)\longrightarrow H(\mathcal{T}_4)\Bigr),\qquad [\alpha]\mapsto[d_5\alpha];\]
that is, a class of the quotient survives exactly when a lift can be corrected by an element of $\mathcal{T}_4$ into a $d_5$-cycle.

On $\eta'\mathcal{A}\langle 0\rangle[\eta',e]$, (iii) says that modulo $\mathcal{T}_4$ the elements $\zeta_1^2$ and $\zeta_2$ act as $\eta'$ and $e$:
\begin{align*}
\bar d_5(\eta'^ae^kc\,\zeta_1^2)&=\eta'^{a+1}e^kc,\qquad \bar d_5(\eta'^ae^kc\,\zeta_2)=\eta'^ae^{k+1}c,\\
\bar d_5(\eta'^ae^kc\,\zeta_1^2\zeta_2)&=\eta'^{a+1}e^kc\,\zeta_2+\eta'^ae^{k+1}c\,\zeta_1^2.
\end{align*}
With $\zeta_1^2,\zeta_2$ in the role of $\kappa_1,\kappa_2$ this is the Koszul complex of $(\eta',e)$ on $\eta'\mathcal{A}\langle 1\rangle[\eta',e]$. That is a regular sequence, so the homology is $\eta'\mathcal{A}\langle 1\rangle[\eta',e]/(\eta'^2,\eta'e)=\eta'\mathcal{A}\langle 1\rangle$. These classes $\eta'c$ are $d_5$-cycles already in $E_4$, and they are killed by $\eta'$, $e$ and $w$, since $\eta'^2c=\bar d_5(\eta'c\,\zeta_1^2)$, $e\eta'c=\bar d_5(\eta'c\,\zeta_2)$ and $w\eta'c=\eta'(\zeta_1^2e+\zeta_2\eta')c=\bar d_5(\eta'c\,\zeta_1^2\zeta_2)$ modulo $\mathcal{T}_4$.

On $\mathcal{A}\langle 1\rangle[e,w]$ the differential is zero, and $[\alpha]$ maps to the class of $d_5\alpha=\nu\,\partial\alpha$ by \eqref{eq:d5free}. This vanishes in $H(\mathcal{T}_4)$ exactly when $\nu\,\partial\alpha\in d_5(\mathcal{T}_4)=\nu d\,\mathcal{A}\langle 1\rangle[e,w]$, by (i). Multiplication by $\nu$ is injective on $\mathcal{A}\langle 0\rangle[e,w]$, and $\nu d=\nu(w+\zeta_1^2e)$, so the condition is
\[\partial\alpha\in(w+\zeta_1^2e)\,\mathcal{A}\langle 1\rangle[e,w].\]
We write $\alpha=\alpha_0+\alpha_1\zeta_3$ with $\alpha_0,\alpha_1\in E(\zeta_1^4,\zeta_2^2)[e,w]$. Since $\partial$ sends $\zeta_1^4$ and $\zeta_2^2$ into $\mathcal{A}\langle 1\rangle\{1,\zeta_1^2\}$ and $\zeta_3$ to $\zeta_2$, the component of $\partial\alpha$ along $\zeta_2$ is $\alpha_1$, which must therefore vanish. So $\alpha=p_0+p_1\zeta_1^4+p_2\zeta_2^2+p_3\zeta_1^4\zeta_2^2$ with $p_i\in\F_2[e,w]$, and
\[\partial\alpha=(p_1+p_3\zeta_2^2)+(p_2+p_3\zeta_1^4)\zeta_1^2.\]
The condition now reads $e(p_1+p_3\zeta_2^2)=w(p_2+p_3\zeta_1^4)$. Comparing coefficients of $\zeta_2^2$ gives $p_3=0$, and comparing coefficients of $1$ gives $ep_1=wp_2$, so $p_1=wp$ and $p_2=ep$ for some $p\in\F_2[e,w]$. Thus $\alpha$ is an $\F_2[e,w]$-combination of $1$ and $g$, and \eqref{eq:d5free} and \eqref{eq:d5nu} give $d_5(g+\nu\zeta_1^2\zeta_2)=\nu d+\nu d=0$. This proves (1).

For (2), $\mathcal{T}_6$ is the image of the $d_5$-cycles of $\mathcal{T}_4$. By (i), $d_5(v+\nu\zeta_1^2\zeta_2\beta)=\nu d\beta$ is nonzero for $\beta\ne0$. So the $d_5$-cycles of $\mathcal{T}_4$ are the image of the transfer $E_4^{C_2}\to E_4$, and $E_6^{C_2}=E_4^{C_2}$ because the $C_2$-level spectral sequence has no $d_5$.
\end{proof}

\begin{proposition}\label{prop:O6}
The part of $E_6$ of slice degree $t\equiv2\bmod4$ agrees with that of $E_4$, and is the free $\F_2[e,w]$-module of rank $8$ on the classes $\mathrm{tr}(c\,xyz)$ for $c\in\mathcal{A}\langle 1\rangle$, in bidegrees $(7+|c|,-1-|c|)$.
\end{proposition}
\begin{proof}
In slice degree $\equiv2$ every class of $E_2$ is a transfer, and by Proposition \ref{prop:T4} every class of $E_4$ is the transfer of a class of $E_4^{C_2}$, the classes $\nu\zeta_1^2\zeta_2\,c\,e^kw^n$ lying in slice degree $t\equiv0$. So $d_5$ vanishes on this part, by naturality; and as $d_5$ raises the slice degree by $4$, every $d_5$ into it starts in it. Hence $E_6=E_4$ there.

By Corollary \ref{cor:E4pageC2}, the part of $E_4^{C_2}$ in these slice degrees is that of odd word length in $x,y,z$, which is free over $\mathcal{A}\langle 1\rangle[x^2,y^2]$ on $x,y,z,xyz$. The transfer carries $x^{2k}y^{2n}$ to $e^kw^n$ by Frobenius reciprocity and kills $x$, $y$ and $z$ by Proposition \ref{prop:tronE4}. The classes $c\,e^kw^n\mathrm{tr}(x)=d_3(c\,e^kw^n\zeta_1^2)$, $c\,e^kw^n\mathrm{tr}(y)=d_3(c\,e^kw^n\zeta_2)$ and $c\,e^kw^n\mathrm{tr}(z)=d_3(c\,e^kw^n\zeta_1^2\zeta_2)$ are the only $d_3$-boundaries in these slice degrees coming from the norm part. What is left is free on $\mathcal{A}\langle 1\rangle[e,w]\{\mathrm{tr}(xyz)\}$.
\end{proof}

\subsection{The $d_7$ differentials at $C_4$} We first note that $e$ and $w$ survive.

\begin{lemma}\label{lem:ew}
$e$ and $w$ are permanent cycles, so every $E_r$ with $r\ge6$ is an $\F_2[e,w]$-module and every $d_r$ with $r\ge7$ is $\F_2[e,w]$-linear.
\end{lemma}
\begin{proof}
$e=Nx_\sigma a_\lambda$ is a product of permanent cycles, and $w$ is a $d_3$- and a $d_5$-cycle by Proposition \ref{prop:d3atC4} and Lemma \ref{lem:wlambda}. For $r\ge7$ the target of $d_r(w)$ lies in stem $3$, where $E_6$ vanishes by Theorem \ref{thm:E6C4} and Proposition \ref{prop:O6}.
\end{proof}

\begin{lemma}\label{lem:imd7}
$d_7(E_6)\subseteq\mathcal{T}_6$.
\end{lemma}
\begin{proof}
$d_7$ raises the slice degree by $6$, hence exchanges the two residues of the slice degree modulo $4$. A class with $t\equiv0$ therefore has its image in slice degree $t\equiv2$, which consists of transfers, and a class with $t\equiv2$ is itself a transfer, so its image is one by naturality.
\end{proof}

Consequently $d_7$ induces an $\F_2[e,w]$-linear map
\[\overline{d}_7\colon E_6/\mathcal{T}_6\longrightarrow \mathcal{T}_6/d_7(\mathcal{T}_6),\]
which is linear because $e$ and $w$ are permanent cycles (Lemma \ref{lem:ew}) and $\F_2[e,w]\,\mathcal{T}_6\subseteq\mathcal{T}_6$.

\begin{proposition}\label{prop:d7nu}
For $c\in\mathcal{A}\langle 1\rangle$,
\[\nu\,c=\mathrm{tr}(c\,x\,\gamma y),\qquad \nu\,c\,\zeta_1^2=\mathrm{tr}(c\,x\,\gamma z),\qquad \nu\,c\,\zeta_2=\mathrm{tr}(c\,z\,\gamma y),\]
and consequently
\[d_7(\nu c)=\mathrm{tr}\bigl(d_7(c)\,x\,\gamma y\bigr),\quad
d_7(\nu c\zeta_1^2)=\mathrm{tr}\bigl(d_7(c)\,x\,\gamma z\bigr),\quad
d_7(\nu c\zeta_2)=\mathrm{tr}\bigl(d_7(c)\,z\,\gamma y\bigr).\]
\end{proposition}
\begin{proof}
The first identity follows from $\nu=\mathrm{tr}(x\gamma y)$. For the second, $z=\zeta_1^2y+\zeta_2x$ gives $\zeta_1^2\gamma y=\gamma z+\zeta_2\gamma x$, so
\[\nu c\zeta_1^2=\mathrm{tr}(c\,x\,\zeta_1^2\gamma y)=\mathrm{tr}(c\,x\,\gamma z)+\zeta_2\,\mathrm{tr}(c\,x\,\gamma x),\]
and $\mathrm{tr}(x\gamma x)=\mathrm{tr}(\Res(e))=0$. The third is the same argument with $\zeta_2x=z+\zeta_1^2y$ and $\mathrm{tr}(c\,y\,\gamma y)=0$. In each case the argument of the transfer is a product of $c$ with permanent cycles of the $C_2$-page, and the differentials follow from naturality of the transfer and the Leibniz rule.
\end{proof}

\begin{proposition}\label{prop:d7bar}
$\overline{d}_7=0$.
\end{proposition}
\begin{proof}
By $\F_2[e,w]$-linearity it is enough to check the ten module generators of Theorem \ref{thm:E6C4}(1), namely $1$ at $(0,0)$, $g$ at $(8,-4)$, and the eight classes $\eta'c$. Their targets lie in the nine bidegrees
\[(-1,7),\ (7,3),\ (0,10),\ (4,6),\ (6,4),\ (10,0),\ (11,-1),\ (13,-3),\ (17,-7),\]
all of slice degree $t\equiv2\bmod4$, hence all inside the free module of Proposition \ref{prop:O6}. A class $e^kw^n\mathrm{tr}(c\,xyz)$ there has bidegree $(7+|c|+2k+4n,\,-1-|c|+2k)$, and solving for each of the first six bidegrees above gives no solution, so $E_6$ vanishes there. At the last three there is exactly one solution each, and the groups are spanned by
\[w\,\mathrm{tr}(xyz),\qquad e\,\mathrm{tr}(\zeta_1^4xyz),\qquad w\,\mathrm{tr}(\zeta_2^2xyz)\]
respectively. Each is $d_7$ of a class of $\mathcal{T}_6$ in the corresponding source bidegree: by Proposition \ref{prop:d7nu} and Proposition \ref{prop:C2ZSS},
\[d_7(\nu\zeta_2^3)=w\,\mathrm{tr}(xyz),\qquad d_7(\nu\zeta_1^4\zeta_3)=e\,\mathrm{tr}(\zeta_1^4xyz),\qquad d_7(\nu\zeta_1^2\zeta_2^2\zeta_3)=w\,\mathrm{tr}(\zeta_2^2xyz).\qedhere\]
\end{proof}

\begin{sseqdata}[name = Esix, hsss, x range = {0}{20}, y range = {-17}{21}, y tick step = 4, differentials = {dseven}]
\class(0,0)
\class(1,3)
\class(2,2)
\class(4,0)
\class(4,4)
\class(5,-1)
\class(6,2)
\class(6,6)
\class(7,-3)
\class(7,-3)
\class(7,-1)
\class(8,-4)
\class(8,-4)
\class(8,0)
\class(8,4)
\class(8,8)
\class(9,-5)
\class(9,-1)
\class(9,1)
\class(10,-6)
\class(10,-2)
\class(10,2)
\class(10,6)
\class(10,10)
\class(11,-7)
\class(11,-7)
\class(11,-5)
\class(11,-3)
\class(11,-1)
\class(11,1)
\class(11,3)
\class(12,-8)
\class(12,-8)
\class(12,-4)
\class(12,-4)
\class(12,0)
\class(12,0)
\class(12,4)
\class(12,8)
\class(12,12)
\class(13,-9)
\class(13,-9)
\class(13,-7)
\class(13,-5)
\class(13,-3)
\class(13,-1)
\class(13,1)
\class(13,3)
\class(13,5)
\class(14,-10)
\class(14,-10)
\class(14,-8)
\class(14,-6)
\class(14,-2)
\class(14,-2)
\class(14,2)
\class(14,2)
\class(14,6)
\class(14,10)
\class(14,14)
\class(15,-11)
\class(15,-7)
\class(15,-7)
\class(15,-7)
\class(15,-5)
\class(15,-5)
\class(15,-3)
\class(15,-1)
\class(15,-1)
\class(15,1)
\class(15,3)
\class(15,5)
\class(15,7)
\class(16,-12)
\class(16,-12)
\class(16,-8)
\class(16,-8)
\class(16,-6)
\class(16,-4)
\class(16,-4)
\class(16,0)
\class(16,0)
\class(16,0)
\class(16,4)
\class(16,4)
\class(16,8)
\class(16,12)
\class(16,16)
\class(17,-13)
\class(17,-11)
\class(17,-9)
\class(17,-9)
\class(17,-7)
\class(17,-5)
\class(17,-5)
\class(17,-5)
\class(17,-3)
\class(17,-3)
\class(17,-1)
\class(17,1)
\class(17,1)
\class(17,3)
\class(17,5)
\class(17,7)
\class(17,9)
\class(18,-14)
\class(18,-14)
\class(18,-12)
\class(18,-10)
\class(18,-10)
\class(18,-8)
\class(18,-6)
\class(18,-6)
\class(18,-4)
\class(18,-2)
\class(18,-2)
\class(18,2)
\class(18,2)
\class(18,2)
\class(18,6)
\class(18,6)
\class(18,10)
\class(18,14)
\class(18,18)
\class(19,-15)
\class(19,-11)
\class(19,-11)
\class(19,-9)
\class(19,-7)
\class(19,-7)
\class(19,-7)
\class(19,-5)
\class(19,-5)
\class(19,-3)
\class(19,-3)
\class(19,-3)
\class(19,-1)
\class(19,-1)
\class(19,-1)
\class(19,1)
\class(19,3)
\class(19,3)
\class(19,5)
\class(19,7)
\class(19,9)
\class(19,11)
\class(20,-16)
\class(20,-14)
\class(20,-12)
\class(20,-12)
\class(20,-10)
\class(20,-8)
\class(20,-8)
\class(20,-8)
\class(20,-6)
\class(20,-4)
\class(20,-4)
\class(20,-4)
\class(20,-2)
\class(20,0)
\class(20,0)
\class(20,0)
\class(20,4)
\class(20,4)
\class(20,4)
\class(20,8)
\class(20,8)
\class(20,12)
\class(20,16)
\class(20,20)
\d7(10,-6,1,1)
\d7(12,-8,1,1)
\d7(12,-4,1,1)
\d7(14,-10,1,1)
\d7(14,-6,1,1)
\d7(14,-2,1,1)
\d7(16,-12,1,1)
\d7(16,-12,2,2)
\d7(16,-8,1,1)
\d7(16,-8,2,2)
\d7(16,-4,1,1)
\d7(16,0,1,1)
\d7(17,-13,1,1)
\d7(18,-14,1,1)
\d7(18,-12,1,1)
\d7(18,-10,1,1)
\d7(18,-10,2,2)
\d7(18,-6,1,1)
\d7(18,-6,2,2)
\d7(18,-2,1,1)
\d7(18,2,1,1)
\d7(19,-15,1,1)
\d7(19,-11,1,1)
\d7(20,-16,1,1)
\d7(20,-14,1,1)
\d7(20,-12,1,1)
\d7(20,-12,2,2)
\d7(20,-10,1,1)
\d7(20,-8,1,1)
\d7(20,-8,2,2)
\d7(20,-8,3,3)
\d7(20,-4,1,1)
\d7(20,-4,2,2)
\d7(20,0,1,1)
\d7(20,4,1,1)

\classoptions["{\scriptscriptstyle 1}" left, label distance = -1pt](0,0,1)
\classoptions["{\scriptscriptstyle \eta^{\prime}}" left, label distance = -1pt](1,3,1)
\classoptions["{\scriptscriptstyle e}" left, label distance = -1pt](2,2,1)
\classoptions["{\scriptscriptstyle w}" above, label distance = -1pt](4,0,1)
\classoptions["{\scriptscriptstyle \eta^{\prime}\zeta_1^4}" left, label distance = -1pt](5,-1,1)
\classoptions["{\scriptscriptstyle \mathrm{tr}(xyz)}" below, label distance = -1pt](7,-1,1)
\classoptions["{\scriptscriptstyle g,\,\eta^{\prime}\zeta_3}" below left, label distance = -1pt](8,-4,1)
\end{sseqdata}

\begin{figure}[htbp]
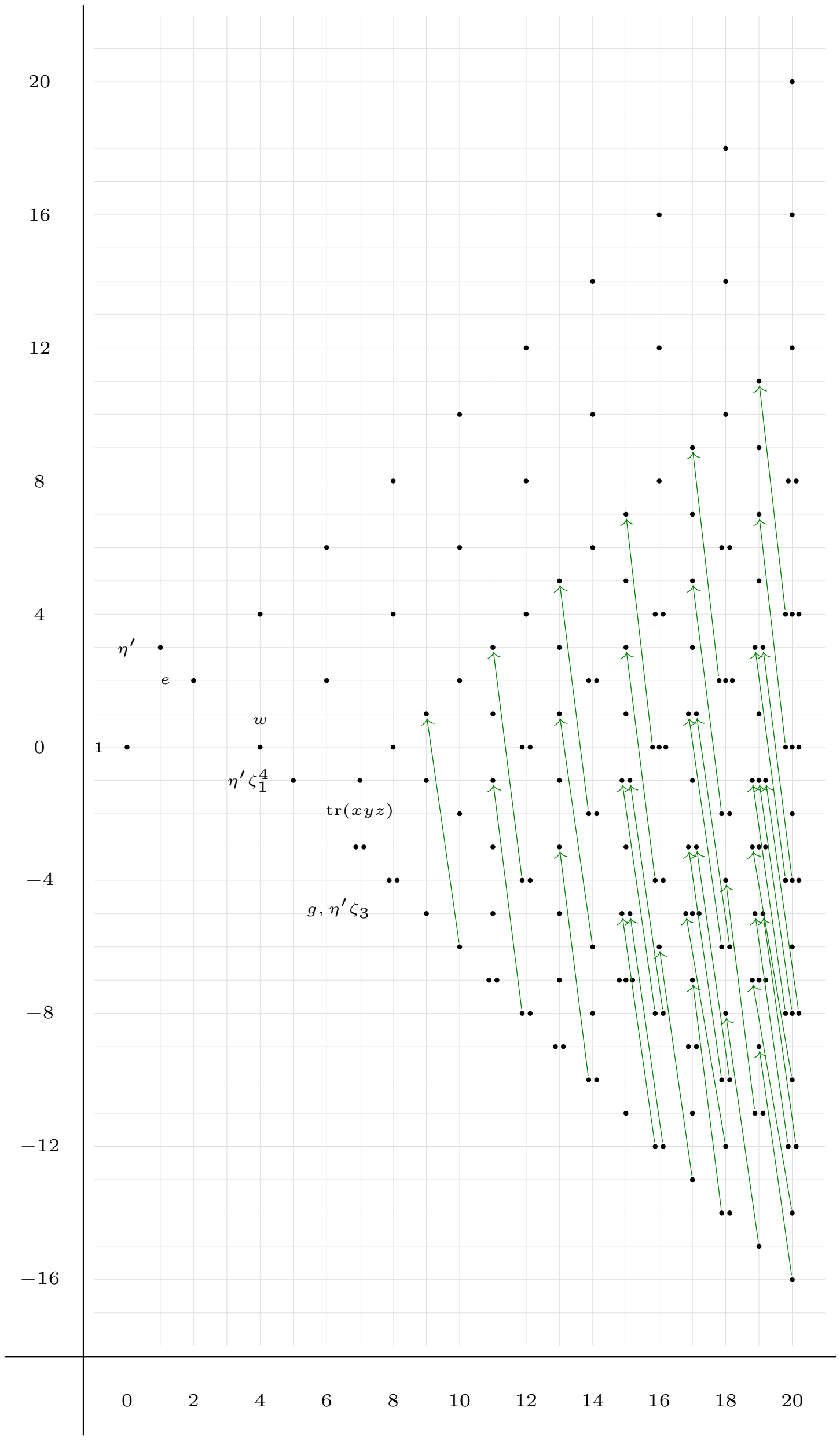
\centering
\printpage[name = Esix, page = 7]
\caption{The page $E_6=E_7$ in stems $0$ to $20$, with the $d_7$-differentials.}
\label{fig:E6d7}
\end{figure}

\subsection{The $E_8$-page at $C_4$} We can now assemble the page from the two halves. Lemma \ref{lem:imd7} says $d_7$ never leaves $\mathcal{T}_6$ and Proposition \ref{prop:d7bar} says it is zero on the quotient, so all of the new information sits in $\mathcal{T}_6$, where it is a $C_2$-level computation.

\begin{theorem}\label{thm:E8C4}
$E_8=E_9$ sits in a short exact sequence
\[0\longrightarrow\mathcal{T}_8\longrightarrow E_8\longrightarrow E_6/\mathcal{T}_6\longrightarrow0,\]
with $E_6/\mathcal{T}_6$ as in Theorem \ref{thm:E6C4}, and $\mathcal{T}_8=H(\mathcal{T}_6,d_7)$, where $\mathcal{T}_6=\mathrm{tr}(E_6^{C_2})$ and $d_7$ is the transfer of the $C_2$-level $d_7$.
\end{theorem}
\begin{proof}
By Lemma \ref{lem:imd7} the image of $d_7$ lies in $\mathcal{T}_6$, so
\[E_8/\mathcal{T}_8=\ker d_7\big/\bigl(\ker d_7\cap\mathcal{T}_6+\operatorname{im}d_7\bigr)=\ker\overline{d}_7,\]
which is all of $E_6/\mathcal{T}_6$ by Proposition \ref{prop:d7bar}. Again because $\operatorname{im}d_7\subseteq\mathcal{T}_6$, we have $\mathcal{T}_6\cap\operatorname{im}d_7=d_7(\mathcal{T}_6)$, so $\mathcal{T}_8=H(\mathcal{T}_6,d_7)$. The rest is Theorem \ref{thm:E6C4}(2) and naturality of the transfer.
\end{proof}

\begin{corollary}\label{cor:fg}
$E_8$ is a finitely generated $\F_2[e,w]$-module, with minimal generators in $18$ bidegrees, all of stem at most $21$.
\end{corollary}
\begin{proof}
Finite generation follows from Theorem \ref{thm:E8C4}, as $\mathcal{T}_6$ is a quotient of the finitely generated $\F_2[x^2,y^2]$-module $E_6^{C_2}$, on which $x^2$ and $y^2$ act through $e$ and $w$. The generators are found by a direct computation from Theorems \ref{thm:E6C4} and \ref{thm:E8C4}.
\end{proof}

\section{The $d_9$-differential}\label{sec:d9}
By Lemma \ref{lem:ew}, the differentials $d_r$ with $r\ge7$ are $\F_2[e,w]$-linear. Along the diagonal, the page carries the tower $e^k$ in bidegrees $(2k,2k)$, a free $\F_2[e]$-summand of $E_8/\mathcal{T}_8$. As we saw in Section \ref{sec:C4E2}, $e^4$ must be hit by a differential.

\begin{lemma}\label{lem:weld}
If $d_9(\nu\zeta_1^4)=e^3$ then $d_9(e^{k-3}\nu\zeta_1^4)=e^k$ for every $k\ge3$. Conversely, if $d_9$ is nonzero out of $(9,-1)$ then it is nonzero out of $(7,-3)$.
\end{lemma}
\begin{proof}
The first is $\F_2[e,w]$-linearity. For the second, $e\colon E_8^{7,-3}\to E_8^{9,-1}$ is surjective, so a nonzero value on $(9,-1)$ is $e$ times one on $(7,-3)$.
\end{proof}

\begin{theorem}\label{thm:bit}
The differential $d_9$ on $E_8$ is determined by $\F_2[e,w]$-linearity, naturality, the coaction and the value
\[d_9(\nu\zeta_1^4)=e^3,\qquad (7,-3)\longrightarrow(6,6).\]
\end{theorem}
\begin{proof}
We first show that $d_9$ is determined up to the single value $d_9(\nu\zeta_1^4)\in\{0,e^3\}$. On $E_8/\mathcal{T}_8$, which is spanned over $\F_2[e,w]$ by $1$, $g$ and the eight $\eta'c$ classes, the differential vanishes: $d_9(1)=0$, $d_9(g)$ lies in $(7,5)$ where $E_8$ vanishes, and $e\,d_9(\eta'c)=d_9(e\eta'c)=0$ because $e\eta'=d_5(\eta'\zeta_2)$ vanishes from $E_6$ on, while $e$ is injective on $E_8$ in the target bidegrees $(|c|,12-|c|)$, as one checks directly. On the image of the transfer $d_9$ vanishes, the $C_2$-level having collapsed.

What remains is the part of $\mathcal{T}_8$ spanned over $\F_2[e,w]$ by classes $\nu\,c\,b$ with $c\in\mathcal{A}\langle 1\rangle$ and $b\in\{1,\zeta_1^2,\zeta_2\}$. Every $c\ne1$ has a coaction term $c'\otimes\zeta_1^4$ with $c'$ of lower degree, so $\psi(\nu cb)$ involves $\nu\zeta_1^4b$ and classes of lower stem. As $d_9$ is a comodule map, its value on $\nu cb$ is determined by its values on classes of lower stem up to a comodule primitive in the target bidegree. By induction on the stem, all of these values are therefore determined by $d_9(\nu\zeta_1^4)$, which lies in $E_8^{6,6}=\F_2\{e^3\}$. A finite check shows that the primitives in the other target bidegrees are accounted for by $\F_2[e,w]$-linearity and Lemma \ref{lem:weld}.

It remains to show that $d_9(\nu\zeta_1^4)=e^3$. The class $e^4$ in $(8,8)$ must die. Differentials of even length vanish, so a source lies in $(9,8-r)$ with $r\ge9$ odd, and of these $E_8$ is nonzero only in $(9,-1)$ and $(9,-5)$. The second does not survive to $E_{13}$, as it supports a $d_9$ determined by the first part of the proof. So $e^4$ is hit by a $d_9$ from $(9,-1)$, and Lemma \ref{lem:weld} gives $d_9(\nu\zeta_1^4)=e^3$.
\end{proof}

\begin{sseqdata}[name = Eeight, hsss, x range = {0}{20}, y range = {-15}{21}, y tick step = 4, differentials = {dnine}]
\class(0,0)
\class(1,3)
\class(2,2)
\class(4,0)
\class(4,4)
\class(5,-1)
\class(6,2)
\class(6,6)
\class(7,-3)
\class(7,-3)
\class(7,-1)
\class(8,-4)
\class(8,-4)
\class(8,0)
\class(8,4)
\class(8,8)
\class(9,-5)
\class(9,-1)
\class(10,-2)
\class(10,2)
\class(10,6)
\class(10,10)
\class(11,-7)
\class(11,-7)
\class(11,-5)
\class(11,-3)
\class(11,1)
\class(12,-8)
\class(12,-4)
\class(12,0)
\class(12,0)
\class(12,4)
\class(12,8)
\class(12,12)
\class(13,-9)
\class(13,-9)
\class(13,-7)
\class(13,-5)
\class(13,-1)
\class(13,3)
\class(14,-10)
\class(14,-8)
\class(14,-2)
\class(14,2)
\class(14,2)
\class(14,6)
\class(14,10)
\class(14,14)
\class(15,-11)
\class(15,-7)
\class(15,-7)
\class(15,-7)
\class(15,-3)
\class(15,1)
\class(15,5)
\class(16,-4)
\class(16,0)
\class(16,0)
\class(16,4)
\class(16,4)
\class(16,8)
\class(16,12)
\class(16,16)
\class(17,-11)
\class(17,-9)
\class(17,-9)
\class(17,-5)
\class(17,-5)
\class(17,-1)
\class(17,3)
\class(17,7)
\class(18,-14)
\class(18,-2)
\class(18,2)
\class(18,2)
\class(18,6)
\class(18,6)
\class(18,10)
\class(18,14)
\class(18,18)
\class(19,-11)
\class(19,-7)
\class(19,-7)
\class(19,-3)
\class(19,-3)
\class(19,1)
\class(19,5)
\class(19,9)
\class(20,-4)
\class(20,0)
\class(20,0)
\class(20,4)
\class(20,4)
\class(20,8)
\class(20,8)
\class(20,12)
\class(20,16)
\class(20,20)
\d9(7,-3,1,1)
\d9(9,-5,1,1)
\d9(9,-1,1,1)
\d9(11,-7,1,1)
\d9(11,-3,1,1)
\d9(11,1,1,1)
\d9(13,-9,1,1)
\d9(13,-9,2,2)
\d9(13,-5,1,1)
\d9(13,-1,1,1)
\d9(13,3,1,1)
\d9(15,-11,1,1)
\d9(15,-7,1,1)
\d9(15,-7,2,2)
\d9(15,-3,1,1)
\d9(15,1,1,1)
\d9(15,5,1,1)
\d9(17,-9,1,1)
\d9(17,-9,2,2)
\d9(17,-5,1,1)
\d9(17,-5,2,2)
\d9(17,-1,1,1)
\d9(17,3,1,1)
\d9(17,7,1,1)
\d9(19,-11,1,1)
\d9(19,-7,1,1)
\d9(19,-7,2,2)
\d9(19,-3,1,1)
\d9(19,-3,2,2)
\d9(19,1,1,1)
\d9(19,5,1,1)
\d9(19,9,1,1)

\classoptions["{\scriptscriptstyle 1}" left, label distance = -1pt](0,0,1)
\classoptions["{\scriptscriptstyle \eta^{\prime}}" left, label distance = -1pt](1,3,1)
\classoptions["{\scriptscriptstyle e}" left, label distance = -1pt](2,2,1)
\classoptions["{\scriptscriptstyle w}" above, label distance = -1pt](4,0,1)
\classoptions["{\scriptscriptstyle \eta^{\prime}\zeta_1^4}" left, label distance = -1pt](5,-1,1)
\classoptions["{\scriptscriptstyle \mathrm{tr}(xyz)}" below, label distance = -1pt](7,-1,1)
\classoptions["{\scriptscriptstyle g,\,\eta^{\prime}\zeta_3}" below left, label distance = -1pt](8,-4,1)
\classoptions["{\scriptscriptstyle \eta^{\prime}\zeta_1^4\zeta_3}" left, label distance = -1pt](12,-8,1)
\classoptions["{\scriptscriptstyle \eta^{\prime}\zeta_2^2\zeta_3}" left, label distance = -1pt](14,-10,1)
\classoptions["{\scriptscriptstyle \eta^{\prime}\zeta_1^4\zeta_2^2\zeta_3}" left, label distance = -1pt](18,-14,1)
\end{sseqdata}

\begin{figure}[htbp]
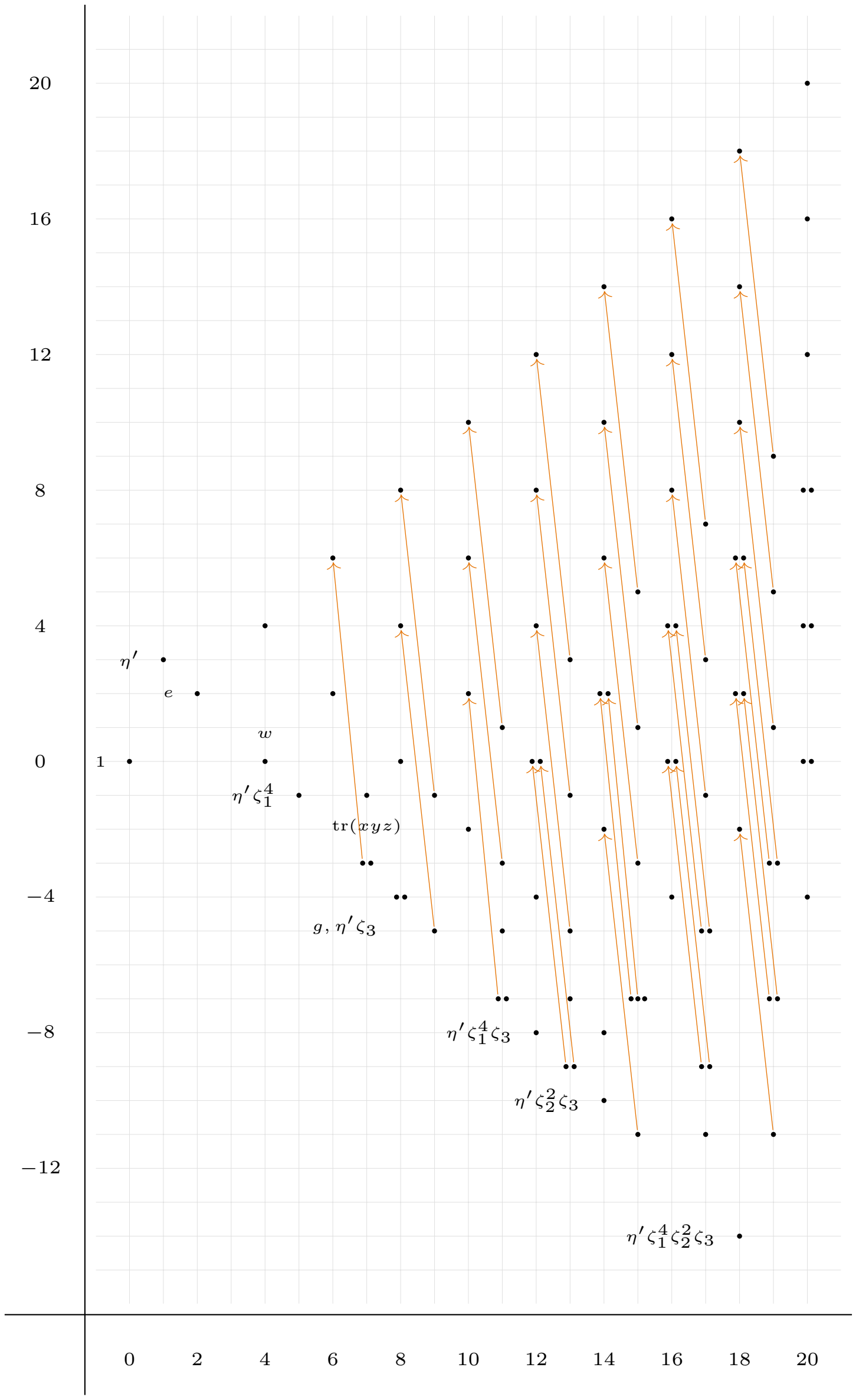
\centering
\printpage[name = Eeight, page = 9]
\caption{The page $E_8=E_9$ in stems $0$ to $20$, with the $d_9$-differentials. The tower $e^k$ on the diagonal is killed by $d_9(e^{k-3}\nu\zeta_1^4)=e^k$ for $k\ge3$.}
\label{fig:E8d9}
\end{figure}

\section{Collapse}\label{sec:collapse}
By Theorem \ref{thm:bit}, $e^3=0$ in $E_{10}$. If $e$ also acts injectively on $E_8$ in high stems, then $E_{10}$ in high stems is controlled by the homology of $E_8/e$, which is $w$-periodic and acyclic.

\subsection{Modules over $\F_2[e,w]$} The remaining work is to show that $E_{10}$ vanishes in all stems above $18$, and since $E_8$ is nonzero in infinitely many stems, this cannot be done stem by stem. Instead we use that $E_8$ is a finitely generated module over $R=\F_2[e,w]$ (Corollary \ref{cor:fg}) and that $d_9$ is $R$-linear.

The dimensions of a finitely generated $R$-module are governed by its Hilbert series, which is a rational function with denominator $(1-x^2y^2)(1-x^4)$; this gives a linear recursion for the dimensions of $E_8$ in high stems. Combined with the $w$-periodicity of $E_8/e$, the recursion shows that $e$ acts injectively on $E_8$ in high stems, and then the vanishing of $E_{10}$ there reduces to a computation in finitely many stems.

Let $R=\F_2[e,w]$, bigraded by $|e|=(2,2)$ and $|w|=(4,0)$. A finitely generated bigraded $R$-module $M$ is finite-dimensional in each bidegree, and we write
\[h_M(x,y)=\sum_{s,f}\dim M^{s,f}\,x^sy^f\]
for its \emph{Hilbert series}.

\begin{proposition}\label{prop:Rmodules}
Let $M$ be a finitely generated bigraded $R$-module.
\begin{enumerate}
\item If $M$ has a finite resolution $0\to F_m\to\cdots\to F_0\to M\to0$ by finitely generated free $R$-modules on generators in stems $s\le N$, then for every $(s,f)$ with $s>N$
\begin{equation}\label{eq:rec}
\dim M^{s,f}-\dim M^{s-2,f-2}-\dim M^{s-4,f}+\dim M^{s-6,f-2}=0 .
\end{equation}
\item If $M$ is generated in stems $\le N$ and $M^{s,*}=0$ for $N<s\le N+4$, then $M^{s,*}=0$ for every $s>N$.
\item Let $d$ be an $R$-linear differential on $M$ of bidegree $(-1,r)$, and suppose that $e^n$ acts by zero on $H(M,d)$. If $e\colon M^{s,*}\to M^{s+2,*}$ is injective for every $s\ge N$, and $H(M/e,d)$ vanishes in stems $\ge N'$, then $H(M,d)$ vanishes in stems $\ge\max(N+3,N')+2n-2$.
\end{enumerate}
\end{proposition}
\begin{proof}
(1) A free $R$-module on generators in bidegrees $(s_i,f_i)$ has Hilbert series $\sum_ix^{s_i}y^{f_i}/(1-x^2y^2)(1-x^4)$. Hilbert series are additive along exact sequences, so $h_M\cdot(1-x^2y^2)(1-x^4)$ is a Laurent polynomial supported in stems $\le N$, and its coefficient of $x^sy^f$ is the left-hand side of \eqref{eq:rec}.

(2) In a stem $s>N$ there are no generators, so $M^{s,*}=eM^{s-2,*}+wM^{s-4,*}$. For $s>N+4$ both terms vanish by induction on $s$.

(3) Let $z$ be a $d$-cycle in a stem $s\ge\max(N+3,N')$. Its image in $M/e$ is a cycle, hence a boundary, so $z-dy=eu$ for some $y$ and some $u$ in stem $s-2$. Then $e\,du=d(eu)=0$, and $du$ lies in stem $s-3\ge N$, so $du=0$ by injectivity and $[z]=e[u]$. Repeating the argument with $u$ in place of $z$, $n$ times in all, gives $[z]\in e^nH(M,d)=0$.
\end{proof}

\subsection{The page $E_8$ in high stems} We apply Proposition \ref{prop:Rmodules} to $M=E_8$, with $d=d_9$.

\begin{lemma}\label{lem:window}
$w\colon(E_8/e)^{s-4,f}\to(E_8/e)^{s,f}$ is surjective for every $s\ge22$.
\end{lemma}
\begin{proof}
This follows from Proposition \ref{prop:Rmodules}(2) applied to $M=E_8/(e,w)$, which is generated in stems $\le21$ by Corollary \ref{cor:fg} and vanishes in stems $22$ to $25$ by direct computation.
\end{proof}

The next proposition combines the bound on the dimensions of $E_8/e$ from above by \eqref{eq:rec} with the bound from below of Lemma \ref{lem:window}.

\begin{proposition}\label{prop:period}
The page $E_8$ is periodic in high stems, in the following sense.
\begin{enumerate}
\item $e\colon E_8^{s,f}\to E_8^{s+2,f+2}$ is injective for every $s\ge19$.
\item $w\colon(E_8/e)^{s,f}\to(E_8/e)^{s+4,f}$ is an isomorphism for every $s\ge20$.
\item Above stem $19$, $E_8/e$ is $4$-periodic of rank $4$, supported in the bidegrees $(4k,0)$, $(4k,-4)$, $(4k+1,-9)$, $(4k+1,-13)$.
\end{enumerate}
\end{proposition}
\begin{proof}
From the generators of Corollary \ref{cor:fg} one computes a resolution of $E_8$ as in Proposition \ref{prop:Rmodules}(1), with all generators in stems $\le24$, so \eqref{eq:rec} holds for $E_8$ in every stem $\ge25$.

Writing $\alpha(s,f)=\dim E_8^{s,f}$ and $\beta(s,f)=\dim(E_8/e)^{s,f}$, one has
\[\beta(t,f)=\alpha(t,f)-\mathrm{rank}\bigl(e\colon E_8^{t-2,f-2}\to E_8^{t,f}\bigr)\ \ge\ \alpha(t,f)-\alpha(t-2,f-2),\]
with equality if and only if $e$ is injective on $E_8^{t-2,f-2}$. We induct on $t\ge25$, assuming that $e$ is injective on all of stem $t-6$. Then
\[\beta(t,f)\le\beta(t-4,f)=\alpha(t-4,f)-\alpha(t-6,f-2)=\alpha(t,f)-\alpha(t-2,f-2)\le\beta(t,f),\]
the first step by Lemma \ref{lem:window}, the second by the inductive hypothesis, and the third by \eqref{eq:rec}. So all are equalities. The last one says that $e$ is injective on stem $t-2$, and the first one says that the surjection $w\colon(E_8/e)^{t-4,f}\to(E_8/e)^{t,f}$ of Lemma \ref{lem:window} is between vector spaces of the same dimension, hence an isomorphism. The base cases $t=25,\dots,28$ need injectivity in stems $19$ to $22$, and this, as well as the case $s=20$ of (2), is a direct computation. Part (3) follows from (2) and the computation of $E_8/e$ in stems $20$ to $24$.
\end{proof}

\begin{proposition}\label{prop:acyclic}
$H(E_8/e,d_9)$ vanishes in every stem $\ge20$.
\end{proposition}
\begin{proof}
The $d_9$-cycle $w$ acts invertibly on $E_8/e$ above stem $19$ by Proposition \ref{prop:period}(2), so $H(E_8/e,d_9)$ is $4$-periodic there and it suffices to treat $20\le s\le23$. There $E_8/e$ is supported in the four families of Proposition \ref{prop:period}(3), and $d_9$ carries $(E_8/e)^{4k+1,-13}$ isomorphically onto $(E_8/e)^{4k,-4}$, and $(E_8/e)^{4k+1,-9}$ onto $(E_8/e)^{4k,0}$.
\end{proof}

\begin{theorem}\label{thm:vanish}
$E_{10}$ vanishes in every stem $>18$.
\end{theorem}
\begin{proof}
We apply Proposition \ref{prop:Rmodules}(3) to $M=E_8$ and $d=d_9$. Here $n=3$, since $e^3=d_9(\nu\zeta_1^4)$ by Theorem \ref{thm:bit}, and hence $e^3z=d_9(\nu\zeta_1^4z)$ for every $d_9$-cycle $z$; $N=19$ by Proposition \ref{prop:period}(1); and $N'=20$ by Proposition \ref{prop:acyclic}. This gives vanishing in stems $\ge26$, and stems $19$ to $25$ are a direct computation.
\end{proof}

Now that the page is confined to stems $\le18$, the following holds for degree reasons.

\begin{corollary}\label{cor:collapse}
$d_r=0$ for every $r\ge11$, so $E_{10}=E_\infty$.
\end{corollary}

The page $E_\infty=E_{10}$ is therefore $24$-dimensional, concentrated in stems $0$ to $18$, and zero in stems $3$ and $9$. Five of its classes have slice degree $t\equiv2\bmod4$: of the eight generators $\mathrm{tr}(c\,xyz)$ of Proposition \ref{prop:O6}, five survive. The other $19$ classes are the ten classes $1,e,w,e^2,ew,w^2,g,ge,gw,gw^2$ of $\F_2[e,w]\{1,g\}$, the eight classes $\eta'c$, and one class of the form $\nu d\,\zeta_1^2\zeta_2^2+\cdots$ in $(15,-7)$. Table \ref{tab:Einf} lists them.

\begin{sseqdata}[name = Einf, hsss, x range = {0}{19}, y range = {-14}{5}, y tick step = 2]
\class(0,0)
\class(1,3)
\class(2,2)
\class(4,0)
\class(4,4)
\class(5,-1)
\class(6,2)
\class(7,-3)
\class(7,-1)
\class(8,-4)
\class(8,-4)
\class(8,0)
\class(10,-2)
\class(11,-7)
\class(11,-5)
\class(12,-8)
\class(12,-4)
\class(13,-7)
\class(14,-10)
\class(14,-8)
\class(15,-7)
\class(16,-4)
\class(17,-11)
\class(18,-14)
\classoptions["\scriptscriptstyle m_{1}" right, label distance = -1pt](0,0,1)
\classoptions["\scriptscriptstyle m_{2}" right, label distance = -1pt](1,3,1)
\classoptions["\scriptscriptstyle m_{3}" right, label distance = -1pt](2,2,1)
\classoptions["\scriptscriptstyle m_{4}" right, label distance = -1pt](4,0,1)
\classoptions["\scriptscriptstyle m_{5}" right, label distance = -1pt](4,4,1)
\classoptions["\scriptscriptstyle m_{6}" right, label distance = -1pt](5,-1,1)
\classoptions["\scriptscriptstyle m_{7}" right, label distance = -1pt](6,2,1)
\classoptions["\scriptscriptstyle m_{8}" right, label distance = -1pt](7,-3,1)
\classoptions["\scriptscriptstyle m_{9}" right, label distance = -1pt](7,-1,1)
\classoptions["\scriptscriptstyle m_{10}" right, label distance = -1pt](8,-4,1)
\classoptions["\scriptscriptstyle m_{11}" below right, label distance = -1pt](8,-4,2)
\classoptions["\scriptscriptstyle m_{12}" right, label distance = -1pt](8,0,1)
\classoptions["\scriptscriptstyle m_{13}" right, label distance = -1pt](10,-2,1)
\classoptions["\scriptscriptstyle m_{14}" right, label distance = -1pt](11,-7,1)
\classoptions["\scriptscriptstyle m_{15}" right, label distance = -1pt](11,-5,1)
\classoptions["\scriptscriptstyle m_{16}" right, label distance = -1pt](12,-8,1)
\classoptions["\scriptscriptstyle m_{17}" right, label distance = -1pt](12,-4,1)
\classoptions["\scriptscriptstyle m_{18}" right, label distance = -1pt](13,-7,1)
\classoptions["\scriptscriptstyle m_{19}" right, label distance = -1pt](14,-10,1)
\classoptions["\scriptscriptstyle m_{20}" right, label distance = -1pt](14,-8,1)
\classoptions["\scriptscriptstyle m_{21}" right, label distance = -1pt](15,-7,1)
\classoptions["\scriptscriptstyle m_{22}" right, label distance = -1pt](16,-4,1)
\classoptions["\scriptscriptstyle m_{23}" right, label distance = -1pt](17,-11,1)
\classoptions["\scriptscriptstyle m_{24}" right, label distance = -1pt](18,-14,1)

\end{sseqdata}

\begin{figure}[htbp]
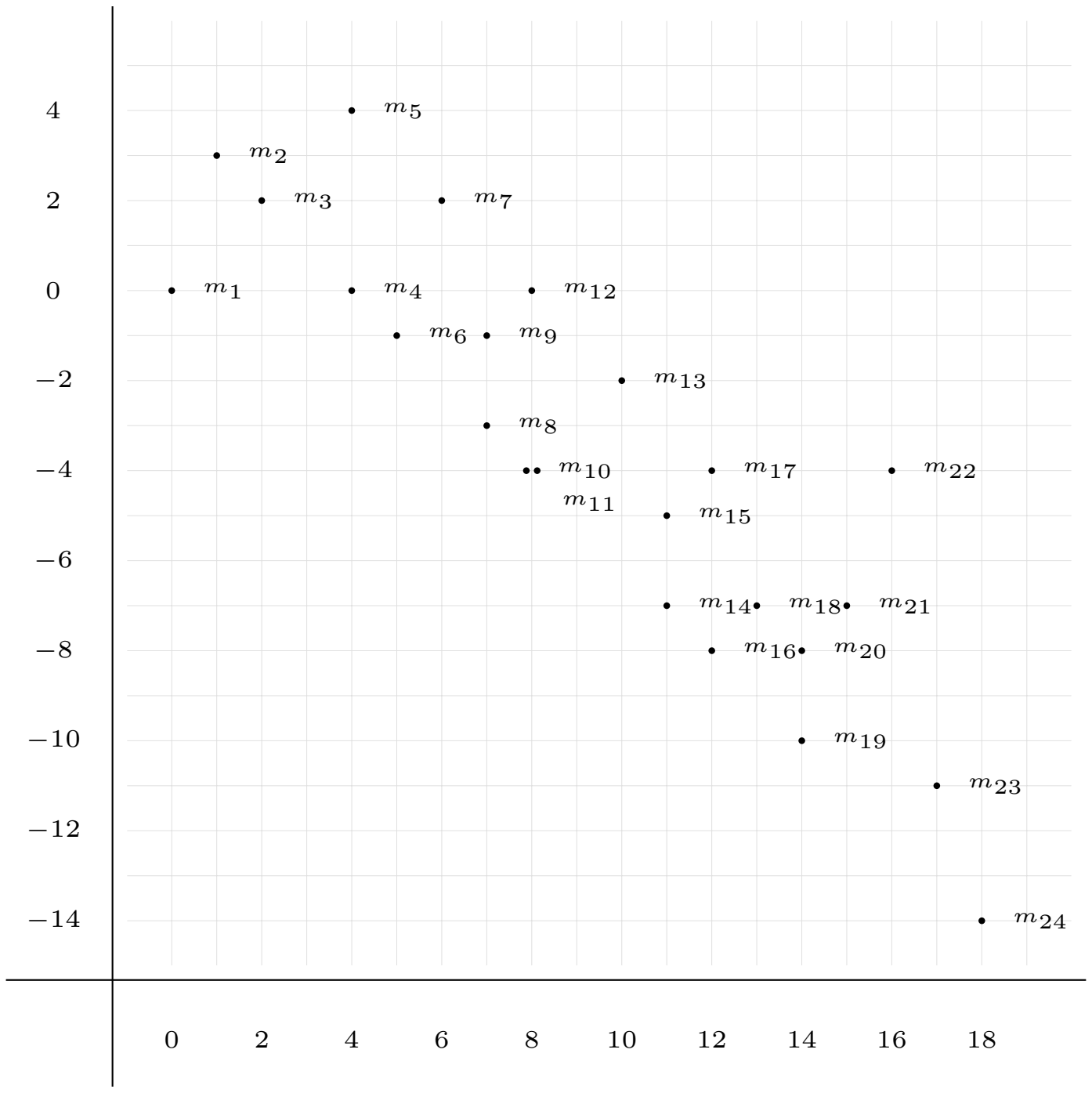
\centering
\printpage[name = Einf, page = 2]
\caption{The $E_\infty$-page at $C_4$, with its $24$ classes labelled as in Table \ref{tab:Einf}.}
\label{fig:Einf}
\end{figure}

\begin{table}[!ht]\centering\footnotesize
\caption{The $24$ classes of $E_\infty=E_{10}$, with their slice degrees $t$.}
\label{tab:Einf}
\renewcommand{\arraystretch}{1.15}
\begin{tabular}{@{}llll@{\hskip 2em}llll@{}}
\toprule
class & bidegree & $t$ & description & class & bidegree & $t$ & description\\
\midrule
$m_{1}$ & $(0,0)$ & $0$ & $1$                    & $m_{13}$ & $(10,-2)$ & $8$ & $ge$\\
$m_{2}$ & $(1,3)$ & $4$ & $\eta'$                & $m_{14}$ & $(11,-7)$ & $4$ & $\eta'\zeta_1^4\zeta_2^2$\\
$m_{3}$ & $(2,2)$ & $4$ & $e$                    & $m_{15}$ & $(11,-5)$ & $6$ & $\mathrm{tr}(\zeta_1^4xyz)$\\
$m_{4}$ & $(4,0)$ & $4$ & $w$                    & $m_{16}$ & $(12,-8)$ & $4$ & $\eta'\zeta_1^4\zeta_3$\\
$m_{5}$ & $(4,4)$ & $8$ & $e^2$                  & $m_{17}$ & $(12,-4)$ & $8$ & $gw$\\
$m_{6}$ & $(5,-1)$ & $4$ & $\eta'\zeta_1^4$      & $m_{18}$ & $(13,-7)$ & $6$ & $\mathrm{tr}(\zeta_2^2xyz)$\\
$m_{7}$ & $(6,2)$ & $8$ & $ew$                   & $m_{19}$ & $(14,-10)$ & $4$ & $\eta'\zeta_2^2\zeta_3$\\
$m_{8}$ & $(7,-3)$ & $4$ & $\eta'\zeta_2^2$      & $m_{20}$ & $(14,-8)$ & $6$ & $\mathrm{tr}(\zeta_3xyz)$\\
$m_{9}$ & $(7,-1)$ & $6$ & $\mathrm{tr}(xyz)$    & $m_{21}$ & $(15,-7)$ & $8$ & $\nu d\,\zeta_1^2\zeta_2^2+\cdots$\\
$m_{10}$ & $(8,-4)$ & $4$ & $g$                  & $m_{22}$ & $(16,-4)$ & $12$ & $gw^2$\\
$m_{11}$ & $(8,-4)$ & $4$ & $\eta'\zeta_3$       & $m_{23}$ & $(17,-11)$ & $6$ & $\mathrm{tr}(\zeta_1^4\zeta_2^2xyz)$\\
$m_{12}$ & $(8,0)$ & $8$ & $w^2$                 & $m_{24}$ & $(18,-14)$ & $4$ & $\eta'\zeta_1^4\zeta_2^2\zeta_3$\\
\bottomrule
\end{tabular}
\end{table}

\section{Extension problems and $\tmf_0(5)$}\label{sec:tmf05}

We now determine $H^{C_4}_*$ from its associated graded $E_\infty$ of Table \ref{tab:Einf}. Throughout this section $R=\inf(\F_2)\otimes_{\inf(\tmf)}\tmf_1(5)$ and $\pi_\bigstar=\pi^{C_4}_\bigstar R$, so that $\pi_n=H^{C_4}_n$, $\pi^{C_2}_nR=H^{C_2}_n$ and $\pi^e_nR=H^e_n$. We choose lifts $m_1,\dots,m_{24}\in H^{C_4}_*$ of the classes of Table \ref{tab:Einf}. The coaction and the products of $H^{C_4}_*$ then agree with those of $E_\infty$ up to terms of larger slice degree; for a product, these are the terms of higher filtration in its stem.

\begin{theorem}\label{thm:Mprime}
There are lifts $m_1,\dots,m_{24}$ with the following properties.
\begin{enumerate}
\item The coaction of $H^{C_4}_*$ is the coaction of $E_\infty$ plus the terms of Table \ref{tab:coaction}. Equivalently, the dual $\mathcal{A}(2)$-module is the one of Table \ref{tab:Mstar}.
\item Products in $H^{C_4}_*$ agree with those in $E_\infty$, except that
\[
m_3m_{19}=m_3m_{20}=m_4m_{16}=m_{10}m_{11}=m_{10}^2=m_{22}.
\]
\end{enumerate}
Moreover, any lifts satisfying (1) satisfy (2).
\end{theorem}

The proof occupies the rest of the section. Sections \ref{subsec:C2input} to \ref{subsec:lastext} determine the coaction, and Sections \ref{subsec:products} and \ref{subsec:beps} the products; the last product $m_{10}^2$ requires an integral computation.

\begin{table}[ht]\centering\footnotesize
\caption{The terms of the coaction of $H^{C_4}_*$ that raise the slice degree: the value of $\psi(m_i)$ minus the coaction of $m_i$ in $E_\infty$, for each $i$ for which it is nonzero.}
\label{tab:coaction}
\renewcommand{\arraystretch}{1.25}
\begin{tabular}{@{}ll@{}}
\toprule
class & terms\\
\midrule
$m_6$ & $\xi_1\otimes m_5$\\
$m_8$ & $\xi_2\otimes m_5$\\
$m_9$ & $\xi_1\otimes m_7+(\xi_1^3+\xi_2)\otimes m_5$\\
$m_{10}$, $m_{11}$ & $\xi_1^2\otimes m_7+\xi_1\xi_2\otimes m_5$\\
$m_{14}$ & $\xi_3\otimes m_5$\\
$m_{15}$ & $\xi_1\otimes m_{13}+\xi_1^5\otimes m_7+(\xi_3+\xi_1\xi_2^2)\otimes m_5$\\
$m_{16}$ & $\xi_2^2\otimes m_7+\xi_1\xi_3\otimes m_5$\\
$m_{18}$ & $\xi_2\otimes m_{13}+(\xi_3+\xi_1^4\xi_2)\otimes m_7+(\xi_2^3+\xi_1^2\xi_3)\otimes m_5$\\
$m_{19}$ & $\xi_1^4\otimes m_{13}+\xi_1^2\xi_2^2\otimes m_7+\xi_2\xi_3\otimes m_5$\\
$m_{20}$ & $\xi_1^2\otimes m_{17}+(\xi_1^4+\xi_1\xi_2)\otimes m_{13}+(\xi_1^6+\xi_2^2)\otimes m_{12}$\\
 & $\quad+(\xi_1\xi_3+\xi_1^2\xi_2^2+\xi_1^5\xi_2)\otimes m_7+(\xi_2\xi_3+\xi_1\xi_2^3+\xi_1^3\xi_3)\otimes m_5$\\
$m_{23}$ & $\xi_1^2\otimes m_{21}+\xi_3\otimes m_{13}+\xi_1^4\xi_3\otimes m_7+\xi_2^2\xi_3\otimes m_5$\\
$m_{24}$ & $\xi_1\otimes m_{23}+\xi_1^2\otimes m_{22}+(\xi_1^3+\xi_2)\otimes m_{21}+\xi_1^4\otimes m_{20}+\xi_1^5\otimes m_{18}+\xi_2^2\otimes m_{17}$\\
 & $\quad+(\xi_3+\xi_1\xi_2^2)\otimes m_{15}+\xi_1\xi_3\otimes m_{13}+\xi_1^4\xi_2^2\otimes m_{12}+(\xi_1^4\xi_3+\xi_1^5\xi_2^2)\otimes m_9$\\
 & $\quad+\xi_1^5\xi_3\otimes m_7+\xi_1\xi_2^2\xi_3\otimes m_5$\\
\bottomrule
\end{tabular}
\end{table}

\subsection{The $C_2$-level input}\label{subsec:C2input}
We first record the room for extensions in the products, and the values of restriction and transfer that we need.

\begin{lemma}\label{lem:room}
Let $2\le i\le j$. The product $m_im_j$ agrees with the product in $E_\infty$, unless $\{i,j\}$ is one of
\[
\{3,19\},\ \{3,20\},\ \{4,16\},\ \{6,14\},\ \{6,15\},\ \{10,10\},\ \{10,11\},\ \{11,11\}.
\]
In these eight cases the product in $E_\infty$ is zero, and $m_im_j\in\F_2\{m_{22}\}$.
\end{lemma}
\begin{proof}
The product $m_im_j$ differs from the product in $E_\infty$ by a class of higher filtration in stem $|m_i|+|m_j|$, that is, by a class of slice degree larger than $t_i+t_j$, where $t_i$ is the slice degree of $m_i$. Every class of positive degree has $t_i\ge4$, and by Table \ref{tab:Einf} the only class of slice degree larger than $8$ is $m_{22}$, in stem $16$ with $t=12$. The pairs in stem $16$ with $t_i+t_j<12$ are the eight listed, and their products in $E_\infty$ lie in bidegrees of stem $16$ and slice degree $8$ or $10$, where $E_\infty$ is zero.
\end{proof}

\begin{lemma}\label{lem:trvalues}
We have $\mathrm{tr}(y)=m_3$, $\mathrm{tr}(z^3)=m_{17}$ and $\mathrm{tr}(y^2z^3)=m_{22}$, and $\Res(m_{21})=\Res(m_{22})=0$.
\end{lemma}
\begin{proof}
We recall that $\Res\circ\mathrm{tr}=1+\gamma$ on $H^{C_2}_*$. By Proposition \ref{prop:C2gamma}, $(1+\gamma)(y)=x^2$ and $(1+\gamma)(z^3)=y^2z^2$. Since $H^{C_4}_2=\F_2\{m_3\}$ and $x^2\neq0$, we get $\mathrm{tr}(y)=m_3$. The class $z^3$ has filtration $-6$, and the only class of $H^{C_4}_{12}$ of filtration at least $-6$ is $m_{17}$; as $\Res\,\mathrm{tr}(z^3)=y^2z^2\ne0$, we get $\mathrm{tr}(z^3)=m_{17}$. On $E_\infty$ one has $\Res(w)=y\gamma y=y^2$, by Proposition \ref{prop:imtr}. So $\Res(m_4)=y^2+\epsilon\,x^2y$ for some $\epsilon\in\F_2$, and by Frobenius reciprocity
\[
\mathrm{tr}(y^2z^3)=\mathrm{tr}\bigl(\Res(m_4)\,z^3\bigr)=m_4\,\mathrm{tr}(z^3)=m_4m_{17}=m_{22},
\]
since $x^2yz^3=0$ and $m_4m_{17}=gw^2$ by Lemma \ref{lem:room}. Finally $\Res(m_{22})$ lies in filtration at least $-4$ of $H^{C_2}_{16}=\F_2\{y^2z^3\}$, and $y^2z^3$ has filtration $-6$. Every term of the representative $N^2a_\sigma x_\sigma u_\lambda^2\zeta_1^2\zeta_2^2+N^2a_\sigma x_\sigma a_\lambda u_\lambda\zeta_2\zeta_3$ of $m_{21}$ is divisible by $a_\sigma$, and $\Res(a_\sigma)=0$, so $\Res(m_{21})$ lies in filtration above $-7$ of $H^{C_2}_{15}=\F_2\{xyzq\}$, and $xyzq$ has filtration $-7$.
\end{proof}

\subsection{The coaction up to two coefficients}\label{subsec:four}
We write $\psi=\psi_\infty+\varphi$, where $\psi_\infty$ is the coaction of $E_\infty$ and $\varphi$ consists of the terms that raise the slice degree. A term $\xi^I\otimes m_j$ can occur in $\varphi(m_i)$ only if $m_j$ lies in stem $|m_i|-|\xi^I|$ and has larger slice degree than $m_i$; over all $i$, there are $151$ triples $(i,I,j)$ satisfying these conditions. Coassociativity of $\psi$ and of $\psi_\infty$ gives the equation
\[
(\Delta\otimes1)\varphi+(1\otimes\psi_\infty)\varphi+(1\otimes\varphi)\psi_\infty=(1\otimes\varphi)\varphi
\]
on $\varphi$. Changing the lift of $m_i$ by a class of the same stem and higher filtration changes $\varphi$, and there are seven such changes of lift, one for each pair $m_i,m_j$ in the same stem with $m_j$ of higher filtration.

To cut down the solutions we use the integral homology $H^{\tmf}_*(-;\Z_{(2)})$ of $\tmf_1(5)$, $\tmf_1(5)^{C_2}$ and $\tmf_0(5)$. These groups are the values of the Mackey functor $\underline{\pi}_n$ of the $C_4$-spectrum $\inf(H\Z)\otimes_{\inf(\tmf)}\tmf_1(5)$, and we write $H^{e}_n(\Z)$, $H^{C_2}_n(\Z)$ and $H^{C_4}_n(\Z)$ for them. This Mackey functor is a module over $\underline{\pi}_0\tmf_1(5)$, which is the constant Mackey functor $\underline{\Z}_{(2)}$ by Lemma \ref{lem:fournu}, and the same is true of $\underline{\pi}_n\tmf_1(5)$. So in both cases $\mathrm{tr}\circ\Res$ is multiplication by $2$, for $e\subset C_2$ and for $C_2\subset C_4$; and since $\Res\circ\mathrm{tr}$ is the sum over the Weyl group, the restriction to the trivial group is an isomorphism, after inverting $2$, onto the invariants of the Weyl group.

\begin{lemma}\label{lem:molien}
Let $L=H^{\tmf}_*(\tmf_1(5);\Z_{(2)})$.
\begin{enumerate}
\item $\pi_*\tmf_0(5)\otimes\mathbb{Q}$ is the ring of invariants of $\pi_*\tmf_1(5)\otimes\mathbb{Q}=\mathbb{Q}[a,b]$ under $\gamma\colon a\mapsto b\mapsto-a$.
\item $L$ is free of rank $24$ with $L/2=H^e_*$, and $L\otimes\mathbb{Q}=\mathbb{Q}[a,b]/(c_4,c_6)$.
\item $H^{C_2}_*(\Z)\otimes\mathbb{Q}$ and $H^{C_4}_*(\Z)\otimes\mathbb{Q}$ are the invariants of $L\otimes\mathbb{Q}$ under $\gamma^2$ and under $\gamma$. Both vanish in degrees not divisible by $4$, the former has dimension $1,3,4,3,1$ in degrees $0,4,8,12,16$, and the latter has dimension $1$ in degree $4$.
\end{enumerate}
\end{lemma}
\begin{proof}
Part (1) and the first claim of (3) were just explained, with the action of $\gamma$ given by Proposition \ref{prop:tate}. Since $H^e_*$ is concentrated in even degrees, $L$ is torsion-free, so it is free with $L/2=H^e_*$. The classes $c_4$ and $c_6$ of Proposition \ref{prop:tate} are coprime in $\mathbb{Q}[a,b]$, as a common factor would divide $c_4^3-c_6^2=1728\,a^5b^5(a^2-11ab-b^2)$, and none of $a$, $b$ and $a^2-11ab-b^2$ divides $c_4$. So $c_4,c_6$ is a regular sequence on $\mathbb{Q}[a,b]$, the K\"unneth spectral sequence $\Tor^{\pi_*\tmf}(\mathbb{Q},\pi_*\tmf_1(5))\Rightarrow L\otimes\mathbb{Q}$ collapses, and $L\otimes\mathbb{Q}=\mathbb{Q}[a,b]/(c_4,c_6)$, which has Poincar\'e series
\[(1-t^8)(1-t^{12})/(1-t^2)^2=1+2t^2+3t^4+4t^6+4t^8+4t^{10}+3t^{12}+2t^{14}+t^{16}.\]
This proves (2). Since $\gamma^2$ acts by $-1$ on $a$ and $b$, it acts by $(-1)^{n/2}$ on $(L\otimes\mathbb{Q})_n$, which gives the dimensions of the $\gamma^2$-invariants. The $\gamma$-invariants of $(L\otimes\mathbb{Q})_4=\mathbb{Q}\{a^2,ab,b^2\}$ are spanned by $a^2+b^2$.
\end{proof}

\begin{lemma}\label{lem:exponent}
Every torsion class of $H^{C_2}_*(\Z)$ is annihilated by $2$, and every torsion class of $H^{C_4}_*(\Z)$ by $4$. For $n\equiv0$ or $3\bmod4$, the group $H^{C_2}_n(\Z)$ is torsion-free and every torsion class of $H^{C_4}_n(\Z)$ is annihilated by $2$.
\end{lemma}
\begin{proof}
A torsion class $x\in H^{C_2}_n(\Z)$ restricts to a torsion class of $L$, hence to zero, so $2x=\mathrm{tr}(\Res(x))=0$. A torsion class $x\in H^{C_4}_n(\Z)$ restricts to a torsion class of $H^{C_2}_n(\Z)$, so $2\Res(x)=0$ and $4x=\mathrm{tr}(2\Res(x))=0$. For $n\equiv0\bmod4$, the dimension of $H^{C_2}_n$ is $1,3,4,3,1$ for $n=0,4,8,12,16$ and zero otherwise, by Theorem \ref{thm:C2abutment}, which is the rank of $H^{C_2}_n(\Z)$ by Lemma \ref{lem:molien}(3). By the universal coefficient theorem, the dimension of $H^{C_2}_n$ is this rank plus the number of cyclic summands in the torsion of $H^{C_2}_n(\Z)$ and of $H^{C_2}_{n-1}(\Z)$, so both groups are torsion-free. Hence $\Res(x)=0$ for a torsion class $x\in H^{C_4}_n(\Z)$ with $n\equiv0,3\bmod4$, and $2x=0$.
\end{proof}

We write $H(H^{C_4}_*;Q_0)$ for the Margolis homology of $Q_0=\Sq^1_*$. It is the $E_2$-page of the Bockstein spectral sequence of $H^{C_4}_*(\Z)$, so its dimension in degree $n$ is the rank of $H^{C_4}_n(\Z)$ plus the number of cyclic summands of order at least $4$ in the torsion of $H^{C_4}_n(\Z)$ and of $H^{C_4}_{n-1}(\Z)$.

\begin{corollary}\label{cor:Q0}
The Margolis homology $H(H^{C_4}_*;Q_0)$ is one-dimensional in degree $4$. In particular $\Sq^1_*m_6=m_5$.
\end{corollary}
\begin{proof}
The torsion of $H^{C_4}_4(\Z)$ and of $H^{C_4}_3(\Z)$ has exponent $2$ by Lemma \ref{lem:exponent}, so the dimension of $H(H^{C_4}_*;Q_0)_4$ is the rank of $H^{C_4}_4(\Z)$, which is $1$ by Lemma \ref{lem:molien}(3). As $H^{C_4}_3=0$, both $m_4$ and $m_5$ are $Q_0$-cycles, and the only possible boundary is $Q_0(m_6)$, which must therefore be nonzero. In $E_\infty$ the coaction of $m_6=\eta'\zeta_1^4$ has no term $\xi_1\otimes m_4$, and $m_5$ is the only class of degree $4$ of larger slice degree than $m_6$, so $\Sq^1_*m_6=m_5$.
\end{proof}

\begin{remark}
In degrees $14$ and $15$ the Margolis homology is one-dimensional, spanned by $m_{19}$ and $m_{21}$, while the rational homology vanishes by Lemma \ref{lem:molien}(3). So $H^{C_4}_{14}(\Z)$ has a cyclic summand of order at least $4$, which is $\Z/4$ by Lemma \ref{lem:exponent}. In the Adams spectral sequence of Section \ref{sec:adamsss}, this summand is the differential $d_2(x_{15,0})=x_{14,2}$ between the $h_0$-towers on $m_{21}$ and $m_{19}$.
\end{remark}

\begin{theorem}\label{thm:four}
There are lifts $m_1,\dots,m_{24}$ and coefficients $c,c'\in\F_2$ such that the coaction of $H^{C_4}_*$ is the coaction of $E_\infty$ plus the terms of Table \ref{tab:coaction}, with the following two changes: in $\psi(m_{24})$ the term $\xi_1^2\otimes m_{22}$ has coefficient $c$, and in $\psi(m_{23})$ and $\psi(m_{24})$ the terms $\xi^I\otimes m_{21}$ have coefficient $c'$.
\end{theorem}
\begin{proof}
The solutions $\varphi$ of the above equation form a finite set, which we have enumerated in Sage, solving the equation over $\F_2$ by increasing difference of slice degrees. There are $2304$ solutions, falling into $36$ classes under the changes of lift. The same computation gives $\Sq^1_*m_6$ for each solution; it is nonzero in exactly four of the $36$ classes, and zero in the other $32$. By Corollary \ref{cor:Q0}, the coaction of $H^{C_4}_*$ lies in one of the four, and these are the four values of $(c,c')$.
\end{proof}

Since $m_{21},\dots,m_{24}$ are the only classes in their stems, the two coefficients
\begin{equation}\label{eq:ABbits}
\Sq^2_*\,m_{24}=c\,m_{22},\qquad \Sq^2_*\,m_{23}=c'\,m_{21}
\end{equation}
do not depend on the choice of lifts. They are not visible in integer degrees. We show $c=c'=1$ using the $RO(C_4)$-graded degrees $n+\sigma$ and $n+\sigma-\lambda$.

\subsection{The coefficient $c$}\label{subsec:A}
Since $m_{21},\dots,m_{24}$ restrict to zero, they are divisible by $a_\sigma$. We use this to move into the degrees $n+\sigma$.

\begin{lemma}\label{lem:sigmadeg}
For every $n$ there is an exact sequence
\[
H^{C_2}_{n+1}\xrightarrow{\ \mathrm{tr}\ }\pi_{n+\sigma}\xrightarrow{\ a_\sigma\ }H^{C_4}_n\xrightarrow{\ \Res\ }H^{C_2}_n .
\]
Consequently there are unique classes $\tilde m_{22}$, $\tilde m_{23}$, $\tilde m_{24}$ with $a_\sigma\tilde m_i=m_i$, and $\pi_{n+\sigma}=\F_2\{\tilde m_{n+6}\}$ for $n=16,17,18$. There is a class $\tilde m_{21}\in\pi_{15+\sigma}$ with $a_\sigma\tilde m_{21}=m_{21}$, unique modulo $\ell=\mathrm{tr}(y^2z^3)$, and $\ell\neq0$ spans the kernel of $a_\sigma$ on $\pi_{15+\sigma}$.
\end{lemma}
\begin{proof}
The exact sequence comes from the cofiber sequence $C_4/C_{2+}\to S^0\to S^\sigma$, and continues to the left with $\Res\colon H^{C_4}_{n+1}\to H^{C_2}_{n+1}$. We have $H^{C_2}_n=0$ for $n\ge17$ and $H^{C_2}_{16}=\F_2\{y^2z^3\}$, and $\Res$ vanishes on $H^{C_4}_{15}$ and $H^{C_4}_{16}$ by Lemma \ref{lem:trvalues}. In particular the transfer $H^{C_2}_{16}\to\pi_{15+\sigma}$ is injective.
\end{proof}

The class $x_\sigma$ of Proposition \ref{prop:positiveconehomologyC4} is a permanent cycle (see the proof of Lemma \ref{lemma:permanentcyclesE2}), and we write $x_\sigma\in\pi_{1-\sigma}$ also for a class it detects. Then $\psi(x_\sigma)=1\otimes x_\sigma+\xi_1\otimes a_\sigma$, since in the degrees $k-\sigma$ with $k\le0$ the $E_2$-page consists of $a_\sigma$ alone.

\begin{theorem}\label{thm:Abit}
$\Sq^2_*\,m_{24}=m_{22}$, that is, $c=1$. Moreover $x_\sigma\tilde m_{23}=m_{24}$ and $x_\sigma\tilde m_{22}=m_{23}$.
\end{theorem}
\begin{proof}
The product $x_\sigma\tilde m_{24}$ lies in $\pi_{19}=H^{C_4}_{19}=0$. We write
\[
\psi(\tilde m_{24})=1\otimes\tilde m_{24}+\xi_1\otimes u_1+\xi_1^2\otimes u_2+\cdots,
\]
with $u_1\in\pi_{17+\sigma}$ and $u_2\in\pi_{16+\sigma}$. The coaction is multiplicative, so by the coaction of $x_\sigma$ the $\xi_1$- and $\xi_1^2$-components of $0=\psi(x_\sigma)\psi(\tilde m_{24})$ give
\[
x_\sigma u_1+a_\sigma\tilde m_{24}=0,\qquad x_\sigma u_2+a_\sigma u_1=0.
\]
The first gives $x_\sigma u_1=m_{24}\neq0$, so $u_1\neq0$, and $u_1=\tilde m_{23}$ by Lemma \ref{lem:sigmadeg}. The second now gives $x_\sigma u_2=m_{23}\neq0$, so $u_2=\tilde m_{22}$. Finally $a_\sigma$ is a comodule primitive, so $\Sq^2_*\,m_{24}=a_\sigma u_2=m_{22}$.
\end{proof}

The same argument, applied to $\tilde m_{23}$, shows that $\Sq^1_*\tilde m_{23}=0$ and $x_\sigma\,\Sq^2_*\tilde m_{23}=m_{22}$. Since $x_\sigma\ell=m_{22}$ by Frobenius reciprocity, this does not determine the $\tilde m_{21}$-component of $\Sq^2_*\tilde m_{23}$, which is $c'$.

\subsection{The coefficient $c'$}\label{subsec:lastext}
After multiplication by $a_\lambda$, the class $\tilde m_{23}$ becomes a transfer from $C_2$, where the coaction is given by Proposition \ref{prop:C2coactions}. We can therefore work out its coaction using the coaction on the $C_2$-level. We use $a_\lambda$-linearity then to prove the following.

\begin{theorem}\label{thm:Bbit}
$\Sq^2_*\,m_{23}=m_{21}$, that is, $c'=1$.
\end{theorem}

At $C_2$ we use the notation of Proposition \ref{prop:C2ROSS}, writing $a_{\sigma_2}$ for the class $a_\sigma$ there, so that $x_2=a_{\sigma_2}\zeta_2+x_1\zeta_1^2$ and $x_3=a_{\sigma_2}\zeta_3+x_1\zeta_2^2+x_2\zeta_1^4$.

\begin{proposition}\label{prop:Bdetector}
$a_\lambda\tilde m_{21}\neq0$, and $\Sq^2_*(a_\lambda\tilde m_{23})=c'\,a_\lambda\tilde m_{21}$.
\end{proposition}
\begin{proof}
Since $a_\sigma$ is a comodule primitive, \eqref{eq:ABbits} and Lemma \ref{lem:sigmadeg} give $\Sq^2_*\tilde m_{23}\in c'\tilde m_{21}+\F_2\{\ell\}$. The cofiber sequence $S(\lambda)_+\to S^0\to S^\lambda$ gives an exact sequence
\[
\pi_{15+\sigma}(S(\lambda)_+\otimes R)\longrightarrow\pi_{15+\sigma}\xrightarrow{\ a_\lambda\ }\pi_{15+\sigma-\lambda} .
\]
The circle $S(\lambda)$ has one free $C_4$-cell in each of the dimensions $0$ and $1$. Since $H^e_{15}=0$ and $H^e_{16}=\F_2\{\zeta_1^6\zeta_2^2t_1^2\}$, the first group is the coinvariants of $H^e_{16}$, and it maps onto the span of $\mathrm{tr}^{C_4}_e(\zeta_1^6\zeta_2^2t_1^2)=\mathrm{tr}(y^2z^3)=\ell$, by Proposition \ref{prop:C2res}. Hence the kernel of $a_\lambda$ is $\F_2\{\ell\}$, so $a_\lambda\tilde m_{21}\neq0$. As $a_\lambda$ is a comodule primitive, $\Sq^2_*(a_\lambda\tilde m_{23})=c'\,a_\lambda\tilde m_{21}$.
\end{proof}

\begin{lemma}\label{lem:freecircle}
The groups $\pi^{C_4}_{16+\sigma}(S(\lambda)_+\otimes R)$ and $\pi^{C_2}_{17}(S(\lambda)_+\otimes R)$ are both $\F_2$, and the transfer $\mathrm{tr}^{C_4}_{C_2}$ between them is zero.
\end{lemma}
\begin{proof}
The cell structure of $S(\lambda)$, with one free $C_4$-cell in each of the dimensions $0$ and $1$, gives a cofiber sequence $C_{4+}\otimes R\xrightarrow{1-\gamma}C_{4+}\otimes R\to S(\lambda)_+\otimes R$. For $G\subseteq C_4$ we have $\pi^G_V(C_{4+}\otimes R)=\Z[C_4/G]\otimes H^e_{|V|}$, with the transfer $\mathrm{tr}^{C_4}_{C_2}$ given by the augmentation $\Z[C_4/C_2]\to\Z$. Since $H^e_{17}=0$, both groups in question are the kernel of $1-\gamma$ on $\Z[C_4/G]\otimes H^e_{16}$, where $H^e_{16}=\F_2$ with trivial action. For $G=C_4$ this is $\F_2$, and for $G=C_2$ it is the span of $1+\gamma$ in $\F_2\{1,\gamma\}$, which the augmentation sends to $2=0$.
\end{proof}

We consider the classes
\[
\theta=\bar t_1^{\,2}x_2^3x_3\quad\text{and}\quad\theta'=\bar t_1^{\,2}x_1x_2^2x_3 ,
\]
of degrees $18-2\sigma_2$ and $16-2\sigma_2$ and slice degree $t=4$. These are the restrictions of the degrees $17+\sigma-\lambda$ and $15+\sigma-\lambda$, and $\theta,\theta'$ are permanent cycles by Proposition \ref{prop:C2ROSS}. Formally, $\theta=z^3q/\bar t_1^{\,2}$ and $\theta'=yz^2q/\bar t_1^{\,2}$, and the following lemma is the identity $\Sq^2_*(z^3q)=yz^2q$ in $H^{C_2}_*$, divided by $\bar t_1^{\,2}$.

\begin{lemma}\label{lem:Bcartan}
$\Sq^2_*\theta=\theta'$ on the $E_2$-page.
\end{lemma}
\begin{proof}
By Proposition \ref{prop:C2coactions} we have $\Sq^1_*x_2=\Sq^1_*x_3=0$, $\Sq^2_*x_2=x_1$ and $\Sq^2_*x_3=0$, and $\bar t_1$ is primitive. Applying the Cartan formula then gives the result.
\end{proof}

\begin{proposition}\label{prop:Btransfer}
Let $\bar\theta\in\pi^{C_2}_{18-2\sigma_2}R$ be detected by $\theta$. Then $\mathrm{tr}(\bar\theta)\in\{0,\ a_\lambda\tilde m_{23}\}$.
\end{proposition}
\begin{proof}
In the exact sequence
\[
\pi_{17+\sigma}\xrightarrow{\ a_\lambda\ }\pi_{17+\sigma-\lambda}\xrightarrow{\ \partial\ }\pi_{16+\sigma}(S(\lambda)_+\otimes R),
\]
we have $\partial\,\mathrm{tr}(\bar\theta)=\mathrm{tr}(\partial\bar\theta)=0$ by Lemma \ref{lem:freecircle}. So $\mathrm{tr}(\bar\theta)$ lies in $a_\lambda\pi_{17+\sigma}=\F_2\{a_\lambda\tilde m_{23}\}$, by Lemma \ref{lem:sigmadeg}.
\end{proof}

It remains to show that $\mathrm{tr}(\theta')$ survives to $E_\infty$ in degree $15+\sigma-\lambda$. It is a permanent cycle, being the transfer of one. Since slice degrees are even and nonnegative, and $d_r$ raises slice degree by $r-1$, the only differentials that could hit it are a $d_3$ from slice degree $t=2$ and a $d_5$ from slice degree $t=0$.

The degrees $n+\sigma-\lambda$ lie outside the positive cone of Proposition \ref{prop:positiveconehomologyC4}. Every slice of positive slice degree, norm or induced, shifts them back into the positive cone. The zero slice, however, does not.

\begin{lemma}\label{lem:mixedzero}
$\bigl(S^{\lambda-\sigma}\otimes H\underline{\Z}\bigr)^{C_4}\simeq H\Z/2$. Consequently the zero slice contributes $\mathcal A(2)_*=\mathcal{A}\langle0\rangle\{k,k'\}$ in degrees $n+\sigma-\lambda$. Here $|k|=\sigma-\lambda$, the class $k'$ represents $\zeta_1k$, and $\Res(k)=0$. All $d_3(kb)$ vanish, and $d_3(k')=\mathrm{tr}(a_{\sigma_2}^3\bar t_1)$.
\end{lemma}
\begin{proof}
There is a fiber sequence
\[
(S^{\lambda-\sigma}\otimes H\underline{\Z})^{C_4}\to(S^\lambda\otimes H\underline{\Z})^{C_4}\xrightarrow{\Res}(S^{2\sigma_2}\otimes H\underline{\Z})^{C_2}.
\]
The middle term has $\Z/4\{a_\lambda\}$ in degree $0$ and $\Z\{u_\lambda\}$ in degree $2$, and the right term has $\Z/2\{a_{\sigma_2}^2\}$ and $\Z\{u_{2\sigma_2}\}$. The restriction is the quotient map in degree $0$ and an isomorphism in degree $2$. Hence the fiber is $H\Z/2$, whose homology relative to $\tmf$ is $\mathcal A(2)_*$. Since $\Res(k)=0$, Frobenius reciprocity kills the product of $k$ with any transfer. The formula for $d_3(k')$ follows from $x_\sigma k'=u_\lambda$ and $d_3(u_\lambda)=\mathrm{tr}(a_{\sigma_2}^3\bar t_1)$ \cite[Theorem 11.13]{hhrc4}, \cite[Theorem 3.1]{hswx}.
\end{proof}

For a class $v$ of the $C_2$-level $E_2$-page in degree $16-2\sigma_2$ and slice degree $t=4$, we write $\phi(v)\in\F_2$ for the coefficient of $a_{\sigma_2}^3x_1\bar t_1^{\,2}\zeta_2^2\zeta_3$ in the polynomial $v|_{\gamma\bar t_1=\bar t_1}$ obtained from $v$ by substituting $\bar t_1$ for $\gamma\bar t_1$. Similarly, for a class $v$ of the $C_4$-level $E_2$-page in degree $15+\sigma-\lambda$ and slice degree $t=4$, we write $f(v)\in\F_2$ for the sum of the coefficients of the two basis elements $\mathrm{tr}(a_{\sigma_2}^3x_1\bar t_1^{\,2}\zeta_2^2\zeta_3)$ and $Na_\lambda x_\lambda\zeta_2^2\zeta_3$ in $v$. Both are linear in $v$.

\begin{lemma}\label{lem:specialize}
For every $v$ as above, $f(\mathrm{tr}(v))=\phi(v)$. Consequently $f$ vanishes on all $d_3$-boundaries.
\end{lemma}
\begin{proof}
Exactly three monomials of the $C_2$-level $E_2$-page become $a_{\sigma_2}^3x_1\bar t_1^{\,2}\zeta_2^2\zeta_3$ under the substitution $\gamma\bar t_1=\bar t_1$. Two of them, namely this monomial and its $\gamma$-conjugate, transfer to the first basis element, since $\mathrm{tr}\circ\gamma=\mathrm{tr}$. The third, $a_{\sigma_2}^3x_1\bar t_1\gamma\bar t_1\zeta_2^2\zeta_3$, transfers to the second basis element. Indeed, by Frobenius reciprocity and Proposition \ref{prop:positiveconehomologyC4},
\[
\mathrm{tr}(a_{\sigma_2}^3x_1\bar t_1\gamma\bar t_1)=N\,\mathrm{tr}(a_{\sigma_2}^3x_1)=Na_\lambda x_\lambda .
\]
All other monomials transfer into the span of the remaining basis elements.

For the second claim, note that classes of slice degree $t=2$ come from the induced slices on $\bar t_1$ and $\gamma\bar t_1$, so every incoming $d_3$ has the form $\mathrm{tr}(d_3v)$. By Proposition \ref{prop:C2ROSS}, $d_3$ at $C_2$ takes values in the ideal $(\bar t_1+\gamma\bar t_1)$, which the substitution $\gamma\bar t_1=\bar t_1$ kills.
\end{proof}

\begin{lemma}\label{lem:mixedd5}
In degree $16+\sigma-\lambda$, the slice degree $t=0$ part of $E_4$ has basis $k\zeta_2^3\zeta_3$ and $k\zeta_1^6\zeta_2\zeta_3$. Modulo $d_3$-boundaries,
\[
d_5(k\zeta_2^3\zeta_3)=0,\qquad
d_5(k\zeta_1^6\zeta_2\zeta_3)=\mathrm{tr}\bigl(\bar t_1^{\,2}(a_{\sigma_2}x_1^3\zeta_1^4\zeta_3+a_{\sigma_2}^3x_1\zeta_1^6\zeta_3+a_{\sigma_2}^4\zeta_1^4\zeta_2\zeta_3)\bigr).
\]
In particular $f$ vanishes on both values.
\end{lemma}
\begin{proof}
The two classes $k'\zeta_1^2\zeta_2^2\zeta_3$ and $k'\zeta_1^6\zeta_2^3$ have linearly independent $d_3$ by Lemma \ref{lem:mixedzero}. Multiplication by $Na_\sigma$ and by $Nx_\sigma$ carries degree $15+\sigma-\lambda$ to the degrees $16+\sigma$ and $17+\sigma$. These degrees are covered by the positive-cone model of Proposition \ref{prop:positiveconehomologyC4}, and $d_5$ there is computed by the same method as in Theorem \ref{thm:d5}. A direct computation shows that the pair $(Na_\sigma,Nx_\sigma)$ is jointly injective on the relevant $E_4$-group, which determines $d_5$. The last claim is Lemma \ref{lem:specialize}, since no specialized monomial on the right equals $a_{\sigma_2}^3x_1\bar t_1^{\,2}\zeta_2^2\zeta_3$.
\end{proof}

\begin{proposition}\label{prop:thetasurvives}
The class $\mathrm{tr}(\theta')$ is nonzero in $E_\infty$.
\end{proposition}
\begin{proof}
Expanding with $\zeta_1^8=\zeta_2^4=0$ gives
\[
\theta'=\bar t_1^{\,2}\bigl(x_1^4\zeta_1^4\zeta_2^2+a_{\sigma_2}x_1^3\zeta_1^4\zeta_3+a_{\sigma_2}^2x_1^2\zeta_1^6\zeta_2^2+a_{\sigma_2}^3x_1\zeta_2^2\zeta_3+a_{\sigma_2}^3x_1\zeta_1^4\zeta_2^3\bigr),
\]
so $f(\mathrm{tr}\,\theta')=\phi(\theta')=1$. By Lemmas \ref{lem:specialize} and \ref{lem:mixedd5}, $f$ vanishes on everything that can hit $\mathrm{tr}(\theta')$.
\end{proof}

\begin{proof}[proof of Theorem \ref{thm:Bbit}]
The spectral sequence is one of comodules and the transfer is a comodule map. By Lemma \ref{lem:Bcartan} the class $\Sq^2_*\,\mathrm{tr}(\bar\theta)=\mathrm{tr}(\Sq^2_*\bar\theta)$ is therefore detected by $\mathrm{tr}(\theta')$, which is nonzero by Proposition \ref{prop:thetasurvives}. Hence $\mathrm{tr}(\bar\theta)\neq0$, so $\mathrm{tr}(\bar\theta)=a_\lambda\tilde m_{23}$ by Proposition \ref{prop:Btransfer}. Proposition \ref{prop:Bdetector} now gives $c'\,a_\lambda\tilde m_{21}\neq0$, so $c'=1$.
\end{proof}

Theorems \ref{thm:four}, \ref{thm:Abit} and \ref{thm:Bbit} prove part (1) of Theorem \ref{thm:Mprime}. We record the result as a module over $\mathcal{A}(2)$. Let $m_i^*$ be the basis of $\Hom(H^{C_4}_*,\F_2)$ dual to the $m_i$, so that $\Sq^r(m_j^*)=\sum_i\langle\Sq^r_*m_i,m_j\rangle\,m_i^*$.

\begin{table}[ht]\centering\footnotesize
\caption{The $\mathcal{A}(2)$-module dual to $H^{C_4}_*$. The class $m_i^*$ lies in degree $|m_i|$. Listed are the nonzero values of $\Sq^1$, $\Sq^2$ and $\Sq^4$; all classes not listed are annihilated by all three.}
\label{tab:Mstar}
\renewcommand{\arraystretch}{1.15}
\begin{tabular}{@{}llll@{\hskip 2em}llll@{}}
\toprule
 & $\Sq^1$ & $\Sq^2$ & $\Sq^4$ & & $\Sq^1$ & $\Sq^2$ & $\Sq^4$\\
\midrule
$m_2^*$ & $m_3^*$ & & $m_6^*$ & $m_{12}^*$ & & & $m_{17}^*$\\
$m_3^*$ & & $m_4^*$ & & $m_{13}^*$ & $m_{15}^*$ & $m_{17}^*$ & $m_{19}^*+m_{20}^*$\\
$m_4^*$ & & & $m_{10}^*$ & $m_{14}^*$ & $m_{16}^*$ & & \\
$m_5^*$ & $m_6^*$ & $m_7^*$ & $m_{12}^*$ & $m_{15}^*$ & & $m_{18}^*$ & \\
$m_6^*$ & & $m_8^*$ & & $m_{16}^*$ & & $m_{19}^*$ & \\
$m_7^*$ & $m_9^*$ & $m_{10}^*+m_{11}^*$ & $m_{13}^*$ & $m_{17}^*$ & & $m_{20}^*$ & \\
$m_8^*$ & $m_{10}^*+m_{11}^*$ & & $m_{14}^*$ & $m_{18}^*$ & $m_{20}^*$ & & $m_{23}^*$\\
$m_9^*$ & & & $m_{15}^*$ & $m_{19}^*$, $m_{20}^*$ & & & $m_{24}^*$\\
$m_{11}^*$ & & & $m_{16}^*$ & $m_{21}^*$ & & $m_{23}^*$ & \\
 & & & & $m_{22}^*$ & & $m_{24}^*$ & \\
 & & & & $m_{23}^*$ & $m_{24}^*$ & & \\
\bottomrule
\end{tabular}
\end{table}

\subsection{Products}\label{subsec:products}
By Lemma \ref{lem:room}, the ring structure of $H^{C_4}_*$ is that of $E_\infty$ up to the coefficients of $m_{22}$ in eight products. Seven of them follow from the comodule algebra structure. The idea is that a product into stems $17$ and above vanishes, while the coaction of such a product involves products in stem $16$.

\begin{proposition}\label{prop:products}
For any lifts as in Theorem \ref{thm:Mprime}(1),
\[
m_3m_{19}=m_3m_{20}=m_4m_{16}=m_{10}m_{11}=m_{22},\qquad m_6m_{14}=m_6m_{15}=m_{11}^2=0.
\]
\end{proposition}
\begin{proof}
We use four products that agree with those of $E_\infty$, by Lemma \ref{lem:room}: $m_4m_{17}=m_{22}$, as $w\cdot gw=gw^2$; $m_5m_{16}=0$, as $e\eta'c=0$ from $E_6$ on; $m_2m_{21}=0$, as every term of the representative of $m_{21}$ is divisible by $a_\sigma u_\lambda$, and $\eta'a_\sigma u_\lambda=Na_\lambda\cdot a_\sigma^2u_\lambda=0$; and $m_7m_{13}=0$, as $E_\infty$ vanishes in bidegree $(16,0)$. In each case below we apply a component $\psi_K$ to a product that vanishes, either because it lies above stem $18$ or because it lies in stem $17$ or $18$ in higher filtration than $m_{23}$ or $m_{24}$. We have
\[
\psi_K(ab)=\sum_{I+J=K}\psi_I(a)\,\psi_J(b),
\]
and we list the nonzero terms, read off from the coaction of $E_\infty$ and Table \ref{tab:coaction}.
\begin{enumerate}
\item $\psi_{\xi_2^2}(m_4m_{24})=m_4(m_{16}+m_{17})+m_2m_{21}$. Hence $m_4m_{16}=m_4m_{17}=m_{22}$.
\item $\psi_{\xi_1^2}(m_4m_{19})=m_4m_{16}+m_3m_{19}$. Hence $m_3m_{19}=m_{22}$.
\item $\psi_{\xi_1^2}(m_4m_{20})=m_4m_{17}+m_3m_{20}$. Hence $m_3m_{20}=m_{22}$.
\item $\psi_{\xi_1^6}(m_{10}m_{19})=m_{10}m_{11}+m_4m_{16}+m_7m_{13}$. Hence $m_{10}m_{11}=m_{22}$.
\item $\psi_{\xi_1^6}(m_{11}m_{19})=m_{11}^2+m_7m_{13}$. Hence $m_{11}^2=0$.
\item $\psi_{\xi_1}(m_6m_{16})=m_6m_{14}+m_5m_{16}$. Hence $m_6m_{14}=0$.
\item $\psi_{\xi_1^2}(m_6m_{18})=m_6m_{15}$. Hence $m_6m_{15}=0$.
\end{enumerate}
The terms used are $\Sq^2_*m_4=m_3$, $\psi_{\xi_2}(m_4)=m_2$, $\Sq^4_*m_{10}=m_4$, $\Sq^2_*m_{19}=m_{16}$, $\psi_{\xi_1^6}(m_{19})=m_{11}$ and $\Sq^1_*m_{16}=m_{14}$ from $E_\infty$, and $\psi_{\xi_2^2}(m_{24})=m_{16}+m_{17}$, $\psi_{\xi_2}(m_{24})=m_{21}$, $\Sq^2_*m_{20}=m_{17}$, $\Sq^2_*m_{10}=\Sq^2_*m_{11}=m_7$, $\Sq^4_*m_{19}=m_{13}$, $\Sq^1_*m_6=m_5$ and $\Sq^2_*m_{18}=m_{15}$ from Table \ref{tab:coaction}.
\end{proof}

The last product $m_{10}^2$ is a square, so the coaction cannot see it: $\psi(m_{10})^2=\sum_I\xi^{2I}\otimes\psi_I(m_{10})^2$, and $\psi_I(m_{10})^2=0$ for $I\neq0$. Instead we compare with the $C_2$-level, where the answer is encoded in the coefficients $c_2$ of Theorem \ref{thm:C2abutment} and $\varepsilon$ of Proposition \ref{prop:C2gamma}.

\begin{proposition}\label{thm:x9sq}
$m_{10}^2=(c_2+\varepsilon)\,m_{22}$.
\end{proposition}
\begin{proof}
First, $\mathrm{tr}(q)=m_{10}+\epsilon_1m_{11}+\epsilon_2m_{12}$ for some $\epsilon_1,\epsilon_2\in\F_2$. Indeed, on $E_\infty$ we have $\Res(m_{10})=\Res(w)\zeta_1^4+\Res(e)\zeta_2^2=y^2\zeta_1^4+x^2\zeta_2^2=z^2$ by Corollary \ref{cor:E4pageC2}, while $\Res(m_{11})$ and $\Res(m_{12})$ lie in filtration above $-4$, as $\Res(\eta')=0$ and $\Res(w^2)=y^4$. So the coefficient of $m_{10}$ in $\mathrm{tr}(q)$ is the coefficient of $z^2$ in $\Res\,\mathrm{tr}(q)=(1+\gamma)(q)=z^2+y^2z+\varepsilon y^4$, which is $1$.

Since $m_{11}^2=0$ by Proposition \ref{prop:products}, and $m_{12}^2=0$ by Lemma \ref{lem:room}, we get $\mathrm{tr}(q)^2=m_{10}^2$. On the other hand, by Frobenius reciprocity,
\[
\mathrm{tr}(q)^2=\mathrm{tr}\bigl(q\,\Res\,\mathrm{tr}(q)\bigr)=\mathrm{tr}\bigl(q^2+q\gamma(q)\bigr)=\mathrm{tr}\bigl(z^2q+y^2zq+\varepsilon y^4q\bigr).
\]
In $H^{C_2}_*$ we have $z^2q=c_2\,y^2z^3$, $y^2zq=z^4=0$ and $y^4q=y^2z^3$, by Theorem \ref{thm:C2abutment}. So $\mathrm{tr}(q)^2=(c_2+\varepsilon)\,\mathrm{tr}(y^2z^3)=(c_2+\varepsilon)\,m_{22}$ by Lemma \ref{lem:trvalues}.
\end{proof}

\subsection{The coefficient $c_2+\varepsilon$}\label{subsec:beps}
It remains to show that $c_2+\varepsilon=1$. The coefficients $c_2$ and $\varepsilon$ depend on the choice of lifts, but their sum does not, by Proposition \ref{thm:x9sq}. They are statements about mod $2$ reductions of torsion-free integral classes, and can be decided by an integral computation, once we know the $\tmf$-algebra $\tmf_1(5)$ and its $C_4$-action explicitly. By Quinn's theorem, these are given by the classical description of $\Gamma_1(5)$-modular forms, which we now recall.

\begin{proposition}\label{prop:tate}
Two-locally, $\pi_*\tmf_1(5)=\Z_{(2)}[a,b]$ with $|a|=|b|=2$, and the universal curve is
\[Y^2+(a+b)XY+a^2bY=X^3+abX^2.\]
The unit $\tmf\to\tmf_1(5)$ sends $c_4$ and $c_6$ to
\begin{align*}
c_4&=a^4-12a^3b+14a^2b^2+12ab^3+b^4,\\
c_6&=-a^6+18a^5b-75a^4b^2-75a^2b^4-18ab^5-b^6,
\end{align*}
and the generator $\gamma$ of $C_4$ acts by $a\mapsto b\mapsto-a$. The classes $a$ and $b$ are the classes $t_1$ and $\gamma t_1$ of Proposition \ref{prop:underlyingE2}.
\end{proposition}
\begin{proof}
This is \cite[Remark 3.2.2, Lemma 4.3.3]{quinn}, with $a=x$ and $b=y$ in the notation there; the curve is the Tate normal form of Behrens--Ormsby \cite[Corollary 1.1.10]{BO} with $u=a$ and $a_1=a+b$. The generator $\bar t_1$ of the form of $\bpcm{1}$ that we use maps to the lift $\bar x$ of $x$ \cite[Construction 4.2.1]{quinn}, so that $t_1=x$ and $\gamma t_1=y$. The values of $c_4$ and $c_6$ are the standard expressions in the Weierstrass coefficients $a_1=a+b$, $a_2=ab$, $a_3=a^2b$ and $a_4=a_6=0$ \cite[Section III.1]{silverman}.
\end{proof}

Recall from Lemma \ref{lem:molien}(2) that $L=H^{\tmf}_*(\tmf_1(5);\Z_{(2)})$ is free of rank $24$ with $L/2=H^e_*$ and $L\otimes\mathbb{Q}=\mathbb{Q}[a,b]/(c_4,c_6)$. By Proposition \ref{prop:tate} one checks that
\[
L=\Z_{(2)}\bigl[a,\ \tfrac{a+b}{2},\ \tfrac{a^3}{2}\bigr]\subset\mathbb{Q}[a,b]/(c_4,c_6).
\]
Indeed, the Hurewicz images of $a+b$ and $a^3$ vanish mod $2$, which gives the inclusion. Since $a=t_1$ and $b=\gamma t_1$ by Proposition \ref{prop:tate}, the classes $a+b$ and $a^3$ are the images of $v_1$ modulo $2$ and of $v_2$ modulo $(2,v_1)$ \cite[Theorem 1.1]{BHSZ}, so the reductions of $(a+b)/2$ and $a^3/2$ have nonzero $\zeta_1^2$- and $\zeta_2^2$-components, so the inclusion is onto mod $2$. The class $\hat Z=a^5(a+b)^3/16$ generates $L_{16}$, and $\gamma\hat Z=\hat Z$.

\begin{proposition}\label{prop:abeps}
$c_2+\varepsilon=1$.
\end{proposition}
\begin{proof}
By Lemma \ref{lem:exponent}, the groups $H^{C_2}_n(\Z)$ are torsion-free for $n=4,8,16$. So $y^2$, $z$, $q$ and $y^2z^3$ are reductions of torsion-free integral classes $W$, $Z$, $Q$ and $V$, and $H^{\tmf}_{16}(\tmf_1(5)^{C_2};\Z_{(2)})=\Z_{(2)}\{V\}$, as its reduction mod $2$ injects into $H^{C_2}_{16}=\F_2$. By Proposition \ref{prop:C2res}, $\mathrm{tr}(\hat Z)=V$ up to a unit, and as $\gamma^2$ acts trivially on $L_{16}$, $\Res(V)=2\hat Z$ up to a unit.

By Proposition \ref{prop:C2gamma}, $\gamma(Q)=Q+D+2H$ with $D=Z^2+WZ+\varepsilon W^2$ and $H$ integral. Since $\gamma V=V$, we have $\gamma(Q)^2=\gamma(Q^2)=Q^2$, so modulo $4V$
\[
0\equiv2QD+D^2\equiv2QD+Z^4+W^2Z^2+\varepsilon W^4+2Z^3W,
\]
the remaining terms $2\varepsilon W^2Z^2$ and $2\varepsilon W^3Z$ of $D^2$ being divisible by $4V$, as $y^4z^2=y^6z=0$. Modulo $2$ we have $QD\equiv(c_2+\varepsilon)V$ and $Z^3W\equiv V$, by the relations used in the proof of Theorem \ref{thm:x9sq}. Finally $Z^4$, $W^2Z^2$ and $W^4$ lie in $4V\Z_{(2)}$: since $\Res(V)$ is $2\hat Z$ up to a unit, this says that $\Res(Z)^4$, $\Res(W)^2\Res(Z)^2$ and $\Res(W)^4$ are divisible by $8\hat Z$ in $L$, which one computes, for any choice of lifts. Hence $2(c_2+\varepsilon+1)V\equiv0$ modulo $4V$.
\end{proof}

With Proposition \ref{thm:x9sq}, this gives $m_{10}^2=m_{22}$ and completes the proof of Theorem \ref{thm:Mprime}.

\section{$\tmf_0(5)$ as a $\tmf$-module}\label{sec:tmfmodule}

In this section we prove Theorem \ref{thm:introB}, and we describe $\tmf_1(5)$ and $\tmf_1(5)^{C_2}$ as $\tmf$-modules using Meier's splitting. We call a $\tmf$-module of the form $\tmf\otimes K$ for a spectrum $K$ \emph{induced}. For a spectrum $K$, $H^{\tmf}_*(\tmf\otimes K)=H_*(K)$, with the $\mathcal{A}(2)_*$-coaction obtained from the $\mathcal{A}_*$-coaction along the quotient $\mathcal{A}_*\to\mathcal{A}(2)_*$. A $\tmf$-module $Y$ is \emph{connective of finite type} if $\pi_kY=0$ for $k\ll0$ and each $\pi_kY$ is a finitely generated $\Z_{(2)}$-module, and for $Y\neq0$ connective, its \emph{bottom degree} is the least $b$ with $\pi_bY\neq0$.

\subsection{Cell structures} We first record the range in which the Hurewicz map of $\tmf$ is an isomorphism.

\begin{lemma}\label{lem:tmfhur}
Let $W$ be a spectrum with $\pi_kW=0$ for $k<b$. Then the map $W\to\tmf\otimes W$ induces an isomorphism on $\pi_k$ for $k\le b+6$ and a surjection for $k=b+7$.
\end{lemma}
\begin{proof}
The unit $\mathbb S\to\tmf$ induces an isomorphism on $\pi_k$ for $k\le6$ and a surjection for $k=7$, as $\pi_7\tmf=0$ \cite{bauer}\cite{bmq}. So its cofiber $\tmf/\mathbb S$ is $7$-connected, and the cofiber $\tmf/\mathbb S\otimes W$ of $W\to\tmf\otimes W$ is $(b+7)$-connected.
\end{proof}

Cell structures of $\tmf$-modules are built exactly as minimal cell structures of spectra, with $H^{\tmf}_*$ in place of mod $2$ homology and the K\"unneth spectral sequence $\Tor^{\pi_*\tmf}(\F_2,\pi_*Y)\Rightarrow H^{\tmf}_*(Y)$ \cite[Theorem IV.4.1]{ekmm} in place of the Hurewicz theorem.

\begin{lemma}\label{lem:mincells}
Let $Y$ be a connective $\tmf$-module of finite type.
\begin{enumerate}
\item If $Y\neq0$, then $H^{\tmf}_*(Y)$ is nonzero, its lowest nonzero degree is the bottom degree $b$ of $Y$, and the Hurewicz map induces $\pi_bY\otimes\F_2\cong H^{\tmf}_b(Y)$.
\item If $H^{\tmf}_*(Y)$ is finite with homogeneous basis $y_1,\dots,y_N$ of degrees $d_1\le\dots\le d_N$, then there are $\tmf$-modules $0=Y_0\to Y_1\to\dots\to Y_N=Y$ over $Y$ and cofiber sequences of $\tmf$-modules
\[
\Sigma^{d_{n+1}-1}\tmf\xrightarrow{\ \alpha_n\ }Y_n\longrightarrow Y_{n+1}
\]
such that $H^{\tmf}_*(Y_n)$ maps injectively to $H^{\tmf}_*(Y)$, with image spanned by $y_1,\dots,y_n$.
\end{enumerate}
\end{lemma}
\begin{proof}
In the K\"unneth spectral sequence, $\Tor_0$ in degree $b$ is $\pi_bY\otimes\F_2$, which is nonzero by the Nakayama lemma, and nothing else contributes in degrees $\le b$; this gives (1). For (2), suppose $Y_n\to Y$ has been constructed. Then $H^{\tmf}_*(Y/Y_n)$ is spanned by the images of $y_{n+1},\dots,y_N$, so by (1) the class $y_{n+1}$ lifts to $\pi_{d_{n+1}}(Y/Y_n)$. The composite $\Sigma^{d_{n+1}}\tmf\to Y/Y_n\to\Sigma Y_n$ is $\Sigma\alpha_n$, and $Y_{n+1}$ is the cofiber of $\alpha_n$, which maps to $Y$ with the required homology.
\end{proof}

When $y_1,\dots,y_n$ are the basis elements of degree at most $d$, we write $Y_{\le d}$ for $Y_n$. Then $H^{\tmf}_*(Y_{\le d})$ is the part of $H^{\tmf}_*(Y)$ in degrees at most $d$, which is a subcomodule as the coaction lowers degree, and $H^{\tmf}_*(Y/Y_{\le d})$ is the quotient comodule.

\begin{lemma}\label{lem:narrow}
Let $Y$ be a connective $\tmf$-module of finite type such that $H^{\tmf}_*(Y)$ is finite and concentrated in degrees $b,b+1,\dots,b+8$. Then $Y\simeq\tmf\otimes K$ for a finite spectrum $K$ with one cell of dimension $d$ for each basis element of $H^{\tmf}_d(Y)$. In particular $H_*(K)\cong H^{\tmf}_*(Y)$ as $\mathcal{A}(2)_*$-comodules.
\end{lemma}
\begin{proof}
We take the cell structure of Lemma \ref{lem:mincells}. We show by induction on $n$ that $Y_n\simeq\tmf\otimes K_n$, where $K_n$ is a finite spectrum with cells in dimensions $d_1,\dots,d_n$. For $n=0$ take $K_0=*$. Suppose $Y_n\simeq\tmf\otimes K_n$. The attaching map $\alpha_n$ is an element of
\[
[\Sigma^{d_{n+1}-1}\tmf,\ \tmf\otimes K_n]_{\tmf}=\pi_{d_{n+1}-1}(\tmf\otimes K_n).
\]
Since $K_n$ has cells in dimensions at least $b$ and $d_{n+1}-1\le b+7$, Lemma \ref{lem:tmfhur} gives $\beta\in\pi_{d_{n+1}-1}K_n$ whose image is $\alpha_n$. We set $K_{n+1}=K_n\cup_\beta e^{d_{n+1}}$. Then $\tmf\otimes\beta$ is the $\tmf$-module map corresponding to $\alpha_n$, so
\[
\tmf\otimes K_{n+1}\simeq\mathrm{cofib}(\tmf\otimes\beta)\simeq\mathrm{cofib}(\alpha_n)=Y_{n+1}.\qedhere
\]
\end{proof}

We will apply this to the fixed points of $\tmf_1(5)$, which are connective of finite type \cite{hhrc4}. Finally, we need the following fact about the homotopy of the fixed points.

\begin{lemma}\label{lem:fournu}
The Mackey functor $\underline{\pi}_0\tmf_1(5)$ is the constant Mackey functor $\underline{\Z}_{(2)}$, generated by the unit. Consequently every torsion class of $\pi_*\tmf_0(5)$ is annihilated by $4$; in particular $4\nu=0$ in $\pi_3\tmf_0(5)$.
\end{lemma}
\begin{proof}
The $C_4$-spectrum $\tmf_1(5)$ is strongly even \cite[Lemma 4.3.3]{quinn}: the restrictions $\pi_0\tmf_0(5)\to\pi_0\tmf_1(5)^{C_2}\to\pi_0\tmf_1(5)=\Z_{(2)}$ are isomorphisms and the Weyl actions are trivial in degree $0$, so each transfer is multiplication by the index. In particular $\mathrm{tr}^{C_4}_e(1)=4$, so that $\mathrm{tr}^{C_4}_e(\Res^{C_4}_e(x))=4x$ for every $x\in\pi_*\tmf_0(5)$ by Frobenius reciprocity. As $\pi_*\tmf_1(5)$ is concentrated in even degrees and torsion-free, $\Res^{C_4}_e(x)=0$ for every torsion class $x$, and $4x=0$.
\end{proof}

\subsection{Meier's splitting} By \cite{Mei23}, $\tmf_1(n)$ carries a genuine $C_2$-structure induced by $(E,P)\mapsto(E,-P)$ on $\overline{\mathcal M}_1(n)$, that it is strongly even in the sense of \cite{HM17}, and that for $n=3$ it is the $C_2$-spectrum of \cite{HM17}, a form of $\bprn{2}$. For $n=5$ the subgroup $C_2\subseteq C_4\cong(\Z/5)^\times$ is $\{\pm1\}$, and $\Res^{C_4}_{C_2}\bpcm{1}$ is strongly even with underlying spectrum $\tmf_1(5)$ \cite[Lemma 4.3.3]{quinn}, so the two $C_2$-spectra agree \cite[Lemma 3.4]{HM17}.

\begin{theorem}[Meier {\cite[Theorem 5.6]{Mei23}}, see also {\cite[Example 3.7, Theorem 6.19]{Mei22}}]
\label{thm:meier}
Two-locally there is an equivalence of genuine $C_2$-spectra, compatible with the $\tmf$-module structures,
\[
  \tmf_1(5)\ \simeq\ \tmf_1(3)\ \oplus\ \Sigma^{\rho}\tmf_1(3)\ \oplus\ \Sigma^{2\rho}\tmf_1(3).
\]
\end{theorem}

\begin{corollary}\label{cor:meier}
\begin{enumerate}
\item There is an equivalence of $\tmf$-modules $\tmf_1(5)\simeq\tmf\otimes\bigl(DA(1)\oplus\Sigma^2DA(1)\oplus\Sigma^4DA(1)\bigr)$.
\item There is an equivalence of $\tmf$-modules
\[\tmf_1(5)^{C_2}\simeq\tmf\otimes X\oplus T_1\oplus T_2,\]
where $T_k=(\Sigma^{k\rho}\tmf_1(3))^{C_2}$ and $X$ is the $10$-cell complex of \cite[Corollary 5.13]{CHR}.
\item For $k=1,2$ there is a cofiber sequence of $\tmf$-modules $\tmf\otimes K_k'\to T_k\to\tmf\otimes K_k''$, where $K_k'$ and $K_k''$ are finite complexes with cells in dimensions $1$ to $9$ and $10$ to $16$, respectively.
\end{enumerate}
\end{corollary}
\begin{proof}
Part (1) follows from Theorem \ref{thm:meier} and $\tmf_1(3)\simeq\tmf\otimes DA(1)$ \cite{HM}\cite[\S3]{Mat}. For (2), take $C_2$-fixed points in Theorem \ref{thm:meier}, and use $\tmf_1(3)^{C_2}=\tmf_0(3)\simeq\tmf\otimes X$ \cite{HM17}\cite[Corollary 5.13]{CHR}. For (3), $T_1$ and $T_2$ are summands of $\tmf_1(5)^{C_2}$, so they are connective of finite type, and their $\tmf$-homology lies in degrees at most $16$ by Theorem \ref{thm:C2abutment}. By Lemma \ref{lem:fournu}(1), $\pi_0\tmf_1(5)^{C_2}\cong\Z_{(2)}$, which is indecomposable, and since $\pi_0\tmf_0(3)\neq0$, we get $\pi_0T_1=\pi_0T_2=0$. As $H^{C_2}_*$ vanishes in negative degrees, Lemma \ref{lem:mincells} shows that $T_k$ has bottom degree at least $1$, so that $H^{\tmf}_*(T_k)$ lies in degrees $1$ to $16$, and we apply Lemma \ref{lem:narrow} to $(T_k)_{\le9}$ and $T_k/(T_k)_{\le9}$.
\end{proof}

\subsection{The two layers} We write $H^{C_4}_{\le8}$ for the span of $m_1,\dots,m_{12}$ and $H^{C_4}_{\ge10}$ for the span of $m_{13},\dots,m_{24}$. Since the coaction lowers degree, $H^{C_4}_{\le8}$ is a subcomodule with quotient $H^{C_4}_{\ge10}$.

\begin{theorem}\label{thm:tmfmodule}
The $\tmf$-module $\tmf_0(5)$ has the following properties.
\begin{enumerate}
\item There are finite spectra $K_1$ and $K_2$, with cells in dimensions
\[
K_1\colon\ 0,1,2,4,4,5,6,7,7,8,8,8,\qquad K_2\colon\ 10,11,11,12,12,13,14,14,15,16,17,18,
\]
and a cofiber sequence of $\tmf$-modules
\begin{equation}\label{eq:twolayer}
\tmf\otimes K_1\longrightarrow\tmf_0(5)\longrightarrow\tmf\otimes K_2
\end{equation}
which induces the short exact sequence $0\to H^{C_4}_{\le8}\to H^{C_4}_*\to H^{C_4}_{\ge10}\to0$ on $H^{\tmf}_*$.
\item There is no spectrum $K$ with $\tmf_0(5)\simeq\tmf\otimes K$ as $\tmf$-modules. In particular \eqref{eq:twolayer} does not split, and its boundary map $\tmf\otimes K_2\to\Sigma\,\tmf\otimes K_1$ is not of the form $\tmf\otimes h$ for a map of spectra $h\colon K_2\to\Sigma K_1$.
\item $\tmf_0(5)$ is indecomposable as a $\tmf$-module.
\end{enumerate}
\end{theorem}

\begin{proof}[proof of Theorem \ref{thm:tmfmodule}(1)]
Let $Y=\tmf_0(5)$. Then $H^{\tmf}_*(Y_{\le8})=H^{C_4}_{\le8}$ is concentrated in degrees $0$ to $8$, so Lemma \ref{lem:narrow} gives $Y_{\le8}\simeq\tmf\otimes K_1$. The module $Y/Y_{\le8}$ is connective of finite type, and its $\tmf$-homology $H^{C_4}_{\ge10}$ is concentrated in degrees $10$ to $18$, as $H^{C_4}_9=0$. So its bottom degree is $10$ by Lemma \ref{lem:mincells}, and Lemma \ref{lem:narrow} gives $Y/Y_{\le8}\simeq\tmf\otimes K_2$.
\end{proof}

The cohomology of $K_1$ and $K_2$ is the $\mathcal{A}(2)$-module dual to $H^{C_4}_{\le8}$ and to $H^{C_4}_{\ge10}$, given by Table \ref{tab:Mstar}. Their cell diagrams are drawn in Figure \ref{fig:cells} in the introduction. They show the action of $\Sq^1$, $\Sq^2$ and $\Sq^4$, which is all that $H^{\tmf}_*$ records.

Part (2) follows from the next proposition. Let $m_i^*$ be the basis of $\Hom(H^{C_4}_*,\F_2)$ dual to the $m_i$, as in Table \ref{tab:Mstar}.

\begin{proposition}\label{prop:noAlift}
There is no $\mathcal{A}$-module structure on $\Hom(H^{C_4}_*,\F_2)$ whose restriction to $\mathcal{A}(2)$ is the given one.
\end{proposition}
\begin{proof}
Suppose there were such a structure. In $\mathcal{A}$ we have the Adem relations
\[
\Sq^4\Sq^6=\Sq^{10}+\Sq^8\Sq^2,\qquad\Sq^2\Sq^8=\Sq^{10}+\Sq^9\Sq^1,\qquad\Sq^1\Sq^8=\Sq^9,
\]
and hence
\[
\Sq^4\Sq^6=\Sq^2\Sq^8+\Sq^1\Sq^8\Sq^1+\Sq^8\Sq^2.
\]
We apply this to $m_2^*$, in degree $1$. Since $H^{C_4}_*$ vanishes in degrees $3$ and $9$, we have $\Sq^2m_2^*=0$ and $\Sq^8m_2^*=0$, so $\Sq^4\Sq^6m_2^*=\Sq^1\Sq^8\Sq^1m_2^*$ lies in $\Sq^1$ of degree $10$, which is spanned by $\Sq^1m_{13}^*=m_{15}^*$. On the other hand, using $\Sq^6=\Sq^2\Sq^4+\Sq^1\Sq^4\Sq^1$ in $\mathcal{A}(2)$ and Table \ref{tab:Mstar}, we get
\[
\Sq^6m_2^*=\Sq^2m_6^*+\Sq^1\Sq^4m_3^*=m_8^*,\qquad\Sq^4\Sq^6m_2^*=\Sq^4m_8^*=m_{14}^*.
\]
Since $m_{14}^*\notin\F_2\{m_{15}^*\}$, this is a contradiction.
\end{proof}

\begin{proof}[proof of Theorem \ref{thm:tmfmodule}(2)]
If $\tmf_0(5)\simeq\tmf\otimes K$, then $H^{C_4}_*\cong H_*(K)$ as $\mathcal{A}(2)_*$-comodules, and $H^*(K)$ is an $\mathcal{A}$-module extending the $\mathcal{A}(2)$-module $\Hom(H^{C_4}_*,\F_2)$, contradicting Proposition \ref{prop:noAlift}. If \eqref{eq:twolayer} split, then $\tmf_0(5)\simeq\tmf\otimes(K_1\oplus K_2)$, and if $\delta=\tmf\otimes h$, then $\tmf_0(5)\simeq\tmf\otimes\mathrm{fib}(h)$.
\end{proof}

For part (3), there are two inputs: the indecomposability of the comodule $\F_2\{m_2,\dots,m_{24}\}$, which can be read off Table \ref{tab:Mstar}, and the relation $4\nu=0$ of Lemma \ref{lem:fournu}.

\begin{lemma}\label{lem:Mbarindec}
As $\mathcal{A}(2)$-modules, $\Hom(H^{C_4}_*,\F_2)=\F_2\{m_1^*\}\oplus\overline{H}$, where $\overline{H}$ is spanned by $m_2^*,\dots,m_{24}^*$, and both summands are indecomposable.
\end{lemma}
\begin{proof}
By Table \ref{tab:Mstar}, $\Sq^1$, $\Sq^2$ and $\Sq^4$ vanish on $m_1^*$, and $\overline{H}$ is a submodule because it consists of all classes of positive degree. This gives the splitting.

Suppose $\overline{H}=N_1\oplus N_2$ as $\mathcal{A}(2)$-modules. We make two observations. First, $\overline{H}$ is one-dimensional in degree $1$, so we may assume $m_2^*\in N_1$. Second, if $\overline{H}$ is one-dimensional in the degree of a basis class $v$, then $v\in N_1$ or $v\in N_2$, and in that case, if $\theta v$ is a nonzero element of $N_1$ for some $\theta\in\mathcal{A}(2)$, then $v\in N_1$.

All references below are to Table \ref{tab:Mstar}. The chain
\[
m_2^*\xrightarrow{\Sq^4}m_6^*\xrightarrow{\Sq^2}m_8^*\xrightarrow{\Sq^4}m_{14}^*\xrightarrow{\Sq^1}m_{16}^*\xrightarrow{\Sq^2}m_{19}^*\xrightarrow{\Sq^4}m_{24}^*
\]
puts all these classes in $N_1$. The degrees of $m_{23}^*$, $m_{22}^*$, $m_{21}^*$, $m_{18}^*$, $m_{13}^*$ and $m_7^*$ are one-dimensional. Applying the second observation in turn to
\begin{align*}
\Sq^1m_{23}^*&=m_{24}^*, & \Sq^2m_{22}^*&=m_{24}^*, & \Sq^2m_{21}^*&=m_{23}^*,\\
\Sq^4m_{18}^*&=m_{23}^*, & \Sq^2\Sq^1m_{13}^*=\Sq^2m_{15}^*&=m_{18}^*, & \Sq^4m_7^*&=m_{13}^*
\end{align*}
puts all six classes in $N_1$. Applying $\mathcal{A}(2)$ to the classes already in $N_1$ gives
\begin{align*}
m_3^*&=\Sq^1m_2^*, & m_4^*&=\Sq^2m_3^*, & m_{10}^*&=\Sq^4m_4^*, & m_{11}^*&=\Sq^1m_8^*+m_{10}^*,\\
m_9^*&=\Sq^1m_7^*, & m_{15}^*&=\Sq^1m_{13}^*, & m_{17}^*&=\Sq^2m_{13}^*, & m_{20}^*&=\Sq^2m_{17}^*,
\end{align*}
all in $N_1$. Hence $N_1$ contains every $m_j^*$ with $2\le j\le24$, except possibly $m_5^*$ and $m_{12}^*$. So $N_2$ lies in degrees $4$ and $8$, and is at most one-dimensional in each.

A nonzero element of $N_2$ in degree $4$ has the form $m_5^*+\epsilon\,m_4^*$, and its $\Sq^1$ is $m_6^*\neq0$, which must lie in $N_2$; but $N_2$ is zero in degree $5$. A nonzero element of $N_2$ in degree $8$ has the form $m_{12}^*+\epsilon_1m_{10}^*+\epsilon_2m_{11}^*$, and its $\Sq^4$ is $m_{17}^*+\epsilon_2m_{16}^*\neq0$; but $N_2$ is zero in degree $12$. Hence $N_2=0$.
\end{proof}

\begin{proof}[proof of Theorem \ref{thm:tmfmodule}(3)]
Suppose $\tmf_0(5)\simeq Y'\oplus Y''$ as $\tmf$-modules, with $Y'$ and $Y''$ nonzero. Both are connective of finite type, so their $\tmf$-homologies are nonzero by Lemma \ref{lem:mincells}, and dually $\Hom(H^{C_4}_*,\F_2)=N'\oplus N''$ as $\mathcal{A}(2)$-modules, with $N'$ and $N''$ nonzero. One of the two, say $N'$, contains $m_1^*$, and then $\overline{H}=(N'\cap\overline{H})\oplus N''$, since $\overline{H}$ is the part of positive degree. By Lemma \ref{lem:Mbarindec}, $N'\cap\overline{H}=0$, so $N'=\F_2\{m_1^*\}$ and $Y'\simeq\tmf$ by Lemma \ref{lem:mincells}. Thus $\tmf$ is a retract of $\tmf_0(5)$ as a $\tmf$-module, by a split injection $i$, and $i(4\nu)=4\nu\cdot i(1)=0$ by Lemma \ref{lem:fournu}, while $4\nu\neq0$ in $\pi_3\tmf=\Z/8\{\nu\}$.
\end{proof}

\section{The Adams spectral sequence of $\tmf_0(5)$}\label{sec:adamsss}
In this section we compute the Adams spectral sequence
\[E_2=\Ext_{\mathcal{A}(2)_*}(\F_2,H^{C_4}_*)\implies\pi_*\tmf_0(5)\]
in stems at most $32$. The periodic theory $D^{-1}\bpcm{1}$ is $32$-periodic \cite{hhrc4}, and this is the range in which the homotopy of $\tmf_0(5)$ displays its basic pattern. We follow the method of Bruner--Rognes for $\tmf$ \cite{BR}: we compute $E_2$ by machine as a module over $\Ext_{\mathcal{A}(2)_*}(\F_2,\F_2)$, and we determine the differentials using the Leibniz rule together with what is known about $\pi_*\tmf_0(5)$ from the slice spectral sequence \cite{hhrc4}\cite{hswx}.

\subsection{The $E_2$-page} By Lemma \ref{lem:Mbarindec}, the $\mathcal{A}(2)$-module dual to $H^{C_4}_*$ splits as $\F_2\{m_1^*\}\oplus\overline{H}$. Hence
\[E_2=\Ext_{\mathcal{A}(2)_*}(\F_2,\F_2)\oplus\Ext_{\mathcal{A}(2)}(\overline{H},\F_2)\]
as a module over $\Ext_{\mathcal{A}(2)_*}(\F_2,\F_2)$. The first summand is the image of the Adams spectral sequence of $\tmf$ under the unit $\tmf\to\tmf_0(5)$, and we write its classes as elements of $\Ext_{\mathcal{A}(2)_*}(\F_2,\F_2)$. Theorem \ref{thm:tmfmodule}(3) says that this splitting is not realized by $\tmf$-modules, and indeed the differentials below connect the two summands; the first instance is the differential that kills $4\nu$.

We use the names of Bruner--Rognes \cite[Chapter 3]{BR} for the algebra generators $h_0,h_1,h_2,c_0,w_1,\alpha,\beta,d_0,e_0,g,\gamma,\delta$ of $\Ext_{\mathcal{A}(2)_*}(\F_2,\F_2)$ in stems below $48$. We computed both summands of $E_2$, together with the action of these generators, through stem $48$ and filtration $26$, using the program \texttt{ext} of the \texttt{sseq} package \cite{sseq} applied to the module of Table \ref{tab:Mstar}.

\begin{proposition}\label{prop:adamsE2}
In stems at most $34$, the second summand $\Ext_{\mathcal{A}(2)}(\overline{H},\F_2)$ is generated as a module over $\Ext_{\mathcal{A}(2)_*}(\F_2,\F_2)$ by one class $x_{n,s}$ in bidegree $(n,s)$ for each of the $25$ bidegrees
\begin{align*}
&(1,0),\ (4,0),\ (4,2),\ (8,1),\ (8,3),\ (9,1),\ (10,2),\ (12,3),\ (14,2),\\
&(15,0),\ (15,2),\ (16,0),\ (16,3),\ (16,4),\ (18,2),\ (20,3),\ (22,4),\\
&(24,5),\ (25,6),\ (26,4),\ (28,5),\ (32,4),\ (32,5),\ (32,7),\ (34,6).
\end{align*}
The classes in filtration $0$ are the primitives $x_{1,0}=m_2$, $x_{4,0}=m_5$, $x_{15,0}=m_{21}$ and $x_{16,0}=m_{22}$ of $H^{C_4}_*$.
\end{proposition}

To name the generators we use the cell filtration of $\overline{H}$. Filtering $\overline{H}$ by degree gives an algebraic Atiyah--Hirzebruch spectral sequence
\[\bigoplus_j\Ext_{\mathcal{A}(2)_*}(\F_2,\F_2)[j]\implies\Ext_{\mathcal{A}(2)}(\overline{H},\F_2),\]
with one copy of $\Ext_{\mathcal{A}(2)_*}(\F_2,\F_2)$ for each class $m_j$ of $\overline{H}$, shifted by its degree. For $a\in\Ext_{\mathcal{A}(2)_*}(\F_2,\F_2)$ we write $a[j]$ for a class detected by $a$ on the cell $m_j$, so that the classes $[j]$ in filtration $0$ are the primitives $m_j$. Table \ref{tab:cellnames} lists the names of all the generators, and the charts are labelled by them. We computed the maps on $\Ext$ induced by the quotients $\overline{H}^*\to\overline{H}^*_{\le k}$ and by the inclusions of the top cells of $\overline{H}^*_{\le k}$ with \texttt{ext}.

\begin{table}[ht]\centering\footnotesize
\caption{The generators of Proposition \ref{prop:adamsE2} and their names in the cell filtration of $\overline{H}$: $x_{n,s}$ is detected by $a[j]$, and $a[i+j]$ stands for $a[i]+a[j]$.}
\label{tab:cellnames}
\renewcommand{\arraystretch}{1.15}
\begin{tabular}{@{}ll@{\hskip 3em}ll@{}}
\toprule
class & name & class & name\\
\midrule
$x_{1,0}$ & $m_{2}$ & $x_{16,4}$ & $w_1[12]$\\
$x_{4,0}$ & $m_{5}$ & $x_{18,2}$ & $h_2^{2}[17]$\\
$x_{4,2}$ & $h_0^{2}[4]$ & $x_{20,3}$ & $c_0[16]$\\
$x_{8,1}$ & $h_0[12]$ & $x_{22,4}$ & $d_0[12]$\\
$x_{8,3}$ & $h_0^{3}[10+11+12]$ & $x_{24,5}$ & $h_0^{2}\alpha[17]$\\
$x_{9,1}$ & $h_1[12]$ & $x_{25,6}$ & $h_0h_2w_1[19]$\\
$x_{10,2}$ & $h_1^{2}[10+12]$ & $x_{26,4}$ & $h_1c_0[23]$\\
$x_{12,3}$ & $h_0^{3}[17]$ & $x_{28,5}$ & $h_2w_1[23]$\\
$x_{14,2}$ & $h_0^{2}[19]$ & $x_{32,4}$ & $e_0[21]$\\
$x_{15,0}$ & $m_{21}$ & $x_{32,5}$ & $h_0e_0[21]$\\
$x_{15,2}$ & $h_0h_2[17]$ & $x_{32,7}$ & $h_0\alpha^{2}[10+11]$\\
$x_{16,0}$ & $m_{22}$ & $x_{34,6}$ & $c_0\alpha[19]$\\
$x_{16,3}$ & $c_0[10+11]$ &  & \\
\bottomrule
\end{tabular}
\end{table}

We write $r_n$ for the rank of $\pi_n\tmf_0(5)$. By Lemma \ref{lem:molien}(1) it is the dimension of the $\gamma$-invariants of $\mathbb{Q}[a,b]$ in degree $n$. Over $\mathbb{Q}(i)$ the span of $a$ and $b$ has a basis $u,v$ of eigenvectors of $\gamma$ with eigenvalues $i$ and $-i$, and the monomial $u^pv^q$ is invariant exactly when $p\equiv q\bmod4$. Counting these monomials gives
\[
r_n=\begin{cases}2\lfloor n/8\rfloor+1,& n\equiv0\bmod4,\\ 0,&\text{otherwise,}\end{cases}
\]
so $r_n$ is $1,1,3,3,5,5,7,7,9$ for $n=0,4,\dots,32$. The infinite $h_0$-towers of $E_\infty$ detect the free part of $\pi_*\tmf_0(5)$, so there are $r_n$ of them in stem $n$. On $E_2$ there can be more, and the difference is accounted for as follows.

\begin{lemma}\label{lem:towers}
Let $T(n)$ be the number of infinite $h_0$-towers in stem $n$ of $E_2$. Then there are integers $\delta_n\ge0$ with $\delta_0=0$ and
\[
T(n)-r_n=\delta_n+\delta_{n+1}\qquad\text{for all }n\ge0,
\]
where $\delta_n$ is the number of differentials carrying an infinite $h_0$-tower in stem $n$ onto one in stem $n-1$.
\end{lemma}
\begin{proof}
Suppose $d_r(x)=y$. If $y$ lies in an infinite $h_0$-tower, then $d_r(h_0^kx)=h_0^ky\neq0$ for all $k$, so $x$ lies in an infinite $h_0$-tower as well. If on the other hand $x$ lies in an infinite $h_0$-tower and $h_0^ky=0$ for some $k$, then $h_0^kx$ is a $d_r$-cycle, and the part of the tower through $x$ from $h_0^kx$ on survives to $E_{r+1}$. So the number of infinite towers changes only through differentials carrying one infinite tower onto another, necessarily in adjacent stems, and each such differential removes one tower from each of the two stems.
\end{proof}

\begin{figure}[htbp]
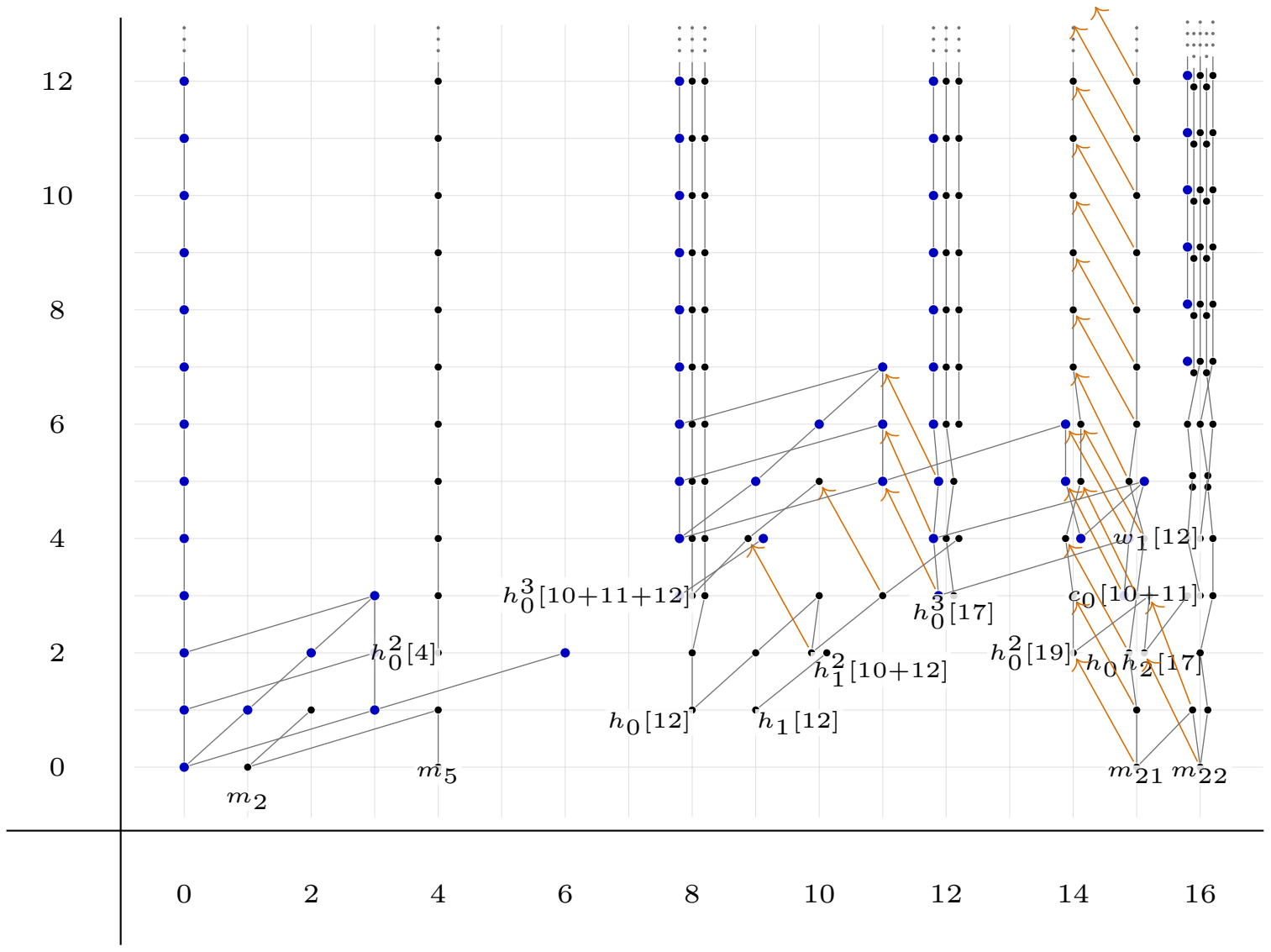
\centering
\adjustbox{max width=\textwidth}{\printpage[name = AdamsEtwoA, page = 2]}
\caption{The Adams $E_2$-page of $\tmf_0(5)$ in stems $0$ to $16$, with the $d_2$-differentials. Blue dots are the classes of the summand $\Ext_{\mathcal{A}(2)_*}(\F_2,\F_2)$, and lines are multiplication by $h_0$, $h_1$ and $h_2$. The module generators of Proposition \ref{prop:adamsE2} are labelled by their names in Table \ref{tab:cellnames}.}
\label{fig:adamsE2a}
\end{figure}

The page $E_2$ is shown in Figures \ref{fig:adamsE2a}--\ref{fig:adamsE2c}. Each $x_{n,s}$ is determined only up to decomposable classes, and we fix our choices in Theorem \ref{thm:adamsdiff} below. In stems at most $34$, the infinite $h_0$-towers of $E_2$ lie in the stems $4k$, where there are exactly $r_{4k}$ of them, and in the pairs of stems $(14,15)$, $(22,23)$, $(26,27)$, $(30,31)$, each of which contains one tower in each stem. So $\delta_n=0$ except for $n=15,23,27,31$, where $\delta_n=1$, and each of these four pairs of towers must be connected by a differential.

\begin{sseqdata}[name = AdamsEtwoA, adams, x range = {0}{16}, y range = {0}{12}]
\class[tmfclass](0,0)
\class[tmfclass](0,1)
\class[tmfclass](0,2)
\class[tmfclass](0,3)
\class[tmfclass](0,4)
\class[tmfclass](0,5)
\class[tmfclass](0,6)
\class[tmfclass](0,7)
\class[tmfclass](0,8)
\class[tmfclass](0,9)
\class[tmfclass](0,10)
\class[tmfclass](0,11)
\class[tmfclass](0,12)
\class[tmfclass](0,13)
\class[tmfclass](0,14)
\class[tmfclass](0,15)
\class[tmfclass](0,16)
\class[tmfclass](0,17)
\class[tmfclass](0,18)
\class[tmfclass](0,19)
\class[tmfclass](0,20)
\class[tmfclass](0,21)
\class[tmfclass](0,22)
\class[tmfclass](0,23)
\class[tmfclass](0,24)
\class(1,0)
\classoptions["{\scriptscriptstyle m_{2}}" below, label distance = -1.5pt](1,0,1)
\class[tmfclass](1,1)
\class(2,1)
\class[tmfclass](2,2)
\class[tmfclass](3,1)
\class[tmfclass](3,2)
\class[tmfclass](3,3)
\class(4,0)
\classoptions["{\scriptscriptstyle m_{5}}" below, label distance = -1.5pt](4,0,1)
\class(4,1)
\class(4,2)
\classoptions["{\scriptscriptstyle h_0^{2}[4]}" left, label distance = -1.5pt](4,2,1)
\class(4,3)
\class(4,4)
\class(4,5)
\class(4,6)
\class(4,7)
\class(4,8)
\class(4,9)
\class(4,10)
\class(4,11)
\class(4,12)
\class(4,13)
\class(4,14)
\class(4,15)
\class(4,16)
\class(4,17)
\class(4,18)
\class(4,19)
\class(4,20)
\class(4,21)
\class(4,22)
\class(4,23)
\class(4,24)
\class[tmfclass](6,2)
\class(8,1)
\classoptions["{\scriptscriptstyle h_0[12]}" below left, label distance = -1.5pt](8,1,1)
\class(8,2)
\class[tmfclass](8,3)
\class(8,3)
\classoptions["{\scriptscriptstyle h_0^{3}[10+11+12]}" left, label distance = -1.5pt](8,3,2)
\class(8,3)
\class[tmfclass](8,4)
\class(8,4)
\class(8,4)
\class[tmfclass](8,5)
\class(8,5)
\class(8,5)
\class[tmfclass](8,6)
\class(8,6)
\class(8,6)
\class[tmfclass](8,7)
\class(8,7)
\class(8,7)
\class[tmfclass](8,8)
\class(8,8)
\class(8,8)
\class[tmfclass](8,9)
\class(8,9)
\class(8,9)
\class[tmfclass](8,10)
\class(8,10)
\class(8,10)
\class[tmfclass](8,11)
\class(8,11)
\class(8,11)
\class[tmfclass](8,12)
\class(8,12)
\class(8,12)
\class[tmfclass](8,13)
\class(8,13)
\class(8,13)
\class[tmfclass](8,14)
\class(8,14)
\class(8,14)
\class[tmfclass](8,15)
\class(8,15)
\class(8,15)
\class[tmfclass](8,16)
\class(8,16)
\class(8,16)
\class[tmfclass](8,17)
\class(8,17)
\class(8,17)
\class[tmfclass](8,18)
\class(8,18)
\class(8,18)
\class[tmfclass](8,19)
\class(8,19)
\class(8,19)
\class[tmfclass](8,20)
\class(8,20)
\class(8,20)
\class[tmfclass](8,21)
\class(8,21)
\class(8,21)
\class[tmfclass](8,22)
\class(8,22)
\class(8,22)
\class[tmfclass](8,23)
\class(8,23)
\class(8,23)
\class[tmfclass](8,24)
\class(8,24)
\class(8,24)
\class(9,1)
\classoptions["{\scriptscriptstyle h_1[12]}" below right, label distance = -1.5pt](9,1,1)
\class(9,2)
\class(9,4)
\class[tmfclass](9,4)
\class[tmfclass](9,5)
\class(10,2)
\classoptions["{\scriptscriptstyle h_1^{2}[10+12]}" below right, label distance = -1.5pt](10,2,1)
\class(10,2)
\class(10,3)
\class(10,5)
\class[tmfclass](10,6)
\class(11,3)
\class[tmfclass](11,5)
\class[tmfclass](11,6)
\class[tmfclass](11,7)
\class[tmfclass](12,3)
\class(12,3)
\classoptions["{\scriptscriptstyle h_0^{3}[17]}" below, label distance = -1.5pt](12,3,2)
\class[tmfclass](12,4)
\class(12,4)
\class(12,4)
\class[tmfclass](12,5)
\class(12,5)
\class[tmfclass](12,6)
\class(12,6)
\class(12,6)
\class[tmfclass](12,7)
\class(12,7)
\class(12,7)
\class[tmfclass](12,8)
\class(12,8)
\class(12,8)
\class[tmfclass](12,9)
\class(12,9)
\class(12,9)
\class[tmfclass](12,10)
\class(12,10)
\class(12,10)
\class[tmfclass](12,11)
\class(12,11)
\class(12,11)
\class[tmfclass](12,12)
\class(12,12)
\class(12,12)
\class[tmfclass](12,13)
\class(12,13)
\class(12,13)
\class[tmfclass](12,14)
\class(12,14)
\class(12,14)
\class[tmfclass](12,15)
\class(12,15)
\class(12,15)
\class[tmfclass](12,16)
\class(12,16)
\class(12,16)
\class[tmfclass](12,17)
\class(12,17)
\class(12,17)
\class[tmfclass](12,18)
\class(12,18)
\class(12,18)
\class[tmfclass](12,19)
\class(12,19)
\class(12,19)
\class[tmfclass](12,20)
\class(12,20)
\class(12,20)
\class[tmfclass](12,21)
\class(12,21)
\class(12,21)
\class[tmfclass](12,22)
\class(12,22)
\class(12,22)
\class[tmfclass](12,23)
\class(12,23)
\class(12,23)
\class[tmfclass](12,24)
\class(12,24)
\class(12,24)
\class(14,2)
\classoptions["{\scriptscriptstyle h_0^{2}[19]}" left, label distance = -1.5pt](14,2,1)
\class(14,3)
\class(14,4)
\class[tmfclass](14,4)
\class[tmfclass](14,5)
\class(14,5)
\class[tmfclass](14,6)
\class(14,6)
\class(14,7)
\class(14,8)
\class(14,9)
\class(14,10)
\class(14,11)
\class(14,12)
\class(14,13)
\class(14,14)
\class(14,15)
\class(14,16)
\class(14,17)
\class(14,18)
\class(14,19)
\class(14,20)
\class(14,21)
\class(14,22)
\class(14,23)
\class(14,24)
\class(15,0)
\classoptions["{\scriptscriptstyle m_{21}}" below, label distance = -1.5pt](15,0,1)
\class(15,1)
\class(15,2)
\class(15,2)
\classoptions["{\scriptscriptstyle h_0h_2[17]}" below, label distance = -1.5pt](15,2,2)
\class[tmfclass](15,3)
\class(15,3)
\class(15,3)
\class[tmfclass](15,4)
\class(15,4)
\class(15,5)
\class[tmfclass](15,5)
\class(15,6)
\class(15,7)
\class(15,8)
\class(15,9)
\class(15,10)
\class(15,11)
\class(15,12)
\class(15,13)
\class(15,14)
\class(15,15)
\class(15,16)
\class(15,17)
\class(15,18)
\class(15,19)
\class(15,20)
\class(15,21)
\class(15,22)
\class(15,23)
\class(15,24)
\class(16,0)
\classoptions["{\scriptscriptstyle m_{22}}" below, label distance = -1.5pt](16,0,1)
\class(16,1)
\class(16,1)
\class(16,2)
\class(16,3)
\class(16,3)
\classoptions["{\scriptscriptstyle c_0[10+11]}" left, label distance = -1.5pt](16,3,2)
\class(16,3)
\class(16,4)
\class(16,4)
\classoptions["{\scriptscriptstyle w_1[12]}" left, label distance = -1.5pt](16,4,2)
\class(16,4)
\class(16,5)
\class(16,5)
\class(16,5)
\class(16,5)
\class(16,6)
\class(16,6)
\class(16,6)
\class[tmfclass](16,7)
\class(16,7)
\class(16,7)
\class(16,7)
\class(16,7)
\class[tmfclass](16,8)
\class(16,8)
\class(16,8)
\class(16,8)
\class(16,8)
\class[tmfclass](16,9)
\class(16,9)
\class(16,9)
\class(16,9)
\class(16,9)
\class[tmfclass](16,10)
\class(16,10)
\class(16,10)
\class(16,10)
\class(16,10)
\class[tmfclass](16,11)
\class(16,11)
\class(16,11)
\class(16,11)
\class(16,11)
\class[tmfclass](16,12)
\class(16,12)
\class(16,12)
\class(16,12)
\class(16,12)
\class[tmfclass](16,13)
\class(16,13)
\class(16,13)
\class(16,13)
\class(16,13)
\class[tmfclass](16,14)
\class(16,14)
\class(16,14)
\class(16,14)
\class(16,14)
\class[tmfclass](16,15)
\class(16,15)
\class(16,15)
\class(16,15)
\class(16,15)
\class[tmfclass](16,16)
\class(16,16)
\class(16,16)
\class(16,16)
\class(16,16)
\class[tmfclass](16,17)
\class(16,17)
\class(16,17)
\class(16,17)
\class(16,17)
\class[tmfclass](16,18)
\class(16,18)
\class(16,18)
\class(16,18)
\class(16,18)
\class[tmfclass](16,19)
\class(16,19)
\class(16,19)
\class(16,19)
\class(16,19)
\class[tmfclass](16,20)
\class(16,20)
\class(16,20)
\class(16,20)
\class(16,20)
\class[tmfclass](16,21)
\class(16,21)
\class(16,21)
\class(16,21)
\class(16,21)
\class[tmfclass](16,22)
\class(16,22)
\class(16,22)
\class(16,22)
\class(16,22)
\class[tmfclass](16,23)
\class(16,23)
\class(16,23)
\class(16,23)
\class(16,23)
\class[tmfclass](16,24)
\class(16,24)
\class(16,24)
\class(16,24)
\class(16,24)
\d[dtwo]2(10,2,1,1)
\d[dtwo]2(11,3,1,1)
\d[dtwo]2(12,3,1,1)
\d[dtwo]2(12,4,1,1)
\d[dtwo]2(12,5,1,1)
\d[dtwo]2(15,0,1,1)
\d[dtwo]2(15,1,1,1)
\d[dtwo]2(15,2,1,1)
\d[dtwo]2(15,3,1,1)
\d[dtwo]2(15,3,2,2)
\d[dtwo]2(15,4,1,1)
\d[dtwo]2(15,4,2,2)
\d[dtwo]2(15,5,1,1)
\d[dtwo]2(15,6,1,1)
\d[dtwo]2(15,7,1,1)
\d[dtwo]2(15,8,1,1)
\d[dtwo]2(15,9,1,1)
\d[dtwo]2(15,10,1,1)
\d[dtwo]2(15,11,1,1)
\d[dtwo]2(15,12,1,1)
\d[dtwo]2(15,13,1,1)
\d[dtwo]2(15,14,1,1)
\d[dtwo]2(15,15,1,1)
\d[dtwo]2(15,16,1,1)
\d[dtwo]2(15,17,1,1)
\d[dtwo]2(15,18,1,1)
\d[dtwo]2(15,19,1,1)
\d[dtwo]2(15,20,1,1)
\d[dtwo]2(15,21,1,1)
\d[dtwo]2(15,22,1,1)
\d[dtwo]2(16,0,1,2)
\d[dtwo]2(16,1,1,3)
\d[dthr]3(4,0,1,1)
\d[dfour]4(16,1,2,2)
\structline(0,0,1)(0,1,1)
\structline(0,0,1)(1,1,1)
\structline(0,0,1)(3,1,1)
\structline(0,1,1)(0,2,1)
\structline(0,1,1)(3,2,1)
\structline(0,2,1)(0,3,1)
\structline(0,2,1)(3,3,1)
\structline(0,3,1)(0,4,1)
\structline(0,4,1)(0,5,1)
\structline(0,5,1)(0,6,1)
\structline(0,6,1)(0,7,1)
\structline(0,7,1)(0,8,1)
\structline(0,8,1)(0,9,1)
\structline(0,9,1)(0,10,1)
\structline(0,10,1)(0,11,1)
\structline(0,11,1)(0,12,1)
\structline(0,12,1)(0,13,1)
\structline(0,13,1)(0,14,1)
\structline(0,14,1)(0,15,1)
\structline(0,15,1)(0,16,1)
\structline(0,16,1)(0,17,1)
\structline(0,17,1)(0,18,1)
\structline(0,18,1)(0,19,1)
\structline(0,19,1)(0,20,1)
\structline(0,20,1)(0,21,1)
\structline(0,21,1)(0,22,1)
\structline(0,22,1)(0,23,1)
\structline(0,23,1)(0,24,1)
\structline(1,0,1)(2,1,1)
\structline(1,0,1)(4,1,1)
\structline(1,1,1)(2,2,1)
\structline(2,2,1)(3,3,1)
\structline(3,1,1)(3,2,1)
\structline(3,1,1)(6,2,1)
\structline(3,2,1)(3,3,1)
\structline(4,0,1)(4,1,1)
\structline(4,2,1)(4,3,1)
\structline(4,3,1)(4,4,1)
\structline(4,4,1)(4,5,1)
\structline(4,5,1)(4,6,1)
\structline(4,6,1)(4,7,1)
\structline(4,7,1)(4,8,1)
\structline(4,8,1)(4,9,1)
\structline(4,9,1)(4,10,1)
\structline(4,10,1)(4,11,1)
\structline(4,11,1)(4,12,1)
\structline(4,12,1)(4,13,1)
\structline(4,13,1)(4,14,1)
\structline(4,14,1)(4,15,1)
\structline(4,15,1)(4,16,1)
\structline(4,16,1)(4,17,1)
\structline(4,17,1)(4,18,1)
\structline(4,18,1)(4,19,1)
\structline(4,19,1)(4,20,1)
\structline(4,20,1)(4,21,1)
\structline(4,21,1)(4,22,1)
\structline(4,22,1)(4,23,1)
\structline(4,23,1)(4,24,1)
\structline(8,1,1)(8,2,1)
\structline(8,1,1)(9,2,1)
\structline(8,2,1)(8,3,3)
\structline(8,3,2)(8,4,2)
\structline(8,3,3)(8,4,3)
\structline(8,3,1)(9,4,2)
\structline(8,3,2)(9,4,1)
\structline(8,4,1)(8,5,1)
\structline(8,4,2)(8,5,2)
\structline(8,4,3)(8,5,3)
\structline(8,4,1)(9,5,1)
\structline(8,4,1)(11,5,1)
\structline(8,5,1)(8,6,1)
\structline(8,5,2)(8,6,2)
\structline(8,5,3)(8,6,3)
\structline(8,5,1)(11,6,1)
\structline(8,6,1)(8,7,1)
\structline(8,6,2)(8,7,2)
\structline(8,6,3)(8,7,3)
\structline(8,6,1)(11,7,1)
\structline(8,7,1)(8,8,1)
\structline(8,7,2)(8,8,2)
\structline(8,7,3)(8,8,3)
\structline(8,8,1)(8,9,1)
\structline(8,8,2)(8,9,2)
\structline(8,8,3)(8,9,3)
\structline(8,9,1)(8,10,1)
\structline(8,9,2)(8,10,2)
\structline(8,9,3)(8,10,3)
\structline(8,10,1)(8,11,1)
\structline(8,10,2)(8,11,2)
\structline(8,10,3)(8,11,3)
\structline(8,11,1)(8,12,1)
\structline(8,11,2)(8,12,2)
\structline(8,11,3)(8,12,3)
\structline(8,12,1)(8,13,1)
\structline(8,12,2)(8,13,2)
\structline(8,12,3)(8,13,3)
\structline(8,13,1)(8,14,1)
\structline(8,13,2)(8,14,2)
\structline(8,13,3)(8,14,3)
\structline(8,14,1)(8,15,1)
\structline(8,14,2)(8,15,2)
\structline(8,14,3)(8,15,3)
\structline(8,15,1)(8,16,1)
\structline(8,15,2)(8,16,2)
\structline(8,15,3)(8,16,3)
\structline(8,16,1)(8,17,1)
\structline(8,16,2)(8,17,2)
\structline(8,16,3)(8,17,3)
\structline(8,17,1)(8,18,1)
\structline(8,17,2)(8,18,2)
\structline(8,17,3)(8,18,3)
\structline(8,18,1)(8,19,1)
\structline(8,18,2)(8,19,2)
\structline(8,18,3)(8,19,3)
\structline(8,19,1)(8,20,1)
\structline(8,19,2)(8,20,2)
\structline(8,19,3)(8,20,3)
\structline(8,20,1)(8,21,1)
\structline(8,20,2)(8,21,2)
\structline(8,20,3)(8,21,3)
\structline(8,21,1)(8,22,1)
\structline(8,21,2)(8,22,2)
\structline(8,21,3)(8,22,3)
\structline(8,22,1)(8,23,1)
\structline(8,22,2)(8,23,2)
\structline(8,22,3)(8,23,3)
\structline(8,23,1)(8,24,1)
\structline(8,23,2)(8,24,2)
\structline(8,23,3)(8,24,3)
\structline(9,1,1)(10,2,2)
\structline(9,2,1)(10,3,1)
\structline(9,4,1)(10,5,1)
\structline(9,5,1)(10,6,1)
\structline(10,2,1)(10,3,1)
\structline(10,2,1)(11,3,1)
\structline(10,6,1)(11,7,1)
\structline(11,3,1)(12,4,3)
\structline(11,5,1)(11,6,1)
\structline(11,5,1)(14,6,1)
\structline(11,6,1)(11,7,1)
\structline(12,3,1)(12,4,1)
\structline(12,3,2)(12,4,2)
\structline(12,3,1)(15,4,1)
\structline(12,4,1)(12,5,1)
\structline(12,4,2)(12,5,2)
\structline(12,4,1)(15,5,2)
\structline(12,5,1)(12,6,1)
\structline(12,5,2)(12,6,2)
\structline(12,6,1)(12,7,1)
\structline(12,6,2)(12,7,2)
\structline(12,6,3)(12,7,3)
\structline(12,7,1)(12,8,1)
\structline(12,7,2)(12,8,2)
\structline(12,7,3)(12,8,3)
\structline(12,8,1)(12,9,1)
\structline(12,8,2)(12,9,2)
\structline(12,8,3)(12,9,3)
\structline(12,9,1)(12,10,1)
\structline(12,9,2)(12,10,2)
\structline(12,9,3)(12,10,3)
\structline(12,10,1)(12,11,1)
\structline(12,10,2)(12,11,2)
\structline(12,10,3)(12,11,3)
\structline(12,11,1)(12,12,1)
\structline(12,11,2)(12,12,2)
\structline(12,11,3)(12,12,3)
\structline(12,12,1)(12,13,1)
\structline(12,12,2)(12,13,2)
\structline(12,12,3)(12,13,3)
\structline(12,13,1)(12,14,1)
\structline(12,13,2)(12,14,2)
\structline(12,13,3)(12,14,3)
\structline(12,14,1)(12,15,1)
\structline(12,14,2)(12,15,2)
\structline(12,14,3)(12,15,3)
\structline(12,15,1)(12,16,1)
\structline(12,15,2)(12,16,2)
\structline(12,15,3)(12,16,3)
\structline(12,16,1)(12,17,1)
\structline(12,16,2)(12,17,2)
\structline(12,16,3)(12,17,3)
\structline(12,17,1)(12,18,1)
\structline(12,17,2)(12,18,2)
\structline(12,17,3)(12,18,3)
\structline(12,18,1)(12,19,1)
\structline(12,18,2)(12,19,2)
\structline(12,18,3)(12,19,3)
\structline(12,19,1)(12,20,1)
\structline(12,19,2)(12,20,2)
\structline(12,19,3)(12,20,3)
\structline(12,20,1)(12,21,1)
\structline(12,20,2)(12,21,2)
\structline(12,20,3)(12,21,3)
\structline(12,21,1)(12,22,1)
\structline(12,21,2)(12,22,2)
\structline(12,21,3)(12,22,3)
\structline(12,22,1)(12,23,1)
\structline(12,22,2)(12,23,2)
\structline(12,22,3)(12,23,3)
\structline(12,23,1)(12,24,1)
\structline(12,23,2)(12,24,2)
\structline(12,23,3)(12,24,3)
\structline(14,2,1)(14,3,1)
\structline(14,2,1)(15,3,3)
\structline(14,3,1)(14,4,1)
\structline(14,4,1)(14,5,2)
\structline(14,4,2)(14,5,1)
\structline(14,4,2)(15,5,2)
\structline(14,5,1)(14,6,1)
\structline(14,5,2)(14,6,2)
\structline(14,6,2)(14,7,1)
\structline(14,7,1)(14,8,1)
\structline(14,8,1)(14,9,1)
\structline(14,9,1)(14,10,1)
\structline(14,10,1)(14,11,1)
\structline(14,11,1)(14,12,1)
\structline(14,12,1)(14,13,1)
\structline(14,13,1)(14,14,1)
\structline(14,14,1)(14,15,1)
\structline(14,15,1)(14,16,1)
\structline(14,16,1)(14,17,1)
\structline(14,17,1)(14,18,1)
\structline(14,18,1)(14,19,1)
\structline(14,19,1)(14,20,1)
\structline(14,20,1)(14,21,1)
\structline(14,21,1)(14,22,1)
\structline(14,22,1)(14,23,1)
\structline(14,23,1)(14,24,1)
\structline(15,0,1)(15,1,1)
\structline(15,0,1)(16,1,1)
\structline(15,1,1)(15,2,1)
\structline(15,2,1)(15,3,2)
\structline(15,2,2)(15,3,3)
\structline(15,2,2)(16,3,1)
\structline(15,3,1)(15,4,1)
\structline(15,3,2)(15,4,2)
\structline(15,4,1)(15,5,2)
\structline(15,4,2)(15,5,1)
\structline(15,5,1)(15,6,1)
\structline(15,6,1)(15,7,1)
\structline(15,7,1)(15,8,1)
\structline(15,8,1)(15,9,1)
\structline(15,9,1)(15,10,1)
\structline(15,10,1)(15,11,1)
\structline(15,11,1)(15,12,1)
\structline(15,12,1)(15,13,1)
\structline(15,13,1)(15,14,1)
\structline(15,14,1)(15,15,1)
\structline(15,15,1)(15,16,1)
\structline(15,16,1)(15,17,1)
\structline(15,17,1)(15,18,1)
\structline(15,18,1)(15,19,1)
\structline(15,19,1)(15,20,1)
\structline(15,20,1)(15,21,1)
\structline(15,21,1)(15,22,1)
\structline(15,22,1)(15,23,1)
\structline(15,23,1)(15,24,1)
\structline(16,0,1)(16,1,1)
\structline(16,0,1)(16,1,2)
\structline(16,1,2)(16,2,1)
\structline(16,2,1)(16,3,3)
\structline(16,3,2)(16,4,1)
\structline(16,3,3)(16,4,3)
\structline(16,4,1)(16,5,1)
\structline(16,4,2)(16,5,2)
\structline(16,4,3)(16,5,4)
\structline(16,5,1)(16,6,1)
\structline(16,5,2)(16,6,2)
\structline(16,5,4)(16,6,3)
\structline(16,6,1)(16,7,2)
\structline(16,6,2)(16,7,3)
\structline(16,6,3)(16,7,5)
\structline(16,7,2)(16,8,2)
\structline(16,7,3)(16,8,3)
\structline(16,7,4)(16,8,4)
\structline(16,7,5)(16,8,5)
\structline(16,8,1)(16,9,1)
\structline(16,8,2)(16,9,2)
\structline(16,8,3)(16,9,3)
\structline(16,8,4)(16,9,4)
\structline(16,8,5)(16,9,5)
\structline(16,9,1)(16,10,1)
\structline(16,9,2)(16,10,2)
\structline(16,9,3)(16,10,3)
\structline(16,9,4)(16,10,4)
\structline(16,9,5)(16,10,5)
\structline(16,10,1)(16,11,1)
\structline(16,10,2)(16,11,2)
\structline(16,10,3)(16,11,3)
\structline(16,10,4)(16,11,4)
\structline(16,10,5)(16,11,5)
\structline(16,11,1)(16,12,1)
\structline(16,11,2)(16,12,2)
\structline(16,11,3)(16,12,3)
\structline(16,11,4)(16,12,4)
\structline(16,11,5)(16,12,5)
\structline(16,12,1)(16,13,1)
\structline(16,12,2)(16,13,2)
\structline(16,12,3)(16,13,3)
\structline(16,12,4)(16,13,4)
\structline(16,12,5)(16,13,5)
\structline(16,13,1)(16,14,1)
\structline(16,13,2)(16,14,2)
\structline(16,13,3)(16,14,3)
\structline(16,13,4)(16,14,4)
\structline(16,13,5)(16,14,5)
\structline(16,14,1)(16,15,1)
\structline(16,14,2)(16,15,2)
\structline(16,14,3)(16,15,3)
\structline(16,14,4)(16,15,4)
\structline(16,14,5)(16,15,5)
\structline(16,15,1)(16,16,1)
\structline(16,15,2)(16,16,2)
\structline(16,15,3)(16,16,3)
\structline(16,15,4)(16,16,4)
\structline(16,15,5)(16,16,5)
\structline(16,16,1)(16,17,1)
\structline(16,16,2)(16,17,2)
\structline(16,16,3)(16,17,3)
\structline(16,16,4)(16,17,4)
\structline(16,16,5)(16,17,5)
\structline(16,17,1)(16,18,1)
\structline(16,17,2)(16,18,2)
\structline(16,17,3)(16,18,3)
\structline(16,17,4)(16,18,4)
\structline(16,17,5)(16,18,5)
\structline(16,18,1)(16,19,1)
\structline(16,18,2)(16,19,2)
\structline(16,18,3)(16,19,3)
\structline(16,18,4)(16,19,4)
\structline(16,18,5)(16,19,5)
\structline(16,19,1)(16,20,1)
\structline(16,19,2)(16,20,2)
\structline(16,19,3)(16,20,3)
\structline(16,19,4)(16,20,4)
\structline(16,19,5)(16,20,5)
\structline(16,20,1)(16,21,1)
\structline(16,20,2)(16,21,2)
\structline(16,20,3)(16,21,3)
\structline(16,20,4)(16,21,4)
\structline(16,20,5)(16,21,5)
\structline(16,21,1)(16,22,1)
\structline(16,21,2)(16,22,2)
\structline(16,21,3)(16,22,3)
\structline(16,21,4)(16,22,4)
\structline(16,21,5)(16,22,5)
\structline(16,22,1)(16,23,1)
\structline(16,22,2)(16,23,2)
\structline(16,22,3)(16,23,3)
\structline(16,22,4)(16,23,4)
\structline(16,22,5)(16,23,5)
\structline(16,23,1)(16,24,1)
\structline(16,23,2)(16,24,2)
\structline(16,23,3)(16,24,3)
\structline(16,23,4)(16,24,4)
\structline(16,23,5)(16,24,5)
\end{sseqdata}

\begin{sseqdata}[name = AdamsEtwoB, adamswide, x range = {16}{25}, y range = {0}{14}]
\class(16,0)
\classoptions["{\scriptscriptstyle m_{22}}" below, label distance = -1.5pt](16,0,1)
\class(16,1)
\class(16,1)
\class(16,2)
\class(16,3)
\class(16,3)
\classoptions["{\scriptscriptstyle c_0[10+11]}" below right, label distance = -1.5pt](16,3,2)
\class(16,3)
\class(16,4)
\class(16,4)
\classoptions["{\scriptscriptstyle w_1[12]}" above right, label distance = -1.5pt](16,4,2)
\class(16,4)
\class(16,5)
\class(16,5)
\class(16,5)
\class(16,5)
\class(16,6)
\class(16,6)
\class(16,6)
\class[tmfclass](16,7)
\class(16,7)
\class(16,7)
\class(16,7)
\class(16,7)
\class[tmfclass](16,8)
\class(16,8)
\class(16,8)
\class(16,8)
\class(16,8)
\class[tmfclass](16,9)
\class(16,9)
\class(16,9)
\class(16,9)
\class(16,9)
\class[tmfclass](16,10)
\class(16,10)
\class(16,10)
\class(16,10)
\class(16,10)
\class[tmfclass](16,11)
\class(16,11)
\class(16,11)
\class(16,11)
\class(16,11)
\class[tmfclass](16,12)
\class(16,12)
\class(16,12)
\class(16,12)
\class(16,12)
\class[tmfclass](16,13)
\class(16,13)
\class(16,13)
\class(16,13)
\class(16,13)
\class[tmfclass](16,14)
\class(16,14)
\class(16,14)
\class(16,14)
\class(16,14)
\class[tmfclass](16,15)
\class(16,15)
\class(16,15)
\class(16,15)
\class(16,15)
\class[tmfclass](16,16)
\class(16,16)
\class(16,16)
\class(16,16)
\class(16,16)
\class[tmfclass](16,17)
\class(16,17)
\class(16,17)
\class(16,17)
\class(16,17)
\class[tmfclass](16,18)
\class(16,18)
\class(16,18)
\class(16,18)
\class(16,18)
\class[tmfclass](16,19)
\class(16,19)
\class(16,19)
\class(16,19)
\class(16,19)
\class[tmfclass](16,20)
\class(16,20)
\class(16,20)
\class(16,20)
\class(16,20)
\class[tmfclass](16,21)
\class(16,21)
\class(16,21)
\class(16,21)
\class(16,21)
\class[tmfclass](16,22)
\class(16,22)
\class(16,22)
\class(16,22)
\class(16,22)
\class[tmfclass](16,23)
\class(16,23)
\class(16,23)
\class(16,23)
\class(16,23)
\class[tmfclass](16,24)
\class(16,24)
\class(16,24)
\class(16,24)
\class(16,24)
\class(17,1)
\class(17,2)
\class(17,3)
\class(17,4)
\class(17,4)
\class[tmfclass](17,4)
\class(17,4)
\class[tmfclass](17,5)
\class(17,5)
\class[tmfclass](17,6)
\class(17,6)
\class(17,8)
\class[tmfclass](17,8)
\class[tmfclass](17,9)
\class(18,1)
\class(18,2)
\class(18,2)
\class(18,2)
\classoptions["{\scriptscriptstyle h_2^{2}[17]}" below, label distance = -1.5pt](18,2,3)
\class(18,3)
\class(18,3)
\class[tmfclass](18,4)
\class[tmfclass](18,5)
\class(18,5)
\class(18,6)
\class(18,6)
\class(18,7)
\class(18,9)
\class[tmfclass](18,10)
\class(19,1)
\class(19,2)
\class(19,3)
\class(19,5)
\class(19,7)
\class[tmfclass](19,9)
\class[tmfclass](19,10)
\class[tmfclass](19,11)
\class(20,3)
\classoptions["{\scriptscriptstyle c_0[16]}" below, label distance = -1.5pt](20,3,1)
\class(20,4)
\class[tmfclass](20,4)
\class(20,4)
\class[tmfclass](20,5)
\class(20,5)
\class[tmfclass](20,6)
\class(20,6)
\class(20,6)
\class[tmfclass](20,7)
\class(20,7)
\class(20,7)
\class(20,7)
\class[tmfclass](20,8)
\class(20,8)
\class(20,8)
\class(20,8)
\class(20,8)
\class[tmfclass](20,9)
\class(20,9)
\class(20,9)
\class(20,9)
\class[tmfclass](20,10)
\class(20,10)
\class(20,10)
\class(20,10)
\class(20,10)
\class[tmfclass](20,11)
\class(20,11)
\class(20,11)
\class(20,11)
\class(20,11)
\class[tmfclass](20,12)
\class(20,12)
\class(20,12)
\class(20,12)
\class(20,12)
\class[tmfclass](20,13)
\class(20,13)
\class(20,13)
\class(20,13)
\class(20,13)
\class[tmfclass](20,14)
\class(20,14)
\class(20,14)
\class(20,14)
\class(20,14)
\class[tmfclass](20,15)
\class(20,15)
\class(20,15)
\class(20,15)
\class(20,15)
\class[tmfclass](20,16)
\class(20,16)
\class(20,16)
\class(20,16)
\class(20,16)
\class[tmfclass](20,17)
\class(20,17)
\class(20,17)
\class(20,17)
\class(20,17)
\class[tmfclass](20,18)
\class(20,18)
\class(20,18)
\class(20,18)
\class(20,18)
\class[tmfclass](20,19)
\class(20,19)
\class(20,19)
\class(20,19)
\class(20,19)
\class[tmfclass](20,20)
\class(20,20)
\class(20,20)
\class(20,20)
\class(20,20)
\class[tmfclass](20,21)
\class(20,21)
\class(20,21)
\class(20,21)
\class(20,21)
\class[tmfclass](20,22)
\class(20,22)
\class(20,22)
\class(20,22)
\class(20,22)
\class[tmfclass](20,23)
\class(20,23)
\class(20,23)
\class(20,23)
\class(20,23)
\class[tmfclass](20,24)
\class(20,24)
\class(20,24)
\class(20,24)
\class(20,24)
\class(21,2)
\class(21,3)
\class(21,4)
\class[tmfclass](21,5)
\class(22,2)
\class(22,4)
\classoptions["{\scriptscriptstyle d_0[12]}" below, label distance = -1.5pt](22,4,1)
\class(22,5)
\class(22,6)
\class(22,7)
\class(22,8)
\class[tmfclass](22,8)
\class[tmfclass](22,9)
\class(22,9)
\class[tmfclass](22,10)
\class(22,10)
\class(22,11)
\class(22,12)
\class(22,13)
\class(22,14)
\class(22,15)
\class(22,16)
\class(22,17)
\class(22,18)
\class(22,19)
\class(22,20)
\class(22,21)
\class(22,22)
\class(22,23)
\class(22,24)
\class(23,3)
\class(23,4)
\class(23,5)
\class(23,6)
\class(23,6)
\class[tmfclass](23,7)
\class(23,7)
\class(23,7)
\class[tmfclass](23,8)
\class(23,8)
\class(23,9)
\class[tmfclass](23,9)
\class(23,10)
\class(23,11)
\class(23,12)
\class(23,13)
\class(23,14)
\class(23,15)
\class(23,16)
\class(23,17)
\class(23,18)
\class(23,19)
\class(23,20)
\class(23,21)
\class(23,22)
\class(23,23)
\class(23,24)
\class(24,3)
\class(24,4)
\class(24,4)
\class(24,5)
\classoptions["{\scriptscriptstyle h_0^{2}\alpha[17]}" below left, label distance = -1.5pt](24,5,1)
\class(24,5)
\class(24,5)
\class[tmfclass](24,6)
\class(24,6)
\class(24,6)
\class(24,7)
\class[tmfclass](24,7)
\class(24,7)
\class(24,7)
\class(24,7)
\class[tmfclass](24,8)
\class(24,8)
\class(24,8)
\class(24,8)
\class(24,8)
\class[tmfclass](24,9)
\class(24,9)
\class(24,9)
\class(24,9)
\class(24,9)
\class(24,9)
\class[tmfclass](24,10)
\class(24,10)
\class(24,10)
\class(24,10)
\class(24,10)
\class[tmfclass](24,11)
\class[tmfclass](24,11)
\class(24,11)
\class(24,11)
\class(24,11)
\class(24,11)
\class(24,11)
\class[tmfclass](24,12)
\class[tmfclass](24,12)
\class(24,12)
\class(24,12)
\class(24,12)
\class(24,12)
\class(24,12)
\class[tmfclass](24,13)
\class[tmfclass](24,13)
\class(24,13)
\class(24,13)
\class(24,13)
\class(24,13)
\class(24,13)
\class[tmfclass](24,14)
\class[tmfclass](24,14)
\class(24,14)
\class(24,14)
\class(24,14)
\class(24,14)
\class(24,14)
\class[tmfclass](24,15)
\class[tmfclass](24,15)
\class(24,15)
\class(24,15)
\class(24,15)
\class(24,15)
\class(24,15)
\class[tmfclass](24,16)
\class[tmfclass](24,16)
\class(24,16)
\class(24,16)
\class(24,16)
\class(24,16)
\class(24,16)
\class[tmfclass](24,17)
\class[tmfclass](24,17)
\class(24,17)
\class(24,17)
\class(24,17)
\class(24,17)
\class(24,17)
\class[tmfclass](24,18)
\class[tmfclass](24,18)
\class(24,18)
\class(24,18)
\class(24,18)
\class(24,18)
\class(24,18)
\class[tmfclass](24,19)
\class[tmfclass](24,19)
\class(24,19)
\class(24,19)
\class(24,19)
\class(24,19)
\class(24,19)
\class[tmfclass](24,20)
\class[tmfclass](24,20)
\class(24,20)
\class(24,20)
\class(24,20)
\class(24,20)
\class(24,20)
\class[tmfclass](24,21)
\class[tmfclass](24,21)
\class(24,21)
\class(24,21)
\class(24,21)
\class(24,21)
\class(24,21)
\class[tmfclass](24,22)
\class[tmfclass](24,22)
\class(24,22)
\class(24,22)
\class(24,22)
\class(24,22)
\class(24,22)
\class[tmfclass](24,23)
\class[tmfclass](24,23)
\class(24,23)
\class(24,23)
\class(24,23)
\class(24,23)
\class(24,23)
\class[tmfclass](24,24)
\class[tmfclass](24,24)
\class(24,24)
\class(24,24)
\class(24,24)
\class(24,24)
\class(24,24)
\class(25,4)
\class(25,5)
\class[tmfclass](25,5)
\class(25,6)
\classoptions["{\scriptscriptstyle h_0h_2w_1[19]}" above left, label distance = -1.5pt](25,6,1)
\class(25,6)
\class(25,7)
\class(25,8)
\class(25,8)
\class[tmfclass](25,8)
\class(25,8)
\class[tmfclass](25,9)
\class(25,9)
\class[tmfclass](25,10)
\class(25,10)
\class(25,12)
\class[tmfclass](25,12)
\class[tmfclass](25,13)
\d[dtwo]2(17,1,1,1)
\d[dtwo]2(18,1,1,1)
\d[dtwo]2(18,2,1,1)
\d[dtwo]2(18,2,2,2)
\d[dtwo]2(18,4,1,1)
\d[dtwo]2(18,6,1,1)
\d[dtwo]2(19,1,1,1)
\d[dtwo]2(19,7,1,1)
\d[dtwo]2(20,3,1,1)
\d[dtwo]2(20,7,1,1)
\d[dtwo]2(20,8,1,1)
\d[dtwo]2(20,9,1,1)
\d[dtwo]2(21,2,1,1)
\d[dtwo]2(22,2,1,1)
\d[dtwo]2(23,3,1,1)
\d[dtwo]2(23,4,1,1)
\d[dtwo]2(23,5,1,1)
\d[dtwo]2(23,6,1,1)
\d[dtwo]2(23,7,1,1)
\d[dtwo]2(23,7,2,2)
\d[dtwo]2(23,8,1,1)
\d[dtwo]2(23,8,2,2)
\d[dtwo]2(23,9,1,1)
\d[dtwo]2(23,10,1,1)
\d[dtwo]2(23,11,1,1)
\d[dtwo]2(23,12,1,1)
\d[dtwo]2(23,13,1,1)
\d[dtwo]2(23,14,1,1)
\d[dtwo]2(23,15,1,1)
\d[dtwo]2(23,16,1,1)
\d[dtwo]2(23,17,1,1)
\d[dtwo]2(23,18,1,1)
\d[dtwo]2(23,19,1,1)
\d[dtwo]2(23,20,1,1)
\d[dtwo]2(23,21,1,1)
\d[dtwo]2(23,22,1,1)
\d[dtwo]2(24,4,1,2)
\d[dtwo]2(24,5,1,3)
\d[dtwo]2(25,5,1,1)
\d[dthr]3(17,4,3,1)
\d[dthr]3(18,2,3,1)
\d[dthr]3(18,5,1,2)
\d[dthr]3(21,3,1,1)
\d[dthr]3(24,6,1,2)
\d[dthr]3(25,8,3,1)
\structline(16,0,1)(16,1,1)
\structline(16,0,1)(16,1,2)
\structline(16,0,1)(17,1,1)
\structline(16,0,1)(19,1,1)
\structline(16,1,2)(16,2,1)
\structline(16,1,1)(17,2,1)
\structline(16,1,2)(17,2,1)
\structline(16,1,2)(19,2,1)
\structline(16,2,1)(16,3,3)
\structline(16,2,1)(19,3,1)
\structline(16,3,2)(16,4,1)
\structline(16,3,3)(16,4,3)
\structline(16,3,1)(17,4,1)
\structline(16,3,1)(17,4,2)
\structline(16,3,2)(17,4,1)
\structline(16,3,2)(17,4,4)
\structline(16,4,1)(16,5,1)
\structline(16,4,2)(16,5,2)
\structline(16,4,3)(16,5,4)
\structline(16,4,2)(17,5,2)
\structline(16,4,2)(19,5,1)
\structline(16,5,1)(16,6,1)
\structline(16,5,2)(16,6,2)
\structline(16,5,4)(16,6,3)
\structline(16,5,3)(17,6,2)
\structline(16,6,1)(16,7,2)
\structline(16,6,2)(16,7,3)
\structline(16,6,3)(16,7,5)
\structline(16,7,2)(16,8,2)
\structline(16,7,3)(16,8,3)
\structline(16,7,4)(16,8,4)
\structline(16,7,5)(16,8,5)
\structline(16,7,1)(17,8,2)
\structline(16,7,4)(17,8,1)
\structline(16,8,1)(16,9,1)
\structline(16,8,2)(16,9,2)
\structline(16,8,3)(16,9,3)
\structline(16,8,4)(16,9,4)
\structline(16,8,5)(16,9,5)
\structline(16,8,1)(17,9,1)
\structline(16,8,1)(19,9,1)
\structline(16,9,1)(16,10,1)
\structline(16,9,2)(16,10,2)
\structline(16,9,3)(16,10,3)
\structline(16,9,4)(16,10,4)
\structline(16,9,5)(16,10,5)
\structline(16,9,1)(19,10,1)
\structline(16,10,1)(16,11,1)
\structline(16,10,2)(16,11,2)
\structline(16,10,3)(16,11,3)
\structline(16,10,4)(16,11,4)
\structline(16,10,5)(16,11,5)
\structline(16,10,1)(19,11,1)
\structline(16,11,1)(16,12,1)
\structline(16,11,2)(16,12,2)
\structline(16,11,3)(16,12,3)
\structline(16,11,4)(16,12,4)
\structline(16,11,5)(16,12,5)
\structline(16,12,1)(16,13,1)
\structline(16,12,2)(16,13,2)
\structline(16,12,3)(16,13,3)
\structline(16,12,4)(16,13,4)
\structline(16,12,5)(16,13,5)
\structline(16,13,1)(16,14,1)
\structline(16,13,2)(16,14,2)
\structline(16,13,3)(16,14,3)
\structline(16,13,4)(16,14,4)
\structline(16,13,5)(16,14,5)
\structline(16,14,1)(16,15,1)
\structline(16,14,2)(16,15,2)
\structline(16,14,3)(16,15,3)
\structline(16,14,4)(16,15,4)
\structline(16,14,5)(16,15,5)
\structline(16,15,1)(16,16,1)
\structline(16,15,2)(16,16,2)
\structline(16,15,3)(16,16,3)
\structline(16,15,4)(16,16,4)
\structline(16,15,5)(16,16,5)
\structline(16,16,1)(16,17,1)
\structline(16,16,2)(16,17,2)
\structline(16,16,3)(16,17,3)
\structline(16,16,4)(16,17,4)
\structline(16,16,5)(16,17,5)
\structline(16,17,1)(16,18,1)
\structline(16,17,2)(16,18,2)
\structline(16,17,3)(16,18,3)
\structline(16,17,4)(16,18,4)
\structline(16,17,5)(16,18,5)
\structline(16,18,1)(16,19,1)
\structline(16,18,2)(16,19,2)
\structline(16,18,3)(16,19,3)
\structline(16,18,4)(16,19,4)
\structline(16,18,5)(16,19,5)
\structline(16,19,1)(16,20,1)
\structline(16,19,2)(16,20,2)
\structline(16,19,3)(16,20,3)
\structline(16,19,4)(16,20,4)
\structline(16,19,5)(16,20,5)
\structline(16,20,1)(16,21,1)
\structline(16,20,2)(16,21,2)
\structline(16,20,3)(16,21,3)
\structline(16,20,4)(16,21,4)
\structline(16,20,5)(16,21,5)
\structline(16,21,1)(16,22,1)
\structline(16,21,2)(16,22,2)
\structline(16,21,3)(16,22,3)
\structline(16,21,4)(16,22,4)
\structline(16,21,5)(16,22,5)
\structline(16,22,1)(16,23,1)
\structline(16,22,2)(16,23,2)
\structline(16,22,3)(16,23,3)
\structline(16,22,4)(16,23,4)
\structline(16,22,5)(16,23,5)
\structline(16,23,1)(16,24,1)
\structline(16,23,2)(16,24,2)
\structline(16,23,3)(16,24,3)
\structline(16,23,4)(16,24,4)
\structline(16,23,5)(16,24,5)
\structline(17,1,1)(18,2,1)
\structline(17,1,1)(18,2,2)
\structline(17,2,1)(18,3,2)
\structline(17,3,1)(17,4,2)
\structline(17,3,1)(20,4,1)
\structline(17,4,3)(17,5,1)
\structline(17,4,3)(18,5,1)
\structline(17,4,4)(18,5,2)
\structline(17,4,3)(20,5,1)
\structline(17,5,1)(17,6,1)
\structline(17,5,2)(18,6,2)
\structline(17,5,1)(20,6,1)
\structline(17,6,2)(18,7,1)
\structline(17,8,1)(18,9,1)
\structline(17,9,1)(18,10,1)
\structline(18,1,1)(18,2,2)
\structline(18,1,1)(21,2,1)
\structline(18,2,1)(18,3,2)
\structline(18,2,2)(18,3,2)
\structline(18,2,3)(18,3,1)
\structline(18,2,1)(19,3,1)
\structline(18,2,3)(21,3,1)
\structline(18,3,1)(21,4,1)
\structline(18,4,1)(18,5,1)
\structline(18,4,1)(21,5,1)
\structline(18,6,1)(18,7,1)
\structline(18,6,1)(19,7,1)
\structline(18,10,1)(19,11,1)
\structline(19,1,1)(19,2,1)
\structline(19,1,1)(22,2,1)
\structline(19,2,1)(19,3,1)
\structline(19,7,1)(20,8,5)
\structline(19,9,1)(19,10,1)
\structline(19,9,1)(22,10,1)
\structline(19,10,1)(19,11,1)
\structline(20,3,1)(20,4,3)
\structline(20,3,1)(21,4,1)
\structline(20,4,2)(20,5,1)
\structline(20,4,3)(20,5,2)
\structline(20,4,2)(21,5,1)
\structline(20,5,1)(20,6,1)
\structline(20,5,2)(20,6,2)
\structline(20,5,2)(20,6,3)
\structline(20,6,2)(20,7,2)
\structline(20,6,3)(20,7,3)
\structline(20,7,1)(20,8,1)
\structline(20,7,2)(20,8,2)
\structline(20,7,3)(20,8,3)
\structline(20,7,4)(20,8,4)
\structline(20,7,1)(23,8,1)
\structline(20,8,1)(20,9,1)
\structline(20,8,2)(20,9,2)
\structline(20,8,3)(20,9,3)
\structline(20,8,4)(20,9,4)
\structline(20,8,1)(23,9,2)
\structline(20,9,1)(20,10,1)
\structline(20,9,2)(20,10,2)
\structline(20,9,3)(20,10,3)
\structline(20,9,4)(20,10,4)
\structline(20,10,1)(20,11,1)
\structline(20,10,2)(20,11,2)
\structline(20,10,3)(20,11,3)
\structline(20,10,4)(20,11,4)
\structline(20,10,5)(20,11,5)
\structline(20,11,1)(20,12,1)
\structline(20,11,2)(20,12,2)
\structline(20,11,3)(20,12,3)
\structline(20,11,4)(20,12,4)
\structline(20,11,5)(20,12,5)
\structline(20,12,1)(20,13,1)
\structline(20,12,2)(20,13,2)
\structline(20,12,3)(20,13,3)
\structline(20,12,4)(20,13,4)
\structline(20,12,5)(20,13,5)
\structline(20,13,1)(20,14,1)
\structline(20,13,2)(20,14,2)
\structline(20,13,3)(20,14,3)
\structline(20,13,4)(20,14,4)
\structline(20,13,5)(20,14,5)
\structline(20,14,1)(20,15,1)
\structline(20,14,2)(20,15,2)
\structline(20,14,3)(20,15,3)
\structline(20,14,4)(20,15,4)
\structline(20,14,5)(20,15,5)
\structline(20,15,1)(20,16,1)
\structline(20,15,2)(20,16,2)
\structline(20,15,3)(20,16,3)
\structline(20,15,4)(20,16,4)
\structline(20,15,5)(20,16,5)
\structline(20,16,1)(20,17,1)
\structline(20,16,2)(20,17,2)
\structline(20,16,3)(20,17,3)
\structline(20,16,4)(20,17,4)
\structline(20,16,5)(20,17,5)
\structline(20,17,1)(20,18,1)
\structline(20,17,2)(20,18,2)
\structline(20,17,3)(20,18,3)
\structline(20,17,4)(20,18,4)
\structline(20,17,5)(20,18,5)
\structline(20,18,1)(20,19,1)
\structline(20,18,2)(20,19,2)
\structline(20,18,3)(20,19,3)
\structline(20,18,4)(20,19,4)
\structline(20,18,5)(20,19,5)
\structline(20,19,1)(20,20,1)
\structline(20,19,2)(20,20,2)
\structline(20,19,3)(20,20,3)
\structline(20,19,4)(20,20,4)
\structline(20,19,5)(20,20,5)
\structline(20,20,1)(20,21,1)
\structline(20,20,2)(20,21,2)
\structline(20,20,3)(20,21,3)
\structline(20,20,4)(20,21,4)
\structline(20,20,5)(20,21,5)
\structline(20,21,1)(20,22,1)
\structline(20,21,2)(20,22,2)
\structline(20,21,3)(20,22,3)
\structline(20,21,4)(20,22,4)
\structline(20,21,5)(20,22,5)
\structline(20,22,1)(20,23,1)
\structline(20,22,2)(20,23,2)
\structline(20,22,3)(20,23,3)
\structline(20,22,4)(20,23,4)
\structline(20,22,5)(20,23,5)
\structline(20,23,1)(20,24,1)
\structline(20,23,2)(20,24,2)
\structline(20,23,3)(20,24,3)
\structline(20,23,4)(20,24,4)
\structline(20,23,5)(20,24,5)
\structline(21,3,1)(21,4,1)
\structline(21,3,1)(24,4,2)
\structline(22,4,1)(22,5,1)
\structline(22,6,1)(22,7,1)
\structline(22,6,1)(23,7,3)
\structline(22,6,1)(25,7,1)
\structline(22,7,1)(22,8,1)
\structline(22,7,1)(25,8,1)
\structline(22,8,1)(22,9,2)
\structline(22,8,2)(22,9,1)
\structline(22,8,2)(23,9,2)
\structline(22,8,2)(25,9,1)
\structline(22,9,1)(22,10,1)
\structline(22,9,2)(22,10,2)
\structline(22,9,1)(25,10,1)
\structline(22,10,2)(22,11,1)
\structline(22,11,1)(22,12,1)
\structline(22,12,1)(22,13,1)
\structline(22,13,1)(22,14,1)
\structline(22,14,1)(22,15,1)
\structline(22,15,1)(22,16,1)
\structline(22,16,1)(22,17,1)
\structline(22,17,1)(22,18,1)
\structline(22,18,1)(22,19,1)
\structline(22,19,1)(22,20,1)
\structline(22,20,1)(22,21,1)
\structline(22,21,1)(22,22,1)
\structline(22,22,1)(22,23,1)
\structline(22,23,1)(22,24,1)
\structline(23,3,1)(24,4,2)
\structline(23,4,1)(23,5,1)
\structline(23,4,1)(24,5,1)
\structline(23,4,1)(24,5,2)
\structline(23,5,1)(23,6,1)
\structline(23,6,1)(23,7,2)
\structline(23,6,2)(23,7,3)
\structline(23,6,2)(24,7,1)
\structline(23,7,1)(23,8,1)
\structline(23,7,2)(23,8,2)
\structline(23,8,1)(23,9,2)
\structline(23,8,2)(23,9,1)
\structline(23,9,1)(23,10,1)
\structline(23,10,1)(23,11,1)
\structline(23,11,1)(23,12,1)
\structline(23,12,1)(23,13,1)
\structline(23,13,1)(23,14,1)
\structline(23,14,1)(23,15,1)
\structline(23,15,1)(23,16,1)
\structline(23,16,1)(23,17,1)
\structline(23,17,1)(23,18,1)
\structline(23,18,1)(23,19,1)
\structline(23,19,1)(23,20,1)
\structline(23,20,1)(23,21,1)
\structline(23,21,1)(23,22,1)
\structline(23,22,1)(23,23,1)
\structline(23,23,1)(23,24,1)
\structline(24,3,1)(25,4,1)
\structline(24,4,1)(24,5,1)
\structline(24,4,1)(24,5,3)
\structline(24,4,1)(25,5,1)
\structline(24,5,1)(24,6,2)
\structline(24,5,2)(24,6,2)
\structline(24,5,3)(24,6,2)
\structline(24,5,3)(24,6,3)
\structline(24,5,2)(25,6,1)
\structline(24,5,2)(25,6,2)
\structline(24,6,1)(24,7,2)
\structline(24,6,2)(24,7,3)
\structline(24,6,3)(24,7,3)
\structline(24,6,3)(24,7,4)
\structline(24,6,3)(24,7,5)
\structline(24,7,2)(24,8,1)
\structline(24,7,3)(24,8,2)
\structline(24,7,4)(24,8,3)
\structline(24,7,5)(24,8,4)
\structline(24,7,1)(25,8,1)
\structline(24,7,1)(25,8,2)
\structline(24,7,4)(25,8,2)
\structline(24,7,4)(25,8,4)
\structline(24,7,5)(25,8,2)
\structline(24,7,5)(25,8,4)
\structline(24,8,1)(24,9,1)
\structline(24,8,2)(24,9,2)
\structline(24,8,3)(24,9,3)
\structline(24,8,4)(24,9,4)
\structline(24,8,5)(24,9,5)
\structline(24,8,5)(25,9,2)
\structline(24,9,1)(24,10,1)
\structline(24,9,2)(24,10,2)
\structline(24,9,3)(24,10,3)
\structline(24,9,4)(24,10,4)
\structline(24,9,5)(24,10,5)
\structline(24,9,6)(25,10,2)
\structline(24,10,1)(24,11,2)
\structline(24,10,2)(24,11,3)
\structline(24,10,3)(24,11,4)
\structline(24,10,4)(24,11,5)
\structline(24,10,5)(24,11,6)
\structline(24,11,2)(24,12,1)
\structline(24,11,3)(24,12,3)
\structline(24,11,4)(24,12,4)
\structline(24,11,5)(24,12,5)
\structline(24,11,6)(24,12,6)
\structline(24,11,7)(24,12,7)
\structline(24,11,1)(25,12,2)
\structline(24,11,7)(25,12,1)
\structline(24,12,1)(24,13,1)
\structline(24,12,2)(24,13,2)
\structline(24,12,3)(24,13,3)
\structline(24,12,4)(24,13,4)
\structline(24,12,5)(24,13,5)
\structline(24,12,6)(24,13,6)
\structline(24,12,7)(24,13,7)
\structline(24,12,2)(25,13,1)
\structline(24,13,1)(24,14,1)
\structline(24,13,2)(24,14,2)
\structline(24,13,3)(24,14,3)
\structline(24,13,4)(24,14,4)
\structline(24,13,5)(24,14,5)
\structline(24,13,6)(24,14,6)
\structline(24,13,7)(24,14,7)
\structline(24,14,1)(24,15,1)
\structline(24,14,2)(24,15,2)
\structline(24,14,3)(24,15,3)
\structline(24,14,4)(24,15,4)
\structline(24,14,5)(24,15,5)
\structline(24,14,6)(24,15,6)
\structline(24,14,7)(24,15,7)
\structline(24,15,1)(24,16,1)
\structline(24,15,2)(24,16,2)
\structline(24,15,3)(24,16,3)
\structline(24,15,4)(24,16,4)
\structline(24,15,5)(24,16,5)
\structline(24,15,6)(24,16,6)
\structline(24,15,7)(24,16,7)
\structline(24,16,1)(24,17,1)
\structline(24,16,2)(24,17,2)
\structline(24,16,3)(24,17,3)
\structline(24,16,4)(24,17,4)
\structline(24,16,5)(24,17,5)
\structline(24,16,6)(24,17,6)
\structline(24,16,7)(24,17,7)
\structline(24,17,1)(24,18,1)
\structline(24,17,2)(24,18,2)
\structline(24,17,3)(24,18,3)
\structline(24,17,4)(24,18,4)
\structline(24,17,5)(24,18,5)
\structline(24,17,6)(24,18,6)
\structline(24,17,7)(24,18,7)
\structline(24,18,1)(24,19,1)
\structline(24,18,2)(24,19,2)
\structline(24,18,3)(24,19,3)
\structline(24,18,4)(24,19,4)
\structline(24,18,5)(24,19,5)
\structline(24,18,6)(24,19,6)
\structline(24,18,7)(24,19,7)
\structline(24,19,1)(24,20,1)
\structline(24,19,2)(24,20,2)
\structline(24,19,3)(24,20,3)
\structline(24,19,4)(24,20,4)
\structline(24,19,5)(24,20,5)
\structline(24,19,6)(24,20,6)
\structline(24,19,7)(24,20,7)
\structline(24,20,1)(24,21,1)
\structline(24,20,2)(24,21,2)
\structline(24,20,3)(24,21,3)
\structline(24,20,4)(24,21,4)
\structline(24,20,5)(24,21,5)
\structline(24,20,6)(24,21,6)
\structline(24,20,7)(24,21,7)
\structline(24,21,1)(24,22,1)
\structline(24,21,2)(24,22,2)
\structline(24,21,3)(24,22,3)
\structline(24,21,4)(24,22,4)
\structline(24,21,5)(24,22,5)
\structline(24,21,6)(24,22,6)
\structline(24,21,7)(24,22,7)
\structline(24,22,1)(24,23,1)
\structline(24,22,2)(24,23,2)
\structline(24,22,3)(24,23,3)
\structline(24,22,4)(24,23,4)
\structline(24,22,5)(24,23,5)
\structline(24,22,6)(24,23,6)
\structline(24,22,7)(24,23,7)
\structline(24,23,1)(24,24,1)
\structline(24,23,2)(24,24,2)
\structline(24,23,3)(24,24,3)
\structline(24,23,4)(24,24,4)
\structline(24,23,5)(24,24,5)
\structline(24,23,6)(24,24,6)
\structline(24,23,7)(24,24,7)
\structline(25,6,1)(25,7,1)
\structline(25,6,2)(25,7,1)
\structline(25,7,1)(25,8,1)
\structline(25,8,3)(25,9,1)
\structline(25,9,1)(25,10,1)
\end{sseqdata}
\begin{figure}[htbp]\centering
\adjustbox{max width=\textwidth}{\printpage[name = AdamsEtwoB, page = 2]}
\caption{The Adams $E_2$-page of $\tmf_0(5)$ in stems $16$ to $25$, with the $d_2$-differentials.}
\label{fig:adamsE2b}
\end{figure}

\begin{sseqdata}[name = AdamsEtwoC, adamswide, x range = {24}{33}, y range = {0}{14}]
\class(24,3)
\class(24,4)
\class(24,4)
\class(24,5)
\classoptions["{\scriptscriptstyle h_0^{2}\alpha[17]}" below right, label distance = -1.5pt](24,5,1)
\class(24,5)
\class(24,5)
\class[tmfclass](24,6)
\class(24,6)
\class(24,6)
\class(24,7)
\class[tmfclass](24,7)
\class(24,7)
\class(24,7)
\class(24,7)
\class[tmfclass](24,8)
\class(24,8)
\class(24,8)
\class(24,8)
\class(24,8)
\class[tmfclass](24,9)
\class(24,9)
\class(24,9)
\class(24,9)
\class(24,9)
\class(24,9)
\class[tmfclass](24,10)
\class(24,10)
\class(24,10)
\class(24,10)
\class(24,10)
\class[tmfclass](24,11)
\class[tmfclass](24,11)
\class(24,11)
\class(24,11)
\class(24,11)
\class(24,11)
\class(24,11)
\class[tmfclass](24,12)
\class[tmfclass](24,12)
\class(24,12)
\class(24,12)
\class(24,12)
\class(24,12)
\class(24,12)
\class[tmfclass](24,13)
\class[tmfclass](24,13)
\class(24,13)
\class(24,13)
\class(24,13)
\class(24,13)
\class(24,13)
\class[tmfclass](24,14)
\class[tmfclass](24,14)
\class(24,14)
\class(24,14)
\class(24,14)
\class(24,14)
\class(24,14)
\class[tmfclass](24,15)
\class[tmfclass](24,15)
\class(24,15)
\class(24,15)
\class(24,15)
\class(24,15)
\class(24,15)
\class[tmfclass](24,16)
\class[tmfclass](24,16)
\class(24,16)
\class(24,16)
\class(24,16)
\class(24,16)
\class(24,16)
\class[tmfclass](24,17)
\class[tmfclass](24,17)
\class(24,17)
\class(24,17)
\class(24,17)
\class(24,17)
\class(24,17)
\class[tmfclass](24,18)
\class[tmfclass](24,18)
\class(24,18)
\class(24,18)
\class(24,18)
\class(24,18)
\class(24,18)
\class[tmfclass](24,19)
\class[tmfclass](24,19)
\class(24,19)
\class(24,19)
\class(24,19)
\class(24,19)
\class(24,19)
\class[tmfclass](24,20)
\class[tmfclass](24,20)
\class(24,20)
\class(24,20)
\class(24,20)
\class(24,20)
\class(24,20)
\class[tmfclass](24,21)
\class[tmfclass](24,21)
\class(24,21)
\class(24,21)
\class(24,21)
\class(24,21)
\class(24,21)
\class[tmfclass](24,22)
\class[tmfclass](24,22)
\class(24,22)
\class(24,22)
\class(24,22)
\class(24,22)
\class(24,22)
\class[tmfclass](24,23)
\class[tmfclass](24,23)
\class(24,23)
\class(24,23)
\class(24,23)
\class(24,23)
\class(24,23)
\class[tmfclass](24,24)
\class[tmfclass](24,24)
\class(24,24)
\class(24,24)
\class(24,24)
\class(24,24)
\class(24,24)
\class(25,4)
\class(25,5)
\class[tmfclass](25,5)
\class(25,6)
\classoptions["{\scriptscriptstyle h_0h_2w_1[19]}" above right, label distance = -1.5pt](25,6,1)
\class(25,6)
\class(25,7)
\class(25,8)
\class(25,8)
\class[tmfclass](25,8)
\class(25,8)
\class[tmfclass](25,9)
\class(25,9)
\class[tmfclass](25,10)
\class(25,10)
\class(25,12)
\class[tmfclass](25,12)
\class[tmfclass](25,13)
\class(26,4)
\classoptions["{\scriptscriptstyle h_1c_0[23]}" below, label distance = -1.5pt](26,4,1)
\class(26,5)
\class(26,5)
\class(26,6)
\class(26,6)
\class(26,6)
\class[tmfclass](26,6)
\class(26,6)
\class[tmfclass](26,7)
\class(26,7)
\class(26,7)
\class(26,7)
\class[tmfclass](26,8)
\class(26,8)
\class(26,9)
\class[tmfclass](26,9)
\class(26,9)
\class(26,10)
\class(26,10)
\class(26,10)
\class(26,11)
\class(26,11)
\class(26,12)
\class(26,13)
\class(26,13)
\class(26,14)
\class[tmfclass](26,14)
\class(26,15)
\class(26,16)
\class(26,17)
\class(26,18)
\class(26,19)
\class(26,20)
\class(26,21)
\class(26,22)
\class(26,23)
\class(26,24)
\class(27,3)
\class(27,4)
\class(27,5)
\class(27,5)
\class(27,5)
\class(27,6)
\class(27,6)
\class[tmfclass](27,6)
\class(27,7)
\class[tmfclass](27,7)
\class(27,7)
\class(27,8)
\class(27,9)
\class(27,9)
\class(27,10)
\class(27,11)
\class(27,11)
\class(27,12)
\class(27,13)
\class[tmfclass](27,13)
\class(27,14)
\class[tmfclass](27,14)
\class(27,15)
\class[tmfclass](27,15)
\class(27,16)
\class(27,17)
\class(27,18)
\class(27,19)
\class(27,20)
\class(27,21)
\class(27,22)
\class(27,23)
\class(27,24)
\class(28,3)
\class(28,4)
\class(28,5)
\classoptions["{\scriptscriptstyle h_2w_1[23]}" below right, label distance = -1.5pt](28,5,1)
\class(28,5)
\class(28,6)
\class(28,6)
\class(28,6)
\class(28,7)
\class(28,7)
\class(28,7)
\class(28,7)
\class(28,8)
\class[tmfclass](28,8)
\class(28,8)
\class(28,8)
\class(28,8)
\class[tmfclass](28,9)
\class(28,9)
\class(28,9)
\class(28,9)
\class[tmfclass](28,10)
\class(28,10)
\class(28,10)
\class(28,10)
\class(28,10)
\class[tmfclass](28,11)
\class(28,11)
\class(28,11)
\class(28,11)
\class(28,11)
\class(28,11)
\class[tmfclass](28,12)
\class(28,12)
\class(28,12)
\class(28,12)
\class(28,12)
\class(28,12)
\class(28,12)
\class[tmfclass](28,13)
\class(28,13)
\class(28,13)
\class(28,13)
\class(28,13)
\class(28,13)
\class[tmfclass](28,14)
\class(28,14)
\class(28,14)
\class(28,14)
\class(28,14)
\class(28,14)
\class(28,14)
\class[tmfclass](28,15)
\class(28,15)
\class(28,15)
\class(28,15)
\class(28,15)
\class(28,15)
\class(28,15)
\class[tmfclass](28,16)
\class(28,16)
\class(28,16)
\class(28,16)
\class(28,16)
\class(28,16)
\class(28,16)
\class[tmfclass](28,17)
\class(28,17)
\class(28,17)
\class(28,17)
\class(28,17)
\class(28,17)
\class(28,17)
\class[tmfclass](28,18)
\class(28,18)
\class(28,18)
\class(28,18)
\class(28,18)
\class(28,18)
\class(28,18)
\class[tmfclass](28,19)
\class(28,19)
\class(28,19)
\class(28,19)
\class(28,19)
\class(28,19)
\class(28,19)
\class[tmfclass](28,20)
\class(28,20)
\class(28,20)
\class(28,20)
\class(28,20)
\class(28,20)
\class(28,20)
\class[tmfclass](28,21)
\class(28,21)
\class(28,21)
\class(28,21)
\class(28,21)
\class(28,21)
\class(28,21)
\class[tmfclass](28,22)
\class(28,22)
\class(28,22)
\class(28,22)
\class(28,22)
\class(28,22)
\class(28,22)
\class[tmfclass](28,23)
\class(28,23)
\class(28,23)
\class(28,23)
\class(28,23)
\class(28,23)
\class(28,23)
\class[tmfclass](28,24)
\class(28,24)
\class(28,24)
\class(28,24)
\class(28,24)
\class(28,24)
\class(28,24)
\class(29,4)
\class(29,5)
\class(29,5)
\class(29,6)
\class(29,6)
\class(29,6)
\class[tmfclass](29,7)
\class(29,7)
\class[tmfclass](29,8)
\class[tmfclass](29,9)
\class(30,3)
\class(30,4)
\class(30,4)
\class(30,5)
\class(30,5)
\class(30,6)
\class[tmfclass](30,6)
\class(30,8)
\class(30,9)
\class(30,10)
\class(30,11)
\class(30,12)
\class[tmfclass](30,12)
\class[tmfclass](30,13)
\class(30,13)
\class[tmfclass](30,14)
\class(30,14)
\class(30,15)
\class(30,16)
\class(30,17)
\class(30,18)
\class(30,19)
\class(30,20)
\class(30,21)
\class(30,22)
\class(30,23)
\class(30,24)
\class(31,3)
\class(31,4)
\class(31,5)
\class(31,6)
\class(31,7)
\class(31,7)
\class(31,8)
\class(31,8)
\class[tmfclass](31,8)
\class(31,9)
\class(31,10)
\class(31,10)
\class[tmfclass](31,11)
\class(31,11)
\class(31,11)
\class[tmfclass](31,12)
\class(31,12)
\class(31,13)
\class[tmfclass](31,13)
\class(31,14)
\class(31,15)
\class(31,16)
\class(31,17)
\class(31,18)
\class(31,19)
\class(31,20)
\class(31,21)
\class(31,22)
\class(31,23)
\class(31,24)
\class(32,4)
\class(32,4)
\classoptions["{\scriptscriptstyle e_0[21]}" below, label distance = -1.5pt](32,4,2)
\class(32,5)
\class(32,5)
\classoptions["{\scriptscriptstyle h_0e_0[21]}" above right, label distance = -1.5pt](32,5,2)
\class(32,5)
\class(32,6)
\class(32,6)
\class(32,6)
\class(32,7)
\class[tmfclass](32,7)
\class[tmfclass](32,7)
\class(32,7)
\classoptions["{\scriptscriptstyle h_0\alpha^{2}[10+11]}" above, label distance = -1.5pt](32,7,4)
\class(32,7)
\class(32,7)
\class(32,8)
\class[tmfclass](32,8)
\class(32,8)
\class(32,8)
\class(32,9)
\class[tmfclass](32,9)
\class(32,9)
\class(32,9)
\class(32,9)
\class(32,9)
\class[tmfclass](32,10)
\class(32,10)
\class(32,10)
\class(32,10)
\class(32,10)
\class(32,11)
\class(32,11)
\class(32,11)
\class[tmfclass](32,11)
\class(32,11)
\class(32,11)
\class(32,11)
\class[tmfclass](32,12)
\class(32,12)
\class(32,12)
\class(32,12)
\class(32,12)
\class(32,12)
\class(32,12)
\class[tmfclass](32,13)
\class(32,13)
\class(32,13)
\class(32,13)
\class(32,13)
\class(32,13)
\class(32,13)
\class(32,13)
\class[tmfclass](32,14)
\class(32,14)
\class(32,14)
\class(32,14)
\class(32,14)
\class(32,14)
\class(32,14)
\class[tmfclass](32,15)
\class[tmfclass](32,15)
\class(32,15)
\class(32,15)
\class(32,15)
\class(32,15)
\class(32,15)
\class(32,15)
\class(32,15)
\class[tmfclass](32,16)
\class[tmfclass](32,16)
\class(32,16)
\class(32,16)
\class(32,16)
\class(32,16)
\class(32,16)
\class(32,16)
\class(32,16)
\class[tmfclass](32,17)
\class[tmfclass](32,17)
\class(32,17)
\class(32,17)
\class(32,17)
\class(32,17)
\class(32,17)
\class(32,17)
\class(32,17)
\class[tmfclass](32,18)
\class[tmfclass](32,18)
\class(32,18)
\class(32,18)
\class(32,18)
\class(32,18)
\class(32,18)
\class(32,18)
\class(32,18)
\class[tmfclass](32,19)
\class[tmfclass](32,19)
\class(32,19)
\class(32,19)
\class(32,19)
\class(32,19)
\class(32,19)
\class(32,19)
\class(32,19)
\class[tmfclass](32,20)
\class[tmfclass](32,20)
\class(32,20)
\class(32,20)
\class(32,20)
\class(32,20)
\class(32,20)
\class(32,20)
\class(32,20)
\class[tmfclass](32,21)
\class[tmfclass](32,21)
\class(32,21)
\class(32,21)
\class(32,21)
\class(32,21)
\class(32,21)
\class(32,21)
\class(32,21)
\class[tmfclass](32,22)
\class[tmfclass](32,22)
\class(32,22)
\class(32,22)
\class(32,22)
\class(32,22)
\class(32,22)
\class(32,22)
\class(32,22)
\class[tmfclass](32,23)
\class[tmfclass](32,23)
\class(32,23)
\class(32,23)
\class(32,23)
\class(32,23)
\class(32,23)
\class(32,23)
\class(32,23)
\class[tmfclass](32,24)
\class[tmfclass](32,24)
\class(32,24)
\class(32,24)
\class(32,24)
\class(32,24)
\class(32,24)
\class(32,24)
\class(32,24)
\class(33,4)
\class(33,4)
\class(33,5)
\class(33,5)
\class(33,5)
\class(33,6)
\class(33,6)
\class[tmfclass](33,8)
\class(33,8)
\class(33,8)
\class(33,9)
\class[tmfclass](33,9)
\class(33,10)
\class(33,10)
\class(33,11)
\class(33,12)
\class(33,12)
\class[tmfclass](33,12)
\class(33,12)
\class[tmfclass](33,13)
\class(33,13)
\class[tmfclass](33,14)
\class(33,14)
\class(33,16)
\class[tmfclass](33,16)
\class[tmfclass](33,17)
\d[dtwo]2(25,5,1,1)
\d[dtwo]2(26,4,1,1)
\d[dtwo]2(26,5,1,1)
\d[dtwo]2(26,6,1,1)
\d[dtwo]2(26,6,2,2)
\d[dtwo]2(26,7,1,1)
\d[dtwo]2(26,8,1,1)
\d[dtwo]2(26,10,1,1)
\d[dtwo]2(27,3,1,2)
\d[dtwo]2(27,4,1,3)
\d[dtwo]2(27,5,1,2)
\d[dtwo]2(27,5,2,3)
\d[dtwo]2(27,6,1,2)
\d[dtwo]2(27,7,1,1)
\d[dtwo]2(27,8,1,2)
\d[dtwo]2(27,9,1,1)
\d[dtwo]2(27,10,1,1)
\d[dtwo]2(27,11,1,1)
\d[dtwo]2(27,11,2,2)
\d[dtwo]2(27,12,1,1)
\d[dtwo]2(27,13,1,1)
\d[dtwo]2(27,14,1,1)
\d[dtwo]2(27,15,1,1)
\d[dtwo]2(27,16,1,1)
\d[dtwo]2(27,17,1,1)
\d[dtwo]2(27,18,1,1)
\d[dtwo]2(27,19,1,1)
\d[dtwo]2(27,20,1,1)
\d[dtwo]2(27,21,1,1)
\d[dtwo]2(27,22,1,1)
\d[dtwo]2(28,3,1,3)
\d[dtwo]2(28,4,1,2)
\d[dtwo]2(28,5,1,2)
\d[dtwo]2(28,5,2,3)
\d[dtwo]2(28,7,1,2)
\d[dtwo]2(28,11,1,2)
\d[dtwo]2(28,12,1,2)
\d[dtwo]2(28,13,1,2)
\d[dtwo]2(29,4,1,1)
\d[dtwo]2(29,5,1,2)
\d[dtwo]2(29,6,1,1)
\d[dtwo]2(29,7,1,1)
\d[dtwo]2(29,8,1,1)
\d[dtwo]2(30,3,1,2)
\d[dtwo]2(30,4,1,2)
\d[dtwo]2(30,4,2,3)
\d[dtwo]2(30,5,1,2)
\d[dtwo]2(31,3,1,2)
\d[dtwo]2(31,4,1,1)
\d[dtwo]2(31,7,1,1)
\d[dtwo]2(31,8,1,1)
\d[dtwo]2(31,9,1,1)
\d[dtwo]2(31,10,1,1)
\d[dtwo]2(31,11,1,1)
\d[dtwo]2(31,11,2,2)
\d[dtwo]2(31,12,1,1)
\d[dtwo]2(31,12,2,2)
\d[dtwo]2(31,13,1,1)
\d[dtwo]2(31,14,1,1)
\d[dtwo]2(31,15,1,1)
\d[dtwo]2(31,16,1,1)
\d[dtwo]2(31,17,1,1)
\d[dtwo]2(31,18,1,1)
\d[dtwo]2(31,19,1,1)
\d[dtwo]2(31,20,1,1)
\d[dtwo]2(31,21,1,1)
\d[dtwo]2(31,22,1,1)
\d[dtwo]2(32,4,1,1)
\d[dtwo]2(32,5,1,2)
\d[dtwo]2(32,6,1,2)
\d[dtwo]2(32,8,1,2)
\d[dtwo]2(32,9,1,3)
\d[dtwo]2(33,4,1,2)
\d[dtwo]2(33,5,1,1)
\d[dtwo]2(33,9,1,1)
\d[dthr]3(25,8,3,1)
\d[dthr]3(26,9,2,2)
\d[dthr]3(30,6,2,1)
\d[dthr]3(31,5,1,1)
\d[dthr]3(32,10,1,2)
\d[dthr]3(33,4,2,2)
\d[dthr]3(33,12,3,1)
\d[dfour]4(31,8,3,2)
\structline(24,3,1)(25,4,1)
\structline(24,4,1)(24,5,1)
\structline(24,4,1)(24,5,3)
\structline(24,4,1)(25,5,1)
\structline(24,4,1)(27,5,1)
\structline(24,4,1)(27,5,2)
\structline(24,4,1)(27,5,3)
\structline(24,5,1)(24,6,2)
\structline(24,5,2)(24,6,2)
\structline(24,5,3)(24,6,2)
\structline(24,5,3)(24,6,3)
\structline(24,5,2)(25,6,1)
\structline(24,5,2)(25,6,2)
\structline(24,5,3)(27,6,2)
\structline(24,6,1)(24,7,2)
\structline(24,6,2)(24,7,3)
\structline(24,6,3)(24,7,3)
\structline(24,6,3)(24,7,4)
\structline(24,6,3)(24,7,5)
\structline(24,6,1)(27,7,2)
\structline(24,6,3)(27,7,3)
\structline(24,7,2)(24,8,1)
\structline(24,7,3)(24,8,2)
\structline(24,7,4)(24,8,3)
\structline(24,7,5)(24,8,4)
\structline(24,7,1)(25,8,1)
\structline(24,7,1)(25,8,2)
\structline(24,7,4)(25,8,2)
\structline(24,7,4)(25,8,4)
\structline(24,7,5)(25,8,2)
\structline(24,7,5)(25,8,4)
\structline(24,8,1)(24,9,1)
\structline(24,8,2)(24,9,2)
\structline(24,8,3)(24,9,3)
\structline(24,8,4)(24,9,4)
\structline(24,8,5)(24,9,5)
\structline(24,8,5)(25,9,2)
\structline(24,8,5)(27,9,2)
\structline(24,9,1)(24,10,1)
\structline(24,9,2)(24,10,2)
\structline(24,9,3)(24,10,3)
\structline(24,9,4)(24,10,4)
\structline(24,9,5)(24,10,5)
\structline(24,9,6)(25,10,2)
\structline(24,10,1)(24,11,2)
\structline(24,10,2)(24,11,3)
\structline(24,10,3)(24,11,4)
\structline(24,10,4)(24,11,5)
\structline(24,10,5)(24,11,6)
\structline(24,11,2)(24,12,1)
\structline(24,11,3)(24,12,3)
\structline(24,11,4)(24,12,4)
\structline(24,11,5)(24,12,5)
\structline(24,11,6)(24,12,6)
\structline(24,11,7)(24,12,7)
\structline(24,11,1)(25,12,2)
\structline(24,11,7)(25,12,1)
\structline(24,12,1)(24,13,1)
\structline(24,12,2)(24,13,2)
\structline(24,12,3)(24,13,3)
\structline(24,12,4)(24,13,4)
\structline(24,12,5)(24,13,5)
\structline(24,12,6)(24,13,6)
\structline(24,12,7)(24,13,7)
\structline(24,12,2)(25,13,1)
\structline(24,12,2)(27,13,2)
\structline(24,13,1)(24,14,1)
\structline(24,13,2)(24,14,2)
\structline(24,13,3)(24,14,3)
\structline(24,13,4)(24,14,4)
\structline(24,13,5)(24,14,5)
\structline(24,13,6)(24,14,6)
\structline(24,13,7)(24,14,7)
\structline(24,13,2)(27,14,2)
\structline(24,14,1)(24,15,1)
\structline(24,14,2)(24,15,2)
\structline(24,14,3)(24,15,3)
\structline(24,14,4)(24,15,4)
\structline(24,14,5)(24,15,5)
\structline(24,14,6)(24,15,6)
\structline(24,14,7)(24,15,7)
\structline(24,14,2)(27,15,2)
\structline(24,15,1)(24,16,1)
\structline(24,15,2)(24,16,2)
\structline(24,15,3)(24,16,3)
\structline(24,15,4)(24,16,4)
\structline(24,15,5)(24,16,5)
\structline(24,15,6)(24,16,6)
\structline(24,15,7)(24,16,7)
\structline(24,16,1)(24,17,1)
\structline(24,16,2)(24,17,2)
\structline(24,16,3)(24,17,3)
\structline(24,16,4)(24,17,4)
\structline(24,16,5)(24,17,5)
\structline(24,16,6)(24,17,6)
\structline(24,16,7)(24,17,7)
\structline(24,17,1)(24,18,1)
\structline(24,17,2)(24,18,2)
\structline(24,17,3)(24,18,3)
\structline(24,17,4)(24,18,4)
\structline(24,17,5)(24,18,5)
\structline(24,17,6)(24,18,6)
\structline(24,17,7)(24,18,7)
\structline(24,18,1)(24,19,1)
\structline(24,18,2)(24,19,2)
\structline(24,18,3)(24,19,3)
\structline(24,18,4)(24,19,4)
\structline(24,18,5)(24,19,5)
\structline(24,18,6)(24,19,6)
\structline(24,18,7)(24,19,7)
\structline(24,19,1)(24,20,1)
\structline(24,19,2)(24,20,2)
\structline(24,19,3)(24,20,3)
\structline(24,19,4)(24,20,4)
\structline(24,19,5)(24,20,5)
\structline(24,19,6)(24,20,6)
\structline(24,19,7)(24,20,7)
\structline(24,20,1)(24,21,1)
\structline(24,20,2)(24,21,2)
\structline(24,20,3)(24,21,3)
\structline(24,20,4)(24,21,4)
\structline(24,20,5)(24,21,5)
\structline(24,20,6)(24,21,6)
\structline(24,20,7)(24,21,7)
\structline(24,21,1)(24,22,1)
\structline(24,21,2)(24,22,2)
\structline(24,21,3)(24,22,3)
\structline(24,21,4)(24,22,4)
\structline(24,21,5)(24,22,5)
\structline(24,21,6)(24,22,6)
\structline(24,21,7)(24,22,7)
\structline(24,22,1)(24,23,1)
\structline(24,22,2)(24,23,2)
\structline(24,22,3)(24,23,3)
\structline(24,22,4)(24,23,4)
\structline(24,22,5)(24,23,5)
\structline(24,22,6)(24,23,6)
\structline(24,22,7)(24,23,7)
\structline(24,23,1)(24,24,1)
\structline(24,23,2)(24,24,2)
\structline(24,23,3)(24,24,3)
\structline(24,23,4)(24,24,4)
\structline(24,23,5)(24,24,5)
\structline(24,23,6)(24,24,6)
\structline(24,23,7)(24,24,7)
\structline(25,5,1)(26,6,1)
\structline(25,5,1)(26,6,2)
\structline(25,5,2)(26,6,4)
\structline(25,6,1)(25,7,1)
\structline(25,6,2)(25,7,1)
\structline(25,6,1)(26,7,2)
\structline(25,6,1)(26,7,3)
\structline(25,6,2)(26,7,2)
\structline(25,6,2)(26,7,3)
\structline(25,6,2)(26,7,4)
\structline(25,6,1)(28,7,2)
\structline(25,6,2)(28,7,2)
\structline(25,7,1)(25,8,1)
\structline(25,7,1)(28,8,1)
\structline(25,8,3)(25,9,1)
\structline(25,8,3)(26,9,2)
\structline(25,8,4)(26,9,3)
\structline(25,8,3)(28,9,1)
\structline(25,9,1)(25,10,1)
\structline(25,9,2)(26,10,3)
\structline(25,9,1)(28,10,1)
\structline(25,10,2)(26,11,2)
\structline(25,12,1)(26,13,2)
\structline(25,13,1)(26,14,2)
\structline(26,4,1)(26,5,1)
\structline(26,4,1)(27,5,1)
\structline(26,4,1)(27,5,2)
\structline(26,4,1)(29,5,1)
\structline(26,5,1)(26,6,1)
\structline(26,5,2)(26,6,3)
\structline(26,5,1)(29,6,1)
\structline(26,5,1)(29,6,3)
\structline(26,5,2)(29,6,2)
\structline(26,5,2)(29,6,3)
\structline(26,6,1)(26,7,4)
\structline(26,6,2)(26,7,4)
\structline(26,6,3)(26,7,3)
\structline(26,6,5)(26,7,2)
\structline(26,6,5)(26,7,3)
\structline(26,6,2)(27,7,3)
\structline(26,6,4)(27,7,2)
\structline(26,6,3)(29,7,2)
\structline(26,6,5)(29,7,2)
\structline(26,7,1)(26,8,1)
\structline(26,7,2)(26,8,2)
\structline(26,7,3)(26,8,2)
\structline(26,7,1)(29,8,1)
\structline(26,8,1)(26,9,2)
\structline(26,8,2)(26,9,1)
\structline(26,8,1)(29,9,1)
\structline(26,9,1)(26,10,2)
\structline(26,10,1)(26,11,2)
\structline(26,10,2)(26,11,1)
\structline(26,10,1)(27,11,2)
\structline(26,11,1)(26,12,1)
\structline(26,12,1)(26,13,1)
\structline(26,13,1)(26,14,1)
\structline(26,14,1)(26,15,1)
\structline(26,14,2)(27,15,2)
\structline(26,15,1)(26,16,1)
\structline(26,16,1)(26,17,1)
\structline(26,17,1)(26,18,1)
\structline(26,18,1)(26,19,1)
\structline(26,19,1)(26,20,1)
\structline(26,20,1)(26,21,1)
\structline(26,21,1)(26,22,1)
\structline(26,22,1)(26,23,1)
\structline(26,23,1)(26,24,1)
\structline(27,3,1)(27,4,1)
\structline(27,3,1)(30,4,1)
\structline(27,3,1)(30,4,2)
\structline(27,4,1)(27,5,2)
\structline(27,4,1)(30,5,1)
\structline(27,5,1)(27,6,1)
\structline(27,5,2)(27,6,1)
\structline(27,5,3)(27,6,2)
\structline(27,5,3)(30,6,1)
\structline(27,6,1)(27,7,1)
\structline(27,6,2)(27,7,3)
\structline(27,6,3)(27,7,2)
\structline(27,7,1)(27,8,1)
\structline(27,8,1)(27,9,1)
\structline(27,9,1)(27,10,1)
\structline(27,10,1)(27,11,1)
\structline(27,11,1)(27,12,1)
\structline(27,11,2)(28,12,7)
\structline(27,12,1)(27,13,1)
\structline(27,13,1)(27,14,1)
\structline(27,13,2)(27,14,2)
\structline(27,13,2)(30,14,1)
\structline(27,14,1)(27,15,1)
\structline(27,14,2)(27,15,2)
\structline(27,15,1)(27,16,1)
\structline(27,16,1)(27,17,1)
\structline(27,17,1)(27,18,1)
\structline(27,18,1)(27,19,1)
\structline(27,19,1)(27,20,1)
\structline(27,20,1)(27,21,1)
\structline(27,21,1)(27,22,1)
\structline(27,22,1)(27,23,1)
\structline(27,23,1)(27,24,1)
\structline(28,3,1)(28,4,1)
\structline(28,3,1)(31,4,1)
\structline(28,4,1)(28,5,2)
\structline(28,4,1)(31,5,1)
\structline(28,5,1)(28,6,1)
\structline(28,5,1)(28,6,2)
\structline(28,5,2)(28,6,3)
\structline(28,5,1)(29,6,3)
\structline(28,6,1)(28,7,2)
\structline(28,6,2)(28,7,2)
\structline(28,6,2)(28,7,3)
\structline(28,6,3)(28,7,4)
\structline(28,6,1)(29,7,2)
\structline(28,6,2)(29,7,2)
\structline(28,6,1)(31,7,2)
\structline(28,6,2)(31,7,2)
\structline(28,7,1)(28,8,1)
\structline(28,7,1)(28,8,4)
\structline(28,7,2)(28,8,1)
\structline(28,7,3)(28,8,3)
\structline(28,7,4)(28,8,5)
\structline(28,7,2)(31,8,2)
\structline(28,8,2)(28,9,1)
\structline(28,8,3)(28,9,2)
\structline(28,8,4)(28,9,3)
\structline(28,8,5)(28,9,4)
\structline(28,8,2)(29,9,1)
\structline(28,9,1)(28,10,1)
\structline(28,9,2)(28,10,2)
\structline(28,9,3)(28,10,3)
\structline(28,9,4)(28,10,5)
\structline(28,10,2)(28,11,2)
\structline(28,10,3)(28,11,3)
\structline(28,10,4)(28,11,4)
\structline(28,10,5)(28,11,5)
\structline(28,11,1)(28,12,1)
\structline(28,11,2)(28,12,2)
\structline(28,11,3)(28,12,3)
\structline(28,11,4)(28,12,4)
\structline(28,11,5)(28,12,5)
\structline(28,11,6)(28,12,6)
\structline(28,11,1)(31,12,1)
\structline(28,12,1)(28,13,1)
\structline(28,12,2)(28,13,2)
\structline(28,12,3)(28,13,3)
\structline(28,12,4)(28,13,4)
\structline(28,12,5)(28,13,5)
\structline(28,12,6)(28,13,6)
\structline(28,12,1)(31,13,2)
\structline(28,13,1)(28,14,1)
\structline(28,13,2)(28,14,2)
\structline(28,13,3)(28,14,3)
\structline(28,13,4)(28,14,4)
\structline(28,13,5)(28,14,5)
\structline(28,13,6)(28,14,6)
\structline(28,14,1)(28,15,1)
\structline(28,14,2)(28,15,2)
\structline(28,14,3)(28,15,3)
\structline(28,14,4)(28,15,4)
\structline(28,14,5)(28,15,5)
\structline(28,14,6)(28,15,6)
\structline(28,14,7)(28,15,7)
\structline(28,15,1)(28,16,1)
\structline(28,15,2)(28,16,2)
\structline(28,15,3)(28,16,3)
\structline(28,15,4)(28,16,4)
\structline(28,15,5)(28,16,5)
\structline(28,15,6)(28,16,6)
\structline(28,15,7)(28,16,7)
\structline(28,16,1)(28,17,1)
\structline(28,16,2)(28,17,2)
\structline(28,16,3)(28,17,3)
\structline(28,16,4)(28,17,4)
\structline(28,16,5)(28,17,5)
\structline(28,16,6)(28,17,6)
\structline(28,16,7)(28,17,7)
\structline(28,17,1)(28,18,1)
\structline(28,17,2)(28,18,2)
\structline(28,17,3)(28,18,3)
\structline(28,17,4)(28,18,4)
\structline(28,17,5)(28,18,5)
\structline(28,17,6)(28,18,6)
\structline(28,17,7)(28,18,7)
\structline(28,18,1)(28,19,1)
\structline(28,18,2)(28,19,2)
\structline(28,18,3)(28,19,3)
\structline(28,18,4)(28,19,4)
\structline(28,18,5)(28,19,5)
\structline(28,18,6)(28,19,6)
\structline(28,18,7)(28,19,7)
\structline(28,19,1)(28,20,1)
\structline(28,19,2)(28,20,2)
\structline(28,19,3)(28,20,3)
\structline(28,19,4)(28,20,4)
\structline(28,19,5)(28,20,5)
\structline(28,19,6)(28,20,6)
\structline(28,19,7)(28,20,7)
\structline(28,20,1)(28,21,1)
\structline(28,20,2)(28,21,2)
\structline(28,20,3)(28,21,3)
\structline(28,20,4)(28,21,4)
\structline(28,20,5)(28,21,5)
\structline(28,20,6)(28,21,6)
\structline(28,20,7)(28,21,7)
\structline(28,21,1)(28,22,1)
\structline(28,21,2)(28,22,2)
\structline(28,21,3)(28,22,3)
\structline(28,21,4)(28,22,4)
\structline(28,21,5)(28,22,5)
\structline(28,21,6)(28,22,6)
\structline(28,21,7)(28,22,7)
\structline(28,22,1)(28,23,1)
\structline(28,22,2)(28,23,2)
\structline(28,22,3)(28,23,3)
\structline(28,22,4)(28,23,4)
\structline(28,22,5)(28,23,5)
\structline(28,22,6)(28,23,6)
\structline(28,22,7)(28,23,7)
\structline(28,23,1)(28,24,1)
\structline(28,23,2)(28,24,2)
\structline(28,23,3)(28,24,3)
\structline(28,23,4)(28,24,4)
\structline(28,23,5)(28,24,5)
\structline(28,23,6)(28,24,6)
\structline(28,23,7)(28,24,7)
\structline(29,4,1)(29,5,1)
\structline(29,4,1)(30,5,1)
\structline(29,4,1)(32,5,1)
\structline(29,4,1)(32,5,3)
\structline(29,5,1)(29,6,1)
\structline(29,5,1)(29,6,3)
\structline(29,5,2)(29,6,2)
\structline(29,5,2)(29,6,3)
\structline(29,5,1)(32,6,1)
\structline(29,5,2)(32,6,2)
\structline(29,6,2)(29,7,2)
\structline(29,6,2)(32,7,1)
\structline(29,7,1)(29,8,1)
\structline(29,7,1)(32,8,2)
\structline(29,8,1)(29,9,1)
\structline(29,8,1)(32,9,2)
\structline(30,3,1)(30,4,1)
\structline(30,3,1)(30,4,2)
\structline(30,3,1)(33,4,1)
\structline(30,4,1)(30,5,1)
\structline(30,4,1)(30,5,2)
\structline(30,4,2)(30,5,2)
\structline(30,4,1)(31,5,1)
\structline(30,4,2)(31,5,1)
\structline(30,4,1)(33,5,1)
\structline(30,4,1)(33,5,3)
\structline(30,4,2)(33,5,2)
\structline(30,4,2)(33,5,3)
\structline(30,5,2)(30,6,1)
\structline(30,5,2)(33,6,1)
\structline(30,8,1)(30,9,1)
\structline(30,10,1)(30,11,1)
\structline(30,10,1)(31,11,3)
\structline(30,10,1)(33,11,1)
\structline(30,11,1)(30,12,1)
\structline(30,11,1)(33,12,2)
\structline(30,12,1)(30,13,2)
\structline(30,12,2)(30,13,1)
\structline(30,12,2)(31,13,2)
\structline(30,12,2)(33,13,1)
\structline(30,13,1)(30,14,1)
\structline(30,13,2)(30,14,2)
\structline(30,13,1)(33,14,1)
\structline(30,14,2)(30,15,1)
\structline(30,15,1)(30,16,1)
\structline(30,16,1)(30,17,1)
\structline(30,17,1)(30,18,1)
\structline(30,18,1)(30,19,1)
\structline(30,19,1)(30,20,1)
\structline(30,20,1)(30,21,1)
\structline(30,21,1)(30,22,1)
\structline(30,22,1)(30,23,1)
\structline(30,23,1)(30,24,1)
\structline(31,3,1)(31,4,1)
\structline(31,4,1)(31,5,1)
\structline(31,6,1)(31,7,2)
\structline(31,6,1)(32,7,1)
\structline(31,7,1)(31,8,2)
\structline(31,7,2)(31,8,2)
\structline(31,8,1)(31,9,1)
\structline(31,8,1)(32,9,1)
\structline(31,8,1)(32,9,6)
\structline(31,8,3)(32,9,2)
\structline(31,9,1)(31,10,1)
\structline(31,9,1)(31,10,2)
\structline(31,10,1)(31,11,2)
\structline(31,10,1)(31,11,3)
\structline(31,10,2)(31,11,3)
\structline(31,10,1)(32,11,1)
\structline(31,10,2)(32,11,1)
\structline(31,11,1)(31,12,1)
\structline(31,11,2)(31,12,2)
\structline(31,12,1)(31,13,2)
\structline(31,12,2)(31,13,1)
\structline(31,13,1)(31,14,1)
\structline(31,14,1)(31,15,1)
\structline(31,15,1)(31,16,1)
\structline(31,16,1)(31,17,1)
\structline(31,17,1)(31,18,1)
\structline(31,18,1)(31,19,1)
\structline(31,19,1)(31,20,1)
\structline(31,20,1)(31,21,1)
\structline(31,21,1)(31,22,1)
\structline(31,22,1)(31,23,1)
\structline(31,23,1)(31,24,1)
\structline(32,4,1)(32,5,1)
\structline(32,4,1)(32,5,3)
\structline(32,4,2)(32,5,2)
\structline(32,4,2)(32,5,3)
\structline(32,4,1)(33,5,1)
\structline(32,4,1)(33,5,2)
\structline(32,4,2)(33,5,2)
\structline(32,5,1)(32,6,1)
\structline(32,5,1)(32,6,3)
\structline(32,5,2)(32,6,2)
\structline(32,5,3)(32,6,3)
\structline(32,5,1)(33,6,2)
\structline(32,5,2)(33,6,2)
\structline(32,5,3)(33,6,2)
\structline(32,6,2)(32,7,1)
\structline(32,6,3)(32,7,4)
\structline(32,6,3)(32,7,6)
\structline(32,7,3)(32,8,2)
\structline(32,7,4)(32,8,3)
\structline(32,7,5)(32,8,3)
\structline(32,7,6)(32,8,3)
\structline(32,7,6)(32,8,4)
\structline(32,7,2)(33,8,1)
\structline(32,7,3)(33,8,1)
\structline(32,7,4)(33,8,2)
\structline(32,7,5)(33,8,3)
\structline(32,7,6)(33,8,2)
\structline(32,8,1)(32,9,1)
\structline(32,8,1)(32,9,5)
\structline(32,8,2)(32,9,2)
\structline(32,8,3)(32,9,3)
\structline(32,8,4)(32,9,3)
\structline(32,8,4)(32,9,4)
\structline(32,8,4)(32,9,6)
\structline(32,8,1)(33,9,1)
\structline(32,9,1)(32,10,4)
\structline(32,9,3)(32,10,2)
\structline(32,9,4)(32,10,3)
\structline(32,9,5)(32,10,5)
\structline(32,9,6)(32,10,4)
\structline(32,9,1)(33,10,1)
\structline(32,9,1)(33,10,2)
\structline(32,9,5)(33,10,1)
\structline(32,9,5)(33,10,2)
\structline(32,10,1)(32,11,4)
\structline(32,10,2)(32,11,5)
\structline(32,10,3)(32,11,2)
\structline(32,10,4)(32,11,3)
\structline(32,10,5)(32,11,6)
\structline(32,11,2)(32,12,3)
\structline(32,11,3)(32,12,4)
\structline(32,11,4)(32,12,1)
\structline(32,11,5)(32,12,2)
\structline(32,11,6)(32,12,5)
\structline(32,11,7)(32,12,6)
\structline(32,11,1)(33,12,1)
\structline(32,11,7)(33,12,1)
\structline(32,11,7)(33,12,2)
\structline(32,11,7)(33,12,4)
\structline(32,12,1)(32,13,1)
\structline(32,12,2)(32,13,2)
\structline(32,12,3)(32,13,3)
\structline(32,12,4)(32,13,4)
\structline(32,12,5)(32,13,5)
\structline(32,12,6)(32,13,6)
\structline(32,12,7)(32,13,7)
\structline(32,12,7)(33,13,2)
\structline(32,13,1)(32,14,1)
\structline(32,13,2)(32,14,2)
\structline(32,13,3)(32,14,3)
\structline(32,13,4)(32,14,4)
\structline(32,13,5)(32,14,5)
\structline(32,13,6)(32,14,6)
\structline(32,13,7)(32,14,7)
\structline(32,13,8)(33,14,2)
\structline(32,14,1)(32,15,2)
\structline(32,14,2)(32,15,3)
\structline(32,14,3)(32,15,4)
\structline(32,14,4)(32,15,5)
\structline(32,14,5)(32,15,6)
\structline(32,14,6)(32,15,7)
\structline(32,14,7)(32,15,8)
\structline(32,15,2)(32,16,1)
\structline(32,15,3)(32,16,3)
\structline(32,15,4)(32,16,4)
\structline(32,15,5)(32,16,5)
\structline(32,15,6)(32,16,6)
\structline(32,15,7)(32,16,7)
\structline(32,15,8)(32,16,8)
\structline(32,15,9)(32,16,9)
\structline(32,15,1)(33,16,2)
\structline(32,15,9)(33,16,1)
\structline(32,16,1)(32,17,1)
\structline(32,16,2)(32,17,2)
\structline(32,16,3)(32,17,3)
\structline(32,16,4)(32,17,4)
\structline(32,16,5)(32,17,5)
\structline(32,16,6)(32,17,6)
\structline(32,16,7)(32,17,7)
\structline(32,16,8)(32,17,8)
\structline(32,16,9)(32,17,9)
\structline(32,16,2)(33,17,1)
\structline(32,17,1)(32,18,1)
\structline(32,17,2)(32,18,2)
\structline(32,17,3)(32,18,3)
\structline(32,17,4)(32,18,4)
\structline(32,17,5)(32,18,5)
\structline(32,17,6)(32,18,6)
\structline(32,17,7)(32,18,7)
\structline(32,17,8)(32,18,8)
\structline(32,17,9)(32,18,9)
\structline(32,18,1)(32,19,1)
\structline(32,18,2)(32,19,2)
\structline(32,18,3)(32,19,3)
\structline(32,18,4)(32,19,4)
\structline(32,18,5)(32,19,5)
\structline(32,18,6)(32,19,6)
\structline(32,18,7)(32,19,7)
\structline(32,18,8)(32,19,8)
\structline(32,18,9)(32,19,9)
\structline(32,19,1)(32,20,1)
\structline(32,19,2)(32,20,2)
\structline(32,19,3)(32,20,3)
\structline(32,19,4)(32,20,4)
\structline(32,19,5)(32,20,5)
\structline(32,19,6)(32,20,6)
\structline(32,19,7)(32,20,7)
\structline(32,19,8)(32,20,8)
\structline(32,19,9)(32,20,9)
\structline(32,20,1)(32,21,1)
\structline(32,20,2)(32,21,2)
\structline(32,20,3)(32,21,3)
\structline(32,20,4)(32,21,4)
\structline(32,20,5)(32,21,5)
\structline(32,20,6)(32,21,6)
\structline(32,20,7)(32,21,7)
\structline(32,20,8)(32,21,8)
\structline(32,20,9)(32,21,9)
\structline(32,21,1)(32,22,1)
\structline(32,21,2)(32,22,2)
\structline(32,21,3)(32,22,3)
\structline(32,21,4)(32,22,4)
\structline(32,21,5)(32,22,5)
\structline(32,21,6)(32,22,6)
\structline(32,21,7)(32,22,7)
\structline(32,21,8)(32,22,8)
\structline(32,21,9)(32,22,9)
\structline(32,22,1)(32,23,1)
\structline(32,22,2)(32,23,2)
\structline(32,22,3)(32,23,3)
\structline(32,22,4)(32,23,4)
\structline(32,22,5)(32,23,5)
\structline(32,22,6)(32,23,6)
\structline(32,22,7)(32,23,7)
\structline(32,22,8)(32,23,8)
\structline(32,22,9)(32,23,9)
\structline(32,23,1)(32,24,1)
\structline(32,23,2)(32,24,2)
\structline(32,23,3)(32,24,3)
\structline(32,23,4)(32,24,4)
\structline(32,23,5)(32,24,5)
\structline(32,23,6)(32,24,6)
\structline(32,23,7)(32,24,7)
\structline(32,23,8)(32,24,8)
\structline(32,23,9)(32,24,9)
\structline(33,4,1)(33,5,1)
\structline(33,4,1)(33,5,2)
\structline(33,4,2)(33,5,2)
\structline(33,4,2)(33,5,3)
\structline(33,5,3)(33,6,1)
\structline(33,10,1)(33,11,1)
\structline(33,10,2)(33,11,1)
\structline(33,11,1)(33,12,2)
\structline(33,12,3)(33,13,1)
\structline(33,13,1)(33,14,1)
\end{sseqdata}
\begin{figure}[htbp]\centering
\adjustbox{max width=\textwidth}{\printpage[name = AdamsEtwoC, page = 2]}
\caption{The Adams $E_2$-page of $\tmf_0(5)$ in stems $24$ to $33$, with the $d_2$-differentials.}
\label{fig:adamsE2c}
\end{figure}

\subsection{The input}\label{subsec:adamsinput} We use the following facts.
\begin{enumerate}
\item The differentials of the Adams spectral sequence of $\tmf$ \cite[Chapter 5]{BR}. In stems below $34$ they are generated under the Leibniz rule by
\begin{align*}
d_2(\alpha)&=h_2w_1, & d_2(\beta)&=h_0d_0, & d_3(e_0)&=c_0w_1,\\
d_3(\alpha^2)&=h_1w_1d_0, & d_3(\beta^2)&=h_1w_1g, & d_4(d_0e_0)&=w_1^2d_0,
\end{align*}
and the unit $\tmf\to\tmf_0(5)$ carries them to differentials on the first summand of $E_2$. The spectral sequence of $\tmf_0(5)$ is a module over that of $\tmf$.
\item The ranks $r_n$ of $\pi_n\tmf_0(5)$ computed above.
\item Every torsion class of $\pi_*\tmf_0(5)$ is annihilated by $4$, by Lemma \ref{lem:fournu}.
\item The orders of the torsion subgroups of $\pi_n\tmf_0(5)$ for $n\le32$, listed in Table \ref{tab:torsionorders}. These can be read off from \cite[Figure 12]{hswx}, which shows the $E_\infty$-page of the slice spectral sequence with the classes of \cite[Proposition 13.4]{hhrc4} omitted. The omitted torsion classes are the products of $\eta$ and $\eta^2$ with $2l$ of the classes on the $0$-line in stem $8l$, so they contribute $2l$ classes to each of the stems $8l+1$ and $8l+2$ \cite[Theorem 13.3(v)]{hhrc4}.
\item The ring map $\tmf_0(5)=\tmf_1(5)^{C_4}\to\mathrm{Tmf}_1(5)^{C_4}\simeq\mathrm{Tmf}_0(5)\to\TMF_0(5)$ under $\tmf$ \cite[Section 6.3]{quinn}, and the computation of $\pi_*\TMF_0(5)$ by Behrens--Ormsby \cite[Theorem 2.3.12]{BO}. In the coordinates of Proposition \ref{prop:tate} we have $\pi_*\TMF_1(5)=\Z_{(2)}[a,b,\Delta^{-1}]$, and the restriction lands in the invariants
\[\Z_{(2)}[b_2,b_4,\delta][\Delta^{-1}]/(b_4^2-b_2^2\delta+4\delta^2),\]
where $b_2=a^2+b^2$, $b_4=a^3b-ab^3$ and $\delta=a^2b^2$ are the classes of \cite{BO} with $u=a$ and $a_1=a+b$, and we use three consequences. First, the images of $\epsilon$ and $\kappa$ in $\pi_*\TMF_0(5)$ are nonzero. Second, $c_4=b_2^2-12b_4+12\delta$, so that $\eta c_4=\eta b_2^2\neq0$ in $\pi_9\TMF_0(5)$. Third, the image of $\pi_{16}\TMF_0(5)$ in the invariants is contained in the ideal generated by $b_2$, $b_4$ and $2\delta^2$.
\end{enumerate}

\begin{table}[ht]\centering\footnotesize
\caption{The orders $2^{t_n}$ of the torsion subgroups of $\pi_n\tmf_0(5)$ for $n\le32$. In the stems not listed, $t_n=0$.}
\label{tab:torsionorders}
\renewcommand{\arraystretch}{1.15}
\begin{tabular}{@{}l*{19}{c}@{}}
\toprule
$n$ & $1$ & $2$ & $3$ & $4$ & $6$ & $8$ & $9$ & $10$ & $14$ & $17$ & $18$ & $19$ & $20$ & $21$ & $22$ & $25$ & $26$ & $27$ & $28$\\
\midrule
$t_n$ & $2$ & $2$ & $2$ & $1$ & $1$ & $1$ & $4$ & $3$ & $1$ & $5$ & $5$ & $2$ & $2$ & $1$ & $2$ & $7$ & $7$ & $1$ & $1$\\
\bottomrule
\end{tabular}
\end{table}

\subsection{The differentials} The differentials are determined by their values on the generators $x_{n,s}$ and on the classes that become indecomposable on later pages.

\begin{theorem}\label{thm:adamsdiff}
With suitable choices of the generators $x_{n,s}$, the nonzero differentials with source in stems at most $32$ on the generators $x_{n,s}$, and on the classes that become indecomposable on $E_3$ and $E_4$, are
\begin{align*}
d_2(x_{10,2})&=h_1x_{8,3}, & d_2(x_{15,0})&=x_{14,2},\\
d_2(x_{16,0})&=x_{15,2}, & d_2(x_{20,3})&=h_2x_{16,4},\\
d_2(x_{26,4})&=x_{25,6}, & d_2(x_{28,5})&=h_0\alpha\beta,\\
d_3(x_{4,0})&=h_0^2h_2, & d_3(x_{18,2})&=h_0e_0,\\
d_3(h_1d_0x_{16,0})&=w_1x_{22,4}, & d_3(e_0x_{16,0}+h_2\beta x_{15,0})&=\epsilon_1',\\
d_4(h_0x_{16,0}+h_1x_{15,0})&=h_1d_0,
\end{align*}
where $\epsilon_1'$ is the unique nonzero class in bidegree $(32,7)$ of the first summand of $E_2$ that is annihilated by $h_0$; it detects the image of $\epsilon_1\in\pi_{32}\tmf$. All other differentials with source in stems at most $32$ follow from these, from the differentials of $\tmf$, from the Leibniz rule and for degree reasons, and $E_5=E_\infty$ in stems at most $32$.
\end{theorem}

\begin{figure}[htbp]
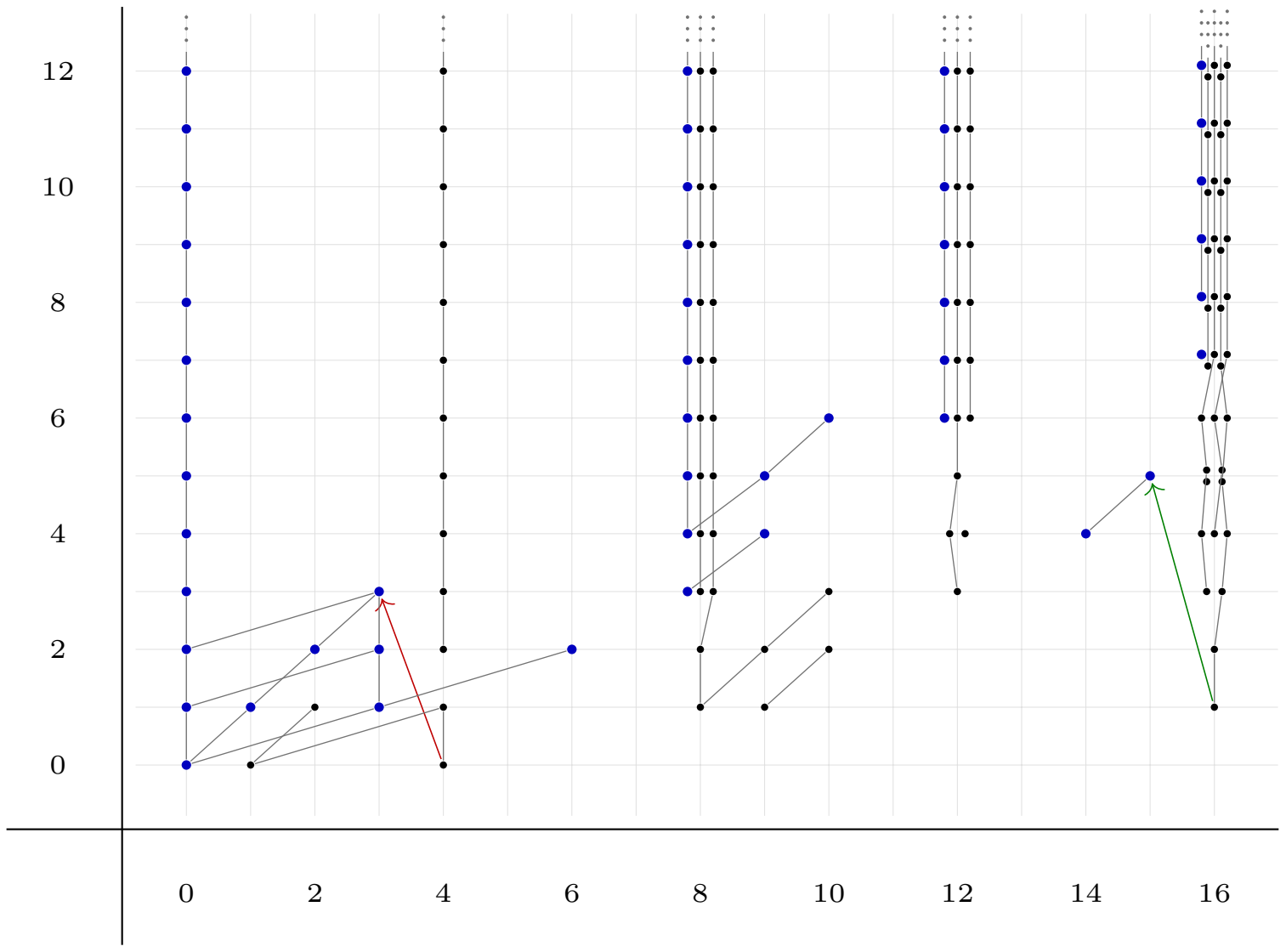
\centering
\adjustbox{max width=\textwidth}{\printpage[name = AdamsEthreeA, page = 3--4]}
\caption{The page $E_3$ in stems $0$ to $16$, with the $d_3$- and $d_4$-differentials.}
\label{fig:adamsE3a}
\end{figure}

\begin{figure}[htbp]
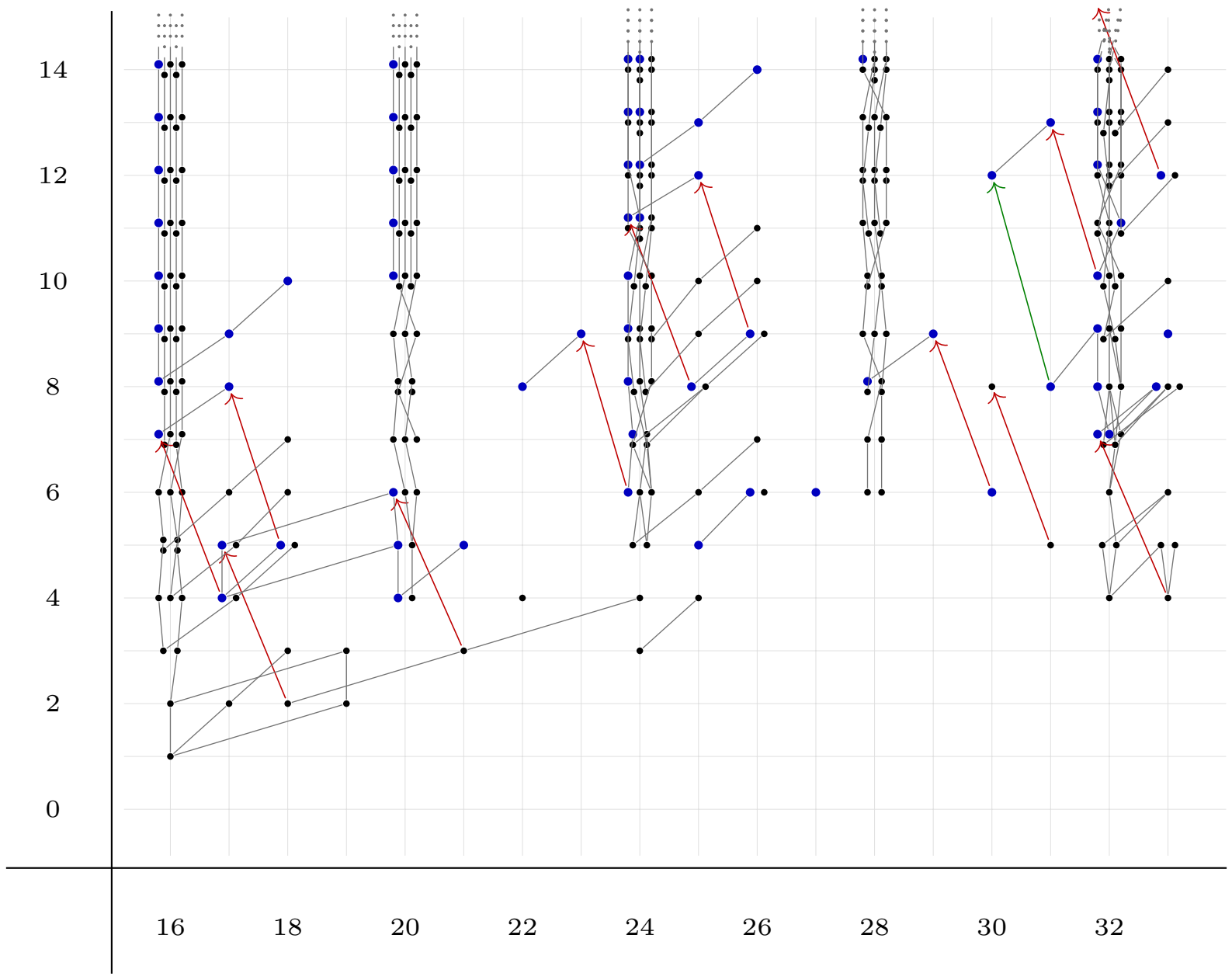
\centering
\adjustbox{max width=\textwidth}{\printpage[name = AdamsEthreeB, page = 3--4]}
\caption{The page $E_3$ in stems $16$ to $33$, with the $d_3$- and $d_4$-differentials.}
\label{fig:adamsE3b}
\end{figure}

The pages $E_3$ and $E_\infty$ are shown in Figures \ref{fig:adamsE3a}--\ref{fig:adamsEinfb}. In particular the classes $x_{1,0}$, $x_{4,2}$, $x_{8,1}$, $x_{8,3}$, $x_{9,1}$, $x_{12,3}$, $x_{16,3}$, $x_{16,4}$, $x_{22,4}$, $x_{24,5}$, $x_{32,4}$, $x_{32,5}$, $x_{32,7}$ and $w_1x_{18,2}$ are permanent cycles, $x_{14,2}$, $x_{15,2}$ and $x_{25,6}$ are boundaries, and the class $x_{1,0}$ detects the element $\eta'$ of \cite{hswx}.

The listed values, and the vanishing of the differentials on the remaining generators, are decided by comparing with the homotopy groups of $\tmf_0(5)$, that is, with the input of Section \ref{subsec:adamsinput}; the Leibniz rule then determines the differential in every bidegree in the range. We give the arguments in the stems up to $16$, which contain all the types of argument used, and only indicate the crucial inputs in the remaining stems.

\begin{proof}
\emph{Stems $1$ to $9$.} The class $x_{1,0}$ is permanent, since a differential on it would truncate the $h_0$-tower on $1$, whose survival is forced by $\pi_0\tmf_0(5)\otimes\mathbb Q\neq0$. The classes $x_{4,2}$, $x_{8,1}$ and $x_{8,3}$ are permanent, as $E_2$ vanishes in stems $3$ above filtration $3$ and in stem $7$. The only possible differential on $x_{4,0}$ is $d_2(x_{4,0})=h_0h_2$ or $d_3(x_{4,0})=h_0^2h_2$. Since $h_0x_{4,0}=h_2x_{1,0}$ is a permanent cycle and $h_0\cdot h_0h_2\neq0$, the first is excluded, and since $4\nu=0$ the class $h_0^2h_2$ detecting $4\nu$ must be hit, which leaves the second. The only possible target of a differential on $x_{9,1}$ annihilated by $h_0$ is $c_0$, which detects $\epsilon$; as $\epsilon$ maps to a nonzero element of $\pi_8\TMF_0(5)$, it is nonzero in $\pi_8\tmf_0(5)$, and so $x_{9,1}$ is permanent.

\emph{Stems $10$ and $11$.} The group $\pi_{11}\tmf_0(5)$ is zero, while $E_2$ contains the class $h_1x_{10,2}$ in stem $11$, which cannot be hit. So $h_1x_{10,2}$ supports a differential. If $d_2(x_{10,2})$ has no $h_1x_{8,3}$-component, then $h_1x_{10,2}$ survives to $E_3$, and the only possible value of $d_3$ on it is $h_1^2w_1$, which detects $\eta^2c_4$. But $\eta^2c_4=\eta^2b_2^2$ is nonzero in $\pi_{10}\TMF_0(5)$, and $E_2$ has no classes in stem $10$ above filtration $6$, so $h_1^2w_1$ survives. Hence $d_2(x_{10,2})\in h_1x_{8,3}+\F_2\{h_1c_0\}$, and replacing $x_{8,3}$ by $x_{8,3}+c_0$ if necessary, $d_2(x_{10,2})=h_1x_{8,3}$.

\emph{Stems $12$ to $15$.} Since $d_2(\alpha)=h_2w_1$, we may choose $x_{12,3}$ to be a $d_2$-cycle, and it is then permanent, as $E_3$ vanishes in stem $11$. The towers in stems $14$ and $15$ must cancel, as $\pi_{14}\tmf_0(5)$ and $\pi_{15}\tmf_0(5)$ are finite, so $d_r(x_{15,0})=h_0^{r-2}x_{14,2}$ for some $r\ge2$. If $r\ge3$, then $x_{14,2}$ would survive alongside $d_0$, which detects $\kappa$ and survives as $\kappa$ maps to a nonzero element of $\pi_{14}\TMF_0(5)$; this contradicts $t_{14}=1$. So $d_2(x_{15,0})=x_{14,2}$. The class $x_{15,2}$ can only support a $d_2$ into $\F_2\{d_0,h_0^2x_{14,2}\}$; the first component vanishes as $\kappa\neq0$, and the second can be removed by changing $x_{15,2}$ by $h_0^2x_{15,0}$. With this choice $x_{15,2}$ is a $d_2$-cycle, and the higher differentials on it vanish for degree reasons.

\emph{Stem $16$.} As $\pi_{15}\tmf_0(5)=0$, the classes $x_{15,2}$ and $h_0x_{15,2}$ must be hit, and the only possible source is $d_2(x_{16,0})=x_{15,2}$. The class $h_1d_0$ detecting $\eta\kappa$ must be hit as well. By the Leibniz rule, the only possible sources are $d_2(x_{16,3})=h_1d_0$ and $d_4(v)=h_1d_0$, where $v=h_0x_{16,0}+h_1x_{15,0}$ is the $d_2$-cycle in bidegree $(16,1)$. We show that $v$ is not a permanent cycle, which gives the second; then $x_{16,3}$ is a $d_2$-cycle, and permanent.

Base change along $\tmf\to H\Z_{(2)}$ gives a map from our spectral sequence to the Adams spectral sequence of the $H\Z_{(2)}$-module $H\Z_{(2)}\otimes_{\tmf}\tmf_0(5)$, whose homotopy is $H^{\tmf}_*(\tmf_0(5);\Z_{(2)})$, and which has $E_2=\Ext_{\mathcal{A}(0)_*}(\F_2,H^{C_4}_*)$; on $E_2$ the map is restriction from $\mathcal{A}(2)_*$ to $\mathcal{A}(0)_*$. The class $m_{22}$ spans the $Q_0$-homology of $H^{C_4}_*$ in degree $16$, and the $Q_0$-homology vanishes in degree $17$, as $Q_0m_{24}=m_{23}$. So $H^{\tmf}_{16}(\tmf_0(5);\Z_{(2)})=\Z_{(2)}\{\bar m_{22}\}$ with $\bar m_{22}$ reducing to $m_{22}$, and in the target spectral sequence the $h_0$-tower on $m_{22}$ survives. The class $v$ maps to $h_0m_{22}$. If $v$ detected an element $z\in\pi_{16}\tmf_0(5)$, then the Hurewicz image of $z$ would be $2\bar m_{22}$ times a unit.

We now compare with the underlying spectrum. In the notation of Section \ref{subsec:beps}, we have \[H^{\tmf}_{16}(\tmf_1(5);\Z_{(2)})=\Z_{(2)}\{\hat Z\},\] and since $\mathrm{tr}(\hat Z)$ reduces to $m_{22}$ by Proposition \ref{prop:C2res}, it is a unit times $\bar m_{22}$; as $\gamma\hat Z=\hat Z$, we get $\Res(\bar m_{22})=4\hat Z$ up to a unit. On the other hand, a direct computation in $L$ shows that the Hurewicz image of every monomial $a^ib^{8-i}$ is $8\hat Z$ times a unit. Hence the Hurewicz image of $y\in\pi_{16}\tmf_1(5)$ is congruent to $8\,y(1,1)\hat Z$ modulo $16\hat Z$, where $y(1,1)$ is the value at $a=b=1$, and therefore
\[\Res(z)(1,1)\equiv1\bmod 2.\]
But $\Res(z)$ lies in the image of $\pi_{16}\TMF_0(5)$, and the ideal generated by $b_2$, $b_4$ and $2\delta^2$ vanishes at $a=b=1$ modulo $2$. This contradiction shows that $v$ supports a differential.

\emph{Stems $17$ to $32$.} The remaining values are decided in the same way. We list the crucial inputs. The differential $d_3(x_{18,2})=h_0e_0$ is forced by $4\bar\kappa=0$ in $\pi_{20}\tmf_0(5)$ (the class $h_0^2g=c_0\alpha$ detecting $4\bar\kappa$ must be hit, and $h_2x_{18,2}$ is the only possible source) together with $t_{17}=5$, which excludes a $w_1x_{9,1}$-component. The values of $d_2(x_{20,3})$ and $d_2(x_{26,4})$ are forced by the Leibniz rule, from $w_1x_{20,3}=\alpha x_{16,4}+h_2x_{25,6}$ and $h_1x_{26,4}=\alpha x_{15,2}$. The vanishing of $d_4$ and $d_7$ on $w_1x_{18,2}$ follows from $t_{25}=t_{26}=7$, and $d_2(x_{28,5})=h_0\alpha\beta$ from $t_{27}=1$. The vanishing of $\pi_{30}\tmf_0(5)$ and $\pi_{31}\tmf_0(5)$ forces $d_3(h_1d_0x_{16,0})=w_1x_{22,4}$, the only differential that can kill these two classes. Finally $\pi_{32}\tmf_0(5)$ is torsion-free, so the image of the class $\epsilon_1$ of order $2$ vanishes, and the class $\epsilon_1'$ detecting it must be hit; the only possible source is the $d_3$ on $e_0x_{16,0}+h_2\beta x_{15,0}$. The remaining generators are permanent cycles for degree reasons or because they may be chosen to be $d_2$-cycles, as for $x_{12,3}$.
\end{proof}

\begin{sseqdata}[name = AdamsEthreeA, adams, x range = {0}{16}, y range = {0}{12}]
\class[tmfclass](0,0)
\class[tmfclass](0,1)
\class[tmfclass](0,2)
\class[tmfclass](0,3)
\class[tmfclass](0,4)
\class[tmfclass](0,5)
\class[tmfclass](0,6)
\class[tmfclass](0,7)
\class[tmfclass](0,8)
\class[tmfclass](0,9)
\class[tmfclass](0,10)
\class[tmfclass](0,11)
\class[tmfclass](0,12)
\class[tmfclass](0,13)
\class[tmfclass](0,14)
\class[tmfclass](0,15)
\class[tmfclass](0,16)
\class[tmfclass](0,17)
\class[tmfclass](0,18)
\class[tmfclass](0,19)
\class[tmfclass](0,20)
\class[tmfclass](0,21)
\class[tmfclass](0,22)
\class[tmfclass](0,23)
\class[tmfclass](0,24)
\class(1,0)
\class[tmfclass](1,1)
\class(2,1)
\class[tmfclass](2,2)
\class[tmfclass](3,1)
\class[tmfclass](3,2)
\class[tmfclass](3,3)
\class(4,0)
\class(4,1)
\class(4,2)
\class(4,3)
\class(4,4)
\class(4,5)
\class(4,6)
\class(4,7)
\class(4,8)
\class(4,9)
\class(4,10)
\class(4,11)
\class(4,12)
\class(4,13)
\class(4,14)
\class(4,15)
\class(4,16)
\class(4,17)
\class(4,18)
\class(4,19)
\class(4,20)
\class(4,21)
\class(4,22)
\class(4,23)
\class(4,24)
\class[tmfclass](6,2)
\class(8,1)
\class(8,2)
\class[tmfclass](8,3)
\class(8,3)
\class(8,3)
\class[tmfclass](8,4)
\class(8,4)
\class(8,4)
\class[tmfclass](8,5)
\class(8,5)
\class(8,5)
\class[tmfclass](8,6)
\class(8,6)
\class(8,6)
\class[tmfclass](8,7)
\class(8,7)
\class(8,7)
\class[tmfclass](8,8)
\class(8,8)
\class(8,8)
\class[tmfclass](8,9)
\class(8,9)
\class(8,9)
\class[tmfclass](8,10)
\class(8,10)
\class(8,10)
\class[tmfclass](8,11)
\class(8,11)
\class(8,11)
\class[tmfclass](8,12)
\class(8,12)
\class(8,12)
\class[tmfclass](8,13)
\class(8,13)
\class(8,13)
\class[tmfclass](8,14)
\class(8,14)
\class(8,14)
\class[tmfclass](8,15)
\class(8,15)
\class(8,15)
\class[tmfclass](8,16)
\class(8,16)
\class(8,16)
\class[tmfclass](8,17)
\class(8,17)
\class(8,17)
\class[tmfclass](8,18)
\class(8,18)
\class(8,18)
\class[tmfclass](8,19)
\class(8,19)
\class(8,19)
\class[tmfclass](8,20)
\class(8,20)
\class(8,20)
\class[tmfclass](8,21)
\class(8,21)
\class(8,21)
\class[tmfclass](8,22)
\class(8,22)
\class(8,22)
\class[tmfclass](8,23)
\class(8,23)
\class(8,23)
\class[tmfclass](8,24)
\class(8,24)
\class(8,24)
\class(9,1)
\class(9,2)
\class[tmfclass](9,4)
\class[tmfclass](9,5)
\class(10,2)
\class(10,3)
\class[tmfclass](10,6)
\class(12,3)
\class(12,4)
\class(12,4)
\class(12,5)
\class[tmfclass](12,6)
\class(12,6)
\class(12,6)
\class[tmfclass](12,7)
\class(12,7)
\class(12,7)
\class[tmfclass](12,8)
\class(12,8)
\class(12,8)
\class[tmfclass](12,9)
\class(12,9)
\class(12,9)
\class[tmfclass](12,10)
\class(12,10)
\class(12,10)
\class[tmfclass](12,11)
\class(12,11)
\class(12,11)
\class[tmfclass](12,12)
\class(12,12)
\class(12,12)
\class[tmfclass](12,13)
\class(12,13)
\class(12,13)
\class[tmfclass](12,14)
\class(12,14)
\class(12,14)
\class[tmfclass](12,15)
\class(12,15)
\class(12,15)
\class[tmfclass](12,16)
\class(12,16)
\class(12,16)
\class[tmfclass](12,17)
\class(12,17)
\class(12,17)
\class[tmfclass](12,18)
\class(12,18)
\class(12,18)
\class[tmfclass](12,19)
\class(12,19)
\class(12,19)
\class[tmfclass](12,20)
\class(12,20)
\class(12,20)
\class[tmfclass](12,21)
\class(12,21)
\class(12,21)
\class[tmfclass](12,22)
\class(12,22)
\class(12,22)
\class[tmfclass](12,23)
\class(12,23)
\class(12,23)
\class[tmfclass](12,24)
\class(12,24)
\class(12,24)
\class[tmfclass](14,4)
\class[tmfclass](15,5)
\class(16,1)
\class(16,2)
\class(16,3)
\class(16,3)
\class(16,4)
\class(16,4)
\class(16,4)
\class(16,5)
\class(16,5)
\class(16,5)
\class(16,5)
\class(16,6)
\class(16,6)
\class(16,6)
\class[tmfclass](16,7)
\class(16,7)
\class(16,7)
\class(16,7)
\class(16,7)
\class[tmfclass](16,8)
\class(16,8)
\class(16,8)
\class(16,8)
\class(16,8)
\class[tmfclass](16,9)
\class(16,9)
\class(16,9)
\class(16,9)
\class(16,9)
\class[tmfclass](16,10)
\class(16,10)
\class(16,10)
\class(16,10)
\class(16,10)
\class[tmfclass](16,11)
\class(16,11)
\class(16,11)
\class(16,11)
\class(16,11)
\class[tmfclass](16,12)
\class(16,12)
\class(16,12)
\class(16,12)
\class(16,12)
\class[tmfclass](16,13)
\class(16,13)
\class(16,13)
\class(16,13)
\class(16,13)
\class[tmfclass](16,14)
\class(16,14)
\class(16,14)
\class(16,14)
\class(16,14)
\class[tmfclass](16,15)
\class(16,15)
\class(16,15)
\class(16,15)
\class(16,15)
\class[tmfclass](16,16)
\class(16,16)
\class(16,16)
\class(16,16)
\class(16,16)
\class[tmfclass](16,17)
\class(16,17)
\class(16,17)
\class(16,17)
\class(16,17)
\class[tmfclass](16,18)
\class(16,18)
\class(16,18)
\class(16,18)
\class(16,18)
\class[tmfclass](16,19)
\class(16,19)
\class(16,19)
\class(16,19)
\class(16,19)
\class[tmfclass](16,20)
\class(16,20)
\class(16,20)
\class(16,20)
\class(16,20)
\class[tmfclass](16,21)
\class(16,21)
\class(16,21)
\class(16,21)
\class(16,21)
\class[tmfclass](16,22)
\class(16,22)
\class(16,22)
\class(16,22)
\class(16,22)
\class[tmfclass](16,23)
\class(16,23)
\class(16,23)
\class(16,23)
\class(16,23)
\class[tmfclass](16,24)
\class(16,24)
\class(16,24)
\class(16,24)
\class(16,24)
\d[dthr]3(4,0,1,1)
\d[dfour]4(16,1,1,1)
\structline(0,0,1)(0,1,1)
\structline(0,0,1)(1,1,1)
\structline(0,0,1)(3,1,1)
\structline(0,1,1)(0,2,1)
\structline(0,1,1)(3,2,1)
\structline(0,2,1)(0,3,1)
\structline(0,2,1)(3,3,1)
\structline(0,3,1)(0,4,1)
\structline(0,4,1)(0,5,1)
\structline(0,5,1)(0,6,1)
\structline(0,6,1)(0,7,1)
\structline(0,7,1)(0,8,1)
\structline(0,8,1)(0,9,1)
\structline(0,9,1)(0,10,1)
\structline(0,10,1)(0,11,1)
\structline(0,11,1)(0,12,1)
\structline(0,12,1)(0,13,1)
\structline(0,13,1)(0,14,1)
\structline(0,14,1)(0,15,1)
\structline(0,15,1)(0,16,1)
\structline(0,16,1)(0,17,1)
\structline(0,17,1)(0,18,1)
\structline(0,18,1)(0,19,1)
\structline(0,19,1)(0,20,1)
\structline(0,20,1)(0,21,1)
\structline(0,21,1)(0,22,1)
\structline(0,22,1)(0,23,1)
\structline(0,23,1)(0,24,1)
\structline(1,0,1)(2,1,1)
\structline(1,0,1)(4,1,1)
\structline(1,1,1)(2,2,1)
\structline(2,2,1)(3,3,1)
\structline(3,1,1)(3,2,1)
\structline(3,1,1)(6,2,1)
\structline(3,2,1)(3,3,1)
\structline(4,0,1)(4,1,1)
\structline(4,2,1)(4,3,1)
\structline(4,3,1)(4,4,1)
\structline(4,4,1)(4,5,1)
\structline(4,5,1)(4,6,1)
\structline(4,6,1)(4,7,1)
\structline(4,7,1)(4,8,1)
\structline(4,8,1)(4,9,1)
\structline(4,9,1)(4,10,1)
\structline(4,10,1)(4,11,1)
\structline(4,11,1)(4,12,1)
\structline(4,12,1)(4,13,1)
\structline(4,13,1)(4,14,1)
\structline(4,14,1)(4,15,1)
\structline(4,15,1)(4,16,1)
\structline(4,16,1)(4,17,1)
\structline(4,17,1)(4,18,1)
\structline(4,18,1)(4,19,1)
\structline(4,19,1)(4,20,1)
\structline(4,20,1)(4,21,1)
\structline(4,21,1)(4,22,1)
\structline(4,22,1)(4,23,1)
\structline(4,23,1)(4,24,1)
\structline(8,1,1)(8,2,1)
\structline(8,1,1)(9,2,1)
\structline(8,2,1)(8,3,3)
\structline(8,3,2)(8,4,2)
\structline(8,3,3)(8,4,3)
\structline(8,3,1)(9,4,1)
\structline(8,4,1)(8,5,1)
\structline(8,4,2)(8,5,2)
\structline(8,4,3)(8,5,3)
\structline(8,4,1)(9,5,1)
\structline(8,5,1)(8,6,1)
\structline(8,5,2)(8,6,2)
\structline(8,5,3)(8,6,3)
\structline(8,6,1)(8,7,1)
\structline(8,6,2)(8,7,2)
\structline(8,6,3)(8,7,3)
\structline(8,7,1)(8,8,1)
\structline(8,7,2)(8,8,2)
\structline(8,7,3)(8,8,3)
\structline(8,8,1)(8,9,1)
\structline(8,8,2)(8,9,2)
\structline(8,8,3)(8,9,3)
\structline(8,9,1)(8,10,1)
\structline(8,9,2)(8,10,2)
\structline(8,9,3)(8,10,3)
\structline(8,10,1)(8,11,1)
\structline(8,10,2)(8,11,2)
\structline(8,10,3)(8,11,3)
\structline(8,11,1)(8,12,1)
\structline(8,11,2)(8,12,2)
\structline(8,11,3)(8,12,3)
\structline(8,12,1)(8,13,1)
\structline(8,12,2)(8,13,2)
\structline(8,12,3)(8,13,3)
\structline(8,13,1)(8,14,1)
\structline(8,13,2)(8,14,2)
\structline(8,13,3)(8,14,3)
\structline(8,14,1)(8,15,1)
\structline(8,14,2)(8,15,2)
\structline(8,14,3)(8,15,3)
\structline(8,15,1)(8,16,1)
\structline(8,15,2)(8,16,2)
\structline(8,15,3)(8,16,3)
\structline(8,16,1)(8,17,1)
\structline(8,16,2)(8,17,2)
\structline(8,16,3)(8,17,3)
\structline(8,17,1)(8,18,1)
\structline(8,17,2)(8,18,2)
\structline(8,17,3)(8,18,3)
\structline(8,18,1)(8,19,1)
\structline(8,18,2)(8,19,2)
\structline(8,18,3)(8,19,3)
\structline(8,19,1)(8,20,1)
\structline(8,19,2)(8,20,2)
\structline(8,19,3)(8,20,3)
\structline(8,20,1)(8,21,1)
\structline(8,20,2)(8,21,2)
\structline(8,20,3)(8,21,3)
\structline(8,21,1)(8,22,1)
\structline(8,21,2)(8,22,2)
\structline(8,21,3)(8,22,3)
\structline(8,22,1)(8,23,1)
\structline(8,22,2)(8,23,2)
\structline(8,22,3)(8,23,3)
\structline(8,23,1)(8,24,1)
\structline(8,23,2)(8,24,2)
\structline(8,23,3)(8,24,3)
\structline(9,1,1)(10,2,1)
\structline(9,2,1)(10,3,1)
\structline(9,5,1)(10,6,1)
\structline(12,3,1)(12,4,1)
\structline(12,4,1)(12,5,1)
\structline(12,5,1)(12,6,2)
\structline(12,6,1)(12,7,1)
\structline(12,6,2)(12,7,2)
\structline(12,6,3)(12,7,3)
\structline(12,7,1)(12,8,1)
\structline(12,7,2)(12,8,2)
\structline(12,7,3)(12,8,3)
\structline(12,8,1)(12,9,1)
\structline(12,8,2)(12,9,2)
\structline(12,8,3)(12,9,3)
\structline(12,9,1)(12,10,1)
\structline(12,9,2)(12,10,2)
\structline(12,9,3)(12,10,3)
\structline(12,10,1)(12,11,1)
\structline(12,10,2)(12,11,2)
\structline(12,10,3)(12,11,3)
\structline(12,11,1)(12,12,1)
\structline(12,11,2)(12,12,2)
\structline(12,11,3)(12,12,3)
\structline(12,12,1)(12,13,1)
\structline(12,12,2)(12,13,2)
\structline(12,12,3)(12,13,3)
\structline(12,13,1)(12,14,1)
\structline(12,13,2)(12,14,2)
\structline(12,13,3)(12,14,3)
\structline(12,14,1)(12,15,1)
\structline(12,14,2)(12,15,2)
\structline(12,14,3)(12,15,3)
\structline(12,15,1)(12,16,1)
\structline(12,15,2)(12,16,2)
\structline(12,15,3)(12,16,3)
\structline(12,16,1)(12,17,1)
\structline(12,16,2)(12,17,2)
\structline(12,16,3)(12,17,3)
\structline(12,17,1)(12,18,1)
\structline(12,17,2)(12,18,2)
\structline(12,17,3)(12,18,3)
\structline(12,18,1)(12,19,1)
\structline(12,18,2)(12,19,2)
\structline(12,18,3)(12,19,3)
\structline(12,19,1)(12,20,1)
\structline(12,19,2)(12,20,2)
\structline(12,19,3)(12,20,3)
\structline(12,20,1)(12,21,1)
\structline(12,20,2)(12,21,2)
\structline(12,20,3)(12,21,3)
\structline(12,21,1)(12,22,1)
\structline(12,21,2)(12,22,2)
\structline(12,21,3)(12,22,3)
\structline(12,22,1)(12,23,1)
\structline(12,22,2)(12,23,2)
\structline(12,22,3)(12,23,3)
\structline(12,23,1)(12,24,1)
\structline(12,23,2)(12,24,2)
\structline(12,23,3)(12,24,3)
\structline(14,4,1)(15,5,1)
\structline(16,1,1)(16,2,1)
\structline(16,2,1)(16,3,2)
\structline(16,3,1)(16,4,1)
\structline(16,3,2)(16,4,3)
\structline(16,4,1)(16,5,1)
\structline(16,4,2)(16,5,2)
\structline(16,4,3)(16,5,4)
\structline(16,5,1)(16,6,1)
\structline(16,5,2)(16,6,2)
\structline(16,5,4)(16,6,3)
\structline(16,6,1)(16,7,2)
\structline(16,6,2)(16,7,3)
\structline(16,6,3)(16,7,5)
\structline(16,7,2)(16,8,2)
\structline(16,7,3)(16,8,3)
\structline(16,7,4)(16,8,4)
\structline(16,7,5)(16,8,5)
\structline(16,8,1)(16,9,1)
\structline(16,8,2)(16,9,2)
\structline(16,8,3)(16,9,3)
\structline(16,8,4)(16,9,4)
\structline(16,8,5)(16,9,5)
\structline(16,9,1)(16,10,1)
\structline(16,9,2)(16,10,2)
\structline(16,9,3)(16,10,3)
\structline(16,9,4)(16,10,4)
\structline(16,9,5)(16,10,5)
\structline(16,10,1)(16,11,1)
\structline(16,10,2)(16,11,2)
\structline(16,10,3)(16,11,3)
\structline(16,10,4)(16,11,4)
\structline(16,10,5)(16,11,5)
\structline(16,11,1)(16,12,1)
\structline(16,11,2)(16,12,2)
\structline(16,11,3)(16,12,3)
\structline(16,11,4)(16,12,4)
\structline(16,11,5)(16,12,5)
\structline(16,12,1)(16,13,1)
\structline(16,12,2)(16,13,2)
\structline(16,12,3)(16,13,3)
\structline(16,12,4)(16,13,4)
\structline(16,12,5)(16,13,5)
\structline(16,13,1)(16,14,1)
\structline(16,13,2)(16,14,2)
\structline(16,13,3)(16,14,3)
\structline(16,13,4)(16,14,4)
\structline(16,13,5)(16,14,5)
\structline(16,14,1)(16,15,1)
\structline(16,14,2)(16,15,2)
\structline(16,14,3)(16,15,3)
\structline(16,14,4)(16,15,4)
\structline(16,14,5)(16,15,5)
\structline(16,15,1)(16,16,1)
\structline(16,15,2)(16,16,2)
\structline(16,15,3)(16,16,3)
\structline(16,15,4)(16,16,4)
\structline(16,15,5)(16,16,5)
\structline(16,16,1)(16,17,1)
\structline(16,16,2)(16,17,2)
\structline(16,16,3)(16,17,3)
\structline(16,16,4)(16,17,4)
\structline(16,16,5)(16,17,5)
\structline(16,17,1)(16,18,1)
\structline(16,17,2)(16,18,2)
\structline(16,17,3)(16,18,3)
\structline(16,17,4)(16,18,4)
\structline(16,17,5)(16,18,5)
\structline(16,18,1)(16,19,1)
\structline(16,18,2)(16,19,2)
\structline(16,18,3)(16,19,3)
\structline(16,18,4)(16,19,4)
\structline(16,18,5)(16,19,5)
\structline(16,19,1)(16,20,1)
\structline(16,19,2)(16,20,2)
\structline(16,19,3)(16,20,3)
\structline(16,19,4)(16,20,4)
\structline(16,19,5)(16,20,5)
\structline(16,20,1)(16,21,1)
\structline(16,20,2)(16,21,2)
\structline(16,20,3)(16,21,3)
\structline(16,20,4)(16,21,4)
\structline(16,20,5)(16,21,5)
\structline(16,21,1)(16,22,1)
\structline(16,21,2)(16,22,2)
\structline(16,21,3)(16,22,3)
\structline(16,21,4)(16,22,4)
\structline(16,21,5)(16,22,5)
\structline(16,22,1)(16,23,1)
\structline(16,22,2)(16,23,2)
\structline(16,22,3)(16,23,3)
\structline(16,22,4)(16,23,4)
\structline(16,22,5)(16,23,5)
\structline(16,23,1)(16,24,1)
\structline(16,23,2)(16,24,2)
\structline(16,23,3)(16,24,3)
\structline(16,23,4)(16,24,4)
\structline(16,23,5)(16,24,5)
\end{sseqdata}

\begin{sseqdata}[name = AdamsEthreeB, adams, x range = {16}{33}, y range = {0}{14}]
\class(16,1)
\class(16,2)
\class(16,3)
\class(16,3)
\class(16,4)
\class(16,4)
\class(16,4)
\class(16,5)
\class(16,5)
\class(16,5)
\class(16,5)
\class(16,6)
\class(16,6)
\class(16,6)
\class[tmfclass](16,7)
\class(16,7)
\class(16,7)
\class(16,7)
\class(16,7)
\class[tmfclass](16,8)
\class(16,8)
\class(16,8)
\class(16,8)
\class(16,8)
\class[tmfclass](16,9)
\class(16,9)
\class(16,9)
\class(16,9)
\class(16,9)
\class[tmfclass](16,10)
\class(16,10)
\class(16,10)
\class(16,10)
\class(16,10)
\class[tmfclass](16,11)
\class(16,11)
\class(16,11)
\class(16,11)
\class(16,11)
\class[tmfclass](16,12)
\class(16,12)
\class(16,12)
\class(16,12)
\class(16,12)
\class[tmfclass](16,13)
\class(16,13)
\class(16,13)
\class(16,13)
\class(16,13)
\class[tmfclass](16,14)
\class(16,14)
\class(16,14)
\class(16,14)
\class(16,14)
\class[tmfclass](16,15)
\class(16,15)
\class(16,15)
\class(16,15)
\class(16,15)
\class[tmfclass](16,16)
\class(16,16)
\class(16,16)
\class(16,16)
\class(16,16)
\class[tmfclass](16,17)
\class(16,17)
\class(16,17)
\class(16,17)
\class(16,17)
\class[tmfclass](16,18)
\class(16,18)
\class(16,18)
\class(16,18)
\class(16,18)
\class[tmfclass](16,19)
\class(16,19)
\class(16,19)
\class(16,19)
\class(16,19)
\class[tmfclass](16,20)
\class(16,20)
\class(16,20)
\class(16,20)
\class(16,20)
\class[tmfclass](16,21)
\class(16,21)
\class(16,21)
\class(16,21)
\class(16,21)
\class[tmfclass](16,22)
\class(16,22)
\class(16,22)
\class(16,22)
\class(16,22)
\class[tmfclass](16,23)
\class(16,23)
\class(16,23)
\class(16,23)
\class(16,23)
\class[tmfclass](16,24)
\class(16,24)
\class(16,24)
\class(16,24)
\class(16,24)
\class(17,2)
\class[tmfclass](17,4)
\class(17,4)
\class[tmfclass](17,5)
\class(17,5)
\class(17,6)
\class[tmfclass](17,8)
\class[tmfclass](17,9)
\class(18,2)
\class(18,3)
\class[tmfclass](18,5)
\class(18,5)
\class(18,6)
\class(18,7)
\class[tmfclass](18,10)
\class(19,2)
\class(19,3)
\class[tmfclass](20,4)
\class(20,4)
\class[tmfclass](20,5)
\class(20,5)
\class[tmfclass](20,6)
\class(20,6)
\class(20,6)
\class(20,7)
\class(20,7)
\class(20,7)
\class(20,8)
\class(20,8)
\class(20,8)
\class(20,8)
\class(20,9)
\class(20,9)
\class(20,9)
\class[tmfclass](20,10)
\class(20,10)
\class(20,10)
\class(20,10)
\class(20,10)
\class[tmfclass](20,11)
\class(20,11)
\class(20,11)
\class(20,11)
\class(20,11)
\class[tmfclass](20,12)
\class(20,12)
\class(20,12)
\class(20,12)
\class(20,12)
\class[tmfclass](20,13)
\class(20,13)
\class(20,13)
\class(20,13)
\class(20,13)
\class[tmfclass](20,14)
\class(20,14)
\class(20,14)
\class(20,14)
\class(20,14)
\class[tmfclass](20,15)
\class(20,15)
\class(20,15)
\class(20,15)
\class(20,15)
\class[tmfclass](20,16)
\class(20,16)
\class(20,16)
\class(20,16)
\class(20,16)
\class[tmfclass](20,17)
\class(20,17)
\class(20,17)
\class(20,17)
\class(20,17)
\class[tmfclass](20,18)
\class(20,18)
\class(20,18)
\class(20,18)
\class(20,18)
\class[tmfclass](20,19)
\class(20,19)
\class(20,19)
\class(20,19)
\class(20,19)
\class[tmfclass](20,20)
\class(20,20)
\class(20,20)
\class(20,20)
\class(20,20)
\class[tmfclass](20,21)
\class(20,21)
\class(20,21)
\class(20,21)
\class(20,21)
\class[tmfclass](20,22)
\class(20,22)
\class(20,22)
\class(20,22)
\class(20,22)
\class[tmfclass](20,23)
\class(20,23)
\class(20,23)
\class(20,23)
\class(20,23)
\class[tmfclass](20,24)
\class(20,24)
\class(20,24)
\class(20,24)
\class(20,24)
\class(21,3)
\class[tmfclass](21,5)
\class(22,4)
\class[tmfclass](22,8)
\class[tmfclass](23,9)
\class(24,3)
\class(24,4)
\class(24,5)
\class(24,5)
\class[tmfclass](24,6)
\class(24,6)
\class(24,6)
\class[tmfclass](24,7)
\class(24,7)
\class(24,7)
\class(24,7)
\class[tmfclass](24,8)
\class(24,8)
\class(24,8)
\class(24,8)
\class(24,8)
\class[tmfclass](24,9)
\class(24,9)
\class(24,9)
\class(24,9)
\class(24,9)
\class(24,9)
\class[tmfclass](24,10)
\class(24,10)
\class(24,10)
\class(24,10)
\class(24,10)
\class[tmfclass](24,11)
\class[tmfclass](24,11)
\class(24,11)
\class(24,11)
\class(24,11)
\class(24,11)
\class(24,11)
\class[tmfclass](24,12)
\class[tmfclass](24,12)
\class(24,12)
\class(24,12)
\class(24,12)
\class(24,12)
\class(24,12)
\class[tmfclass](24,13)
\class[tmfclass](24,13)
\class(24,13)
\class(24,13)
\class(24,13)
\class(24,13)
\class(24,13)
\class[tmfclass](24,14)
\class[tmfclass](24,14)
\class(24,14)
\class(24,14)
\class(24,14)
\class(24,14)
\class(24,14)
\class[tmfclass](24,15)
\class[tmfclass](24,15)
\class(24,15)
\class(24,15)
\class(24,15)
\class(24,15)
\class(24,15)
\class[tmfclass](24,16)
\class[tmfclass](24,16)
\class(24,16)
\class(24,16)
\class(24,16)
\class(24,16)
\class(24,16)
\class[tmfclass](24,17)
\class[tmfclass](24,17)
\class(24,17)
\class(24,17)
\class(24,17)
\class(24,17)
\class(24,17)
\class[tmfclass](24,18)
\class[tmfclass](24,18)
\class(24,18)
\class(24,18)
\class(24,18)
\class(24,18)
\class(24,18)
\class[tmfclass](24,19)
\class[tmfclass](24,19)
\class(24,19)
\class(24,19)
\class(24,19)
\class(24,19)
\class(24,19)
\class[tmfclass](24,20)
\class[tmfclass](24,20)
\class(24,20)
\class(24,20)
\class(24,20)
\class(24,20)
\class(24,20)
\class[tmfclass](24,21)
\class[tmfclass](24,21)
\class(24,21)
\class(24,21)
\class(24,21)
\class(24,21)
\class(24,21)
\class[tmfclass](24,22)
\class[tmfclass](24,22)
\class(24,22)
\class(24,22)
\class(24,22)
\class(24,22)
\class(24,22)
\class[tmfclass](24,23)
\class[tmfclass](24,23)
\class(24,23)
\class(24,23)
\class(24,23)
\class(24,23)
\class(24,23)
\class[tmfclass](24,24)
\class[tmfclass](24,24)
\class(24,24)
\class(24,24)
\class(24,24)
\class(24,24)
\class(24,24)
\class(25,4)
\class[tmfclass](25,5)
\class(25,6)
\class[tmfclass](25,8)
\class(25,8)
\class(25,9)
\class(25,10)
\class[tmfclass](25,12)
\class[tmfclass](25,13)
\class[tmfclass](26,6)
\class(26,6)
\class(26,7)
\class[tmfclass](26,9)
\class(26,9)
\class(26,10)
\class(26,11)
\class[tmfclass](26,14)
\class[tmfclass](27,6)
\class(28,6)
\class(28,6)
\class(28,7)
\class(28,7)
\class[tmfclass](28,8)
\class(28,8)
\class(28,8)
\class(28,8)
\class(28,9)
\class(28,9)
\class(28,9)
\class(28,10)
\class(28,10)
\class(28,10)
\class(28,10)
\class(28,11)
\class(28,11)
\class(28,11)
\class(28,11)
\class(28,11)
\class(28,12)
\class(28,12)
\class(28,12)
\class(28,12)
\class(28,12)
\class(28,12)
\class(28,13)
\class(28,13)
\class(28,13)
\class(28,13)
\class(28,13)
\class[tmfclass](28,14)
\class(28,14)
\class(28,14)
\class(28,14)
\class(28,14)
\class(28,14)
\class(28,14)
\class[tmfclass](28,15)
\class(28,15)
\class(28,15)
\class(28,15)
\class(28,15)
\class(28,15)
\class(28,15)
\class[tmfclass](28,16)
\class(28,16)
\class(28,16)
\class(28,16)
\class(28,16)
\class(28,16)
\class(28,16)
\class[tmfclass](28,17)
\class(28,17)
\class(28,17)
\class(28,17)
\class(28,17)
\class(28,17)
\class(28,17)
\class[tmfclass](28,18)
\class(28,18)
\class(28,18)
\class(28,18)
\class(28,18)
\class(28,18)
\class(28,18)
\class[tmfclass](28,19)
\class(28,19)
\class(28,19)
\class(28,19)
\class(28,19)
\class(28,19)
\class(28,19)
\class[tmfclass](28,20)
\class(28,20)
\class(28,20)
\class(28,20)
\class(28,20)
\class(28,20)
\class(28,20)
\class[tmfclass](28,21)
\class(28,21)
\class(28,21)
\class(28,21)
\class(28,21)
\class(28,21)
\class(28,21)
\class[tmfclass](28,22)
\class(28,22)
\class(28,22)
\class(28,22)
\class(28,22)
\class(28,22)
\class(28,22)
\class[tmfclass](28,23)
\class(28,23)
\class(28,23)
\class(28,23)
\class(28,23)
\class(28,23)
\class(28,23)
\class[tmfclass](28,24)
\class(28,24)
\class(28,24)
\class(28,24)
\class(28,24)
\class(28,24)
\class(28,24)
\class[tmfclass](29,9)
\class[tmfclass](30,6)
\class(30,8)
\class[tmfclass](30,12)
\class(31,5)
\class[tmfclass](31,8)
\class[tmfclass](31,13)
\class(32,4)
\class(32,5)
\class(32,5)
\class(32,6)
\class[tmfclass](32,7)
\class[tmfclass](32,7)
\class(32,7)
\class(32,7)
\class(32,7)
\class[tmfclass](32,8)
\class(32,8)
\class(32,8)
\class[tmfclass](32,9)
\class(32,9)
\class(32,9)
\class(32,9)
\class(32,9)
\class[tmfclass](32,10)
\class(32,10)
\class(32,10)
\class(32,10)
\class(32,10)
\class(32,11)
\class(32,11)
\class[tmfclass](32,11)
\class(32,11)
\class(32,11)
\class(32,11)
\class[tmfclass](32,12)
\class(32,12)
\class(32,12)
\class(32,12)
\class(32,12)
\class(32,12)
\class(32,12)
\class[tmfclass](32,13)
\class(32,13)
\class(32,13)
\class(32,13)
\class(32,13)
\class(32,13)
\class(32,13)
\class(32,13)
\class[tmfclass](32,14)
\class(32,14)
\class(32,14)
\class(32,14)
\class(32,14)
\class(32,14)
\class(32,14)
\class[tmfclass](32,15)
\class[tmfclass](32,15)
\class(32,15)
\class(32,15)
\class(32,15)
\class(32,15)
\class(32,15)
\class(32,15)
\class(32,15)
\class[tmfclass](32,16)
\class[tmfclass](32,16)
\class(32,16)
\class(32,16)
\class(32,16)
\class(32,16)
\class(32,16)
\class(32,16)
\class(32,16)
\class[tmfclass](32,17)
\class[tmfclass](32,17)
\class(32,17)
\class(32,17)
\class(32,17)
\class(32,17)
\class(32,17)
\class(32,17)
\class(32,17)
\class[tmfclass](32,18)
\class[tmfclass](32,18)
\class(32,18)
\class(32,18)
\class(32,18)
\class(32,18)
\class(32,18)
\class(32,18)
\class(32,18)
\class[tmfclass](32,19)
\class[tmfclass](32,19)
\class(32,19)
\class(32,19)
\class(32,19)
\class(32,19)
\class(32,19)
\class(32,19)
\class(32,19)
\class[tmfclass](32,20)
\class[tmfclass](32,20)
\class(32,20)
\class(32,20)
\class(32,20)
\class(32,20)
\class(32,20)
\class(32,20)
\class(32,20)
\class[tmfclass](32,21)
\class[tmfclass](32,21)
\class(32,21)
\class(32,21)
\class(32,21)
\class(32,21)
\class(32,21)
\class(32,21)
\class(32,21)
\class[tmfclass](32,22)
\class[tmfclass](32,22)
\class(32,22)
\class(32,22)
\class(32,22)
\class(32,22)
\class(32,22)
\class(32,22)
\class(32,22)
\class[tmfclass](32,23)
\class[tmfclass](32,23)
\class(32,23)
\class(32,23)
\class(32,23)
\class(32,23)
\class(32,23)
\class(32,23)
\class(32,23)
\class[tmfclass](32,24)
\class[tmfclass](32,24)
\class(32,24)
\class(32,24)
\class(32,24)
\class(32,24)
\class(32,24)
\class(32,24)
\class(32,24)
\class(33,4)
\class(33,5)
\class(33,5)
\class(33,6)
\class[tmfclass](33,8)
\class(33,8)
\class(33,8)
\class[tmfclass](33,9)
\class(33,10)
\class[tmfclass](33,12)
\class(33,12)
\class(33,13)
\class(33,14)
\class[tmfclass](33,16)
\class[tmfclass](33,17)
\d[dthr]3(17,4,1,1)
\d[dthr]3(18,2,1,1)
\d[dthr]3(18,5,1,1)
\d[dthr]3(21,3,1,1)
\d[dthr]3(24,6,1,1)
\d[dthr]3(25,8,1,1)
\d[dthr]3(26,9,1,1)
\d[dthr]3(30,6,1,1)
\d[dthr]3(31,5,1,1)
\d[dthr]3(32,10,1,1)
\d[dthr]3(33,4,1,1)
\d[dthr]3(33,12,1,1)
\d[dfour]4(31,8,1,1)
\structline(16,1,1)(16,2,1)
\structline(16,1,1)(17,2,1)
\structline(16,1,1)(19,2,1)
\structline(16,2,1)(16,3,2)
\structline(16,2,1)(19,3,1)
\structline(16,3,1)(16,4,1)
\structline(16,3,2)(16,4,3)
\structline(16,3,1)(17,4,2)
\structline(16,4,1)(16,5,1)
\structline(16,4,2)(16,5,2)
\structline(16,4,3)(16,5,4)
\structline(16,4,2)(17,5,2)
\structline(16,5,1)(16,6,1)
\structline(16,5,2)(16,6,2)
\structline(16,5,4)(16,6,3)
\structline(16,5,3)(17,6,1)
\structline(16,6,1)(16,7,2)
\structline(16,6,2)(16,7,3)
\structline(16,6,3)(16,7,5)
\structline(16,7,2)(16,8,2)
\structline(16,7,3)(16,8,3)
\structline(16,7,4)(16,8,4)
\structline(16,7,5)(16,8,5)
\structline(16,7,1)(17,8,1)
\structline(16,8,1)(16,9,1)
\structline(16,8,2)(16,9,2)
\structline(16,8,3)(16,9,3)
\structline(16,8,4)(16,9,4)
\structline(16,8,5)(16,9,5)
\structline(16,8,1)(17,9,1)
\structline(16,9,1)(16,10,1)
\structline(16,9,2)(16,10,2)
\structline(16,9,3)(16,10,3)
\structline(16,9,4)(16,10,4)
\structline(16,9,5)(16,10,5)
\structline(16,10,1)(16,11,1)
\structline(16,10,2)(16,11,2)
\structline(16,10,3)(16,11,3)
\structline(16,10,4)(16,11,4)
\structline(16,10,5)(16,11,5)
\structline(16,11,1)(16,12,1)
\structline(16,11,2)(16,12,2)
\structline(16,11,3)(16,12,3)
\structline(16,11,4)(16,12,4)
\structline(16,11,5)(16,12,5)
\structline(16,12,1)(16,13,1)
\structline(16,12,2)(16,13,2)
\structline(16,12,3)(16,13,3)
\structline(16,12,4)(16,13,4)
\structline(16,12,5)(16,13,5)
\structline(16,13,1)(16,14,1)
\structline(16,13,2)(16,14,2)
\structline(16,13,3)(16,14,3)
\structline(16,13,4)(16,14,4)
\structline(16,13,5)(16,14,5)
\structline(16,14,1)(16,15,1)
\structline(16,14,2)(16,15,2)
\structline(16,14,3)(16,15,3)
\structline(16,14,4)(16,15,4)
\structline(16,14,5)(16,15,5)
\structline(16,15,1)(16,16,1)
\structline(16,15,2)(16,16,2)
\structline(16,15,3)(16,16,3)
\structline(16,15,4)(16,16,4)
\structline(16,15,5)(16,16,5)
\structline(16,16,1)(16,17,1)
\structline(16,16,2)(16,17,2)
\structline(16,16,3)(16,17,3)
\structline(16,16,4)(16,17,4)
\structline(16,16,5)(16,17,5)
\structline(16,17,1)(16,18,1)
\structline(16,17,2)(16,18,2)
\structline(16,17,3)(16,18,3)
\structline(16,17,4)(16,18,4)
\structline(16,17,5)(16,18,5)
\structline(16,18,1)(16,19,1)
\structline(16,18,2)(16,19,2)
\structline(16,18,3)(16,19,3)
\structline(16,18,4)(16,19,4)
\structline(16,18,5)(16,19,5)
\structline(16,19,1)(16,20,1)
\structline(16,19,2)(16,20,2)
\structline(16,19,3)(16,20,3)
\structline(16,19,4)(16,20,4)
\structline(16,19,5)(16,20,5)
\structline(16,20,1)(16,21,1)
\structline(16,20,2)(16,21,2)
\structline(16,20,3)(16,21,3)
\structline(16,20,4)(16,21,4)
\structline(16,20,5)(16,21,5)
\structline(16,21,1)(16,22,1)
\structline(16,21,2)(16,22,2)
\structline(16,21,3)(16,22,3)
\structline(16,21,4)(16,22,4)
\structline(16,21,5)(16,22,5)
\structline(16,22,1)(16,23,1)
\structline(16,22,2)(16,23,2)
\structline(16,22,3)(16,23,3)
\structline(16,22,4)(16,23,4)
\structline(16,22,5)(16,23,5)
\structline(16,23,1)(16,24,1)
\structline(16,23,2)(16,24,2)
\structline(16,23,3)(16,24,3)
\structline(16,23,4)(16,24,4)
\structline(16,23,5)(16,24,5)
\structline(17,2,1)(18,3,1)
\structline(17,4,1)(17,5,1)
\structline(17,4,1)(18,5,1)
\structline(17,4,2)(18,5,2)
\structline(17,4,1)(20,5,1)
\structline(17,5,2)(18,6,1)
\structline(17,5,1)(20,6,1)
\structline(17,6,1)(18,7,1)
\structline(17,9,1)(18,10,1)
\structline(18,2,1)(21,3,1)
\structline(19,2,1)(19,3,1)
\structline(20,4,1)(20,5,1)
\structline(20,4,2)(20,5,2)
\structline(20,4,1)(21,5,1)
\structline(20,5,1)(20,6,1)
\structline(20,5,2)(20,6,2)
\structline(20,5,2)(20,6,3)
\structline(20,6,2)(20,7,1)
\structline(20,6,3)(20,7,2)
\structline(20,7,1)(20,8,1)
\structline(20,7,2)(20,8,2)
\structline(20,7,3)(20,8,3)
\structline(20,8,1)(20,9,1)
\structline(20,8,2)(20,9,2)
\structline(20,8,3)(20,9,3)
\structline(20,9,1)(20,10,2)
\structline(20,9,2)(20,10,3)
\structline(20,9,3)(20,10,4)
\structline(20,10,1)(20,11,1)
\structline(20,10,2)(20,11,2)
\structline(20,10,3)(20,11,3)
\structline(20,10,4)(20,11,4)
\structline(20,10,5)(20,11,5)
\structline(20,11,1)(20,12,1)
\structline(20,11,2)(20,12,2)
\structline(20,11,3)(20,12,3)
\structline(20,11,4)(20,12,4)
\structline(20,11,5)(20,12,5)
\structline(20,12,1)(20,13,1)
\structline(20,12,2)(20,13,2)
\structline(20,12,3)(20,13,3)
\structline(20,12,4)(20,13,4)
\structline(20,12,5)(20,13,5)
\structline(20,13,1)(20,14,1)
\structline(20,13,2)(20,14,2)
\structline(20,13,3)(20,14,3)
\structline(20,13,4)(20,14,4)
\structline(20,13,5)(20,14,5)
\structline(20,14,1)(20,15,1)
\structline(20,14,2)(20,15,2)
\structline(20,14,3)(20,15,3)
\structline(20,14,4)(20,15,4)
\structline(20,14,5)(20,15,5)
\structline(20,15,1)(20,16,1)
\structline(20,15,2)(20,16,2)
\structline(20,15,3)(20,16,3)
\structline(20,15,4)(20,16,4)
\structline(20,15,5)(20,16,5)
\structline(20,16,1)(20,17,1)
\structline(20,16,2)(20,17,2)
\structline(20,16,3)(20,17,3)
\structline(20,16,4)(20,17,4)
\structline(20,16,5)(20,17,5)
\structline(20,17,1)(20,18,1)
\structline(20,17,2)(20,18,2)
\structline(20,17,3)(20,18,3)
\structline(20,17,4)(20,18,4)
\structline(20,17,5)(20,18,5)
\structline(20,18,1)(20,19,1)
\structline(20,18,2)(20,19,2)
\structline(20,18,3)(20,19,3)
\structline(20,18,4)(20,19,4)
\structline(20,18,5)(20,19,5)
\structline(20,19,1)(20,20,1)
\structline(20,19,2)(20,20,2)
\structline(20,19,3)(20,20,3)
\structline(20,19,4)(20,20,4)
\structline(20,19,5)(20,20,5)
\structline(20,20,1)(20,21,1)
\structline(20,20,2)(20,21,2)
\structline(20,20,3)(20,21,3)
\structline(20,20,4)(20,21,4)
\structline(20,20,5)(20,21,5)
\structline(20,21,1)(20,22,1)
\structline(20,21,2)(20,22,2)
\structline(20,21,3)(20,22,3)
\structline(20,21,4)(20,22,4)
\structline(20,21,5)(20,22,5)
\structline(20,22,1)(20,23,1)
\structline(20,22,2)(20,23,2)
\structline(20,22,3)(20,23,3)
\structline(20,22,4)(20,23,4)
\structline(20,22,5)(20,23,5)
\structline(20,23,1)(20,24,1)
\structline(20,23,2)(20,24,2)
\structline(20,23,3)(20,24,3)
\structline(20,23,4)(20,24,4)
\structline(20,23,5)(20,24,5)
\structline(21,3,1)(24,4,1)
\structline(22,8,1)(23,9,1)
\structline(24,3,1)(25,4,1)
\structline(24,5,1)(24,6,2)
\structline(24,5,2)(24,6,2)
\structline(24,5,2)(24,6,3)
\structline(24,5,1)(25,6,1)
\structline(24,6,1)(24,7,1)
\structline(24,6,2)(24,7,2)
\structline(24,6,3)(24,7,2)
\structline(24,6,3)(24,7,3)
\structline(24,6,3)(24,7,4)
\structline(24,7,1)(24,8,1)
\structline(24,7,2)(24,8,2)
\structline(24,7,3)(24,8,3)
\structline(24,7,4)(24,8,4)
\structline(24,7,3)(25,8,2)
\structline(24,7,4)(25,8,2)
\structline(24,8,1)(24,9,1)
\structline(24,8,2)(24,9,2)
\structline(24,8,3)(24,9,3)
\structline(24,8,4)(24,9,4)
\structline(24,8,5)(24,9,5)
\structline(24,8,5)(25,9,1)
\structline(24,9,1)(24,10,1)
\structline(24,9,2)(24,10,2)
\structline(24,9,3)(24,10,3)
\structline(24,9,4)(24,10,4)
\structline(24,9,5)(24,10,5)
\structline(24,9,6)(25,10,1)
\structline(24,10,1)(24,11,2)
\structline(24,10,2)(24,11,3)
\structline(24,10,3)(24,11,4)
\structline(24,10,4)(24,11,5)
\structline(24,10,5)(24,11,6)
\structline(24,11,2)(24,12,1)
\structline(24,11,3)(24,12,3)
\structline(24,11,4)(24,12,4)
\structline(24,11,5)(24,12,5)
\structline(24,11,6)(24,12,6)
\structline(24,11,7)(24,12,7)
\structline(24,11,1)(25,12,1)
\structline(24,12,1)(24,13,1)
\structline(24,12,2)(24,13,2)
\structline(24,12,3)(24,13,3)
\structline(24,12,4)(24,13,4)
\structline(24,12,5)(24,13,5)
\structline(24,12,6)(24,13,6)
\structline(24,12,7)(24,13,7)
\structline(24,12,2)(25,13,1)
\structline(24,13,1)(24,14,1)
\structline(24,13,2)(24,14,2)
\structline(24,13,3)(24,14,3)
\structline(24,13,4)(24,14,4)
\structline(24,13,5)(24,14,5)
\structline(24,13,6)(24,14,6)
\structline(24,13,7)(24,14,7)
\structline(24,14,1)(24,15,1)
\structline(24,14,2)(24,15,2)
\structline(24,14,3)(24,15,3)
\structline(24,14,4)(24,15,4)
\structline(24,14,5)(24,15,5)
\structline(24,14,6)(24,15,6)
\structline(24,14,7)(24,15,7)
\structline(24,15,1)(24,16,1)
\structline(24,15,2)(24,16,2)
\structline(24,15,3)(24,16,3)
\structline(24,15,4)(24,16,4)
\structline(24,15,5)(24,16,5)
\structline(24,15,6)(24,16,6)
\structline(24,15,7)(24,16,7)
\structline(24,16,1)(24,17,1)
\structline(24,16,2)(24,17,2)
\structline(24,16,3)(24,17,3)
\structline(24,16,4)(24,17,4)
\structline(24,16,5)(24,17,5)
\structline(24,16,6)(24,17,6)
\structline(24,16,7)(24,17,7)
\structline(24,17,1)(24,18,1)
\structline(24,17,2)(24,18,2)
\structline(24,17,3)(24,18,3)
\structline(24,17,4)(24,18,4)
\structline(24,17,5)(24,18,5)
\structline(24,17,6)(24,18,6)
\structline(24,17,7)(24,18,7)
\structline(24,18,1)(24,19,1)
\structline(24,18,2)(24,19,2)
\structline(24,18,3)(24,19,3)
\structline(24,18,4)(24,19,4)
\structline(24,18,5)(24,19,5)
\structline(24,18,6)(24,19,6)
\structline(24,18,7)(24,19,7)
\structline(24,19,1)(24,20,1)
\structline(24,19,2)(24,20,2)
\structline(24,19,3)(24,20,3)
\structline(24,19,4)(24,20,4)
\structline(24,19,5)(24,20,5)
\structline(24,19,6)(24,20,6)
\structline(24,19,7)(24,20,7)
\structline(24,20,1)(24,21,1)
\structline(24,20,2)(24,21,2)
\structline(24,20,3)(24,21,3)
\structline(24,20,4)(24,21,4)
\structline(24,20,5)(24,21,5)
\structline(24,20,6)(24,21,6)
\structline(24,20,7)(24,21,7)
\structline(24,21,1)(24,22,1)
\structline(24,21,2)(24,22,2)
\structline(24,21,3)(24,22,3)
\structline(24,21,4)(24,22,4)
\structline(24,21,5)(24,22,5)
\structline(24,21,6)(24,22,6)
\structline(24,21,7)(24,22,7)
\structline(24,22,1)(24,23,1)
\structline(24,22,2)(24,23,2)
\structline(24,22,3)(24,23,3)
\structline(24,22,4)(24,23,4)
\structline(24,22,5)(24,23,5)
\structline(24,22,6)(24,23,6)
\structline(24,22,7)(24,23,7)
\structline(24,23,1)(24,24,1)
\structline(24,23,2)(24,24,2)
\structline(24,23,3)(24,24,3)
\structline(24,23,4)(24,24,4)
\structline(24,23,5)(24,24,5)
\structline(24,23,6)(24,24,6)
\structline(24,23,7)(24,24,7)
\structline(25,5,1)(26,6,1)
\structline(25,6,1)(26,7,1)
\structline(25,8,1)(26,9,1)
\structline(25,8,2)(26,9,2)
\structline(25,9,1)(26,10,1)
\structline(25,10,1)(26,11,1)
\structline(25,13,1)(26,14,1)
\structline(28,6,1)(28,7,1)
\structline(28,6,2)(28,7,2)
\structline(28,7,1)(28,8,2)
\structline(28,7,2)(28,8,4)
\structline(28,8,2)(28,9,1)
\structline(28,8,3)(28,9,2)
\structline(28,8,4)(28,9,3)
\structline(28,8,1)(29,9,1)
\structline(28,9,1)(28,10,1)
\structline(28,9,2)(28,10,2)
\structline(28,9,3)(28,10,4)
\structline(28,10,1)(28,11,1)
\structline(28,10,2)(28,11,2)
\structline(28,10,3)(28,11,3)
\structline(28,10,4)(28,11,4)
\structline(28,11,1)(28,12,1)
\structline(28,11,2)(28,12,2)
\structline(28,11,3)(28,12,3)
\structline(28,11,4)(28,12,4)
\structline(28,11,5)(28,12,5)
\structline(28,12,1)(28,13,1)
\structline(28,12,2)(28,13,2)
\structline(28,12,3)(28,13,3)
\structline(28,12,4)(28,13,4)
\structline(28,12,5)(28,13,5)
\structline(28,13,1)(28,14,2)
\structline(28,13,2)(28,14,3)
\structline(28,13,3)(28,14,4)
\structline(28,13,4)(28,14,5)
\structline(28,13,5)(28,14,6)
\structline(28,14,1)(28,15,1)
\structline(28,14,2)(28,15,2)
\structline(28,14,3)(28,15,3)
\structline(28,14,4)(28,15,4)
\structline(28,14,5)(28,15,5)
\structline(28,14,6)(28,15,6)
\structline(28,14,7)(28,15,7)
\structline(28,15,1)(28,16,1)
\structline(28,15,2)(28,16,2)
\structline(28,15,3)(28,16,3)
\structline(28,15,4)(28,16,4)
\structline(28,15,5)(28,16,5)
\structline(28,15,6)(28,16,6)
\structline(28,15,7)(28,16,7)
\structline(28,16,1)(28,17,1)
\structline(28,16,2)(28,17,2)
\structline(28,16,3)(28,17,3)
\structline(28,16,4)(28,17,4)
\structline(28,16,5)(28,17,5)
\structline(28,16,6)(28,17,6)
\structline(28,16,7)(28,17,7)
\structline(28,17,1)(28,18,1)
\structline(28,17,2)(28,18,2)
\structline(28,17,3)(28,18,3)
\structline(28,17,4)(28,18,4)
\structline(28,17,5)(28,18,5)
\structline(28,17,6)(28,18,6)
\structline(28,17,7)(28,18,7)
\structline(28,18,1)(28,19,1)
\structline(28,18,2)(28,19,2)
\structline(28,18,3)(28,19,3)
\structline(28,18,4)(28,19,4)
\structline(28,18,5)(28,19,5)
\structline(28,18,6)(28,19,6)
\structline(28,18,7)(28,19,7)
\structline(28,19,1)(28,20,1)
\structline(28,19,2)(28,20,2)
\structline(28,19,3)(28,20,3)
\structline(28,19,4)(28,20,4)
\structline(28,19,5)(28,20,5)
\structline(28,19,6)(28,20,6)
\structline(28,19,7)(28,20,7)
\structline(28,20,1)(28,21,1)
\structline(28,20,2)(28,21,2)
\structline(28,20,3)(28,21,3)
\structline(28,20,4)(28,21,4)
\structline(28,20,5)(28,21,5)
\structline(28,20,6)(28,21,6)
\structline(28,20,7)(28,21,7)
\structline(28,21,1)(28,22,1)
\structline(28,21,2)(28,22,2)
\structline(28,21,3)(28,22,3)
\structline(28,21,4)(28,22,4)
\structline(28,21,5)(28,22,5)
\structline(28,21,6)(28,22,6)
\structline(28,21,7)(28,22,7)
\structline(28,22,1)(28,23,1)
\structline(28,22,2)(28,23,2)
\structline(28,22,3)(28,23,3)
\structline(28,22,4)(28,23,4)
\structline(28,22,5)(28,23,5)
\structline(28,22,6)(28,23,6)
\structline(28,22,7)(28,23,7)
\structline(28,23,1)(28,24,1)
\structline(28,23,2)(28,24,2)
\structline(28,23,3)(28,24,3)
\structline(28,23,4)(28,24,4)
\structline(28,23,5)(28,24,5)
\structline(28,23,6)(28,24,6)
\structline(28,23,7)(28,24,7)
\structline(30,12,1)(31,13,1)
\structline(31,8,1)(32,9,1)
\structline(32,4,1)(32,5,1)
\structline(32,4,1)(32,5,2)
\structline(32,4,1)(33,5,1)
\structline(32,5,2)(32,6,1)
\structline(32,5,1)(33,6,1)
\structline(32,5,2)(33,6,1)
\structline(32,6,1)(32,7,3)
\structline(32,6,1)(32,7,5)
\structline(32,7,2)(32,8,1)
\structline(32,7,3)(32,8,2)
\structline(32,7,4)(32,8,2)
\structline(32,7,5)(32,8,2)
\structline(32,7,5)(32,8,3)
\structline(32,7,1)(33,8,1)
\structline(32,7,2)(33,8,1)
\structline(32,7,3)(33,8,2)
\structline(32,7,4)(33,8,3)
\structline(32,7,5)(33,8,2)
\structline(32,8,1)(32,9,1)
\structline(32,8,2)(32,9,2)
\structline(32,8,3)(32,9,2)
\structline(32,8,3)(32,9,3)
\structline(32,8,3)(32,9,5)
\structline(32,9,2)(32,10,2)
\structline(32,9,3)(32,10,3)
\structline(32,9,4)(32,10,5)
\structline(32,9,5)(32,10,4)
\structline(32,9,4)(33,10,1)
\structline(32,10,1)(32,11,3)
\structline(32,10,2)(32,11,4)
\structline(32,10,3)(32,11,1)
\structline(32,10,4)(32,11,2)
\structline(32,10,5)(32,11,5)
\structline(32,11,1)(32,12,3)
\structline(32,11,2)(32,12,4)
\structline(32,11,3)(32,12,1)
\structline(32,11,4)(32,12,2)
\structline(32,11,5)(32,12,5)
\structline(32,11,6)(32,12,6)
\structline(32,11,6)(33,12,2)
\structline(32,12,1)(32,13,1)
\structline(32,12,2)(32,13,2)
\structline(32,12,3)(32,13,3)
\structline(32,12,4)(32,13,4)
\structline(32,12,5)(32,13,5)
\structline(32,12,6)(32,13,6)
\structline(32,12,7)(32,13,7)
\structline(32,12,7)(33,13,1)
\structline(32,13,1)(32,14,1)
\structline(32,13,2)(32,14,2)
\structline(32,13,3)(32,14,3)
\structline(32,13,4)(32,14,4)
\structline(32,13,5)(32,14,5)
\structline(32,13,6)(32,14,6)
\structline(32,13,7)(32,14,7)
\structline(32,13,8)(33,14,1)
\structline(32,14,1)(32,15,2)
\structline(32,14,2)(32,15,3)
\structline(32,14,3)(32,15,4)
\structline(32,14,4)(32,15,5)
\structline(32,14,5)(32,15,6)
\structline(32,14,6)(32,15,7)
\structline(32,14,7)(32,15,8)
\structline(32,15,2)(32,16,1)
\structline(32,15,3)(32,16,3)
\structline(32,15,4)(32,16,4)
\structline(32,15,5)(32,16,5)
\structline(32,15,6)(32,16,6)
\structline(32,15,7)(32,16,7)
\structline(32,15,8)(32,16,8)
\structline(32,15,9)(32,16,9)
\structline(32,15,1)(33,16,1)
\structline(32,16,1)(32,17,1)
\structline(32,16,2)(32,17,2)
\structline(32,16,3)(32,17,3)
\structline(32,16,4)(32,17,4)
\structline(32,16,5)(32,17,5)
\structline(32,16,6)(32,17,6)
\structline(32,16,7)(32,17,7)
\structline(32,16,8)(32,17,8)
\structline(32,16,9)(32,17,9)
\structline(32,16,2)(33,17,1)
\structline(32,17,1)(32,18,1)
\structline(32,17,2)(32,18,2)
\structline(32,17,3)(32,18,3)
\structline(32,17,4)(32,18,4)
\structline(32,17,5)(32,18,5)
\structline(32,17,6)(32,18,6)
\structline(32,17,7)(32,18,7)
\structline(32,17,8)(32,18,8)
\structline(32,17,9)(32,18,9)
\structline(32,18,1)(32,19,1)
\structline(32,18,2)(32,19,2)
\structline(32,18,3)(32,19,3)
\structline(32,18,4)(32,19,4)
\structline(32,18,5)(32,19,5)
\structline(32,18,6)(32,19,6)
\structline(32,18,7)(32,19,7)
\structline(32,18,8)(32,19,8)
\structline(32,18,9)(32,19,9)
\structline(32,19,1)(32,20,1)
\structline(32,19,2)(32,20,2)
\structline(32,19,3)(32,20,3)
\structline(32,19,4)(32,20,4)
\structline(32,19,5)(32,20,5)
\structline(32,19,6)(32,20,6)
\structline(32,19,7)(32,20,7)
\structline(32,19,8)(32,20,8)
\structline(32,19,9)(32,20,9)
\structline(32,20,1)(32,21,1)
\structline(32,20,2)(32,21,2)
\structline(32,20,3)(32,21,3)
\structline(32,20,4)(32,21,4)
\structline(32,20,5)(32,21,5)
\structline(32,20,6)(32,21,6)
\structline(32,20,7)(32,21,7)
\structline(32,20,8)(32,21,8)
\structline(32,20,9)(32,21,9)
\structline(32,21,1)(32,22,1)
\structline(32,21,2)(32,22,2)
\structline(32,21,3)(32,22,3)
\structline(32,21,4)(32,22,4)
\structline(32,21,5)(32,22,5)
\structline(32,21,6)(32,22,6)
\structline(32,21,7)(32,22,7)
\structline(32,21,8)(32,22,8)
\structline(32,21,9)(32,22,9)
\structline(32,22,1)(32,23,1)
\structline(32,22,2)(32,23,2)
\structline(32,22,3)(32,23,3)
\structline(32,22,4)(32,23,4)
\structline(32,22,5)(32,23,5)
\structline(32,22,6)(32,23,6)
\structline(32,22,7)(32,23,7)
\structline(32,22,8)(32,23,8)
\structline(32,22,9)(32,23,9)
\structline(32,23,1)(32,24,1)
\structline(32,23,2)(32,24,2)
\structline(32,23,3)(32,24,3)
\structline(32,23,4)(32,24,4)
\structline(32,23,5)(32,24,5)
\structline(32,23,6)(32,24,6)
\structline(32,23,7)(32,24,7)
\structline(32,23,8)(32,24,8)
\structline(32,23,9)(32,24,9)
\structline(33,4,1)(33,5,1)
\structline(33,4,1)(33,5,2)
\end{sseqdata}

\subsection{The $E_\infty$-page and hidden extensions} 
\begin{figure}[htbp]
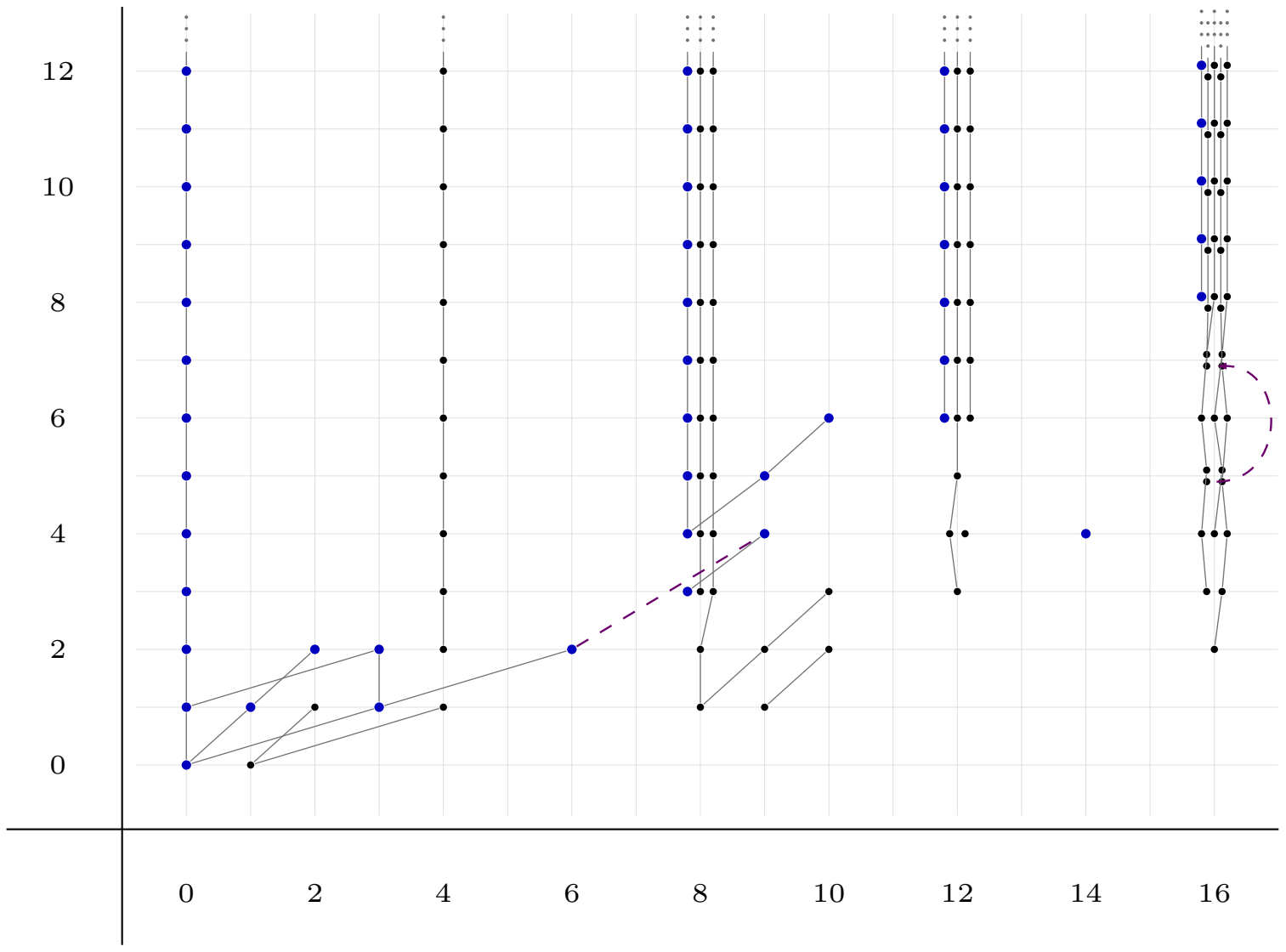
\centering
\adjustbox{max width=\textwidth}{\printpage[name = AdamsEinfA, page = 5]}
\caption{The page $E_5=E_\infty$ in stems $0$ to $16$. The dashed lines are the hidden extensions of Proposition \ref{prop:hidden}: the hidden $\nu$-extension $\nu\cdot\nu^2=\eta\epsilon$ and the hidden extension by $2$ from $y$ to $w_1x_{8,3}$.}
\label{fig:adamsEinfa}
\end{figure}

\begin{figure}[htbp]
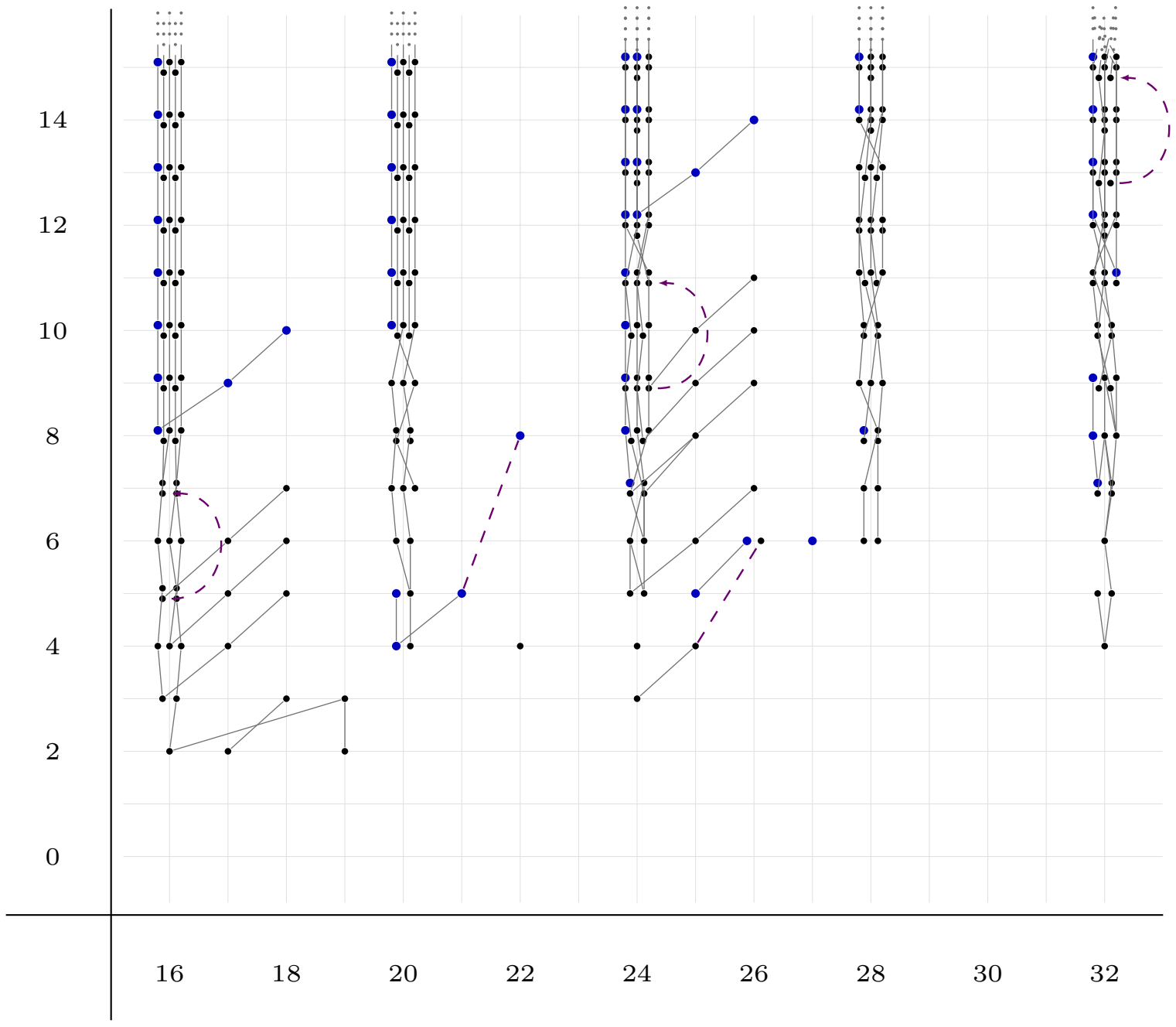
\centering
\adjustbox{max width=\textwidth}{\printpage[name = AdamsEinfB, page = 5]}
\caption{The page $E_5=E_\infty$ in stems $16$ to $32$, with the hidden $\eta$-extensions and the hidden extensions by $2$ of Proposition \ref{prop:hidden} drawn as dashed lines.}
\label{fig:adamsEinfb}
\end{figure}

The result is shown in Figures \ref{fig:adamsEinfa} and \ref{fig:adamsEinfb}.

\begin{corollary}\label{cor:adamsEinf}
In stems $n\le31$, the classes of $E_\infty$ not lying on infinite $h_0$-towers are the following.
\begin{center}\small
\begin{tabular}{@{}r@{\quad}p{0.8\textwidth}@{}}
$n$ & classes\\ \hline
$1$ & $x_{1,0},\ h_1$\\
$2$ & $h_1x_{1,0},\ h_1^2$\\
$3$ & $h_2,\ h_0h_2$\\
$4$ & $h_0x_{4,0}$\\
$6$ & $h_2^2$\\
$8$ & $c_0$\\
$9$ & $x_{9,1},\ h_1x_{8,1},\ h_1c_0,\ h_1w_1$\\
$10$ & $h_1x_{9,1},\ h_0x_{10,2},\ h_1^2w_1$\\
$12$ & $w_1x_{4,0}$\\
$14$ & $d_0$\\
$16$ & $y=w_1x_{8,1}+h_0x_{16,4}$\\
$17$ & $h_1^2x_{15,0},\ h_1x_{16,3},\ w_1x_{9,1},\ h_1w_1x_{8,1},\ h_1w_1^2$\\
$18$ & $h_0^2h_2x_{15,0},\ c_0x_{10,2},\ h_1w_1x_{9,1},\ h_0w_1x_{10,2},\ h_1^2w_1^2$\\
$19$ & $h_0h_2x_{16,0},\ h_0^2h_2x_{16,0}$\\
$20$ & $g,\ h_0g,\ w_1^2x_{4,0}$\\
$21$ & $h_1g$\\
$22$ & $x_{22,4},\ w_1d_0$\\
$24$ & $c_0x_{16,0},\ \beta x_{9,1},\ w_1y$\\
$25$ & $h_1c_0x_{16,0},\ \gamma,\ h_1^2w_1x_{15,0},\ h_1w_1x_{16,3},\ w_1^2x_{9,1},\ h_1w_1^2x_{8,1},\ h_1w_1^3$\\
$26$ & $h_1\gamma,\ w_1x_{18,2},\ h_0^3x_{26,4},\ c_0w_1x_{10,2},\ h_1w_1^2x_{9,1},\ h_0w_1^2x_{10,2},\ h_1^2w_1^3$\\
$27$ & $\alpha\beta$\\
$28$ & $w_1g,\ w_1^3x_{4,0}$\\
\end{tabular}
\end{center}
In stem $32$ they are $x_{32,5}$, $\delta'$, $h_0\delta'$, $w_1^2y$ and, if the $d_4$ on the class in bidegree $(33,5)$ vanishes, $c_0\alpha^2=h_0^2\delta'$. There are $r_n$ infinite towers in stem $n$.
\end{corollary}

So in stems at most $32$, the groups $\pi_n\tmf_0(5)$ are as follows. Their ranks are the $r_n$ above, and their torsion subgroups have the orders of Table \ref{tab:torsionorders}. The torsion classes coming from $\tmf$ are $\eta,\eta^2,\nu,\nu^2,\epsilon,\eta\epsilon,\kappa,\bar\kappa,\eta\bar\kappa,\eta_1,\eta\eta_1,\nu_1$ and $\kappa^2=\epsilon\bar\kappa$, together with the products $\eta c_4$, $\eta^2c_4$ and their $c_4$-multiples, and $\nu$ and $\bar\kappa$ have order $4$. The element $\eta'$ is detected by $x_{1,0}$, and $\eta'\nu$ by $h_0x_{4,0}=h_2x_{1,0}$.

We also enumerate the hidden $2$, $\eta$, and $\nu$ extensions in this range.

\begin{proposition}\label{prop:hidden}
In stems at most $32$ there are the following hidden extensions.
\begin{enumerate}
\item A hidden $\nu$-extension from $h_2^2$ in $(6,2)$ to $h_1c_0$ in $(9,4)$, realizing $\nu^3=\eta\epsilon$.
\item Hidden $\eta$-extensions from $h_1g$ in $(21,5)$ to $w_1d_0$ in $(22,8)$, realizing $\eta^2\bar\kappa$, and from $h_1c_0x_{16,0}$ in $(25,4)$ to $w_1x_{18,2}$ in $(26,6)$.
\item Hidden extensions by $2$ from $y$, $w_1y$ and $w_1^2y$ to $w_1x_{8,3}$, $w_1^2x_{8,3}$ and $w_1^3x_{8,3}$, in stems $16$, $24$ and $32$.
\end{enumerate}
Moreover, each of the classes
\[w_1x_{4,0},\quad w_1^2x_{4,0},\quad w_1^3x_{4,0},\quad c_0x_{16,0},\quad\beta x_{9,1},\quad x_{32,5},\quad\delta'\]
detects an element of infinite order and supports a hidden extension by $2$ into the $h_0$-towers of its stem.
\end{proposition}

\begin{sseqdata}[name = AdamsEinfA, adams, x range = {0}{16}, y range = {0}{12}]
\class[tmfclass](0,0)
\class[tmfclass](0,1)
\class[tmfclass](0,2)
\class[tmfclass](0,3)
\class[tmfclass](0,4)
\class[tmfclass](0,5)
\class[tmfclass](0,6)
\class[tmfclass](0,7)
\class[tmfclass](0,8)
\class[tmfclass](0,9)
\class[tmfclass](0,10)
\class[tmfclass](0,11)
\class[tmfclass](0,12)
\class[tmfclass](0,13)
\class[tmfclass](0,14)
\class[tmfclass](0,15)
\class[tmfclass](0,16)
\class[tmfclass](0,17)
\class[tmfclass](0,18)
\class[tmfclass](0,19)
\class[tmfclass](0,20)
\class[tmfclass](0,21)
\class[tmfclass](0,22)
\class[tmfclass](0,23)
\class[tmfclass](0,24)
\class(1,0)
\class[tmfclass](1,1)
\class(2,1)
\class[tmfclass](2,2)
\class[tmfclass](3,1)
\class[tmfclass](3,2)
\class(4,1)
\class(4,2)
\class(4,3)
\class(4,4)
\class(4,5)
\class(4,6)
\class(4,7)
\class(4,8)
\class(4,9)
\class(4,10)
\class(4,11)
\class(4,12)
\class(4,13)
\class(4,14)
\class(4,15)
\class(4,16)
\class(4,17)
\class(4,18)
\class(4,19)
\class(4,20)
\class(4,21)
\class(4,22)
\class(4,23)
\class(4,24)
\class[tmfclass](6,2)
\class(8,1)
\class(8,2)
\class[tmfclass](8,3)
\class(8,3)
\class(8,3)
\class[tmfclass](8,4)
\class(8,4)
\class(8,4)
\class[tmfclass](8,5)
\class(8,5)
\class(8,5)
\class[tmfclass](8,6)
\class(8,6)
\class(8,6)
\class[tmfclass](8,7)
\class(8,7)
\class(8,7)
\class[tmfclass](8,8)
\class(8,8)
\class(8,8)
\class[tmfclass](8,9)
\class(8,9)
\class(8,9)
\class[tmfclass](8,10)
\class(8,10)
\class(8,10)
\class[tmfclass](8,11)
\class(8,11)
\class(8,11)
\class[tmfclass](8,12)
\class(8,12)
\class(8,12)
\class[tmfclass](8,13)
\class(8,13)
\class(8,13)
\class[tmfclass](8,14)
\class(8,14)
\class(8,14)
\class[tmfclass](8,15)
\class(8,15)
\class(8,15)
\class[tmfclass](8,16)
\class(8,16)
\class(8,16)
\class[tmfclass](8,17)
\class(8,17)
\class(8,17)
\class[tmfclass](8,18)
\class(8,18)
\class(8,18)
\class[tmfclass](8,19)
\class(8,19)
\class(8,19)
\class[tmfclass](8,20)
\class(8,20)
\class(8,20)
\class[tmfclass](8,21)
\class(8,21)
\class(8,21)
\class[tmfclass](8,22)
\class(8,22)
\class(8,22)
\class[tmfclass](8,23)
\class(8,23)
\class(8,23)
\class[tmfclass](8,24)
\class(8,24)
\class(8,24)
\class(9,1)
\class(9,2)
\class[tmfclass](9,4)
\class[tmfclass](9,5)
\class(10,2)
\class(10,3)
\class[tmfclass](10,6)
\class(12,3)
\class(12,4)
\class(12,4)
\class(12,5)
\class[tmfclass](12,6)
\class(12,6)
\class(12,6)
\class[tmfclass](12,7)
\class(12,7)
\class(12,7)
\class[tmfclass](12,8)
\class(12,8)
\class(12,8)
\class[tmfclass](12,9)
\class(12,9)
\class(12,9)
\class[tmfclass](12,10)
\class(12,10)
\class(12,10)
\class[tmfclass](12,11)
\class(12,11)
\class(12,11)
\class[tmfclass](12,12)
\class(12,12)
\class(12,12)
\class[tmfclass](12,13)
\class(12,13)
\class(12,13)
\class[tmfclass](12,14)
\class(12,14)
\class(12,14)
\class[tmfclass](12,15)
\class(12,15)
\class(12,15)
\class[tmfclass](12,16)
\class(12,16)
\class(12,16)
\class[tmfclass](12,17)
\class(12,17)
\class(12,17)
\class[tmfclass](12,18)
\class(12,18)
\class(12,18)
\class[tmfclass](12,19)
\class(12,19)
\class(12,19)
\class[tmfclass](12,20)
\class(12,20)
\class(12,20)
\class[tmfclass](12,21)
\class(12,21)
\class(12,21)
\class[tmfclass](12,22)
\class(12,22)
\class(12,22)
\class[tmfclass](12,23)
\class(12,23)
\class(12,23)
\class[tmfclass](12,24)
\class(12,24)
\class(12,24)
\class[tmfclass](14,4)
\class(16,2)
\class(16,3)
\class(16,3)
\class(16,4)
\class(16,4)
\class(16,4)
\class(16,5)
\class(16,5)
\class(16,5)
\class(16,5)
\class(16,6)
\class(16,6)
\class(16,6)
\class(16,7)
\class(16,7)
\class(16,7)
\class(16,7)
\class[tmfclass](16,8)
\class(16,8)
\class(16,8)
\class(16,8)
\class(16,8)
\class[tmfclass](16,9)
\class(16,9)
\class(16,9)
\class(16,9)
\class(16,9)
\class[tmfclass](16,10)
\class(16,10)
\class(16,10)
\class(16,10)
\class(16,10)
\class[tmfclass](16,11)
\class(16,11)
\class(16,11)
\class(16,11)
\class(16,11)
\class[tmfclass](16,12)
\class(16,12)
\class(16,12)
\class(16,12)
\class(16,12)
\class[tmfclass](16,13)
\class(16,13)
\class(16,13)
\class(16,13)
\class(16,13)
\class[tmfclass](16,14)
\class(16,14)
\class(16,14)
\class(16,14)
\class(16,14)
\class[tmfclass](16,15)
\class(16,15)
\class(16,15)
\class(16,15)
\class(16,15)
\class[tmfclass](16,16)
\class(16,16)
\class(16,16)
\class(16,16)
\class(16,16)
\class[tmfclass](16,17)
\class(16,17)
\class(16,17)
\class(16,17)
\class(16,17)
\class[tmfclass](16,18)
\class(16,18)
\class(16,18)
\class(16,18)
\class(16,18)
\class[tmfclass](16,19)
\class(16,19)
\class(16,19)
\class(16,19)
\class(16,19)
\class[tmfclass](16,20)
\class(16,20)
\class(16,20)
\class(16,20)
\class(16,20)
\class[tmfclass](16,21)
\class(16,21)
\class(16,21)
\class(16,21)
\class(16,21)
\class[tmfclass](16,22)
\class(16,22)
\class(16,22)
\class(16,22)
\class(16,22)
\class[tmfclass](16,23)
\class(16,23)
\class(16,23)
\class(16,23)
\class(16,23)
\class[tmfclass](16,24)
\class(16,24)
\class(16,24)
\class(16,24)
\class(16,24)
\structline(0,0,1)(0,1,1)
\structline(0,0,1)(1,1,1)
\structline(0,0,1)(3,1,1)
\structline(0,1,1)(0,2,1)
\structline(0,1,1)(3,2,1)
\structline(0,2,1)(0,3,1)
\structline(0,3,1)(0,4,1)
\structline(0,4,1)(0,5,1)
\structline(0,5,1)(0,6,1)
\structline(0,6,1)(0,7,1)
\structline(0,7,1)(0,8,1)
\structline(0,8,1)(0,9,1)
\structline(0,9,1)(0,10,1)
\structline(0,10,1)(0,11,1)
\structline(0,11,1)(0,12,1)
\structline(0,12,1)(0,13,1)
\structline(0,13,1)(0,14,1)
\structline(0,14,1)(0,15,1)
\structline(0,15,1)(0,16,1)
\structline(0,16,1)(0,17,1)
\structline(0,17,1)(0,18,1)
\structline(0,18,1)(0,19,1)
\structline(0,19,1)(0,20,1)
\structline(0,20,1)(0,21,1)
\structline(0,21,1)(0,22,1)
\structline(0,22,1)(0,23,1)
\structline(0,23,1)(0,24,1)
\structline(1,0,1)(2,1,1)
\structline(1,0,1)(4,1,1)
\structline(1,1,1)(2,2,1)
\structline(3,1,1)(3,2,1)
\structline(3,1,1)(6,2,1)
\structline(4,2,1)(4,3,1)
\structline(4,3,1)(4,4,1)
\structline(4,4,1)(4,5,1)
\structline(4,5,1)(4,6,1)
\structline(4,6,1)(4,7,1)
\structline(4,7,1)(4,8,1)
\structline(4,8,1)(4,9,1)
\structline(4,9,1)(4,10,1)
\structline(4,10,1)(4,11,1)
\structline(4,11,1)(4,12,1)
\structline(4,12,1)(4,13,1)
\structline(4,13,1)(4,14,1)
\structline(4,14,1)(4,15,1)
\structline(4,15,1)(4,16,1)
\structline(4,16,1)(4,17,1)
\structline(4,17,1)(4,18,1)
\structline(4,18,1)(4,19,1)
\structline(4,19,1)(4,20,1)
\structline(4,20,1)(4,21,1)
\structline(4,21,1)(4,22,1)
\structline(4,22,1)(4,23,1)
\structline(4,23,1)(4,24,1)
\structline(8,1,1)(8,2,1)
\structline(8,1,1)(9,2,1)
\structline(8,2,1)(8,3,3)
\structline(8,3,2)(8,4,2)
\structline(8,3,3)(8,4,3)
\structline(8,3,1)(9,4,1)
\structline(8,4,1)(8,5,1)
\structline(8,4,2)(8,5,2)
\structline(8,4,3)(8,5,3)
\structline(8,4,1)(9,5,1)
\structline(8,5,1)(8,6,1)
\structline(8,5,2)(8,6,2)
\structline(8,5,3)(8,6,3)
\structline(8,6,1)(8,7,1)
\structline(8,6,2)(8,7,2)
\structline(8,6,3)(8,7,3)
\structline(8,7,1)(8,8,1)
\structline(8,7,2)(8,8,2)
\structline(8,7,3)(8,8,3)
\structline(8,8,1)(8,9,1)
\structline(8,8,2)(8,9,2)
\structline(8,8,3)(8,9,3)
\structline(8,9,1)(8,10,1)
\structline(8,9,2)(8,10,2)
\structline(8,9,3)(8,10,3)
\structline(8,10,1)(8,11,1)
\structline(8,10,2)(8,11,2)
\structline(8,10,3)(8,11,3)
\structline(8,11,1)(8,12,1)
\structline(8,11,2)(8,12,2)
\structline(8,11,3)(8,12,3)
\structline(8,12,1)(8,13,1)
\structline(8,12,2)(8,13,2)
\structline(8,12,3)(8,13,3)
\structline(8,13,1)(8,14,1)
\structline(8,13,2)(8,14,2)
\structline(8,13,3)(8,14,3)
\structline(8,14,1)(8,15,1)
\structline(8,14,2)(8,15,2)
\structline(8,14,3)(8,15,3)
\structline(8,15,1)(8,16,1)
\structline(8,15,2)(8,16,2)
\structline(8,15,3)(8,16,3)
\structline(8,16,1)(8,17,1)
\structline(8,16,2)(8,17,2)
\structline(8,16,3)(8,17,3)
\structline(8,17,1)(8,18,1)
\structline(8,17,2)(8,18,2)
\structline(8,17,3)(8,18,3)
\structline(8,18,1)(8,19,1)
\structline(8,18,2)(8,19,2)
\structline(8,18,3)(8,19,3)
\structline(8,19,1)(8,20,1)
\structline(8,19,2)(8,20,2)
\structline(8,19,3)(8,20,3)
\structline(8,20,1)(8,21,1)
\structline(8,20,2)(8,21,2)
\structline(8,20,3)(8,21,3)
\structline(8,21,1)(8,22,1)
\structline(8,21,2)(8,22,2)
\structline(8,21,3)(8,22,3)
\structline(8,22,1)(8,23,1)
\structline(8,22,2)(8,23,2)
\structline(8,22,3)(8,23,3)
\structline(8,23,1)(8,24,1)
\structline(8,23,2)(8,24,2)
\structline(8,23,3)(8,24,3)
\structline(9,1,1)(10,2,1)
\structline(9,2,1)(10,3,1)
\structline(9,5,1)(10,6,1)
\structline(12,3,1)(12,4,1)
\structline(12,4,1)(12,5,1)
\structline(12,5,1)(12,6,2)
\structline(12,6,1)(12,7,1)
\structline(12,6,2)(12,7,2)
\structline(12,6,3)(12,7,3)
\structline(12,7,1)(12,8,1)
\structline(12,7,2)(12,8,2)
\structline(12,7,3)(12,8,3)
\structline(12,8,1)(12,9,1)
\structline(12,8,2)(12,9,2)
\structline(12,8,3)(12,9,3)
\structline(12,9,1)(12,10,1)
\structline(12,9,2)(12,10,2)
\structline(12,9,3)(12,10,3)
\structline(12,10,1)(12,11,1)
\structline(12,10,2)(12,11,2)
\structline(12,10,3)(12,11,3)
\structline(12,11,1)(12,12,1)
\structline(12,11,2)(12,12,2)
\structline(12,11,3)(12,12,3)
\structline(12,12,1)(12,13,1)
\structline(12,12,2)(12,13,2)
\structline(12,12,3)(12,13,3)
\structline(12,13,1)(12,14,1)
\structline(12,13,2)(12,14,2)
\structline(12,13,3)(12,14,3)
\structline(12,14,1)(12,15,1)
\structline(12,14,2)(12,15,2)
\structline(12,14,3)(12,15,3)
\structline(12,15,1)(12,16,1)
\structline(12,15,2)(12,16,2)
\structline(12,15,3)(12,16,3)
\structline(12,16,1)(12,17,1)
\structline(12,16,2)(12,17,2)
\structline(12,16,3)(12,17,3)
\structline(12,17,1)(12,18,1)
\structline(12,17,2)(12,18,2)
\structline(12,17,3)(12,18,3)
\structline(12,18,1)(12,19,1)
\structline(12,18,2)(12,19,2)
\structline(12,18,3)(12,19,3)
\structline(12,19,1)(12,20,1)
\structline(12,19,2)(12,20,2)
\structline(12,19,3)(12,20,3)
\structline(12,20,1)(12,21,1)
\structline(12,20,2)(12,21,2)
\structline(12,20,3)(12,21,3)
\structline(12,21,1)(12,22,1)
\structline(12,21,2)(12,22,2)
\structline(12,21,3)(12,22,3)
\structline(12,22,1)(12,23,1)
\structline(12,22,2)(12,23,2)
\structline(12,22,3)(12,23,3)
\structline(12,23,1)(12,24,1)
\structline(12,23,2)(12,24,2)
\structline(12,23,3)(12,24,3)
\structline(16,2,1)(16,3,2)
\structline(16,3,1)(16,4,1)
\structline(16,3,2)(16,4,3)
\structline(16,4,1)(16,5,1)
\structline(16,4,2)(16,5,2)
\structline(16,4,3)(16,5,4)
\structline(16,5,1)(16,6,1)
\structline(16,5,2)(16,6,2)
\structline(16,5,4)(16,6,3)
\structline(16,6,1)(16,7,1)
\structline(16,6,2)(16,7,2)
\structline(16,6,3)(16,7,4)
\structline(16,7,1)(16,8,2)
\structline(16,7,2)(16,8,3)
\structline(16,7,3)(16,8,4)
\structline(16,7,4)(16,8,5)
\structline(16,8,1)(16,9,1)
\structline(16,8,2)(16,9,2)
\structline(16,8,3)(16,9,3)
\structline(16,8,4)(16,9,4)
\structline(16,8,5)(16,9,5)
\structline(16,9,1)(16,10,1)
\structline(16,9,2)(16,10,2)
\structline(16,9,3)(16,10,3)
\structline(16,9,4)(16,10,4)
\structline(16,9,5)(16,10,5)
\structline(16,10,1)(16,11,1)
\structline(16,10,2)(16,11,2)
\structline(16,10,3)(16,11,3)
\structline(16,10,4)(16,11,4)
\structline(16,10,5)(16,11,5)
\structline(16,11,1)(16,12,1)
\structline(16,11,2)(16,12,2)
\structline(16,11,3)(16,12,3)
\structline(16,11,4)(16,12,4)
\structline(16,11,5)(16,12,5)
\structline(16,12,1)(16,13,1)
\structline(16,12,2)(16,13,2)
\structline(16,12,3)(16,13,3)
\structline(16,12,4)(16,13,4)
\structline(16,12,5)(16,13,5)
\structline(16,13,1)(16,14,1)
\structline(16,13,2)(16,14,2)
\structline(16,13,3)(16,14,3)
\structline(16,13,4)(16,14,4)
\structline(16,13,5)(16,14,5)
\structline(16,14,1)(16,15,1)
\structline(16,14,2)(16,15,2)
\structline(16,14,3)(16,15,3)
\structline(16,14,4)(16,15,4)
\structline(16,14,5)(16,15,5)
\structline(16,15,1)(16,16,1)
\structline(16,15,2)(16,16,2)
\structline(16,15,3)(16,16,3)
\structline(16,15,4)(16,16,4)
\structline(16,15,5)(16,16,5)
\structline(16,16,1)(16,17,1)
\structline(16,16,2)(16,17,2)
\structline(16,16,3)(16,17,3)
\structline(16,16,4)(16,17,4)
\structline(16,16,5)(16,17,5)
\structline(16,17,1)(16,18,1)
\structline(16,17,2)(16,18,2)
\structline(16,17,3)(16,18,3)
\structline(16,17,4)(16,18,4)
\structline(16,17,5)(16,18,5)
\structline(16,18,1)(16,19,1)
\structline(16,18,2)(16,19,2)
\structline(16,18,3)(16,19,3)
\structline(16,18,4)(16,19,4)
\structline(16,18,5)(16,19,5)
\structline(16,19,1)(16,20,1)
\structline(16,19,2)(16,20,2)
\structline(16,19,3)(16,20,3)
\structline(16,19,4)(16,20,4)
\structline(16,19,5)(16,20,5)
\structline(16,20,1)(16,21,1)
\structline(16,20,2)(16,21,2)
\structline(16,20,3)(16,21,3)
\structline(16,20,4)(16,21,4)
\structline(16,20,5)(16,21,5)
\structline(16,21,1)(16,22,1)
\structline(16,21,2)(16,22,2)
\structline(16,21,3)(16,22,3)
\structline(16,21,4)(16,22,4)
\structline(16,21,5)(16,22,5)
\structline(16,22,1)(16,23,1)
\structline(16,22,2)(16,23,2)
\structline(16,22,3)(16,23,3)
\structline(16,22,4)(16,23,4)
\structline(16,22,5)(16,23,5)
\structline(16,23,1)(16,24,1)
\structline(16,23,2)(16,24,2)
\structline(16,23,3)(16,24,3)
\structline(16,23,4)(16,24,4)
\structline(16,23,5)(16,24,5)
\structline[hiddennu](6,2,1)(9,4,1)
\draw[hiddencol, dashed, line width = 0.45pt, -{Latex[length = 2.2pt]}, shorten >= 1.2pt, shorten <= 1.2pt, page constraint = {\page >= 5}] (16,5,3) to[out = 0, in = 0, looseness = 1.7] (16,7,3);
\end{sseqdata}

\begin{sseqdata}[name = AdamsEinfB, adams, x range = {16}{32}, y range = {0}{15}]
\class(16,2)
\class(16,3)
\class(16,3)
\class(16,4)
\class(16,4)
\class(16,4)
\class(16,5)
\class(16,5)
\class(16,5)
\class(16,5)
\class(16,6)
\class(16,6)
\class(16,6)
\class(16,7)
\class(16,7)
\class(16,7)
\class(16,7)
\class[tmfclass](16,8)
\class(16,8)
\class(16,8)
\class(16,8)
\class(16,8)
\class[tmfclass](16,9)
\class(16,9)
\class(16,9)
\class(16,9)
\class(16,9)
\class[tmfclass](16,10)
\class(16,10)
\class(16,10)
\class(16,10)
\class(16,10)
\class[tmfclass](16,11)
\class(16,11)
\class(16,11)
\class(16,11)
\class(16,11)
\class[tmfclass](16,12)
\class(16,12)
\class(16,12)
\class(16,12)
\class(16,12)
\class[tmfclass](16,13)
\class(16,13)
\class(16,13)
\class(16,13)
\class(16,13)
\class[tmfclass](16,14)
\class(16,14)
\class(16,14)
\class(16,14)
\class(16,14)
\class[tmfclass](16,15)
\class(16,15)
\class(16,15)
\class(16,15)
\class(16,15)
\class[tmfclass](16,16)
\class(16,16)
\class(16,16)
\class(16,16)
\class(16,16)
\class[tmfclass](16,17)
\class(16,17)
\class(16,17)
\class(16,17)
\class(16,17)
\class[tmfclass](16,18)
\class(16,18)
\class(16,18)
\class(16,18)
\class(16,18)
\class[tmfclass](16,19)
\class(16,19)
\class(16,19)
\class(16,19)
\class(16,19)
\class[tmfclass](16,20)
\class(16,20)
\class(16,20)
\class(16,20)
\class(16,20)
\class[tmfclass](16,21)
\class(16,21)
\class(16,21)
\class(16,21)
\class(16,21)
\class[tmfclass](16,22)
\class(16,22)
\class(16,22)
\class(16,22)
\class(16,22)
\class[tmfclass](16,23)
\class(16,23)
\class(16,23)
\class(16,23)
\class(16,23)
\class[tmfclass](16,24)
\class(16,24)
\class(16,24)
\class(16,24)
\class(16,24)
\class(17,2)
\class(17,4)
\class(17,5)
\class(17,6)
\class[tmfclass](17,9)
\class(18,3)
\class(18,5)
\class(18,6)
\class(18,7)
\class[tmfclass](18,10)
\class(19,2)
\class(19,3)
\class[tmfclass](20,4)
\class(20,4)
\class[tmfclass](20,5)
\class(20,5)
\class(20,6)
\class(20,6)
\class(20,7)
\class(20,7)
\class(20,7)
\class(20,8)
\class(20,8)
\class(20,8)
\class(20,8)
\class(20,9)
\class(20,9)
\class(20,9)
\class[tmfclass](20,10)
\class(20,10)
\class(20,10)
\class(20,10)
\class(20,10)
\class[tmfclass](20,11)
\class(20,11)
\class(20,11)
\class(20,11)
\class(20,11)
\class[tmfclass](20,12)
\class(20,12)
\class(20,12)
\class(20,12)
\class(20,12)
\class[tmfclass](20,13)
\class(20,13)
\class(20,13)
\class(20,13)
\class(20,13)
\class[tmfclass](20,14)
\class(20,14)
\class(20,14)
\class(20,14)
\class(20,14)
\class[tmfclass](20,15)
\class(20,15)
\class(20,15)
\class(20,15)
\class(20,15)
\class[tmfclass](20,16)
\class(20,16)
\class(20,16)
\class(20,16)
\class(20,16)
\class[tmfclass](20,17)
\class(20,17)
\class(20,17)
\class(20,17)
\class(20,17)
\class[tmfclass](20,18)
\class(20,18)
\class(20,18)
\class(20,18)
\class(20,18)
\class[tmfclass](20,19)
\class(20,19)
\class(20,19)
\class(20,19)
\class(20,19)
\class[tmfclass](20,20)
\class(20,20)
\class(20,20)
\class(20,20)
\class(20,20)
\class[tmfclass](20,21)
\class(20,21)
\class(20,21)
\class(20,21)
\class(20,21)
\class[tmfclass](20,22)
\class(20,22)
\class(20,22)
\class(20,22)
\class(20,22)
\class[tmfclass](20,23)
\class(20,23)
\class(20,23)
\class(20,23)
\class(20,23)
\class[tmfclass](20,24)
\class(20,24)
\class(20,24)
\class(20,24)
\class(20,24)
\class[tmfclass](21,5)
\class(22,4)
\class[tmfclass](22,8)
\class(24,3)
\class(24,4)
\class(24,5)
\class(24,5)
\class(24,6)
\class(24,6)
\class[tmfclass](24,7)
\class(24,7)
\class(24,7)
\class(24,7)
\class[tmfclass](24,8)
\class(24,8)
\class(24,8)
\class(24,8)
\class(24,8)
\class[tmfclass](24,9)
\class(24,9)
\class(24,9)
\class(24,9)
\class(24,9)
\class(24,9)
\class[tmfclass](24,10)
\class(24,10)
\class(24,10)
\class(24,10)
\class(24,10)
\class[tmfclass](24,11)
\class(24,11)
\class(24,11)
\class(24,11)
\class(24,11)
\class(24,11)
\class[tmfclass](24,12)
\class[tmfclass](24,12)
\class(24,12)
\class(24,12)
\class(24,12)
\class(24,12)
\class(24,12)
\class[tmfclass](24,13)
\class[tmfclass](24,13)
\class(24,13)
\class(24,13)
\class(24,13)
\class(24,13)
\class(24,13)
\class[tmfclass](24,14)
\class[tmfclass](24,14)
\class(24,14)
\class(24,14)
\class(24,14)
\class(24,14)
\class(24,14)
\class[tmfclass](24,15)
\class[tmfclass](24,15)
\class(24,15)
\class(24,15)
\class(24,15)
\class(24,15)
\class(24,15)
\class[tmfclass](24,16)
\class[tmfclass](24,16)
\class(24,16)
\class(24,16)
\class(24,16)
\class(24,16)
\class(24,16)
\class[tmfclass](24,17)
\class[tmfclass](24,17)
\class(24,17)
\class(24,17)
\class(24,17)
\class(24,17)
\class(24,17)
\class[tmfclass](24,18)
\class[tmfclass](24,18)
\class(24,18)
\class(24,18)
\class(24,18)
\class(24,18)
\class(24,18)
\class[tmfclass](24,19)
\class[tmfclass](24,19)
\class(24,19)
\class(24,19)
\class(24,19)
\class(24,19)
\class(24,19)
\class[tmfclass](24,20)
\class[tmfclass](24,20)
\class(24,20)
\class(24,20)
\class(24,20)
\class(24,20)
\class(24,20)
\class[tmfclass](24,21)
\class[tmfclass](24,21)
\class(24,21)
\class(24,21)
\class(24,21)
\class(24,21)
\class(24,21)
\class[tmfclass](24,22)
\class[tmfclass](24,22)
\class(24,22)
\class(24,22)
\class(24,22)
\class(24,22)
\class(24,22)
\class[tmfclass](24,23)
\class[tmfclass](24,23)
\class(24,23)
\class(24,23)
\class(24,23)
\class(24,23)
\class(24,23)
\class[tmfclass](24,24)
\class[tmfclass](24,24)
\class(24,24)
\class(24,24)
\class(24,24)
\class(24,24)
\class(24,24)
\class(25,4)
\class[tmfclass](25,5)
\class(25,6)
\class(25,8)
\class(25,9)
\class(25,10)
\class[tmfclass](25,13)
\class[tmfclass](26,6)
\class(26,6)
\class(26,7)
\class(26,9)
\class(26,10)
\class(26,11)
\class[tmfclass](26,14)
\class[tmfclass](27,6)
\class(28,6)
\class(28,6)
\class(28,7)
\class(28,7)
\class[tmfclass](28,8)
\class(28,8)
\class(28,8)
\class(28,8)
\class(28,9)
\class(28,9)
\class(28,9)
\class(28,10)
\class(28,10)
\class(28,10)
\class(28,10)
\class(28,11)
\class(28,11)
\class(28,11)
\class(28,11)
\class(28,11)
\class(28,12)
\class(28,12)
\class(28,12)
\class(28,12)
\class(28,12)
\class(28,12)
\class(28,13)
\class(28,13)
\class(28,13)
\class(28,13)
\class(28,13)
\class[tmfclass](28,14)
\class(28,14)
\class(28,14)
\class(28,14)
\class(28,14)
\class(28,14)
\class(28,14)
\class[tmfclass](28,15)
\class(28,15)
\class(28,15)
\class(28,15)
\class(28,15)
\class(28,15)
\class(28,15)
\class[tmfclass](28,16)
\class(28,16)
\class(28,16)
\class(28,16)
\class(28,16)
\class(28,16)
\class(28,16)
\class[tmfclass](28,17)
\class(28,17)
\class(28,17)
\class(28,17)
\class(28,17)
\class(28,17)
\class(28,17)
\class[tmfclass](28,18)
\class(28,18)
\class(28,18)
\class(28,18)
\class(28,18)
\class(28,18)
\class(28,18)
\class[tmfclass](28,19)
\class(28,19)
\class(28,19)
\class(28,19)
\class(28,19)
\class(28,19)
\class(28,19)
\class[tmfclass](28,20)
\class(28,20)
\class(28,20)
\class(28,20)
\class(28,20)
\class(28,20)
\class(28,20)
\class[tmfclass](28,21)
\class(28,21)
\class(28,21)
\class(28,21)
\class(28,21)
\class(28,21)
\class(28,21)
\class[tmfclass](28,22)
\class(28,22)
\class(28,22)
\class(28,22)
\class(28,22)
\class(28,22)
\class(28,22)
\class[tmfclass](28,23)
\class(28,23)
\class(28,23)
\class(28,23)
\class(28,23)
\class(28,23)
\class(28,23)
\class[tmfclass](28,24)
\class(28,24)
\class(28,24)
\class(28,24)
\class(28,24)
\class(28,24)
\class(28,24)
\class(32,4)
\class(32,5)
\class(32,5)
\class(32,6)
\class[tmfclass](32,7)
\class(32,7)
\class(32,7)
\class(32,7)
\class[tmfclass](32,8)
\class(32,8)
\class(32,8)
\class[tmfclass](32,9)
\class(32,9)
\class(32,9)
\class(32,9)
\class(32,9)
\class(32,10)
\class(32,10)
\class(32,10)
\class(32,10)
\class(32,11)
\class(32,11)
\class[tmfclass](32,11)
\class(32,11)
\class(32,11)
\class(32,11)
\class[tmfclass](32,12)
\class(32,12)
\class(32,12)
\class(32,12)
\class(32,12)
\class(32,12)
\class(32,12)
\class[tmfclass](32,13)
\class(32,13)
\class(32,13)
\class(32,13)
\class(32,13)
\class(32,13)
\class(32,13)
\class(32,13)
\class[tmfclass](32,14)
\class(32,14)
\class(32,14)
\class(32,14)
\class(32,14)
\class(32,14)
\class(32,14)
\class[tmfclass](32,15)
\class(32,15)
\class(32,15)
\class(32,15)
\class(32,15)
\class(32,15)
\class(32,15)
\class(32,15)
\class[tmfclass](32,16)
\class[tmfclass](32,16)
\class(32,16)
\class(32,16)
\class(32,16)
\class(32,16)
\class(32,16)
\class(32,16)
\class(32,16)
\class[tmfclass](32,17)
\class[tmfclass](32,17)
\class(32,17)
\class(32,17)
\class(32,17)
\class(32,17)
\class(32,17)
\class(32,17)
\class(32,17)
\class[tmfclass](32,18)
\class[tmfclass](32,18)
\class(32,18)
\class(32,18)
\class(32,18)
\class(32,18)
\class(32,18)
\class(32,18)
\class(32,18)
\class[tmfclass](32,19)
\class[tmfclass](32,19)
\class(32,19)
\class(32,19)
\class(32,19)
\class(32,19)
\class(32,19)
\class(32,19)
\class(32,19)
\class[tmfclass](32,20)
\class[tmfclass](32,20)
\class(32,20)
\class(32,20)
\class(32,20)
\class(32,20)
\class(32,20)
\class(32,20)
\class(32,20)
\class[tmfclass](32,21)
\class[tmfclass](32,21)
\class(32,21)
\class(32,21)
\class(32,21)
\class(32,21)
\class(32,21)
\class(32,21)
\class(32,21)
\class[tmfclass](32,22)
\class[tmfclass](32,22)
\class(32,22)
\class(32,22)
\class(32,22)
\class(32,22)
\class(32,22)
\class(32,22)
\class(32,22)
\class[tmfclass](32,23)
\class[tmfclass](32,23)
\class(32,23)
\class(32,23)
\class(32,23)
\class(32,23)
\class(32,23)
\class(32,23)
\class(32,23)
\class[tmfclass](32,24)
\class[tmfclass](32,24)
\class(32,24)
\class(32,24)
\class(32,24)
\class(32,24)
\class(32,24)
\class(32,24)
\class(32,24)
\structline(16,2,1)(16,3,2)
\structline(16,2,1)(19,3,1)
\structline(16,3,1)(16,4,1)
\structline(16,3,2)(16,4,3)
\structline(16,3,1)(17,4,1)
\structline(16,4,1)(16,5,1)
\structline(16,4,2)(16,5,2)
\structline(16,4,3)(16,5,4)
\structline(16,4,2)(17,5,1)
\structline(16,5,1)(16,6,1)
\structline(16,5,2)(16,6,2)
\structline(16,5,4)(16,6,3)
\structline(16,5,3)(17,6,1)
\structline(16,6,1)(16,7,1)
\structline(16,6,2)(16,7,2)
\structline(16,6,3)(16,7,4)
\structline(16,7,1)(16,8,2)
\structline(16,7,2)(16,8,3)
\structline(16,7,3)(16,8,4)
\structline(16,7,4)(16,8,5)
\structline(16,8,1)(16,9,1)
\structline(16,8,2)(16,9,2)
\structline(16,8,3)(16,9,3)
\structline(16,8,4)(16,9,4)
\structline(16,8,5)(16,9,5)
\structline(16,8,1)(17,9,1)
\structline(16,9,1)(16,10,1)
\structline(16,9,2)(16,10,2)
\structline(16,9,3)(16,10,3)
\structline(16,9,4)(16,10,4)
\structline(16,9,5)(16,10,5)
\structline(16,10,1)(16,11,1)
\structline(16,10,2)(16,11,2)
\structline(16,10,3)(16,11,3)
\structline(16,10,4)(16,11,4)
\structline(16,10,5)(16,11,5)
\structline(16,11,1)(16,12,1)
\structline(16,11,2)(16,12,2)
\structline(16,11,3)(16,12,3)
\structline(16,11,4)(16,12,4)
\structline(16,11,5)(16,12,5)
\structline(16,12,1)(16,13,1)
\structline(16,12,2)(16,13,2)
\structline(16,12,3)(16,13,3)
\structline(16,12,4)(16,13,4)
\structline(16,12,5)(16,13,5)
\structline(16,13,1)(16,14,1)
\structline(16,13,2)(16,14,2)
\structline(16,13,3)(16,14,3)
\structline(16,13,4)(16,14,4)
\structline(16,13,5)(16,14,5)
\structline(16,14,1)(16,15,1)
\structline(16,14,2)(16,15,2)
\structline(16,14,3)(16,15,3)
\structline(16,14,4)(16,15,4)
\structline(16,14,5)(16,15,5)
\structline(16,15,1)(16,16,1)
\structline(16,15,2)(16,16,2)
\structline(16,15,3)(16,16,3)
\structline(16,15,4)(16,16,4)
\structline(16,15,5)(16,16,5)
\structline(16,16,1)(16,17,1)
\structline(16,16,2)(16,17,2)
\structline(16,16,3)(16,17,3)
\structline(16,16,4)(16,17,4)
\structline(16,16,5)(16,17,5)
\structline(16,17,1)(16,18,1)
\structline(16,17,2)(16,18,2)
\structline(16,17,3)(16,18,3)
\structline(16,17,4)(16,18,4)
\structline(16,17,5)(16,18,5)
\structline(16,18,1)(16,19,1)
\structline(16,18,2)(16,19,2)
\structline(16,18,3)(16,19,3)
\structline(16,18,4)(16,19,4)
\structline(16,18,5)(16,19,5)
\structline(16,19,1)(16,20,1)
\structline(16,19,2)(16,20,2)
\structline(16,19,3)(16,20,3)
\structline(16,19,4)(16,20,4)
\structline(16,19,5)(16,20,5)
\structline(16,20,1)(16,21,1)
\structline(16,20,2)(16,21,2)
\structline(16,20,3)(16,21,3)
\structline(16,20,4)(16,21,4)
\structline(16,20,5)(16,21,5)
\structline(16,21,1)(16,22,1)
\structline(16,21,2)(16,22,2)
\structline(16,21,3)(16,22,3)
\structline(16,21,4)(16,22,4)
\structline(16,21,5)(16,22,5)
\structline(16,22,1)(16,23,1)
\structline(16,22,2)(16,23,2)
\structline(16,22,3)(16,23,3)
\structline(16,22,4)(16,23,4)
\structline(16,22,5)(16,23,5)
\structline(16,23,1)(16,24,1)
\structline(16,23,2)(16,24,2)
\structline(16,23,3)(16,24,3)
\structline(16,23,4)(16,24,4)
\structline(16,23,5)(16,24,5)
\structline(17,2,1)(18,3,1)
\structline(17,4,1)(18,5,1)
\structline(17,5,1)(18,6,1)
\structline(17,6,1)(18,7,1)
\structline(17,9,1)(18,10,1)
\structline(19,2,1)(19,3,1)
\structline(20,4,1)(20,5,1)
\structline(20,4,2)(20,5,2)
\structline(20,4,1)(21,5,1)
\structline(20,5,2)(20,6,1)
\structline(20,5,2)(20,6,2)
\structline(20,6,1)(20,7,1)
\structline(20,6,2)(20,7,2)
\structline(20,7,1)(20,8,1)
\structline(20,7,2)(20,8,2)
\structline(20,7,3)(20,8,3)
\structline(20,8,1)(20,9,1)
\structline(20,8,2)(20,9,2)
\structline(20,8,3)(20,9,3)
\structline(20,9,1)(20,10,2)
\structline(20,9,2)(20,10,3)
\structline(20,9,3)(20,10,4)
\structline(20,10,1)(20,11,1)
\structline(20,10,2)(20,11,2)
\structline(20,10,3)(20,11,3)
\structline(20,10,4)(20,11,4)
\structline(20,10,5)(20,11,5)
\structline(20,11,1)(20,12,1)
\structline(20,11,2)(20,12,2)
\structline(20,11,3)(20,12,3)
\structline(20,11,4)(20,12,4)
\structline(20,11,5)(20,12,5)
\structline(20,12,1)(20,13,1)
\structline(20,12,2)(20,13,2)
\structline(20,12,3)(20,13,3)
\structline(20,12,4)(20,13,4)
\structline(20,12,5)(20,13,5)
\structline(20,13,1)(20,14,1)
\structline(20,13,2)(20,14,2)
\structline(20,13,3)(20,14,3)
\structline(20,13,4)(20,14,4)
\structline(20,13,5)(20,14,5)
\structline(20,14,1)(20,15,1)
\structline(20,14,2)(20,15,2)
\structline(20,14,3)(20,15,3)
\structline(20,14,4)(20,15,4)
\structline(20,14,5)(20,15,5)
\structline(20,15,1)(20,16,1)
\structline(20,15,2)(20,16,2)
\structline(20,15,3)(20,16,3)
\structline(20,15,4)(20,16,4)
\structline(20,15,5)(20,16,5)
\structline(20,16,1)(20,17,1)
\structline(20,16,2)(20,17,2)
\structline(20,16,3)(20,17,3)
\structline(20,16,4)(20,17,4)
\structline(20,16,5)(20,17,5)
\structline(20,17,1)(20,18,1)
\structline(20,17,2)(20,18,2)
\structline(20,17,3)(20,18,3)
\structline(20,17,4)(20,18,4)
\structline(20,17,5)(20,18,5)
\structline(20,18,1)(20,19,1)
\structline(20,18,2)(20,19,2)
\structline(20,18,3)(20,19,3)
\structline(20,18,4)(20,19,4)
\structline(20,18,5)(20,19,5)
\structline(20,19,1)(20,20,1)
\structline(20,19,2)(20,20,2)
\structline(20,19,3)(20,20,3)
\structline(20,19,4)(20,20,4)
\structline(20,19,5)(20,20,5)
\structline(20,20,1)(20,21,1)
\structline(20,20,2)(20,21,2)
\structline(20,20,3)(20,21,3)
\structline(20,20,4)(20,21,4)
\structline(20,20,5)(20,21,5)
\structline(20,21,1)(20,22,1)
\structline(20,21,2)(20,22,2)
\structline(20,21,3)(20,22,3)
\structline(20,21,4)(20,22,4)
\structline(20,21,5)(20,22,5)
\structline(20,22,1)(20,23,1)
\structline(20,22,2)(20,23,2)
\structline(20,22,3)(20,23,3)
\structline(20,22,4)(20,23,4)
\structline(20,22,5)(20,23,5)
\structline(20,23,1)(20,24,1)
\structline(20,23,2)(20,24,2)
\structline(20,23,3)(20,24,3)
\structline(20,23,4)(20,24,4)
\structline(20,23,5)(20,24,5)
\structline(24,3,1)(25,4,1)
\structline(24,5,1)(24,6,1)
\structline(24,5,2)(24,6,1)
\structline(24,5,2)(24,6,2)
\structline(24,5,1)(25,6,1)
\structline(24,6,1)(24,7,2)
\structline(24,6,2)(24,7,2)
\structline(24,6,2)(24,7,3)
\structline(24,6,2)(24,7,4)
\structline(24,7,1)(24,8,1)
\structline(24,7,2)(24,8,2)
\structline(24,7,3)(24,8,3)
\structline(24,7,4)(24,8,4)
\structline(24,7,3)(25,8,1)
\structline(24,7,4)(25,8,1)
\structline(24,8,1)(24,9,1)
\structline(24,8,2)(24,9,2)
\structline(24,8,3)(24,9,3)
\structline(24,8,4)(24,9,4)
\structline(24,8,5)(24,9,5)
\structline(24,8,5)(25,9,1)
\structline(24,9,1)(24,10,1)
\structline(24,9,2)(24,10,2)
\structline(24,9,3)(24,10,3)
\structline(24,9,4)(24,10,4)
\structline(24,9,5)(24,10,5)
\structline(24,9,6)(25,10,1)
\structline(24,10,1)(24,11,1)
\structline(24,10,2)(24,11,2)
\structline(24,10,3)(24,11,3)
\structline(24,10,4)(24,11,4)
\structline(24,10,5)(24,11,5)
\structline(24,11,1)(24,12,1)
\structline(24,11,2)(24,12,3)
\structline(24,11,3)(24,12,4)
\structline(24,11,4)(24,12,5)
\structline(24,11,5)(24,12,6)
\structline(24,11,6)(24,12,7)
\structline(24,12,1)(24,13,1)
\structline(24,12,2)(24,13,2)
\structline(24,12,3)(24,13,3)
\structline(24,12,4)(24,13,4)
\structline(24,12,5)(24,13,5)
\structline(24,12,6)(24,13,6)
\structline(24,12,7)(24,13,7)
\structline(24,12,2)(25,13,1)
\structline(24,13,1)(24,14,1)
\structline(24,13,2)(24,14,2)
\structline(24,13,3)(24,14,3)
\structline(24,13,4)(24,14,4)
\structline(24,13,5)(24,14,5)
\structline(24,13,6)(24,14,6)
\structline(24,13,7)(24,14,7)
\structline(24,14,1)(24,15,1)
\structline(24,14,2)(24,15,2)
\structline(24,14,3)(24,15,3)
\structline(24,14,4)(24,15,4)
\structline(24,14,5)(24,15,5)
\structline(24,14,6)(24,15,6)
\structline(24,14,7)(24,15,7)
\structline(24,15,1)(24,16,1)
\structline(24,15,2)(24,16,2)
\structline(24,15,3)(24,16,3)
\structline(24,15,4)(24,16,4)
\structline(24,15,5)(24,16,5)
\structline(24,15,6)(24,16,6)
\structline(24,15,7)(24,16,7)
\structline(24,16,1)(24,17,1)
\structline(24,16,2)(24,17,2)
\structline(24,16,3)(24,17,3)
\structline(24,16,4)(24,17,4)
\structline(24,16,5)(24,17,5)
\structline(24,16,6)(24,17,6)
\structline(24,16,7)(24,17,7)
\structline(24,17,1)(24,18,1)
\structline(24,17,2)(24,18,2)
\structline(24,17,3)(24,18,3)
\structline(24,17,4)(24,18,4)
\structline(24,17,5)(24,18,5)
\structline(24,17,6)(24,18,6)
\structline(24,17,7)(24,18,7)
\structline(24,18,1)(24,19,1)
\structline(24,18,2)(24,19,2)
\structline(24,18,3)(24,19,3)
\structline(24,18,4)(24,19,4)
\structline(24,18,5)(24,19,5)
\structline(24,18,6)(24,19,6)
\structline(24,18,7)(24,19,7)
\structline(24,19,1)(24,20,1)
\structline(24,19,2)(24,20,2)
\structline(24,19,3)(24,20,3)
\structline(24,19,4)(24,20,4)
\structline(24,19,5)(24,20,5)
\structline(24,19,6)(24,20,6)
\structline(24,19,7)(24,20,7)
\structline(24,20,1)(24,21,1)
\structline(24,20,2)(24,21,2)
\structline(24,20,3)(24,21,3)
\structline(24,20,4)(24,21,4)
\structline(24,20,5)(24,21,5)
\structline(24,20,6)(24,21,6)
\structline(24,20,7)(24,21,7)
\structline(24,21,1)(24,22,1)
\structline(24,21,2)(24,22,2)
\structline(24,21,3)(24,22,3)
\structline(24,21,4)(24,22,4)
\structline(24,21,5)(24,22,5)
\structline(24,21,6)(24,22,6)
\structline(24,21,7)(24,22,7)
\structline(24,22,1)(24,23,1)
\structline(24,22,2)(24,23,2)
\structline(24,22,3)(24,23,3)
\structline(24,22,4)(24,23,4)
\structline(24,22,5)(24,23,5)
\structline(24,22,6)(24,23,6)
\structline(24,22,7)(24,23,7)
\structline(24,23,1)(24,24,1)
\structline(24,23,2)(24,24,2)
\structline(24,23,3)(24,24,3)
\structline(24,23,4)(24,24,4)
\structline(24,23,5)(24,24,5)
\structline(24,23,6)(24,24,6)
\structline(24,23,7)(24,24,7)
\structline(25,5,1)(26,6,1)
\structline(25,6,1)(26,7,1)
\structline(25,8,1)(26,9,1)
\structline(25,9,1)(26,10,1)
\structline(25,10,1)(26,11,1)
\structline(25,13,1)(26,14,1)
\structline(28,6,1)(28,7,1)
\structline(28,6,2)(28,7,2)
\structline(28,7,1)(28,8,2)
\structline(28,7,2)(28,8,4)
\structline(28,8,2)(28,9,1)
\structline(28,8,3)(28,9,2)
\structline(28,8,4)(28,9,3)
\structline(28,9,1)(28,10,1)
\structline(28,9,2)(28,10,2)
\structline(28,9,3)(28,10,4)
\structline(28,10,1)(28,11,1)
\structline(28,10,2)(28,11,2)
\structline(28,10,3)(28,11,3)
\structline(28,10,4)(28,11,4)
\structline(28,11,1)(28,12,1)
\structline(28,11,2)(28,12,2)
\structline(28,11,3)(28,12,3)
\structline(28,11,4)(28,12,4)
\structline(28,11,5)(28,12,5)
\structline(28,12,1)(28,13,1)
\structline(28,12,2)(28,13,2)
\structline(28,12,3)(28,13,3)
\structline(28,12,4)(28,13,4)
\structline(28,12,5)(28,13,5)
\structline(28,13,1)(28,14,2)
\structline(28,13,2)(28,14,3)
\structline(28,13,3)(28,14,4)
\structline(28,13,4)(28,14,5)
\structline(28,13,5)(28,14,6)
\structline(28,14,1)(28,15,1)
\structline(28,14,2)(28,15,2)
\structline(28,14,3)(28,15,3)
\structline(28,14,4)(28,15,4)
\structline(28,14,5)(28,15,5)
\structline(28,14,6)(28,15,6)
\structline(28,14,7)(28,15,7)
\structline(28,15,1)(28,16,1)
\structline(28,15,2)(28,16,2)
\structline(28,15,3)(28,16,3)
\structline(28,15,4)(28,16,4)
\structline(28,15,5)(28,16,5)
\structline(28,15,6)(28,16,6)
\structline(28,15,7)(28,16,7)
\structline(28,16,1)(28,17,1)
\structline(28,16,2)(28,17,2)
\structline(28,16,3)(28,17,3)
\structline(28,16,4)(28,17,4)
\structline(28,16,5)(28,17,5)
\structline(28,16,6)(28,17,6)
\structline(28,16,7)(28,17,7)
\structline(28,17,1)(28,18,1)
\structline(28,17,2)(28,18,2)
\structline(28,17,3)(28,18,3)
\structline(28,17,4)(28,18,4)
\structline(28,17,5)(28,18,5)
\structline(28,17,6)(28,18,6)
\structline(28,17,7)(28,18,7)
\structline(28,18,1)(28,19,1)
\structline(28,18,2)(28,19,2)
\structline(28,18,3)(28,19,3)
\structline(28,18,4)(28,19,4)
\structline(28,18,5)(28,19,5)
\structline(28,18,6)(28,19,6)
\structline(28,18,7)(28,19,7)
\structline(28,19,1)(28,20,1)
\structline(28,19,2)(28,20,2)
\structline(28,19,3)(28,20,3)
\structline(28,19,4)(28,20,4)
\structline(28,19,5)(28,20,5)
\structline(28,19,6)(28,20,6)
\structline(28,19,7)(28,20,7)
\structline(28,20,1)(28,21,1)
\structline(28,20,2)(28,21,2)
\structline(28,20,3)(28,21,3)
\structline(28,20,4)(28,21,4)
\structline(28,20,5)(28,21,5)
\structline(28,20,6)(28,21,6)
\structline(28,20,7)(28,21,7)
\structline(28,21,1)(28,22,1)
\structline(28,21,2)(28,22,2)
\structline(28,21,3)(28,22,3)
\structline(28,21,4)(28,22,4)
\structline(28,21,5)(28,22,5)
\structline(28,21,6)(28,22,6)
\structline(28,21,7)(28,22,7)
\structline(28,22,1)(28,23,1)
\structline(28,22,2)(28,23,2)
\structline(28,22,3)(28,23,3)
\structline(28,22,4)(28,23,4)
\structline(28,22,5)(28,23,5)
\structline(28,22,6)(28,23,6)
\structline(28,22,7)(28,23,7)
\structline(28,23,1)(28,24,1)
\structline(28,23,2)(28,24,2)
\structline(28,23,3)(28,24,3)
\structline(28,23,4)(28,24,4)
\structline(28,23,5)(28,24,5)
\structline(28,23,6)(28,24,6)
\structline(28,23,7)(28,24,7)
\structline(32,4,1)(32,5,1)
\structline(32,4,1)(32,5,2)
\structline(32,5,2)(32,6,1)
\structline(32,6,1)(32,7,2)
\structline(32,6,1)(32,7,4)
\structline(32,7,1)(32,8,1)
\structline(32,7,2)(32,8,2)
\structline(32,7,3)(32,8,2)
\structline(32,7,4)(32,8,2)
\structline(32,7,4)(32,8,3)
\structline(32,8,1)(32,9,1)
\structline(32,8,2)(32,9,2)
\structline(32,8,3)(32,9,2)
\structline(32,8,3)(32,9,3)
\structline(32,8,3)(32,9,5)
\structline(32,9,2)(32,10,1)
\structline(32,9,3)(32,10,2)
\structline(32,9,4)(32,10,4)
\structline(32,9,5)(32,10,3)
\structline(32,10,1)(32,11,4)
\structline(32,10,2)(32,11,1)
\structline(32,10,3)(32,11,2)
\structline(32,10,4)(32,11,5)
\structline(32,11,1)(32,12,3)
\structline(32,11,2)(32,12,4)
\structline(32,11,3)(32,12,1)
\structline(32,11,4)(32,12,2)
\structline(32,11,5)(32,12,5)
\structline(32,11,6)(32,12,6)
\structline(32,12,1)(32,13,1)
\structline(32,12,2)(32,13,2)
\structline(32,12,3)(32,13,3)
\structline(32,12,4)(32,13,4)
\structline(32,12,5)(32,13,5)
\structline(32,12,6)(32,13,6)
\structline(32,12,7)(32,13,7)
\structline(32,13,1)(32,14,1)
\structline(32,13,2)(32,14,2)
\structline(32,13,3)(32,14,3)
\structline(32,13,4)(32,14,4)
\structline(32,13,5)(32,14,5)
\structline(32,13,6)(32,14,6)
\structline(32,13,7)(32,14,7)
\structline(32,14,1)(32,15,1)
\structline(32,14,2)(32,15,2)
\structline(32,14,3)(32,15,3)
\structline(32,14,4)(32,15,4)
\structline(32,14,5)(32,15,5)
\structline(32,14,6)(32,15,6)
\structline(32,14,7)(32,15,7)
\structline(32,15,1)(32,16,1)
\structline(32,15,2)(32,16,3)
\structline(32,15,3)(32,16,4)
\structline(32,15,4)(32,16,5)
\structline(32,15,5)(32,16,6)
\structline(32,15,6)(32,16,7)
\structline(32,15,7)(32,16,8)
\structline(32,15,8)(32,16,9)
\structline(32,16,1)(32,17,1)
\structline(32,16,2)(32,17,2)
\structline(32,16,3)(32,17,3)
\structline(32,16,4)(32,17,4)
\structline(32,16,5)(32,17,5)
\structline(32,16,6)(32,17,6)
\structline(32,16,7)(32,17,7)
\structline(32,16,8)(32,17,8)
\structline(32,16,9)(32,17,9)
\structline(32,17,1)(32,18,1)
\structline(32,17,2)(32,18,2)
\structline(32,17,3)(32,18,3)
\structline(32,17,4)(32,18,4)
\structline(32,17,5)(32,18,5)
\structline(32,17,6)(32,18,6)
\structline(32,17,7)(32,18,7)
\structline(32,17,8)(32,18,8)
\structline(32,17,9)(32,18,9)
\structline(32,18,1)(32,19,1)
\structline(32,18,2)(32,19,2)
\structline(32,18,3)(32,19,3)
\structline(32,18,4)(32,19,4)
\structline(32,18,5)(32,19,5)
\structline(32,18,6)(32,19,6)
\structline(32,18,7)(32,19,7)
\structline(32,18,8)(32,19,8)
\structline(32,18,9)(32,19,9)
\structline(32,19,1)(32,20,1)
\structline(32,19,2)(32,20,2)
\structline(32,19,3)(32,20,3)
\structline(32,19,4)(32,20,4)
\structline(32,19,5)(32,20,5)
\structline(32,19,6)(32,20,6)
\structline(32,19,7)(32,20,7)
\structline(32,19,8)(32,20,8)
\structline(32,19,9)(32,20,9)
\structline(32,20,1)(32,21,1)
\structline(32,20,2)(32,21,2)
\structline(32,20,3)(32,21,3)
\structline(32,20,4)(32,21,4)
\structline(32,20,5)(32,21,5)
\structline(32,20,6)(32,21,6)
\structline(32,20,7)(32,21,7)
\structline(32,20,8)(32,21,8)
\structline(32,20,9)(32,21,9)
\structline(32,21,1)(32,22,1)
\structline(32,21,2)(32,22,2)
\structline(32,21,3)(32,22,3)
\structline(32,21,4)(32,22,4)
\structline(32,21,5)(32,22,5)
\structline(32,21,6)(32,22,6)
\structline(32,21,7)(32,22,7)
\structline(32,21,8)(32,22,8)
\structline(32,21,9)(32,22,9)
\structline(32,22,1)(32,23,1)
\structline(32,22,2)(32,23,2)
\structline(32,22,3)(32,23,3)
\structline(32,22,4)(32,23,4)
\structline(32,22,5)(32,23,5)
\structline(32,22,6)(32,23,6)
\structline(32,22,7)(32,23,7)
\structline(32,22,8)(32,23,8)
\structline(32,22,9)(32,23,9)
\structline(32,23,1)(32,24,1)
\structline(32,23,2)(32,24,2)
\structline(32,23,3)(32,24,3)
\structline(32,23,4)(32,24,4)
\structline(32,23,5)(32,24,5)
\structline(32,23,6)(32,24,6)
\structline(32,23,7)(32,24,7)
\structline(32,23,8)(32,24,8)
\structline(32,23,9)(32,24,9)
\structline[hiddeneta](21,5,1)(22,8,1)
\structline[hiddeneta](25,4,1)(26,6,2)
\draw[hiddencol, dashed, line width = 0.45pt, -{Latex[length = 2.2pt]}, shorten >= 1.2pt, shorten <= 1.2pt, page constraint = {\page >= 5}] (16,5,3) to[out = 0, in = 0, looseness = 1.7] (16,7,3);
\draw[hiddencol, dashed, line width = 0.45pt, -{Latex[length = 2.2pt]}, shorten >= 1.2pt, shorten <= 1.2pt, page constraint = {\page >= 5}] (24,9,6) to[out = 0, in = 0, looseness = 1.7] (24,11,6);
\draw[hiddencol, dashed, line width = 0.45pt, -{Latex[length = 2.2pt]}, shorten >= 1.2pt, shorten <= 1.2pt, page constraint = {\page >= 5}] (32,13,8) to[out = 0, in = 0, looseness = 1.7] (32,15,8);
\end{sseqdata}

\clearpage
\printbibliography

\end{document}